\documentclass[12pt,a4paper]{amsart}
\usepackage{amscd}
\usepackage{amssymb}
\usepackage[centering,text={14.5cm,22cm}]{geometry}
\usepackage{graphicx}
\usepackage[dvipsnames]{xcolor}
\usepackage[all]{xy}
\usepackage{mathrsfs}
\usepackage{marvosym}
\usepackage{stmaryrd}
\usepackage{srcltx}
\usepackage{comment}
\usepackage[colorlinks=true,linkcolor=blue]{hyperref}
\usepackage{url}

\definecolor{shadecolor}{rgb}{1,0.9,0.7}

\newtheorem{theorem}{Theorem}[section]
\newtheorem{lemma}[theorem]{Lemma}
\newtheorem{proposition}[theorem]{Proposition}
\newtheorem{corollary}[theorem]{Corollary}
\newtheorem{conjecture}[theorem]{Conjecture}
\newtheorem{claim}[theorem]{Claim}

\theoremstyle{definition}
\newtheorem{definition}[theorem]{Definition}
\newtheorem{assumptions}[theorem]{Assumptions}
\newtheorem{construction}[theorem]{Construction}

\newtheorem{example}[theorem]{Example}
\newtheorem{examples}[theorem]{Examples}

\newtheorem{basicsetup}[theorem]{Basic Setup}

\theoremstyle{remark}
\newtheorem{remark}[theorem]{Remark}
\newtheorem{remarks}[theorem]{Remarks}

\numberwithin{equation}{section}
\numberwithin{figure}{section}

\newcommand{\NN} {\mathbb{N}}
\newcommand{\ZZ} {\mathbb{Z}}
\newcommand{\QQ} {\mathbb{Q}}
\newcommand{\RR} {\mathbb{R}}
\newcommand{\CC} {\mathbb{C}}

\newcommand{\PP} {\mathbb{P}}
\renewcommand{\AA} {\mathbb{A}}
\newcommand{\GG} {\mathbb{G}}
\newcommand{\LL} {\mathbb{L}}

\newcommand {\shA}  {\mathcal{A}}

\newcommand {\shD}  {\mathcal{D}}

\newcommand {\shI}  {\mathcal{I}}

\newcommand {\shK}  {\mathcal{K}}
\newcommand {\shL}  {\mathcal{L}}
\newcommand {\shM}  {\mathcal{M}}

\newcommand {\shN}  {\mathcal{N}}

\newcommand {\shU}  {\mathcal{U}}

\newcommand {\shP}  {\mathcal{P}}

\newcommand {\shX}  {\mathcal{X}}
\newcommand {\shY}  {\mathcal{Y}}
\newcommand {\shZ}  {\mathcal{Z}}
\newcommand {\bsigma} {\boldsymbol{\sigma}}

\newcommand {\foC}  {\mathfrak{C}}
\newcommand {\foD}  {\mathfrak{D}}

\newcommand {\foM}  {\mathfrak{M}}

\newcommand {\foU}  {\mathfrak{U}}
\newcommand {\foV}  {\mathfrak{V}}

\newcommand {\foX}  {\mathfrak{X}}

\newcommand {\fom}  {\mathfrak{m}}

\newcommand {\bas} {\mathrm{bas}}

\newcommand {\cl}  {\operatorname{cl}}

\newcommand {\Cbf} {{\mathbf{C}}}

\newcommand {\Div}  {\operatorname{Div}}

\newcommand {\ev}  {\operatorname{ev}}

\newcommand {\fine} {\mathrm{fine}}
\newcommand {\fs} {\mathrm{fs}}

\newcommand {\gp}  {{\operatorname{gp}}}
\newcommand {\gl}  {\mathrm{gl}}

\newcommand {\Hom}  {\operatorname{Hom}}

\newcommand {\id}  {\operatorname{id}}

\newcommand {\Int}  {\operatorname{Int}}

\renewcommand {\ker } {\operatorname{ker}}
\newcommand {\kk} {\Bbbk}

\newcommand {\Log} {\mathcal{L}og\hspace{1pt}}

\newcommand {\lra}  {\longrightarrow}
\newcommand {\ls}  {\dagger}

\newcommand {\M} {\mathcal{M}}

\newcommand {\Mbf} {\mathbf{M}}

\renewcommand{\O}  {\mathcal{O}}

\newcommand {\ol} {\overline}

\newcommand {\out}  {\mathrm{out}}

\renewcommand{\P}  {\mathscr{P}}

\newcommand {\Pic}  {\operatorname{Pic}}

\newcommand {\pr}  {\operatorname{pr}}

\newcommand {\Proj} {\operatorname{Proj}}

\newcommand {\rank} {\operatorname{rank}}

\newcommand {\red}  {{\operatorname{red}}}

\newcommand {\rk} {\operatorname{rk}}

\newcommand {\scrM}  {\mathscr{M}}
\newcommand {\scrC}  {\mathscr{C}}

\newcommand {\scrP}  {\mathscr{P}}

\newcommand {\sat}  {{\operatorname{sat}}}

\newcommand {\Spec} {\operatorname{Spec}}
\newcommand {\Spf}  {\operatorname{Spf}}

\newcommand {\spl} {\mathrm{spl}}

\newcommand {\tf} {\mathrm{tf}}

\newcommand {\ul} {\underline}

\newcommand {\virt} {\mathrm{virt}}

\newcommand{\hooklongrightarrow}{\lhook\joinrel\longrightarrow}

\def\mapright#1{\smash{
  \mathop{\longrightarrow}\limits^{#1}}}

\def\mydate{\ifcase\month \or January\or February\or March\or
April\or May\or June\or July\or August\or September\or October\or 
November\or December\fi \space\number\day,\space\number\year}

\begin{document}

%===========================================================
\title{Intrinsic mirror symmetry}

\author{Mark Gross} \address{DPMMS, Centre for Mathematical Sciences,
Wilberforce Road, Cambridge, CB3 0WB, UK}
%\curraddr{}
\email{mgross@dpmms.cam.ac.uk}
\thanks{This work was partially supported by the ERC Advanced Grant
MSAG, 
EPSRC grant EP/N03189X/1 and a Royal Society Wolfson Research
Merit Award.}
\thanks{AMS MSC: 14J33, 14N35, 14A21.}

\author{Bernd Siebert} \address{
Department of Mathematics, RLM 9.160, University of Texas at Austin,
2515 Speedway Stop C1200, Austin TX 78712
}
%\curraddr{}
\email{siebert@math.utexas.edu}

\date{\today}
\maketitle
\tableofcontents

%\bigskip

%===========================================================
%===========================================================
%===========================================================
%===========================================================

\section*{Introduction}
For some time, we have been developing a program for constructing
mirror pairs and understanding mirror symmetry 
\cite{PartI},\cite{PartII},\cite{Annals}. The last reference in particular
gave a construction of mirrors to what we call toric degenerations of
Calabi-Yau varieties. On the other hand, work with Hacking and
Keel in two dimensions,
\cite{GHKI} and \cite{GHKSK3}, has led to constructions of mirrors
to either log Calabi-Yau pairs $(X,D)$ or maximally unipotent
degenerations of K3 surfaces $X\rightarrow S$. In the first case,
$X$ is a non-singular projective surface and $D\in |-K_X|$
is a reduced nodal anti-canonical divisor. In the second case,
$X\rightarrow S$ is a relatively minimal but normal crossings degeneration.
Based on this experience, in \cite{GHS}, we, along with Paul Hacking,
set up the general framework which incorporates known and anticipated
constructions from both of the above points of view.

Roughy put, suppose we wish to construct the mirror to a log Calabi-Yau
pair $(X,D)$ or a maximally unipotent degeneration $X\rightarrow S$
of log Calabi-Yau manifolds. We call these the \emph{absolute}
and \emph{relative} cases respectively. We would carry out the following steps:
\begin{enumerate}
\item
Construct the dual intersection complex $(B,\P)$ of $(X,D)$ or
$X\rightarrow S$. Construct a \emph{wall structure} (also
often called a \emph{scattering diagram}) $\foD$ on $B$.
This is a data structure which controls the deformation theory of
an open subset $\check X_0^{\circ}$ of a singular scheme $\check X_0$
combinatorially determined by $B$. This enables the construction
of flat deformations $\check X_I^\circ\rightarrow \Spec S_I$ of
$\check X_0^{\circ}$ over some Artinian ring $S_I=\kk[P]/I$, 
for some monoid $P$ and co-Artinian monomial ideal $I$.
\item
Assuming the wall-crossing structure is \emph{consistent}, we may
use \emph{broken lines} to construct \emph{theta functions}. 
In the absolute case, these will
be functions on $\check X_I^{\circ}$ that
generate $\Gamma(\check X_I^{\circ},\O_{\check X_I^{\circ}})$ as an
$S_I$-module. In the relative case, there is an ample line bundle
$\shL$ on $\check X_I^{\circ}$ and the theta functions are sections
of various powers of $\shL$ which 
generate $\bigoplus_{d\ge 0} \Gamma(\check X_I^{\circ},\shL^{\otimes d})$.
We then obtain a (partial) compactification of $\check X_I^{\circ}$ by setting
\[
\hbox{$\check X_I = \Spec \Gamma(\check X_I^{\circ},\O_{\check X_I^{\circ}})$
or 
$\check X_I = \Proj \bigoplus_{d\ge 0} \Gamma(\check X_I^{\circ},\shL^{\otimes d})$}
\]
in the two cases. Taking the limit over all $I$ gives the mirror family.
\item 
The multiplication rule for theta functions can then be described
in terms of tropical geometry on $B$, essentially counting tropical
trees with two inputs and one output. Thus the coordinate ring or the
homogeneous coordinate ring of $\check X_I$ can be described directly
from the data of $B$ and $\foD$ without the intervening algebro-geometric
constructions.
\end{enumerate}

The paper \cite{GHS} explains the details of the above steps, but relies
on input a consistent wall-crossing structure $\foD$, without providing
an explanation for the source of such a diagram. Previously, 
such wall-crossing
structures have been constructed for nice toric degenerations
\cite{Annals}, log Calabi-Yau surfaces $(X,D)$ in \cite{GHKI},
and for maximally unipotent degenerations of K3 surfaces in
\cite{GHKSK3}. 

In \cite{Utah}, we announced two breakthroughs which allow us to
complete the above program in all generality. The key new
ingredient which allowed for these breakthroughs is the introduction
of \emph{punctured Gromov-Witten invariants}, \cite{ACGS18}, jointly
with Dan Abramovich and Qile Chen.
This is the first
of two papers devoted to these announced results. The second paper,
\cite{CanScat}, produces consistent wall-crossing structures in
all dimensions analogous to those constructed in \cite{GHKI} and
\cite{GHKSK3} in dimension two. That, by itself, completes the above
program, but \cite{CanScat} requires some assumptions on
$(X,D)$ which are only guaranteed in general
by the existence of good minimal models in the sense of \cite{NiXuYu}.
The existence of good minimal models for log Calabi-Yau varieties
is not in general known given the current state of 
the Minimal Model Program. The current paper, however, short circuits
the first two steps of the program, and describes the product rule
for theta functions by interpreting the tropical counts of
step (3) above in terms of (punctured) Gromov-Witten theory of the
pair $(X,D)$ or the degeneration $X\rightarrow S$. Furthermore,
the results of this paper apply in considerably broader contexts than
those of \cite{CanScat}, producing a ring associated to a pair $(X,D)$
with many fewer assumptions on $D$. See Theorems \ref{mainassociativity1} and
\ref{mainassociativity2} for precise statements. For example, if $D$ is
a smooth divisor with $K_X+D$ non-negative on all effective curve classes
or non-positive on all effective curve classes, our construction will give
an associative ring, a kind of ``relative quantum cohomology." 
However, such a simple example is inaccessible by the techniques of 
\cite{CanScat}.

A bit more precisely, the structure constants for the (homogeneous)
coordinate ring of the mirror are described in terms of actual
holomorphic curve counts of $3$-pointed curves with two ``input''
points and one ``output'' point. In particular, the theta functions
are indexed by integral points $B(\ZZ)$ of the dual intersection complex $B$,
and we write $\vartheta_p$ for the theta function indexed by $p\in B(\ZZ)$.
Then we desire structure constants for the product rule as an $S_I$-algebra
\[
\vartheta_{p}\cdot \vartheta_{q}=\sum_{r\in B(\ZZ)} \alpha_{pqr} \vartheta_r
\]
with $\alpha_{pqr}\in S_I=\kk[P]/I$ as before. 
Typically, we take $P$ a submonoid of $H_2(X,\ZZ)$
which contains the classes of all effective curves in $X$, is
saturated, and whose invertible
elements coincide with the torsion part of $H_2(X,\ZZ)$. See Basic
Setup \ref{setup3} for more details.
We write 
\[
\alpha_{pqr}=\sum_{A \in P\setminus I}
N^A_{pqr} t^{A},
\]
where $t^{A}\in \kk[P]$ denotes the monomial corresponding to
$A\in P$. According to logarithmic Gromov-Witten theory, we
can view points of $B(\ZZ)$ as specifying contact orders for marked points
on stable maps. We shall give, in \S\ref{sec:invariant construction},
a precise definition for the numbers $N^A_{pqr}$ as a count
of \emph{punctured} log maps $f:(C,x_1,x_2,x_{\out})\rightarrow X$
representing the class $A$,
with $C$ of genus zero, with contact orders at the three points 
specified by $p$, $q$, and $-r$,
and satisfying a point constraint at $x_{\out}$.

The negative contact order at $x_{\out}$ is one of the novel features
of our construction, and is a key feature of punctured maps
\cite{ACGS18}. It may appear strange in algebraic geometry,
but it fits well with the construction of symplectic cohomology.
We defer further discussions of these invariants until 
\S\ref{constructionsection} and \S\ref{sec:invariant construction}.

While these invariants can be defined for any simple normal crossings or
toroidal crossings pair $(X,D)$ (subject to hypotheses on the log structure
induced by $D$ stated in \S\ref{constructionsection}), the main issue
is to prove associativity of the resulting product rule. In general
the product rule is not expected to be associative. Philosophically,
that is because we are only building a part of an actual ``relative
quantum cohomology ring", and
are truncating the product rule rather crudely, so that 
associativity may not be preserved.
We show associativity under two situations. If $\Theta_{X/\kk}$ denotes
the logarithmic tangent bundle, then we require that either (1) 
$\pm c_1(\Theta_{X/\kk})$ is nef, or (2) $D$ is simple normal crossings
and $K_X+D$ is numerically equivalent to an effective $\QQ$-divisor supported
on $D$. In the latter case, we do not use the entire dual intersection
complex of the pair $(X,D)$, but only the Kontsevich-Soibelman skeleton,
see \S\ref{logCYcasesection}. This is a crucial point, as it allows
us to avoid the existence of minimal models that
\cite{CanScat} relies on. Case (1) already 
includes the case of a minimal log Calabi-Yau manifold, where $K_X+D=0$,
but Case (2) allows non-minimal log Calabi-Yau manifolds. Thus we
have no restriction on $(X,D)$ arising from the Minimal Model Program.

While we do not wish to dwell on the philosophy behind this construction
(see however \cite{Utah} for more discussion)
let us mention the most important point, namely that this construction
should be viewed, in the log Calabi-Yau case, as an
algebro-geometric version of $SH^0(X\setminus D)$ (symplectic
cohomology). This fits well
with the conjectures stated in the first preprint version of
\cite{GHKI} concerning the relationship between symplectic cohomology
and the mirror construction given there, as well as \cite{GP21}.
For a precise statement, see Conjecture~\ref{conj:SH}.
For degenerations of
Calabi-Yau manifolds, the case of elliptic curves and the relation between
symplectic cohomology, tropical geometry and punctured invariants
has been studied in detail by H\"ulya Arg\"uz \cite{Ar}, lending
support to this philosophical point of view.

Most of the paper is devoted to a proof of the associativity theorems
in these two cases, Theorems~\ref{mainassociativity1} and
\ref{mainassociativity2}. We do not try to include arguments that
we have really constructed the (homogeneous) coordinate ring of the
mirror, and leave this for future work. The main argument here that
we are constructing mirrors is the fact that this construction 
completes step (3) above, and the comparison with the results of 
steps (1) and (2) is given in \cite{CanScat}. It should even
be an interesting problem to confirm directly that the construction given
here agrees with the mirror construction of \cite{Annals}, which
is known to generalize the Batyrev-Borisov construction of mirrors,
\cite{BB}, \cite{GBB}. This is only known in dimensions
one or two (the elliptic curve or K3 surface cases), see \cite{Gonch}.

We defer further technical details of our results until the next
section, which contains precise statements of the setup and
main theorems. We try our best in \S\ref{constructionsection} to avoid
technical issues of log geometry, and that section should be easily
readable to the algebraic geometer unfamiliar with the log world.
After this, the structure of the paper is as follows. In \S\ref{sec:review},
we review facts about logarithmic geometry, stable log maps and
punctured maps as needed for this paper. In 
\S\ref{sec:invariant construction}, we set up the necessary machinery
to define the invariants $N^{\beta}_{pqr}$ described above. 

\S\ref{section:sketch} is both a sketch of the proofs of the main theorems
and the beginning of said proofs. At its heart, the proof is exactly
as one would expect following the proof of the WDVV equations, i.e.,
associativity of quantum cohomology. One considers ways in which
four-pointed punctured maps varying in a virtually one-dimensional moduli
space
can degenerate. However, unfortunately
everything becomes much harder in the log world. We give many examples
of phenomena which need to be taken into account, and reduce the proof
of associativity to three key theorems. Each of the sections 
\ref{sec: modulus invariance}--\ref{section:key comparison}
then proves one of these theorems, while \S\ref{sec:forgetful} is
devoted to the study of forgetful morphisms in the logarithmic context,
enabling us to prove $\vartheta_0$ is a unit in the (homogeneous)
coordinate ring of the mirror. This is analogous to the fundamental class
axiom in ordinary Gromov-Witten theory.

\medskip

Finally, we remark that Sean Keel and Tony Yue Yu 
\cite{KY19}
have obtained results
in the same direction as this paper, using Yu's non-Archimedean approach
to Gromov-Witten theory. However, at this point, their approach
only applies to mirrors to affine log Calabi-Yau manifolds containing
a copy of an algebraic torus as a dense open subset. We have no such
constraints. On the other hand, in the case they consider,
the Gromov-Witten invariants necessary to define structure constants
for multiplication are proved to be enumerative, and hence are positive
integers. Positivity is not true in general: indeed, there exist simple
examples of log Calabi-Yau threefolds $(X,D)$ where $X\setminus D$ is
not affine, for which the structure constants are sometimes negative.
However, it is not known if integrality holds in general.
This positivity is a quite powerful observation in many situations, e.g., it explains 
the positivity in cluster theory proved in \cite{GHKK}. The 
construction of Keel and Yu has now been shown in \cite{J22a}
to agree with that of the
current paper in the cases covered by \cite{KY19}, giving the positivity
and integrality of our structure constants in this case.

\medskip

\emph{Notation:} We follow the convention that if $X$ is a log scheme or
stack, then $\ul{X}$ is the underlying scheme or stack. We almost always
write $\shM_X$ for the sheaf of monoids on $X$ and $\alpha_X:\shM_X
\rightarrow\O_X$ for the structure map.
We write $\Spec\kk^{\dagger}$ for the standard log point with
sheaf of monoids $\kk^{\times}\times \NN$.
 If $P$ is a
monoid, we write $P^{\vee}:=\Hom(P,\NN)$ and $P^*=\Hom(P,\ZZ)$. We also
write $P_{\RR}$ for the real cone spanned by $P\subseteq P^{\gp}\otimes_{\ZZ}
\RR$.
Given a rational polyhedral cone $\sigma$, a 
\emph{tangent vector} to the cone is 
an element of the real span of $\sigma$.

\medskip

\emph{Acknowledgements:} Our various collaborations have played a key
role in the development of the ideas presented here. This paper owes
an obvious debt to the collaboration with Paul Hacking and Sean Keel,
in which the tropical description of the product rule was first proved in
the surface case.
The collaboration with Dan Abramovich and Qile Chen developing
punctured invariants and a gluing formula for log Gromov-Witten theory
lies at the technical heart of this paper. Indeed, it was already in
2011 in conversations with Abramovich and Chen that it became clear
that there would be new types of logarithmic invariants which might
be useful for completing our program. Finally, it was in conversations
of B.S.\ with Mohammed Abouzaid in 2012
and M.G.\ with Daniel Pomerleano in 2015
that the significance of constructing the
ring in this paper became apparent, and the paper owes a great deal
to discussions with them. Andrew Kresch kindly provided technical help
with Chow groups of stacks, and in particular provided the proof of
Lemma~\ref{lem:kresch lemma}. Finally, we would also like to thank
Y.\ Lekili, D.~Ranganathan and J.\ Wise for useful conversations.
Finally, we thank the referees for their careful reading and 
helpful comments.

\section{The construction and its variants}
\label{constructionsection}

\subsection{Algebras associated to pairs}
\label{pairalgebra}

\begin{basicsetup}
\label{setup1}
We work over an algebraically closed field $\kk$ of characteristic zero.
We fix throughout a projective log smooth
morphism $g:X\rightarrow S$ where $X$ carries a Zariski log structure,
and either $S=\Spec\kk$ or $S$ is a non-singular \emph{non-complete}
curve carrying a divisorial log
structure coming from a single closed point $s_0\in S$.
We refer to these
two cases as the \emph{absolute} and \emph{relative} cases. 
We further assume throughout that
$\overline{\shM}^{\gp}_X\otimes_{\ZZ}\QQ$ 
is globally
generated, i.e., the natural map $\Gamma(X,\overline{\shM}_X^{\gp}\otimes_{\ZZ}
\QQ)
\rightarrow \overline{\shM}^{\gp}_{X,x}\otimes_{\ZZ}\QQ$ 
is surjective for all $x\in X$.\footnote{This
assumption is currently necessary for a good theory of punctured
invariants \cite{ACGS18}, see loc.\ cit.\ Thm.~3.12.
It is possible that this assumption could be 
removed from the current paper, but this would not enhance readability,
and at any rate the assumption always holds in the simple normal crossings case.}

\end{basicsetup}

\begin{example}
\label{snccaseexample1}
For readers unfamiliar with logarithmic geometry, the basic example to
keep in mind is that $X$ is determined by an underlying
smooth scheme $\ul{X}$ and a simple normal crossings divisor $D\subseteq
\ul{X}$, and we write this as $X=(\ul{X}, D)$. Here the condition that
$X$ be Zariski means $D$ is simple normal crossings rather than just normal
crossings, so that $\ul{X}=\PP^2$ with $D$ a nodal cubic is not allowed.
This constraint is imposed here since logarithmic Gromov-Witten theory is best
understood in the case of a Zariski target space, but we do not expect this
restriction to be essential to the theory.

In the absolute case, there is no further constraint on $X$. However, if $S$
is a curve, then a morphism $\ul{g}:\ul{X}\rightarrow \ul{S}$ induces
a log morphism $g:X\rightarrow S$ provided that $\ul{g}^{-1}(s_0)\subseteq D$,
and it will be log smooth if further $\ul{g}$ is \'etale locally of the form
$\ul{\AA}^n\rightarrow \ul{\AA}^1$ given by $(x_1,\ldots,x_n)\mapsto 
\prod_{i=1}^r x_i^{a_i}$ for some collection of integers $a_i\in\NN$
 and $D\subseteq \ul{\AA}^n$ is given by
$D=V(\prod_{j\in J} x_j)$ for some index set $J\subseteq \{1,\ldots,n\}$
with $\{1,\ldots,r\}\subseteq J$.

A key example for us is the traditional setting for mirror
symmetry, where  $X\rightarrow S$ is a maximally unipotent degeneration
of Calabi-Yau manifolds, (see e.g., \cite{Norway}, Def.~16.17),
which we can then assume is log smooth
via semi-stable reduction \cite{KKMS}, p.~53. If in fact $K_{X/S}=0$,
then maximal unipotency is equivalent to the existence of a zero-dimensional
stratum of the central fibre $X_{s_0}$.

More generally, if $X$ is log smooth over $S=\Spec\kk$ or $S$ a curve
with divisorial log structure, the log structure on $X$ comes from
a divisor $D\subseteq \ul{X}$ with $(\ul{X}, D)$ \emph{toroidal crossings},
i.e., \'etale locally $\ul{X}$ is a toric variety and $D$ is its toric
boundary. \qed
\end{example}

\begin{basicsetup}
\label{setup2}
Functorially associated to the log scheme $X$ is its \emph{tropicalization},
a generalized cone
complex $\Sigma(X)$ with underlying topological space
$|\Sigma(X)|$ presented as a colimit of a diagram of 
rational polyhedral cones with maps given by integral inclusions of faces,
see \cite{ACGS17}, \S2.1.4. See also \S\ref{sec:tropical curves}
for a brief review. If $S$ is a point, then 
$|\Sigma(S)|$ is a point, while if $S$ is a curve, $|\Sigma(S)|=\RR_{\ge 0}$.
In any event, the structure morphism $X\rightarrow S$ induces a map of
cone complexes $\Sigma(X) \rightarrow \Sigma(S)$. We denote by
$\Sigma(X)(\ZZ)$ the set of points of $|\Sigma(X)|$ which are the image 
of some integral point $p\in \sigma$ for some $\sigma\in \Sigma(X)$.
In particular, we obtain a map $\Sigma(X)(\ZZ)\rightarrow \Sigma(S)(\ZZ)
=\NN$ in the case that $\dim S=1$.

We note that since $X$ is Zariski and log smooth over $\Spec\kk$ in both the
absolute and relative cases, 
there is a one-to-one correspondence between connected logarithmic strata of
$X$ and cones of $\Sigma(X)$, and that in fact $\Sigma(X)$ is a cone complex,
rather than a generalised cone complex. In particular there are no 
self-identifications of faces --- this is the classical case of toroidal
embeddings in \cite{KKMS}. Note this allows distinct cones meeting
along a union of faces. If further $X$ is in fact toric with $D$
its toric boundary, then $\Sigma(X)$ coincides, as an abstract cone complex,
with the fan of $X$.
\end{basicsetup}

\begin{example}
\label{snccaseexample2}
Again, to describe these objects without reference to log geometry, 
suppose that we are in the simple normal crossings case described in 
Example~\ref{snccaseexample1}. Then we can write $D=D_1+\cdots+D_s$
as a sum of irreducible divisors. Let us assume further that for
any index set $I\subseteq \{1,\ldots,s\}$, 
\[
D_I:=\bigcap_{i\in I} D_i
\]
is connected. Then it is useful to describe $\Sigma(X)$ as the dual cone
complex to $D$ in the following way. Let
$\Div_D(X)$ denote the group of divisors of $X$ supported on $D$, and
write $\Div_D(X)^*=\Hom(\Div_D(X),\ZZ)$, $\Div_D(X)^*_{\RR}
=\Hom(\Div_D(X),\RR)$. Write $D_1^*,\ldots,D_s^*$ for the basis of the
latter vector space dual to $D_1,\ldots,D_s$. We then define $\Sigma(X)$
to be the collection of cones
\[
\Sigma(X):=\left\{ \sum_{i\in I} \RR_{\ge 0} D_i^*\,\bigg|\,
\hbox{$I\subseteq \{1,\ldots,s\}$ an index set with
$D_I\not=\emptyset$}\right\}.
\]
We then have
\[
|\Sigma(X)|=\bigcup_{\sigma\in \Sigma(X)} \sigma \subseteq \Div_D(X)^*_{\RR}
\]
and
\[
\Sigma(X)(\ZZ)=\left\{\sum a_i D_i^*\,\bigg|\,
a_i\in \NN, \bigcap_{i: a_i>0} D_i\not=\emptyset\right\}.
\]
In the relative case, the map $\Sigma(g):\Sigma(X)
\rightarrow \Sigma(S)$ is then easily described as the map induced
by the linear map $\Div_D(X)^*_{\RR}\rightarrow \RR$ given by evaluation
on the (not necessarily reduced) divisor $g^*(s_0)$. \qed
\end{example}

\begin{examples}
\label{runningexample1}
We consider two running examples. They are very similar and our construction
is quite simple in these cases, but they exhibit slightly different features.

(1) Consider the log Calabi-Yau $(\PP^2,\bar D=\bar D_1\cup \bar D_2\cup \bar D_3)$ 
where the $\bar D_i$ are lines, and blow up a point on $\bar D_1$ which is not
a double point of $\bar D$ to get $X$.
Take $D_i$ to be the strict transform of $\bar D_i$, getting a pair
$(X,D)$ with $D=\sum D_i$. In this case $\Sigma(X)$ consists of three
two-dimensional cones; as a polyhedral complex in $\Div_D(X)^*_{\RR}$,
it consists of the cones $\RR_{\ge 0}D_i^*+\RR_{\ge 0}D_j^*$ for $i\not=j$
and their faces.
Of course $|\Sigma(X)|$ is homeomorphic to $\RR^2$.

(2) As a second example, consider the case of $\ul{X}=\PP^2$ and
$D=D_1+D_2$, with $D_1$ a line and $D_2$ a conic meeting transversally
at two points. This does not satisfy the connectedness
hypotheses of Example 
\ref{snccaseexample2}, but we can still easily 
describe $\Sigma(X)$, which consists
of five cones: a zero dimensional cone $0$ corresponding to the open
stratum of $X$, two one-dimensional cones $\rho_1$ and $\rho_2$ corresponding
to the open strata of $D_1$, $D_2$, and two two-dimensional cones
$\sigma_1$ and $\sigma_2$ corresponding to the two points of $D_1\cap D_2$.
The boundary of each $\sigma_i$ is $\rho_1\cup\rho_2$, so that $\sigma_1$
and $\sigma_2$ are glued together along their entire boundaries. Again,
$|\Sigma(X)|$ is homeomorphic to $\RR^2$, split into two two-dimensional cones.
See Figure~\ref{fig:RunningExpl2}.
\qed
\end{examples}

\begin{basicsetup}
\label{setup3}
We choose additional data of a finitely generated abelian group $H_2(X)$ 
of possible degree data for curves on $X$. In the absolute case,
this could be $1$-cycles
on $X$ modulo algebraic or numerical equivalence, or it could be
$\Hom(\Pic(X),\ZZ)$. If we work over $\CC$, we can use ordinary singular
homology $H_2(X,\ZZ)$. In the relative case, we can use the group generated
by curves contracted by the morphism $X\rightarrow S$, modulo algebraic
or numerical equivalence, or again, over $\CC$, use $H_2(X,\ZZ)$. 
There are three key requirements this group must satisfy:
\begin{enumerate}
\item Any
stable map $\ul{f}:\ul{C}/\ul{W}\rightarrow \ul{X}$ 
with $\ul{W}$ connected induces a
curve class $f_*[\ul{C}_{\bar w}]\in H_2(X)$ whenever $\bar w\in \ul{W}$ is
a geometric point, and this class is independent of the choice of $
\bar w$.
\item If $f_*[\ul{C}_{\bar w}]$ is
a torsion class in $H_2(X)$, then in fact $\ul{f}$ is constant on
$\ul{C}_{\bar w}$ and $f_*[\ul{C}_{\bar w}]=0$.
\item
The moduli space of (ordinary)
stable maps of class $A\in H_2(X)$, 
genus $g$ and $n$ marked points to $\ul{X}$ is
of finite type. 
\end{enumerate}

We also choose a finitely generated 
monoid $P\subseteq H_2(X)$ such that (1) $P$ contains
the classes of all stable maps to $\ul{X}$; (2) $P$ is saturated; (3) 
the group of invertible elements of $P$ coincides with the torsion part
of $H_2(X)$.
\end{basicsetup}

\begin{example}
In Examples~\ref{runningexample1}, in the first example we take
$H_2(X)=\Pic(X)^*=\ZZ[L-E]\oplus\ZZ[E]$, where $L$ is the total transform of
a line in $\PP^2$ and $E$ is the exceptional curve, and take
$P$ to be the submonoid generated by $[L-E]$ and $[E]$. In the second
case, we take $H_2(X)=\Pic(X)^*=\ZZ$,
and $P=\NN$, generated by the class of a line. \qed
\end{example}

\begin{construction}
\label{const:first}
In what follows, we now choose any monoid ideal $I\subseteq P$
such that $P\setminus I$ is finite, i.e., 
\begin{equation}
\label{eq:SI}
S_I:= \kk[P]/I
\end{equation}
is Artinian. We always confuse a monoid ideal $I$ with the monomial
ideal in $\kk[P]$ it generates.
We use the convention that given an element $A\in P$,
$t^A\in \kk[P]$ is the corresponding element.
We then define 
\begin{equation}
\label{defofRI}
R_I:=\bigoplus_{p\in \Sigma(X)(\ZZ)} S_I \vartheta_p,
\end{equation}
a free $S_I$-module, and our main goal is to define, in certain cases,
an algebra structure on $R_I$ by defining structure constants
\begin{equation}
\label{structureequation}
\vartheta_{p_1}\cdot\vartheta_{p_2}=\sum_{r\in \Sigma(X)(\ZZ)}
\alpha_{p_1p_2r} \vartheta_r,
\end{equation}
with only a finite number of $\alpha_{p_1p_2r}\in S_I$ non-zero.
In turn, we can write
\begin{equation}
\label{structureconstants}
\alpha_{p_1p_2r} = \sum_{A \in P\setminus I} 
N^A_{p_1p_2r} 
t^{A}.
\end{equation}
Here, in fact we will define
$N^A_{p_1p_2r} 
\in \QQ$, so that the construction yields schemes defined over $\QQ$.

The notion of \emph{contact order} specified
by a tangent vector to a cone of $\Sigma(X)$ is reviewed in 
\S\ref{sec:contact orders}. Contact orders are a way of imposing
tangency conditions on stable or punctured curves. For example,
in the situation of Example~\ref{snccaseexample2}, if
$p\in \Div_D(X)^*$ is tangent to a cone of $\Sigma(X)$, 
then we interpret $p$ as a possible
tangency condition for a punctured stable map $f:(C,x)\rightarrow X$ at
$x$, demanding that $f$ make contact with the divisor $D_i$ at $x$
to order $\langle D_i, p\rangle$.

In Definition~\ref{def:Nbeta}, we will define $N^A_{p_1p_2r}$ 
as a punctured Gromov-Witten
invariant, as introduced in \cite{ACGS18}. 
Roughly put, this is defined as the number of genus zero three-pointed
punctured maps representing the curve class $A$,
with marked points $x_1,x_2,x_{\out}$ having contact
orders $p_1,p_2$ and $-r$ respectively and with a certain point
constraint at $x_{\out}$.

The point with
contact order $-r$, being of negative contact order provided $r\not=0$,
is made sense of using the notion of a punctured point in \cite{ACGS18},
reviewed in \S\ref{sec:punctured log curves}.
We further impose a point constraint at $x_{\out}$. This
is done as follows. First, any $r\in \Sigma(X)(\ZZ)$ determines a (locally
closed)
stratum $Z_r^{\circ}\subset X$. Indeed, let $\sigma\in\Sigma(X)$ be the smallest
cone containing $r$. Then $\sigma$ corresponds to a 
locally closed stratum $Z_r^{\circ}$,
whose closure we denote by $Z_r$. The open subset $Z_r^{\circ}\subset
Z_r$ is the complement of the union of lower dimensional strata contained
in $Z_r$.
Second, choose a point $z\in Z_r^{\circ}$. At a schematic level,
we only count punctured maps which map $x_{\out}$ to $z$. However,
this constraint is imposed logarithmically, and to explain this
properly requires more development. Thus we defer this definition
to \S\ref{sec:invariant construction}.
We summarize the discussion there by stating that we define a moduli
space $\scrM(X,\beta,z)$ associated to the choice of data $A,
p_1,p_2,r$ and $z$ whose virtual dimension is $c_1(\Theta_{X/\kk})\cdot
A$. We then define
\[
N^{A}_{p_1p_2r} = 
\begin{cases} \deg [\scrM(X,\beta,z)]^{\virt} & \hbox{if 
$\virt.\dim \scrM(X,\beta,z)
=0$,}\\
0 & \hbox{otherwise.}
\end{cases}
\]
\qed
\end{construction}

The first main theorem of the paper is

\begin{theorem} 
\label{mainassociativity1}
Let $X$ be a Zariski log scheme with a projective log smooth morphism
$g:X\rightarrow S$ with $S\cong\Spec\kk$ or $S$ a non-complete curve
carrying a divisorial log structure coming from a single point $s_0\in S$.
If $c_1(\Theta_{X/\kk})$ is nef or anti-nef,
the structure constants $N^A_{p_1p_2r}$ define, via 
\eqref{structureequation}, a commutative,
associative $S_I$-algebra structure on $R_I$ with unit given by $\vartheta_0$. 
\end{theorem}

The proof of this theorem will be completed in 
\S\ref{section:sketch}--\ref{section:key comparison}.

\begin{remark}
Here $\Theta_{X/\kk}$ denotes the \emph{logarithmic tangent bundle}.
In particular, $c_1(\Theta_{X/\kk})=-(K_X+D)$. 
We recall for more general audiences that a divisor $H$ is \emph{nef} if
$H\cdot A\ge 0$ for all curve classes $A$ represented by an effective
curve class. We say $H$ is \emph{anti-nef} if $-H$ is nef. 
\end{remark}

\subsection{The log Calabi-Yau case}
\label{logCYcasesection}
We continue with the situation of Basic Setups~\ref{setup1}, \ref{setup2},
and \ref{setup3}.

\begin{definition}
Suppose the divisor $D\subseteq \ul{X}$ is simple normal crossings.
We say that the pair $(\ul{X}, D)$ is \emph{log Calabi-Yau} if
$-c_1(\Theta_{X/\kk})=K_X+D$ is numerically equivalent 
to an effective $\QQ$-divisor supported on $D$.
Similarly, if $\dim S=1$, we say $X\rightarrow S$ is \emph{relative
log Calabi-Yau} if $-c_1(\Theta_{X/S})$ is numerically equivalent to 
an effective $\QQ$-divisor supported on $D$.
\end{definition}

\begin{remark}
We note this definition of log Calabi-Yau is much more liberal than is
often the case in the literature. For example, from the point of view
of the Minimal Model Program, we do not insist on log Kodaira dimension
zero, while from the point of view of metric geometry, there is no
reason to expect a complete Ricci-flat metric on $X\setminus D$.
\end{remark}

In either the relative or absolute cases, after an additional
choice, we can define a
sub-cone complex $B\subseteq \Sigma(X)$, known as the \emph{Kontsevich-Soibelman
skeleton}, as follows. Write 
$D=D_1+\cdots+D_s$ with the $D_i$ the irreducible components of
$D$. We fix a choice of $a_i\in \QQ$, $1\le i\le s$, with $a_i\ge 0$
for all $i$ and
\begin{equation}
\label{eq:c1 log CY}
-c_1(\Theta_{X/S})\equiv_{\QQ}\sum_i a_i D_i.
\end{equation}
The construction will depend on this choice, 
see Example~\ref{ex:yanki}.

In the absolute case, we call a divisor $D_i$ \emph{good} if $a_i=0$. 

In the relative case, if
$D_i\subseteq g^{-1}(s_0)$, we define the \emph{weight} of $D_i$ to be $w_i:=a_i/\mu_i$,
where $\mu_i$ is the multiplicity of $D_i$ in the divisor $g^*(s_0)$.
 We call a divisor $D_i$
\emph{good} if either (1) $g|_{D_i}:D_i\rightarrow S$ is surjective and 
$a_i=0$, or (2) $g(D_i)=s_0\in S$, and $w_i=\min\{w_j\,|\, 
g(D_j)=s_0\}$.
Note that if $w=\min\{w_j\,|\, g(D_j)=s_0\}$, 
then as the divisor $g^*(s_0)
=\sum_i \mu_i D_i$
is numerically equivalent to the zero divisor, we can replace the expression
$\sum a_i D_i$ by the 
numerically equivalent expression $\sum a_i D_i - wg^*(s_0)$, and then
any divisor $D_i$ is good if and only if $a_i=0$. 

There is a one-to-one
correspondence between one-dimensional cones of $\Sigma(X)$ and 
the divisors $D_i$. Let $\scrP$ be the collection of cones of $\Sigma(X)$
whose one-dimensional faces correspond to divisors $D_i$ which are good,
and set
\[
B=\bigcup_{\sigma\in\P} \sigma.
\]
This is the Kontsevich-Soibelman skeleton of $X$, see \cite{KS06} and
\cite{NiXu}.

We remark that in this non-minimal log Calabi-Yau case, we 
restrict to the simple normal crossings case, rather than allowing more
general toroidal crossings (log smooth). This is crucial
because some key parts
of the argument involve intersection numbers of curves with individual
irreducible components of $D$, and thus we require these components to be
Cartier. This only occurs in the normal crossings case. In any event,
the normal crossings situation is the natural one from the point of view
of the theory of \cite{NiXu} and \cite{NiXuYu}.

\begin{remark}
In the most relevant case for mirror symmetry, the 
choice of \eqref{eq:c1 log CY} will be unique. This occurs when we begin
with a pair $(X',D')$ which is a \emph{log Calabi-Yau minimal model}.
This means that $(X',D')$ is a $\QQ$-factorial divisorial log terminal
(dlt) pair $(X',D')$ such that $K_{X'}+D'\equiv 0$. One takes $\phi:X\rightarrow
X'$ a resolution of singularities so that $\phi^{-1}(D)$ is an snc divisor.
In this case, we obtain an equality as in \eqref{eq:c1 log CY}
with the $a_i$ uniquely determined as the log discrepancies of the
exceptional divisors. See \cite{CanScat}, \S1.2 for more discussion.

In particular, the ability to work with non-minimal log Calabi-Yau
manifolds but restricting to the Kontsevich-Soibelman skeleton introduced
above gives us a powerful tool for circumventing the non-existence of
minimal models with snc boundary divisor.
\end{remark}

\begin{remark}
In general, there is very little that can be said about the topology
of the topological space $|\Sigma(X)|$. However, in the log Calabi-Yau
case, $B$ is expected to be a strong deformation retract of $\Sigma(X)$,
see \cite{NiXu}. The topological space $B$ is expected to have nicer
properties, see e.g., \cite{KollarXu}, where evidence is given for the
folklore conjecture that if $\dim_{\RR} B=\dim X$, then $B$ 
is the cone over a finite quotient of a sphere. 
\end{remark}

We then define $B(\ZZ)=B\cap \Sigma(X)(\ZZ)$. Using a choice of $P
\subseteq H_2(X)$ as in \S\ref{pairalgebra}, we define
\begin{equation}
\label{eq:logCY RI}
R_I:=\bigoplus_{p\in B(\ZZ)} S_I \vartheta_p.
\end{equation}
Similarly as before we define
\begin{equation}
\label{structureequation2}
\vartheta_{p_1}\cdot\vartheta_{p_2}=\sum_{r\in B(\ZZ)}
\alpha_{p_1p_2r} \vartheta_r,
\end{equation}
with $\alpha_{p_1p_2r}$ defined as in \eqref{structureconstants}.

\begin{theorem}
\label{mainassociativity2}
Let $X$ be a Zariski log scheme with a projective log smooth morphism
$g:X\rightarrow S$ with $S\cong\Spec\kk$ or $S$ a non-complete curve
carrying a divisorial log structure coming from a single point $s_0\in S$.
Suppose further we are given an expression \eqref{eq:c1 log CY} with $a_i\ge 0$
for all $i$. Let $B$ be the Kontsevich-Soibelman skeleton defined by the
choice of expression in \eqref{eq:c1 log CY}. Let $R_I$ be as in
\eqref{eq:logCY RI}. Then the
structure constants $N^A_{p_1p_2r}$ define, via
\eqref{structureequation2}, a
commutative, associative $S_I$-algebra structure on $R_I$ with
unit given by $\vartheta_0$.
\end{theorem}

This will also be proved in \S\ref{section:sketch}-\ref{section:key comparison}
in parallel with 
Theorem~\ref{mainassociativity1}.

\subsection{Mirrors, examples, basic properties, and auxilliary constructions}
\label{sec:mirrors}

The following result is generally useful for understanding when punctured
stable maps might exist, and has many consequences for us. Recall
first that any $s\in\Gamma(X,\overline\shM^{\gp}_X)$ yields an 
$\O_X^{\times}$-torsor, under pull-back of $s$ via the quotient map
$\shM^{\gp}_X\rightarrow\overline{\shM}^{\gp}_X$. We denote this torsor
as $\shL^{\times}_s$, with corresponding line bundle $\shL_s$.
Further, for $\eta\in X$ the generic point of a stratum of $X$, the section 
$s$ can be evaluated 
on any element $u\in \Hom(\overline\shM_{X, \eta},\ZZ)$, the space
of integral tangent vectors to the cone $\sigma_{\eta}\in \Sigma(X)$
corresponding to the stratum. We write this evaluation as $\langle u,s\rangle$.
We can now recall the following proposition from \cite{ACGS18}, 
Cor.~2.30:

\begin{proposition}
\label{intersectionnumbers}
Let $X\rightarrow S$ be as in \S\ref{pairalgebra}, and
suppose given a punctured map $f:(C,x_1,\ldots,x_n)\rightarrow X$
with contact orders
given by $u_{x_1},\ldots,u_{x_n}$,
where $u_{x_i}$ is an integral tangent vector to a cone $\sigma_i\in
\Sigma(X)$. 
Then 
we must have, for any $s\in \Gamma(X,\overline{\shM}_X^{\gp})$,
\[
\deg \ul{f}^*(\shL_s)=-\sum_{i=1}^n \langle u_{x_i},s\rangle.
\]
\end{proposition}

In the situation of Example~\ref{snccaseexample2}, we can rephrase the
above condition, bearing in mind that any tangent vector to a cone
in $\Sigma(X)$ can be viewed as an element of $\Div_D(X)^*$:

\begin{corollary}
\label{cor:intersectionnumbers}
Suppose that $D=D_1+\cdots+D_s$ is simple normal crossings
as in Example~\ref{snccaseexample2}. 
Then in the situation
of Proposition~\ref{intersectionnumbers}, 
for any divisor $D'$ supported on $D$, we have
\[
\deg \ul{f}^*\O_X(D') = \sum_i \langle u_{x_i}, D' \rangle.
\]
\end{corollary}

\begin{proof}
This is a restatement of the previous proposition, bearing in mind that
in this case $\overline{\shM}_X = \bigoplus_{1\le i \le s} \NN_{D_i}$,
where $\NN_{D_i}$ is the constant sheaf on $D_i$ with stalk $\NN$.
A section $s$ given by $1 \in \NN_{D_i}$ yields the line bundle
$\O_X(-D_i)$, and the claim follows.
\end{proof}

In what follows, write $B=\Sigma(X)$ in the general case of \S\ref{pairalgebra}, and
$B\subseteq \Sigma(X)$ the Kontsevich-Soibelman skeleton
in the log Calabi-Yau case.
In both cases, we write $\P$ for the cones making up $B$. 

We will prove the following lemma in \S\ref{subsec:the invariants}.

\begin{lemma}
\label{lem:constant-maps}
If $A=0$, then 
$N^A_{p_1p_2r} = 1$ if there exists a cone $\sigma\in\P$
containing $p_1,p_2$ and $r$ with $p_1+p_2=r$ inside $\sigma$. Otherwise
$N^A_{p_1p_2r}=0$.
\end{lemma}

\begin{proposition}
\label{prop:degree 0 mult}
Let $\fom=P\setminus P^{\times}$ be the maximal monomial ideal; 
recall
by assumption, the invertible elements $P^{\times}$ of $P$ coincide with
the torsion part of $H_2(X)$.
Then the $S_{\fom}:=S/\fom=\kk[P^{\times}]$-algebra structure on $R_{\fom}$ is given by
\[
\vartheta_{p_1}\cdot\vartheta_{p_2}=
\sum_{\sigma} \vartheta_{r_{\sigma}} 
\]
where the sum is over all minimal cones $\sigma$ of $\P$ containing both
$p_1$ and $p_2$ and $r_{\sigma}$ is the sum of $p_1$ and $p_2$ inside the cone 
$\sigma$. In particular, if $p_1,p_2$ are not contained in a common
cone of $\P$, the product is zero.
\end{proposition}

\begin{proof}
Recall from Basic Setup~\ref{setup3} that if a curve class $A\in P^{\times}$
is represented by a stable map, then the stable map must be constant, and
in particular $A=0$. Thus we have the structure constants
\[
\alpha_{p_1p_2r} = N^0_{p_1p_2r}t^0 = N^0_{p_1p_2r}.
\]
Thus the result follows immediately from Lemma~\ref{lem:constant-maps}.
\end{proof}

\begin{example}
\label{runningexample3}
Continuing with Examples~\ref{runningexample1}, in the first case,
let $v_i\in B(\ZZ)$ be the primitive generator of the ray $\rho_i$
corresponding to $D_i$. Then one immediately obtains
\[
R_{\fom}\cong \kk[\vartheta_{v_1},\vartheta_{v_2},\vartheta_{v_3}]
/(\vartheta_{v_1}\vartheta_{v_2}\vartheta_{v_3}).
\]
Indeed, from the proposition, any expression of the form
$\vartheta_{v_i}^a\vartheta_{v_j}^b$ agrees with $\vartheta_{av_i+bv_j}$,
where $av_i+bv_j$ lies in the cone $\RR_{\ge 0}v_i+\RR_{\ge 0}v_j$.
On the other hand, $\vartheta_{v_1}\vartheta_{v_2}\vartheta_{v_3}
=\vartheta_{v_1+v_2}\vartheta_{v_3}=0$ since there is no cone containing
both $v_1+v_2$ and $v_3$.

In the second case,
\begin{figure}
\input{RunningExplCones.pspdftex}
\caption{Cones and vectors in the second case of Example~\ref{runningexample3}.
The shaded areas are the two cones, and $B$ identifies the rays $\rho_1$,
$\rho'_1$ and the tangent vectors $v'_1$ and $v_1$.}
\label{fig:RunningExpl2}
\end{figure}
let $v_1,v_2\in B(\ZZ)$ be the
primitive generators of the cones $\rho_1$, $\rho_2$ and let $v_{\sigma_1}, v_{\sigma_2}$
be the sum $v_1+v_2$ inside either $\sigma_1$ or $\sigma_2$
respectively. Then
from the above proposition, we have an isomorphism of $S_{\fom}=\kk$-algebras
\[
R_{\fom} \cong \kk[\vartheta_{v_1},\vartheta_{v_2},
\vartheta_{v_{\sigma_1}},\vartheta_{v_{\sigma_2}}]/(\vartheta_{v_1}
\vartheta_{v_2}-
\vartheta_{v_{\sigma_1}}-\vartheta_{v_{\sigma_2}},
\vartheta_{v_{\sigma_1}}\vartheta_{v_{\sigma_2}}).
\]
Note $\Spec R_{\fom}$ is a union of two planes intersecting along the union of
coordinate axes on the planes. See Figure~\ref{fig:RunningExpl2}.
\end{example}

\begin{construction}[Mirrors to log Calabi-Yau manifolds]
\label{mirrorconstruction1}
Morally, in the case that $X$ is log Calabi-Yau (in particular the
absolute case), one should view $R_I$ as the
ring of regular functions on a thickened neighbourhood of the large complex
structure limit of the mirror family to $X\setminus D$. More explicitly, in
the case that $B$ is pure-dimensional with $\dim_{\RR} B = \dim X$, we define
the mirror family to $X\setminus D$ by taking the direct limit of 
families of schemes $\check X_I:=\Spec R_I \rightarrow \Spec S_I$, getting a 
formal flat family of affine schemes
\[
\check \foX\rightarrow \Spf \widehat{\kk[P]},
\]
where $\widehat{\kk[P]}$ is the completion of $\kk[P]$ with respect to the maximal
ideal $\fom$. Note flatness is obvious since $R_I$ is a free, hence flat, $S_I$-module.

In the case that $\dim X=2$, we have $(X,D)$ a Looijenga pair, i.e., $X$ is a projective
non-singular rational surface and $D\in |-K_X|$ is a cycle of rational curves.
In this case, mirrors have already been constructed in \cite{GHKI}, and 
the construction of the present paper 
agrees with that of \cite{GHKI}, as follows from \cite{CanScat}, Ex.~3.14 and
Thm.~6.1. 
Note that in any event, the
definition of the product given here is a logarithmic version of the tropical
product
rule for theta functions on mirrors to Looijenga pairs described in 
\cite{GHKI}, Thm.~2.34, and generalized to all dimensions in 
\cite{GHS}, Thm.~3.24. The fact that the logarithmic and tropical definitions
of the product agree is \cite{CanScat}, Thm.~6.1.

Even in dimension three, however,
it is worth noting that the mirror family constructed
in this way may not be exactly what one anticipates as a mirror.
For example, Pomerleano gives in \cite{Pom}, App.~C, an example
for which the general fibre of our mirror family is a singular affine variety.
In this example, $X$ is a non-singular hypersurface in $\PP^2\times\PP^2$
of bidegree $(1,1)$, with $D=D_1+D_2+D_3$. One takes $D_1$, $D_2$ and
$D_3$ to be the intersections of $X$ with a general hypersurface
of bidegrees $(1,0)$, $(0,1)$ and $(1,1)$ respectively. While it may
be difficult to determine the precise mirror, Pomerleano
makes clever use of Corollary~\ref{cor:intersectionnumbers} 
to give a rough shape of the equations of
the mirror. The particular shape of the equations then
implies the existence of a curve of singularities in the generic fibre of
the mirror family.

This raises a philosophical question as to what the correct mirror
to $(X,D)$ is. The paper \cite{Pom} suggests that 
the wrapped Fukaya category of $X\setminus D$ is in fact a categorical crepant
resolution of the mirror we have constructed here. 
Thus it may be natural to choose a crepant resolution of the
mirror family, if one exists. Conjectures of Bondal and Orlov would
then suggest that the derived category of the resolution is independent
of the choice, and hence this resolution 
allows one to study questions of homological mirror symmetry.
Alternatively, one may feel that as a crepant resolution does not necessarily
exist and in any event is not canonical, 
it is more natural to keep the mirror family unresolved.

Alternatively, in many cases, one may
potentially be able to degenerate the pair $(X,D)$
to obtain a family $(\shX,\shD)\rightarrow S$, and apply 
Construction~\ref{mirrorconstruction2} below to obtain a non-singular
mirror which is the resolution of the one obtained from $(X,D)$. This
procedure has not been carried out in Pomerleano's example described above, 
but we expect this to be possible.
From another point of view, we may expect that for a log
Calabi-Yau pair $(X,D)$, the boundary has to be ``sufficiently degenerate''
in order for the mirror to be non-singular, just as a degeneration of
Calabi-Yau varieties $\shX\rightarrow S$ needs to be sufficiently
degenerate (a ``large complex structure limit'') in order to guarantee
a well-behaved mirror. It is an interesting question to determine
when the mirror is non-singular. Methods applied in \cite{Annals} or
\cite{GHKI} are insufficient to answer this question.

Of course, in general one does not necessarily expect the general fibre
of the mirror family to be non-singular: already in \cite{Bat}, \cite{BB},
mirrors to Calabi-Yau manifolds of dimension $\ge 4$ may be singular
with no possible crepant resolution. However, these mirrors may still
be viewed as non-singular orbifolds, and a natural categorical crepant
resolution would then be the derived category of the non-singular orbifold.
\end{construction}

\begin{remark}
In the above construction, we are claiming to construct a mirror family
to $X\setminus D$. Thus the question arises as to the dependence of the
construction on the particular choice of compactification $X$ of $X\setminus
D$. 
By \cite[\S10]{J22b}, if $(\widetilde X,\widetilde D)\rightarrow (X,D)$ is
a logarithmically \'etale projective birational morphism (e.g., obtained
by blowing up a sequence of strata of $D$), then the mirror family to
$X\setminus D$ we construct above is a base-change of the mirror family
to $\widetilde X\setminus \widetilde D$ induced by the map
$H_2(\widetilde X)\rightarrow H_2(X)$. More general changes to the
compactification are not currently known to produce similarly compatible
families. However, in those cases where the Keel-Yu construction of
\cite{KY19} apply, the comparison result of \cite{J22a} and the
compactification independence of \cite{KY19} imply our construction
is also independent of the compactification, up to base-change.
\end{remark}

\begin{construction}[Mirrors to degenerations of log Calabi-Yau manifolds]
\label{mirrorconstruction2}
In the relative log Calabi-Yau case, instead $R_I$ can be viewed as the homogeneous
coordinate ring of a family of mirrors to the degeneration $X\rightarrow S$.
Explicitly, in this case we note that $R_I$ carries a natural grading, defined
as follows. Let $\rho\in \Gamma(X,\overline\shM_X)$ be the pull-back of
$1\in \Gamma(S,\overline{\shM}_S)\cong\NN$, so that $\shL_\rho\cong \O_X(-g^*(s_0))$.
We can define the \emph{degree function} on $B$ by $\deg p = \langle p, \rho\rangle$
for $p\in B$.
For $p\in B(\ZZ)$, we can then set $\deg \vartheta_p= \deg p$. Noting that
the divisor $g^*(s_0)$ is numerically equivalent to the divisor $0$ on $X$,
by Proposition~\ref{intersectionnumbers} 
it follows that if $N^A_{p_1p_2r}
\not=0$, we must have
\[
\deg\vartheta_r = \deg \vartheta_{p_1}+\deg\vartheta_{p_2}.
\]
Hence the ring structure is graded, and we can define a mirror formal
family 
as the direct limit of families of schemes $\check X_I=\Proj R_I\rightarrow
\Spec S_I$ to get a flat family of formal schemes
\[
\check\foX\rightarrow \Spf \widehat{\kk[P]}
\]
as before. Of course $\Proj R_I$ is projective over $\Spec R_I^0$, where $R_I^0$
is the degree $0$ part of $R_I$. It is easy to check that $R_I^0$ is precisely
the ring obtained by applying our construction to a general fibre of $g$. 

In the case that $g:X\rightarrow S$ is in fact a degeneration of Calabi-Yau manifolds,
so that $D=g^{-1}(s_0)$ set-theoretically, we can in fact do better. We can define the ring
\[
\widehat R= \bigoplus_{p\in B(\ZZ)} \widehat{\kk[P]} \vartheta_p,
\]
and define the multiplication rule formally using \eqref{structureequation} and 
\eqref{structureconstants}. However, now \eqref{structureequation} will always
be a finite sum. Indeed, if $p$ is a primitive generator of a ray of $\P$,
corresponding to a divisor $D_p\subseteq g^{-1}(s_0)$, $\deg p$ is simply the
coefficient of $D_p$ in $g^*(s_0)$, and $\deg$ is linear on each cone of $\P$.
Thus one sees that $\deg^{-1}(d)$ is a compact subspace of $B$, and 
$B(\ZZ)\cap \deg^{-1}(d)$ is thus a finite set. Thus, since the product rule
respects the grading, \eqref{structureequation} only involves a finite number of
$r$. Further, from \eqref{structureconstants}, we can view the $\alpha_{p_1p_2r}
\in \widehat{\kk[P]}$. Thus $\widehat{R}$ carries a $\widehat{\kk[P]}$-algebra
structure which is associative, as it is associative modulo each $I$ by
Theorem~\ref{mainassociativity2}. We then define the mirror family, in the
case that $\dim_{\RR} B=\dim X$, to be the flat family
\[
\check\shX=\Proj \widehat R \rightarrow\Spec\widehat{\kk[P]}.
\]
Alternatively, one can apply Grothendieck existence to the formal family
$\check\foX\rightarrow \Spf\widehat{\kk[P]}$, but the above construction is
more explicit.
\end{construction}

\begin{construction}[The Rees construction]
\label{reesconstruction}
This construction is motivated by \cite{GHKK}, \S8.5, and gives a cheap
way of compactifying mirrors in the log Calabi-Yau case.
Suppose in the log Calabi-Yau case that $D$ supports a nef divisor $D'$.
Then $\O_X(-D')=\shL_s$ for some $s\in \Gamma(X,\overline{\shM}_X^{\gp})$.
Indeed, as in the proof of Corollary~\ref{cor:intersectionnumbers},
if $D=D_1+\cdots+D_s$, then $\overline{\shM}_X^{\gp}
=\bigoplus_i \ul{\ZZ}_{D_i}$, where $\ul{\ZZ}_{D_i}$ is the constant sheaf
with support $D_i$ and stalk $\ZZ$. Then if $D'=\sum_i a_iD_i$, one may take
$s=(a_i)_i\in\bigoplus_i \Gamma(X,\ul{\ZZ}_{D_i})$.
It then follows from Proposition~\ref{intersectionnumbers}
that if $N^A_{p_1p_2r}\not=0$, we have
\[
0\le A\cdot D' = \langle p_1, s\rangle + \langle p_2, s\rangle -
\langle r,s\rangle, 
\]
and thus 
\[
\langle p_1, s\rangle + \langle p_2,s\rangle \ge \langle r, s\rangle.
\]
This gives rise to a filtered ring structure on $R_I$, and we can then take
the corresponding graded Rees algebra. Explicitly, we consider
\[
\widetilde R_I = \bigoplus_{d \ge 0} \bigoplus_{p\in B(\ZZ)\atop \langle p, s\rangle
\le d} u^d S_I \vartheta_p \subseteq R_I[u],
\]
with $S_I$ as in \eqref{eq:SI}.
The above inequality shows that $\widetilde R_I$ is a graded subring of $R_I[u]$,
with grading given by the power of $u$. Thus we similarly obtain 
a family ${\check X}'_I:=\Proj \widetilde R_I \rightarrow \Spec S_I$, and a corresponding
formal family
\[
\check\foX'\rightarrow \Spf \widehat{\kk[P]}. 
\]
If further, the set of $p\in B(\ZZ)$ with $\langle p,s \rangle\le d$ is always finite,
then a similar argument as in Construction~\ref{mirrorconstruction2} constructs a projective
family
\[
\check \shX' \rightarrow \Spec \widehat{\kk[P]}.
\]
Note that $\check X'_I$ and $\check\foX'$ contain $\check X_I$ and
$\check\foX$ respectively as open subschemes, as can be seen
by the fact that $R_I$ is the degree $0$ part of the localization of $R'_I$ at
$u$. These open subschemes are dense as $u$ is a non-zero-divisor.
\end{construction}

\begin{example}
\label{runningexample2}
Returning to Example~\ref{runningexample1}, let us now compute the rings $R_I$
for all ideals $I$. We consider the first example. From the form
of $R_{\fom}$ given in Example~\ref{runningexample3}, it is sufficient
to deform the equation $\vartheta_{v_1}\vartheta_{v_2}\vartheta_{v_3}=0$
by computing this triple product. We first calculate 
$\vartheta_{v_1}\vartheta_{v_2}$ by enumerating possible $A$ and
$r$ for which $N^A_{v_1v_2r}$ may be non-zero.
We can write a curve class $A=a_1(L-E)+a_2E$ and take
$r=\sum b_i v_i$, where at least one of the $b_i$ is zero. Note we
must have all $a_i$, $b_i$ non-negative. 
If $N^A_{p_1p_2r}\not=0$, then by Corollary~\ref{cor:intersectionnumbers},
applied with $\shL_s=\O_X(-D_i)$, $1\le i\le 3$, we have
\begin{align*}
a_2 = {} & A\cdot D_1 = 
\langle v_1,D_1\rangle + \langle v_2,D_1\rangle -\langle r,D_1\rangle
=1-b_1\\
a_1 = {} & A\cdot D_2 =  
\langle v_1,D_2\rangle + \langle v_2,D_2\rangle -\langle r,D_2\rangle
=1-b_2\\
a_1 = {} & A\cdot D_3 =  
\langle v_1,D_3\rangle + \langle v_2,D_3\rangle -\langle r,D_3\rangle
= -b_3
\end{align*}
By the third equation we necessarily have $a_1=b_3=0$, and hence by 
the second equation $b_2=1$. From the first equation we then have the
options $b_1=0$, $a_2=1$ or $b_1=1$, $a_2=0$. In the latter case, we have
a constant map, so by Lemma~\ref{lem:constant-maps} the contribution is 
monomial, namely the
term $\vartheta_{v_1+v_2}$. In the former case, $A=E$. 
However, any stable map with image $E$ is disjoint from $D_2$, and thus
there is no such curve. Thus we see $\vartheta_{v_1}\vartheta_{v_2}=
\vartheta_{v_1+v_2}$. 

We next calculate $\vartheta_{v_1+v_2}\vartheta_{v_3}$. Again, with
$A$ and $r$ represented as above, we now must have
\begin{align*}
a_2 = {} & A\cdot D_1 = 
\langle v_1+v_2,D_1\rangle + \langle v_3,D_1\rangle -\langle r,D_1\rangle
=1-b_1\\
a_1 = {} & A\cdot D_2 =  
\langle v_1+v_2,D_2\rangle + \langle v_3,D_2\rangle -\langle r,D_2\rangle
=1-b_2\\
a_1 = {} & A\cdot D_3 =  
\langle v_1+v_2,D_3\rangle + \langle v_3,D_3\rangle -\langle r,D_3\rangle
= 1-b_3
\end{align*}
Noting $a_i,b_i\ge 0$, 
we have potentially four possibilities $a_1=0,1$ and $a_2=0,1$,
with the $b_i$ being determined by the above equations.
However, if
$a_1=a_2=0$, we obtain $b_1=b_2=b_3=1$, which is not allowed. So we have
possibilities for $A$ being $E$, $L-E$ or $L$. We immediately
rule out the case that $A=E$ as before. If $A=L-E$,
then $r=v_1$. In this case, one picks a point $z\in D_1^{\circ}$. 
One can then check that the stable map $\PP^1\rightarrow \ul{X}$
identifying $\PP^1$ with $D_1$ can be given in a unique way
the structure of a punctured
map with the desired tangency conditions realizing the constraint
$z$, and that this punctured map is unobstructed. Thus 
$N^{[L-E]}_{v_1+v_2,v_3,v_1}=1$.
On the other hand, if $A=L$, then $r=0$ and after fixing
$z\in X\setminus D$, there is a unique line passing through $z$ and
$D_1\cap D_2$, and similarly $N^{[L]}_{v_1+v_2,v_3,0}=1$. Thus we obtain
\[
\vartheta_{v_1+v_2}\vartheta_{v_3}=t^{[L]}\vartheta_0+
t^{[L-E]}\vartheta_{v_1}.
\]
Putting this together, and using that $\vartheta_0$ is the unit element,
we see that
\begin{equation}
\label{eq:first RI desc}
R_I\cong S_I[\vartheta_{v_1},\vartheta_{v_2},\vartheta_{v_3}]
/(\vartheta_{v_1}\vartheta_{v_2}\vartheta_{v_3}-t^{[L]}
-t^{[L-E]}\vartheta_{v_1}).
\end{equation}

Note this does not tell us what all products of theta functions are.
In this case, it is possible to show the following. First,
if $p_1,p_2$ both lie in the same one of the four cones
$\RR_{\ge 0}v_1+\RR_{\ge 0} v_2$, $\RR_{\ge 0}v_1+\RR_{\ge 0}v_3$,
$\RR_{\ge 0}v_2+\RR_{\ge 0}(v_2+v_3)$ or $\RR_{\ge 0}v_3+\RR_{\ge 0}(v_2+v_3)$,
the multiplication is purely monomial\footnote{This
claim is in fact quite non-trivial, as there may be many potential curve
classes which are represented by actual punctured maps with the desired
tangency conditions. However, all the non-constant maps which may contribute
are in fact obstructed and contribute Gromov-Witten invariant $0$. To
actually show this claim, it is easiest to use the canonical scattering
diagram approach to theta functions of \cite{GHKI}, in which the claim
becomes straightforward. Indeed, \cite{CanScat}, Ex.~3.14 and Thm.~6.1 
show the construction of \cite{GHKI} and the construction
of this paper coincide. \cite{CanScat}, Thm.~6.1 further shows that
the canonical scattering diagram formalism constructs the same
mirror as this paper in all dimensions.},
 i.e.,
\begin{equation}
\label{eq:monomial product}
\vartheta_{p_1}\cdot\vartheta_{p_2}=\vartheta_{p_1+p_2}.
\end{equation}
 
Second, one can also calculate as above that
\begin{equation}
\label{eq:thetav2v3 product}
\vartheta_{v_2}\vartheta_{v_3}=\vartheta_{v_2+v_3}+t^{[L-E]}.
\end{equation} 
The first term comes from the constant map and the second
arises as there is a unique curve in the pencil $|L-E|$ passing through
a given point $z\in X\setminus D$. From all products described so far,
all remaining products can be calculated.

\medskip

Turning to (2) of Examples~\ref{runningexample1},
from the description of $R_{\fom}$ in Example~\ref{runningexample3},
it is enough to calculate the products $\vartheta_{v_1}\cdot\vartheta_{v_2}$
and $\vartheta_{v_{\sigma_1}}\cdot \vartheta_{v_{\sigma_2}}$ in all generality.
 This is quite easily done: using Corollary~\ref{cor:intersectionnumbers}, we 
reduce the calculation to the calculation of two of the 
$N^A_{p_1p_2r}$. First consider the
product $\vartheta_{v_1}\cdot\vartheta_{v_2}$. Take 
$r=a_1v_1+a_2v_2$, with the sum calculated in one of $\sigma_1, \sigma_2$.
Then it follows from Corollary~\ref{cor:intersectionnumbers}
that if $N^A_{v_1v_2r}\not=0$, we have
\begin{align*}
A\cdot D_1 = {} &\langle v_1,D_1\rangle +\langle v_2,D_1\rangle
-\langle r,D_1\rangle= 1-a_1\\
A\cdot D_2 = {} &\langle v_1,D_2\rangle +\langle v_2,D_2\rangle
-\langle r,D_2\rangle= 1-a_2
\end{align*}
But since $D_2$ is a conic, $A\cdot D_2$ is always even, 
and hence there
is no choice but for $a_2=1$, $A=0$, and we are reduced to the constant
map case already considered in Example~\ref{runningexample3}. Thus
$\vartheta_{v_1}\cdot\vartheta_{v_2}=\vartheta_{v_{\sigma_1}}
+\vartheta_{v_{\sigma_2}}$ for any ideal $I$. 

On the other hand, suppose $N^A_{v_{\sigma_1}v_{\sigma_2}r}\not=0$,
with $r$ described as above. Then similarly we obtain
\begin{align*}
A\cdot D_1 = {} &\langle v_{\sigma_1},D_1\rangle +\langle 
v_{\sigma_2},D_1\rangle -\langle r,D_1\rangle= 
1+1-a_1\\
A\cdot D_2 = {} & \langle v_{\sigma_1},D_2\rangle +\langle 
v_{\sigma_2},D_2\rangle -\langle r,D_2\rangle= 
1+1-a_2
\end{align*}
and now we have the possibilities that $a_2=2$, in which case $A=0$ 
and we are
again in the situation analyzed in Example~\ref{runningexample3}, or $a_2=0$, in which
case $A$ is the class of a line and $a_1=1$, so that $r=v_1$. Fixing
any point $z\in Z^{\circ}_{v_1}=D_1\setminus D_2$, it is easy to see there
is exactly one punctured map $f:(C,x_1,x_2,x_{\out})\rightarrow X$ in 
$\scrM(X,\beta,z)$. Here
$C\cong\PP^1$ and $\ul{f}$ induces an isomorphism between $\ul{C}$ and $D_1$,
taking $x_i$ to the point of $D_1\cap D_2$ corresponding to $\sigma_i$ and
$x_{\out}$ to $z$. Thus $N^A_{v_{\sigma_1}v_{\sigma_2}v_1}=1$ is the
only non-trivial contribution to the structure constants, and we obtain, with
$I=(t^k)\subseteq\kk[P]=\kk[t]$,
\[
R_I\cong (\kk[t]/(t^k))[\vartheta_{v_1},\vartheta_{v_2},
\vartheta_{v_{\sigma_1}},\vartheta_{v_{\sigma_2}}]/(\vartheta_{v_1}
\vartheta_{v_2}-
\vartheta_{v_{\sigma_1}}-\vartheta_{v_{\sigma_2}},
\vartheta_{v_{\sigma_1}}\vartheta_{v_{\sigma_2}}-t\vartheta_{v_1}).
\]
This result agrees with the one given by the procedure of
\cite{GHKI}.

Note that in neither of these two cases do we need to work over the 
Artinian ring $R_I$ or
the formal setting as the product rule is actually polynomial.
\end{example}

\begin{remark}
The above example is a special case of a phenomenon observed for surfaces already
in \cite{GHKI}, and which generalises as follows. Suppose that there is an
\emph{ample} divisor
$D'$ supported on $D$ with $D'=\sum_i a_iD_i$ with all $a_i>0$. (Here we are assuming
$D$ is normal crossings). Then the multiplication rule is polynomial, and hence defines
a $\kk[P]$-algebra structure on 
\[
R=\bigoplus_{p\in B(\ZZ)} \kk[P] \vartheta_p.
\]
Indeed, let $s\in \Gamma(X,\overline{\shM}_X^{\gp})$ be defined by taking the value $a_i$ on
the summand $\ZZ_{D_i}$ of $\overline{\shM}_X^{\gp}$, so that $\shL_s\cong \O_X(-D')$ as in the proof of Corollary~\ref{cor:intersectionnumbers}.
Then $s$ defines a function $s:B\rightarrow \RR$ by $s(p)
=\langle s,p\rangle$ for $p\in B$.
As $s$ is positive on the generator of each ray of $\P$,
necessarily $s^{-1}([0,d])$ is compact. Hence $B(\ZZ)\cap s^{-1}([0,d])$ is always
finite, from which it follows as in Construction~\ref{reesconstruction} that
for any $p_1,p_2$ the sum on the right-hand side of \eqref{structureequation}
is finite. Further, each $\alpha_{p_1p_2r}$ is necessarily a finite sum. Indeed,
by Proposition~\ref{intersectionnumbers}, $A\cdot D'$ 
is determined by $p_1,p_2$
and $r$, and because $D'$ is ample, by the fact the Hilbert scheme is finite type,
there are only a finite number of effective curve classes $A$ with the given 
value of $A\cdot D'$. Thus $\alpha_{p_1p_2r}\in \kk[P]$, and so 
\eqref{structureconstants} defines a ring structure on $R$, necessarily associative.
\end{remark}

\begin{example}
\label{ex:yanki}
Here is a simple example of how in the log Calabi-Yau case, the ring
we construct depends on the representative of the numerical equivalence
class taken for $c_1(\Theta_{X/\kk})$ as given in
\eqref{eq:c1 log CY}. Let $\ul{X}=\PP^1$, $D=D_1+D_2+D_3$
where the $D_i$'s are distinct points. In this case $c_1(\Theta_{X/\kk})$
is nef, so we may apply our construction to get a ring.
Here $\Sigma(X)$ consists of a union of three half-lines with a common 
origin. As the only
curve class $A$ with $A\cdot c_1(\Theta_{X/\kk})=0$
is $A=0$, in fact the only contributions to the product are
from constant maps, in which case Lemma~\ref{lem:constant-maps} applies.
With $S_I=\kk[t]/(t^k)$ for some $k$, we thus immediately see that
\[
R_I=S_I[\vartheta_1,\vartheta_2,\vartheta_3]/(\vartheta_1\vartheta_2,
\vartheta_1\vartheta_3,\vartheta_2\vartheta_3),
\]
and $\Spec R_I$ is a union of the three
coordinate lines in $\AA^3\times\Spec S_I$.

On the other hand, if we write $c_1(\Theta_{X/\kk})\equiv_{\QQ} D_3$ and pass to
the Kontsevich-Soibelman skeleton, we have $B$ a union of two half-lines
with a common origin. Hence we throw out the theta function
$\vartheta_3$, and obtain
\[
R_I=S_I[\vartheta_1,\vartheta_2]/(\vartheta_1\vartheta_2).
\]
Thus $\Spec R_I$ is a union of the two coordinates lines in
$\AA^2\times\Spec S_I$.

Both of these can be viewed as mirrors to $(\PP^1,D)$, in a certain sense.
\cite{AAEKO} suggested as mirror to $(\PP^1,D)$ a Landau-Ginzburg model
$W:\AA^3\rightarrow \AA^1$ given by $W=x_1x_2x_3$. The critical locus
of this potential function is the union of three affine lines, and it
is expected that this critical locus can also be viewed as a mirror,
see \cite{GKR}, as well as \cite{Ruddat}.

For the choice of mirror which is a union of two affine lines,
the partially wrapped Fukaya category on the
thrice punctured sphere is isomorphic to a quotient of the derived
category of $\Spec \kk[x_1,x_2]/(x_1x_2)$, see \cite{LP18}, \S 1.

From both points of view, we are not actually 
constructing the mirror, but rather
a scheme related to the mirror.

As also explained to us by Y.\ Lekili, a choice of nowhere vanishing
top dimension holomorphic form on $X\setminus D$ induces a grading
on symplectic cohomology. For us, this choice is encoded in the representative
for $c_1(\Theta_{X/\kk})$. Without this choice, $R_I$ as
given in Theorem \ref{mainassociativity1} is a piece of a 
non-graded symplectic cohomology ring, but with the choice, 
symplectic cohomology becomes a graded ring and the $R_I$
constructed in Theorem \ref{mainassociativity2} should be the degree zero part
of this graded symplectic cohomology.
\end{example}

We can turn this into a more precise conjectural statement about 
the relationship between symplectic cohomology and our ring $R_I$ in
a more restrictive log Calabi-Yau context.
Evidence for the following conjecture comes from \cite{GP21}
and \cite{Pom}:

\begin{conjecture}
\label{conj:SH}
Suppose given $D\subseteq \ul{X}$ a simple normal crossings divisor. Suppose
that $D$ supports an effective
ample divisor on $X$ and that there exists a rational
section $\Omega$ of the canonical bundle $\omega_X$ with $\Omega|_{X\setminus
D}$ nowhere vanishing and $\Omega$ having at worst simple poles along
$D$. This data gives an expression in $\Div(X)$
\[
K_X+D \sim (\Omega) + D = \sum_i a_iD_i,
\]
where $(\Omega)$ is the divisor of zeros and poles of $\Omega$ and
$a_i\ge 0$. Using this expression as in \eqref{eq:c1 log CY},
we obtain the ring $R_I$ from 
Theorem~\ref{mainassociativity2}. The form $\Omega$ also induces
a grading on the symplectic cohomology ring $\mathrm{SH}^*(X\setminus D)$.
Then $R_I$ and $\mathrm{SH}^0(X\setminus D)$ are
isomorphic after a suitable base-change.
\end{conjecture}

\begin{example}
\label{ex:relative}
Here is a simple example in the relative setting. We start with the family of quadrics in $\PP^3$ from \cite{Invitation},\S3.2:
\[
\PP^3_{x,y,z,w}\times\AA^1_s\supseteq \overline X:=V\big(xy-sz(z+w)\big)\lra \AA^1.
\]
The fibre over $s_0=0$ is the union of two planes $\overline F_1= V(s,x)$,
$\overline F_2 = V(s,y)$, as sketched in Figure~\ref{fig:relative_example}. We
consider $\ol X$ as a family of log Calabi-Yau pairs with boundary
divisor
\[
\left(V(zw)\cap \overline X\right)\cup \overline X_{s_0} = 
\overline D_1+\overline D_2+\overline D_3+\overline F_1+\overline F_2,
\]
where $\overline D_1 = V(x,z)$, $\overline D_2 = V(y,z)$,
$\overline D_3 = V(w, xy-sz^2)$ all surject onto $S$.

Now $\overline X$ has two singular points, each ordinary double points, at
\[
x_1:=\big((0:0:1:-1),0\big), \quad x_2:=\big((0:0:0:1),0\big),
\]
and $\ol X\to S$ fails to be log smooth at $x_1$. We obtain a small resolution
of both double points by blowing up $\ol F_2$, with strict transform a normal
crossings pair $(X,D)$ and $X\to S$ log smooth. We illustrate the geometry of
$D$ on the left-hand side of Figure~\ref{fig:relative_example}. The morphism
$F_2\rightarrow\overline F_2$ is the blow-up of $\overline F_2$ at the two
points $x_1,x_2$, with exceptional curves $E_1,E_2$. Further, $D_1$ is the
blow-up of $\overline D_1$ at $x_2$, and $E_2=D_1\cap F_2$. The remaining
boundary divisors are unchanged.

%\begin{figure}
%\input{FigureRelative.pspdftex}
%\caption{The geometry of Example \ref{ex:relative}.}
%\label{fig:relative_example}
%\end{figure}

\begin{figure}
\input{FigureRelative_alt.pspdftex}
\caption{The geometry of Example \ref{ex:relative}.}
\label{fig:relative_example}
\end{figure}

Using the method of discrete Legendre transform, wall structures and generalized
theta functions from our previous work \cite{PartI,Annals,GHS}, one can now
easily deduce the mirror algebra of theta functions. For readers familiar with
this theory we provide in Figure~\ref{Figure10} the discrete Legendre transform
of \cite[Fig.\,3.1]{Invitation} describing $B$ for the mirror geometry as a
polyhedral affine manifold with a simple singular point on the bounded edge.

\begin{figure}
\input{Figure10.pspdftex}
\caption{The polyhedral affine manifold for the mirror geometry of \cite[\S3.1]{Invitation}.}
\label{Figure10}
\end{figure}

A finer, but compatible, picture is obtained from the intrinsic mirror ring
$R_I$ in the present paper, which we now compute for all $I$ by similar
methods as in Example~\ref{runningexample2}. We again take $H_2(X)=\Pic(X)^*$
and $P\subseteq H_2(X)$ a suitable submonoid as usual. Let $m$ be the
class of the line $X_s\cap D_2$, while $E_2$ and $L$ denote
the classes of the exceptional curve over $x_2$ and a line in $F_1\cong
\PP^2$, respectively. Note that $L+E_2$ is the class of the line 
$X_s\cap D_1$ and $L+E_1$ is the class $m$ of $X_s\cap D_2$.
Writing $v_i:=D_i^*$, $w_i:=F_i^*$,
$u_i:=\vartheta_{v_i}$, $W_i:= \vartheta_{w_i}$, we obtain
\[
R_I\cong {A_I[u_1,u_2,u_3][W_1,W_2]\over
\big(
u_1u_2u_3-t^{L+E_2}u_1-t^mu_2,
u_2W_1-t^{E_2}u_1W_2, 
u_1u_3W_2-t^LW_1-t^mW_2\big)
}.
\]
Here $u_1,u_2,u_3$ have degree $0$ and $W_1,W_2$ degree $1$, and the equations
are homogeneous with respect to this grading. We explain the origin of each of
the equations. 

For the first one,
we note that any product involving only degree
$0$ theta functions only involves 
contact orders with $D_1,D_2,D_3$, and as such,
we may choose the point constraint $z$ to lie in some fixed $X_s$, $s\not=s_0$.
The products are then precisely the same as the products arising
from the pair $(X_s, X_s\cap D)$. Here, $X_s\cong\PP^1\times\PP^1$,
and $X_s\cap D$ is a union of curves of bidgree $(1,0)$, $(0,1)$ and
$(1,1)$. Arguments as in Example~\ref{runningexample2} then yield
this first equation.

The second defining equation arises from two products. The first is 
$u_1W_2=\vartheta_{v_1+w_2}$, arising from 
Proposition~\ref{prop:degree 0 mult}. The second is
$u_2W_1=t^{E_2}\vartheta_{v_1+w_2}$. This arises
as the exceptional curve $E_2$ can be enhanced to a punctured log curve
which is responsible for this product. In particular, $E_2\cdot D_2=
E_2\cdot F_1=1$ and $E_2\cdot D_1=E_2\cdot F_2=-1$. This shows the existence
of a punctured log enhancement on $E_2$ with contact orders
$v_2,w_1$ and $-(v_1+w_2)$ is compatible with
Corollary~\ref{cor:intersectionnumbers}.

Finally, the last equation arises via
\[
(u_1u_3)W_2
=(\vartheta_{v_1+v_3}+t^m)W_2
=t^LW_1+t^mW_2.
\]
The first equality can be viewed as a product arising from $(X_s,
D\cap X_s)$ as before, while the second product arises from a line
on $F_1$ which passes through the point $D_1\cap D_3\cap F_1$ and
a chosen general point $z\in F_1$. This line then also meets $F_2$
at one point, as depicted in Figure~\ref{fig:relative_example}.

The interested reader may check that the mirror algebra of generalized theta
functions mentioned above, for the scattering diagram from \cite{Annals}, is
the algebra over $\kk[t]/(t^k)$ obtained from the localization of $R_I$ at
$E_1$ by the base change
\begin{eqnarray*}
t^m\mapsto t^{\kappa_1},\quad t^{E_2}\mapsto t^{\kappa_2},\quad
t^L\mapsto at^{\kappa_1},\quad t^{E_1}\mapsto a^{-1}.
\end{eqnarray*}
Here $\kappa_1$ is the kink of the chosen MPL function along the vertical bounded edge,
$\kappa_2$ the kink along the vertical unbounded edge emanating from the
vertex labelled $W_2$, and $a\in\kk^\times$ a parameter determining the gluing
data. The order $k$ can be chosen arbitrarily large by considering appropriate
$I$. More universally one can also work over a base ring generated by
$t_1,t_2,a^{\pm1}$ where $t_i$ pulls back to $t^{\kappa_i}$, see
\cite[\S A2]{GHS}.
\end{example}

\begin{construction}[The torus action]
Suppose we are in the situation of Example \ref{snccaseexample2}.
Then there is a natural torus action on $\Spec R_I$ 
which has proved to be very useful, see e.g., \cite{GHKI}, \S5
or \cite{GHS}, \S4.4. Recall that $\Div_D(X)^*$ is the dual
of the free abelian group generated by the irreducible components of
$D$. One can easily describe an action of the algebraic torus 
$\Spec \kk[\Div_D(X)^*]$
on $\Spec R_I$, as follows.
Indeed, it is sufficient to give a $\Div_D(X)^*$-grading
on $R_I$. We do this by defining the $\Div_D(X)^*$-degree of
$\vartheta_p$ to be $p\in\Div_D(X)^*$, and define the 
$\Div_D(X)^*$-degree of $t^A\in \kk[P]$ to be the functional given by
$D_i\mapsto D_i\cdot A$.
That our product rule
on $R_I$ then respects the grading follows again from 
Corollary~\ref{cor:intersectionnumbers}. Indeed, if the coefficient
of $t^A\vartheta_r$ is non-zero in the product $\vartheta_p\cdot\vartheta_q$,
we have by that corollary that
\[
\langle r, D_i\rangle + D_i \cdot A = \langle p, D_i\rangle
+\langle q, D_i\rangle,
\]
showing degrees are additive.

The same holds in the log Calabi-Yau case. 
\end{construction}

\section{Preliminaries and review}
\label{sec:review}

\subsection{Generalities on log geometry}
We assume the reader is familiar with logarithmic geometry as is necessary
to understand logarithmic Gromov-Witten invariants. However there are
several more specific points we will be needing, and we review these
quickly and give references. All log structures are fine and saturated (fs),
except as explicitly mentioned, such as in the case of punctured log
structures, which are only fine.

\subsubsection{Tropicalization}
\label{sec:tropicalization}
We assume the reader is familiar with the notation of \cite{ACGS17},
\S 2, and in particular use freely the notion of generalized cone complexes
(\S 2.1.2) 
and tropicalization of logarithmic schemes and stacks (\S 2.1.4). 
For more subtle issues about tropicalization, the reader may consult
\cite{ACGS18},~App.~C. However, in general such subtleties will not arise
in this paper. We write $\Sigma$ for the
tropicalization functor from the category of fine or fs log 
algebraic
stacks to the category of generalized cone complexes. Recall
$\mathbf{Cones}$ is the category whose objects are pairs
$(\sigma, N)$ where $N$ is a lattice and $\sigma\subseteq N_{\RR}$
is a rational polyhedral cone with $\sigma^{\gp}=N_{\RR}$, and morphisms
are maps of cones induced by linear maps of the lattices. A generalized
cone complex is (an equivalence class of)
a diagram in $\mathbf{Cones}$, all of whose morphisms
are face morphisms.
The cones of the tropicalization of a log stack $X$ 
are indexed by geometric points $\bar x\in |X|$
with the associated cone being
\[
\sigma_{\bar x}:=\Hom(\overline{\shM}_{X,\bar x},\RR_{\ge 0}),
\]
and lattice
\[
N_{\sigma_{\bar x}}:=\Hom(\overline{\shM}_{X,\bar x},\ZZ).
\]
The face maps of cones in $\Sigma(X)$ are induced by generization maps. 
See \cite{ACGS17}, \S 2.1.4 or \cite{ACGS18}, App.~C for more details.

\subsubsection{Artin fans}
\label{subsubsec:artin fan}
Given a log scheme
$X$, we have its associated \emph{Artin fan} $\shX$, as constructed in 
\cite{ACMW17}, Prop.~3.1.1, (where it is denoted $\shA_X$) 
see \cite{ACGS17}, \S2.2 for a summary. 
There is a factorization
\[
X\rightarrow \shX \rightarrow \Log_{\kk}
\]
of the tautological morphism $X\rightarrow\Log_{\kk}$, where $\Log_{\kk}$
is Olsson's stack
parameterizing all fine saturated log schemes over $\Spec\kk$.
As Olsson uses this notation for the stack of all fine log schemes
over $\kk$, in fact our usage agrees with Olsson's
stack ${\mathcal{T}\!or}_\kk$ of \cite{OlssonENS}. The morphism
$\shX\rightarrow\Log_{\kk}$ is \'etale and representable by algebraic
spaces, and is the initial such factorization. In particular, $X\rightarrow
\shX$ is strict, and one should view $\shX$ as capturing the combinatorics
of the log structure on $X$.

\subsubsection{Idealized log structures}

We recall the notion of \emph{idealized log schemes} from \cite{Ogus},
III,Def.~1.3.1. An idealized log scheme is a log scheme
equipped with a sheaf of monoid ideals $\shK_X\subseteq\shM_X$
with the property that $\alpha_X(\shK_X)=0$. This then gives a category
of idealized log schemes, with a morphism $f:X\rightarrow Y$
satisfying $f^{\flat}(\shK_Y)\subseteq\shK_X$. 

Further, given
a morphism $f:X\rightarrow Y$ of ordinary log schemes, and given
an idealized structure $\shK_Y\subseteq \shM_Y$, we obtain 
an idealized structure on $X$, denoted by $f^{\bullet}\shK_Y$, to
be the sheaf of monoid ideals generated by $f^{\flat}(\shK_Y)$.

The trivial idealized structure on $X$ is given by $\shK_X=\emptyset$,
the empty ideal. 

A standard example of an idealized log scheme is $A_{P,K}:=\Spec\kk[P]/K$,
where $P$ is a fine saturated monoid and $K\subseteq P$ is a monoid
ideal of $P$. Then $\shM_{A_{P,K}}$ is the standard toric log structure
induced by the pre-log structure $P\rightarrow \O_{A_{P,K}}$,
and $\shK$ is the ideal generated by $K$ under the
image of the canonical morphism $P\rightarrow \shM_{A_{P,K}}$.

Here we only consider \emph{coherent} idealized log structures in
fine or fs log schemes. This means that locally on $X$, the ideal
$\shK_X$ is generated by a finite set of sections: 
see \cite{Ogus}, II,Prop.~2.6.1.

\subsubsection{Fine saturated fibre products}

Fibre products in the category of fine or fs log schemes or stacks cause
a great deal of difficulties. Given (fine, fs) log schemes $X,Y,S$
with morphisms $X\rightarrow S$, $Y\rightarrow S$, we distinguish between
$X\times_S Y$, $X\times^{\fine}_S Y$ and $X\times^{\fs}_S Y$ as the
fibre products in the category of log schemes, fine log schemes and fs
log schemes. In general, $\ul{X\times_S Y}=\ul{X}\times_{\ul{S}}\ul{Y}$,
while the same does not hold for the other fibre products. However,
there is an obvious functorial morphism 
$X\times^{\fine}_S Y\rightarrow X\times_SY$ if $X,Y$ and $S$ are fine,
and this morphism is given by \emph{integralization} of the log structure.
This morphism is always a closed
embedding. There is also an obvious functorial morphism
$X\times^{\fs}_S Y\rightarrow X\times^{\fine}_SY$ if $X,Y$ and $S$ are
fs, and this morphism is given by \emph{saturation}, which is always
finite and surjective. See \cite{Ogus}, III,Prop.~2.1.5. 

If one of the morphisms $X\rightarrow S$ or $Y\rightarrow S$ is 
strict, then all three fibre products agree, and being strict is stable
under base-change.

\begin{remark}
\label{rem:sat finite representable}
In the case that $X,Y$ and $S$ are stacks with fine or fs log structures,
the construction of integralization and saturation in
\cite{Ogus}, III,Prop.~2.1.5, in particular implies that 
$\ul{X\times^{\fine}_SY}\rightarrow \ul{X}
\times_{\ul{S}}\ul{Y}$
or $\ul{X\times^{\fs}_SY}\rightarrow \ul{X}\times_{\ul{S}}\ul{Y}$ are finite and
representable morphisms.
\end{remark}

To understand the behaviour of fibre products at a tropical level,
we recall \cite{ACGS17}, Prop.~5.2:

\begin{proposition}
\label{tropicalproduct}
Let $X,Y$ and $S$ be fs log schemes, with morphisms $f_1:X\rightarrow S$,
$f_2:Y\rightarrow S$. Let $Z=X\times^{\fs}_S Y$,
$p_1,p_2$ the projections. Suppose given geometric points 
$\bar z\in Z$, $\bar x= p_1(\bar z)$,
$\bar y=p_2(\bar z)$, and $\bar s=f_1(p_1(\bar z))=f_2(p_2(\bar z))$. Then
\[
\Hom(\overline{\shM}_{Z,\bar z},\NN)=
\Hom(\overline{\shM}_{X,\bar x},\NN)\times_{\Hom(\overline{\shM}_{S,\bar s},\NN)}
\Hom(\overline{\shM}_{Y,\bar y},\NN)
\]
and
\[
\sigma_{\bar z}=\sigma_{\bar x}\times_{\sigma_{\bar s}} \sigma_{\bar y}.
\]
\end{proposition}

\subsubsection{Transversality and flatness}
\label{sec:transversality}

The goal in this subsection is to give an approach for working around
a fundamental problem in log geometry: log smooth or log \'etale
morphisms need not be flat. Indeed, they are only guaranteed to be 
flat if they are also integral
morphisms. For the definition of integral morphism, see \cite{Ogus}, 
III,Def.~2.5.1. More generally, log flat morphisms, defined
in \cite{Ogus},~IV,Def.~4.1.1, are flat if integral:

\begin{proposition}
\label{prop: flatness sorites}
Let $f:X\rightarrow Y$ be a morphism of fs log schemes. Then
\begin{enumerate}
\item If $f$ is log smooth, then $f$ is log flat.
\item If $X$ and $Y$ are locally Noetherian fine log schemes and $\ul{f}$
is locally of finite presentation, then if $f$ is log flat and integral,
$\ul{f}$ is flat.
\end{enumerate}
\end{proposition}

\begin{proof} (1) and (2) are \cite{Ogus}, IV,Thms.~ 4.1.2, 1, and 
4.3.5, 1 respectively. 
\end{proof}

Integrality of morphisms can often be tested combinatorially:

\begin{proposition}
\label{prop:integral morphism}
Let $\theta:\NN^r=P\rightarrow Q$ be a morphism between saturated toric monoids,
dual to a morphism of cones 
\[
\theta^t:\sigma_Q=\Hom(Q,\RR_{\ge 0}) \rightarrow \sigma_P=\Hom(\NN^r,\RR_{\ge 0}).
\]
Then $\theta$ is integral if and only if the
restriction of $\theta^t$ to every face of $\sigma_Q$ surjects onto a face
of $\sigma_P$.
\end{proposition}

\begin{proof}
By \cite{Ogus}, I,Prop.~4.7.5, $\theta$ being integral is equivalent to
$\theta$ being $\QQ$-integral. This is where the assumption that
$P=\NN^r$ is used. The $\QQ$-integrality in turn is equivalent to $\theta$ being
locally exact by \cite{Ogus}, I,Prop.~4.7.7, which means that for
each face $G\subset Q$, the induced morphism of localized monoids
$P_{\theta^{-1}(G)}\rightarrow Q_G$ is exact. By \cite{Ogus}, I,Prop.~4.3.7,1, 
this exactness is equivalent to surjectivity of the map of cones
$\Hom(Q_G, \RR_{\ge 0})\rightarrow \Hom(P_{\theta^{-1}(G)},\RR_{\ge 0})$.
However, these cones can be described as faces of $\sigma_Q$ and
$\sigma_P$ respectively, namely $G^{\perp}\cap \sigma_Q$ and 
$(\theta^{-1}(G))^\perp\cap \sigma_P$. Note that the latter face is the smallest
face of $\sigma_P$ containing $\theta^t(G^{\perp}\cap\sigma_Q)$. Hence
the hypothesis on $\theta^t$ is equivalent to the desired exactness, 
completing the proof.
\end{proof}

\begin{remark}
Alternatively, one can prove the above proposition by observing that
the combinatorial condition on $\theta^t$ implies that all fibres of 
the induced morphism of toric varieties $\Spec\kk[Q]\rightarrow \Spec\kk[P]$
are equi-dimensional of dimension $\rk Q^{\gp}-\rk P^{\gp}$. Thus
by ``miracle flatness,'' \cite{Mats}, Thm.~23.1, this induced
morphism is flat, which is equivalent to $\theta$ being integral.
Since miracle flatness requires a non-singular target space, this
gives another explanation for why we must take $P$ to be the free monoid.
\end{remark}

\begin{definition}
\label{def: transverse}
 Let $f:X\rightarrow Y$ be a morphism between fs log stacks. We say a 
morphism $g:W\rightarrow Y$ is \emph{transverse to $f$} if
the first projection $W\times_Y^{\fs} X\rightarrow W$
is an integral morphism.
\end{definition}

The idea is that if $f$ is log smooth, $\ul{f}$ may not be flat,
but it may become flat after a judicious choice of base change by some $g$.
In certain cases, there is a simple combinatorial test for transversality,
which is a consequence of the following lemma.

\begin{lemma}
\label{lem: product of cones}
Suppose given rational polyhedral cones $\sigma_1,\sigma_2,\tau$ with
linear maps $f_i:\sigma_i\rightarrow \tau$, and consider the cone
$\sigma:=\sigma_1\times_{\tau}\sigma_2$ with projections $p_i:\sigma\rightarrow
\sigma_i$. Suppose that for each face
$F_2$ of $\sigma_2$, $f_1^{-1}(f_2(F_2))$ is a face of $\sigma_1$.
Then $p_1$ gives a surjection of each face of $\sigma$ onto a face of
$\sigma_1$.
\end{lemma}

\begin{proof}
Note that $\sigma_1\times_{\tau}\sigma_2$ is the intersection of
$\sigma_1\times\sigma_2\subset \sigma_1^{\gp}\times\sigma_2^{\gp}$
with the linear space $\ker f$, where $f:\sigma_1^{\gp}\times\sigma_2^{\gp}
\rightarrow \tau^{\gp}$ is given by $f(a,b)=f_1(a)-f_2(b)$, confusing
$f_i$ with its linear extension. Thus, as each face
of $\sigma_1\times\sigma_2$ is of the form $F_1\times F_2$ for $F_i$ a
face of $\sigma_i$, each face of $\sigma$ is of the form
$F_1\times_{\tau} F_2$. We show that $p_1$ maps this face surjectively
onto the face $F:=f_1^{-1}(f_2(F_2))\cap F_1$ of $\sigma_1$.

Indeed, if $(a_1,a_2)\in F_1\times_{\tau} F_2$, then $f_1(a_1)=f_2(a_2)$
so $a_1\in f_1^{-1}(f_2(a_2))\cap F_1\subseteq F$. Thus $p_1(a_1,a_2)\in F$ and
$p_1$ maps $F_1\times_{\tau}F_2$ into $F$. 
Further, given $a_1\in F$, there exists an $a_2\in F_2$
such that $f_1(a_1)=f_2(a_2)$. Then $(a_1,a_2)\in F_1\times_{\tau}F_2$,
showing $p_1$ maps this face surjectively onto $F$.
\end{proof}

\begin{definition}
An fs log stack $W$ is \emph{free} if all stalks of $\overline{\shM}_W$
are free monoids, i.e., isomorphic to $\NN^r$ for some $r$.
\end{definition}

The crucial observation for us is:

\begin{theorem}
\label{theorem: transversality}
Let $f:X\rightarrow Y$, $g:W\rightarrow Y$ be log morphisms between
fs log stacks. Suppose further that $W$ is free. Then $g$ is transverse
to $f$ if for any geometric points $\bar x\in X$, $\bar y\in Y$ and
$\bar w\in W$ with $f(\bar x)=\bar y = g(\bar w)$, and for any face
$F$ of $\sigma_{\bar x}$, we have that $\Sigma(g)^{-1}(\Sigma(f)(F))
\subset \sigma_{\bar w}$ is a face of $\sigma_{\bar w}$.
\end{theorem}

\begin{proof}
If $\bar z$ is a geometric point of $W\times_Y^{\fs} X$ with images
$\bar w$ and $\bar x$ under the two projections,
then $\sigma_{\bar z}=\sigma_{\bar w}\times_{\sigma_{\bar y}}
\sigma_{\bar x}$ by Proposition~\ref{tropicalproduct}. Integrality
of the projection $W\times^{\fs}_Y X\rightarrow W$ now follows from
Proposition~\ref{prop:integral morphism} and Lemma~\ref{lem: product of cones}.
\end{proof}

\begin{remark}
\label{rem:transverse}
This notion of transversality is a key tool in this paper. Here is a simple
example of how this notion is useful. Given a finite \'etale morphism
of ordinary schemes $\ul{f}:\ul{X}\rightarrow\ul{Y}$, with $\ul{Y}$ connected,
all fibres are of finite length, and the length of any fibre is an invariant 
of the morphism, i.e., the degree. If one instead considers a log \'etale
morphism $f:X\rightarrow Y$, finiteness is not a particularly useful property,
as interesting log \'etale morphisms are seldom finite. Instead, it is 
natural to replace the finiteness condition with properness.

In this case, scheme-theoretic fibres are not necessarily of finite length.
However, if $g:W\rightarrow Y$ is a morphism from a log point $W$ which
is transverse to $f$, then $W\times_Y X$ is of finite length over $W$.
One can then prove this length is independent of the particular choice of
$W$, hence allowing the definition of a log degree. We do not need this
particular result in this paper, so we do not give a proof. However, the
basic philosophy is used repeatedly: we use a choice of transverse morphism
to control the dimension of a log fibre of a given morphism. The above
results show that transverse morphisms can be chosen via a tropical analysis.
\end{remark}

\subsection{Stable punctured maps}
\label{sec:punctured log curves}

Punctured maps, developed in \cite{ACGS18}, 
are a crucial generalization of stable log
maps of \cite{JAMS}, \cite{AC14}, \cite{Chen14}.
The main point is they allow marked points with negative
contact order. 

We will review the most salient points.

\subsubsection{Punctures and punctured maps}
\label{subsubsec:punctured maps}
The general definition of a puncturing is given in \cite{ACGS18},
Def.~2.1. For the purposes of this paper, we only consider
the situation of an fs log scheme (or stack) $Y=(\ul{Y},\shM_Y)$
where $\shM_Y=\shM\oplus_{\O_Y^{\times}} \shP$ for some log structures
$\shM$ and $\shP$ on $Y$. We further assume that $\shP$ is a DF(1)
log structure, i.e., there is a surjective homomorphism
$\ul{\NN}\rightarrow \overline{\shP}$ from the constant sheaf $\ul{\NN}$
with stalk $\NN$. Then a \emph{puncturing}
 $Y^{\circ}$ of $Y$ is a fine logarithmic structure 
\[
\shM_{Y^{\circ}}\subset \shM\oplus_{\O_Y^{\times}} \shP^{\gp}
\]
containing $\shM_Y$ such that 
\begin{enumerate}
\item The inclusion $\shM_Y\rightarrow \shM_{Y^{\circ}}$ is a morphism of
fine logarithmic structures on $\ul{Y}$.
\item 
For any geometric point $\bar y$ of $\ul{Y}$ let $s_{\bar y}\in 
\shM_{Y^\circ,\bar y}$ be
such that $s_{\bar y}\not \in \shM_{\bar{y}}\oplus_{\O^{\times}}\shP_{\bar{y}}$.
Representing $s_{\bar y}=(m_{\bar y},p_{\bar y})\in \shM_{\bar y}
\oplus_{\O^{\times}_Y}
\shP_{\bar y}^{\gp}$,
we have $\alpha_{\shM_{Y^\circ}}(s_{\bar y})=\alpha_{\shM}(m_{\bar y})=0$ in $\O_{Y,\bar y}$.
\end{enumerate}
Note that we do not require the punctured log structure to be saturated,
but only fine. For some concrete relevant examples, see
\cite{ACGS18},~Ex.~2.12,2.19.

Puncturings suffer the problem that there may be many different choices
of the puncturing log structure $\shM_{Y^{\circ}}$. However, this 
non-uniqueness can be removed in the presence of a morphism $f:Y^{\circ}
\rightarrow X$, by demanding that $\shM_{Y^{\circ}}$ be as small as possible.
Formally: 

\begin{definition}
\label{def:prestable}
A morphism $f:Y^{\circ}\rightarrow X$ from a punctured log scheme
$Y^{\circ}$  is \emph{pre-stable} if $\shM_{Y^{\circ}}$
is generated as a sheaf of fine monoids by $\shM_Y$ and the image of
$f^*\shM_X$ under $f^{\flat}$.
\end{definition}

The only cases that we need to consider in this paper is when we are given
a log smooth family of curves $\pi:C\rightarrow W$ with marked
points $x_1,\ldots,x_n:\ul{W}\rightarrow \ul{C}$. We use the convention
as in \cite{ACGS18} that the log structure $\shM_C$ is written as
$\shM\oplus_{\O_C^{\times}}
\shP$ where $\shP$ is the divisorial log structure associated to the
marked points $x_1,\ldots,x_n$, and $\shM$ is the log structure on $C$
for which the points $x_1,\ldots,x_n$ are not marked.

A punctured map is then a morphism $f:C^{\circ}\rightarrow X$
for some choice of puncturing $C^{\circ}$. 
Throughout the entire paper, all such morphisms are
assumed to be pre-stable in the sense just defined,
and this hypothesis will generally not be mentioned.

Punctured maps allow negative contact orders as follows. Fix
$\bar w\in W$ and one of the marked points $x=x_i\in C_{\bar w}$. Then
$f$ induces a composed morphism, denoted $u_x$,
\begin{equation}
\label{eq:ux def}
P_{x}:=\overline{\shM}_{X,f(x)}\mapright{\bar f^{\flat}}
\overline{\shM}_{C^{\circ},x}\subset
\overline{\shM}_{W,\bar w}\oplus \overline\shP_x^{\gp}=\overline{\shM}_{W,\bar w}
\oplus \ZZ\mapright{\pr_2} \ZZ
\end{equation}
where $\pr_2$ denotes the second projection. Here $u_x$ is called the
\emph{contact order} of the point $x$, and we have $u_x\in P_x^*=
\Hom(P_x,\ZZ)$, whereas in more traditional logarithmic Gromov-Witten
theory the contact order $u_x$ lies in $P_x^{\vee}=\Hom(P_x,\NN)$.

We use the convention of \cite{ACGS18} that if $u_x\in P_x^{\vee}$,
we call $x$ a \emph{marked point}. However, if we call a point a punctured
point, we do not exclude the possibility that it is a marked point.
The assumption of pre-stability guarantees that marked points
behave exactly as marked points in the original logarithmic
Gromov-Witten theory, and only punctured points with $u_x\not\in P_x^{\vee}$
present new phenomena.

\subsubsection{Contact orders for punctured points}
\label{sec:contact orders}
The issue of specifying contact orders at punctured points is rather
subtle, and for an exhaustive discussion, we refer the reader to \cite{ACGS18},
\S2.4 and Appendix A. However, for this paper we generally only need
a relatively simple form of this discussion. We define the
set of \emph{global contact orders}
\[
\foC(X):=\operatorname{colim}_{\sigma\in\Sigma(X)} N_{\sigma},
\]
where $N_{\sigma}$, as before, is the integral lattice
associated to $\sigma$.  Here the colimit is defined in
the category of sets rather than abelian groups.\footnote{This
is in fact a special case of the general setup in \cite{ACGS18},\S2.4,
where this set would be written as $\foC_0(X)$.}

Thus, given a punctured map $f:C^\circ/W \rightarrow X$, each
punctured section $x:\ul{W}\rightarrow \ul{C}$ and each geometric point 
$\bar w\in W$
induces via \eqref{eq:ux def} $u_{x(\bar w)}\in N_{\sigma_{f(x(\bar w))}}$
and hence an element of $\foC(X)$, which we write
as $\bar u_{x(\bar w)}$. Further, $\bar u_{x(\bar w)}$ is now locally
constant on $W$, see \cite{ACGS18}, Lem.~2.41. If $\bar u_{x(\bar w)}$ 
is constant
on $W$ (e.g., if $W$ is connected), then we call $\bar u_{x(\bar w)}$ 
the \emph{global contact order} of the puncture $x$, writing it
as $\bar u_x$.

Note for any $\sigma\in\Sigma(X)$, there is a canonical map
\[
\iota_\sigma:N_{\sigma}\rightarrow \foC(X).
\]
We say that $\bar u\in \foC(X)$ has \emph{finite monodromy}
if for any $\sigma\in\Sigma(X)$, $\iota^{-1}(\bar u)\subseteq
N_{\sigma}$ is finite. We say $\bar u$ is \emph{monodromy-free}
if $\iota_{\sigma}^{-1}(\bar u)$ has at most one element for each
$\sigma\in\Sigma(X)$.

As we review shortly, when we consider
moduli spaces of punctured maps, we will need to specify families of 
contact orders for each punctured point. Infinite monodromy
may well prevent these moduli spaces being finite type. Hence,
as in \cite{ACGS18}, we must for now restrict to the case that
$\overline{\shM}^{\gp}_X\otimes\QQ$ is generated by 
its global sections as in Basic
Setup~\ref{setup1}. We do not believe that this restriction is essential
to the theory, but it does not seem to us to be worth the extra pain
at this point to remove this hypothesis. We remark that with
this hypothesis, all contact orders are monodromy free, see
\cite{ACGS18}, Rem.~3.13. Recent work of Johnston \cite{J22b}
allows one to remove this generatedness hypothesis in the theory of
punctured invariants. However, we will not work in that degree of
generality here, as this hypothesis also plays a role in other ways.

\begin{remarks}
\label{rem:contact order remarks}
(1) In \cite{ACGS18},~App.~A introduced the notion of
a family of contact orders for $X$. This is a slightly complicated notion
which is rarely needed in the degree of generality stated in that Appendix.
Instead, we will summarize here the results necessary in the simple case
that the contact order is given by
$u \in N_{\sigma}$ with $u \in \Int(\sigma)$ or $-u \in\Int(\sigma)$, 
representing
a global contact order $\bar u$. 
Note first that we may then view $u \in \Sigma(X)(\ZZ)$ or
$-u \in \Sigma(X)(\ZZ)$, and a point  $r\in \Sigma(X)(\ZZ)$ can
be viewed as determining two contact orders, one positive and one negative.
Take $r=u$ or $-u$ in the two cases.
\cite{ACGS18},~App.~A 
constructs a stack or scheme $\shZ_{\bar u}$ and $Z_{\bar u}$
equipped with strict morphisms to $\shX$ and $X$ respectively.
In the case at hand, these are particularly simple: they are simply
the strata $\shZ_r$ or $Z_r$ introduced in Construction~\ref{const:first}.
Equivalently, these are the strata corresponding to the cone $\sigma$.

Per \cite{ACGS18},~App.~A, these now satisfy the following property.
Given a punctured map $f:C^\circ/W \rightarrow X$ or 
$f:C^{\circ}/W\rightarrow \shX$, and a puncture $x$ given by a section
$x:\ul{W}\rightarrow \ul{C}$ with global contact order $\bar u$, 
then the morphism
\[
f\circ x:(\ul{W},x^*\shM_{C^{\circ}})\longrightarrow C^{\circ}\longrightarrow 
\hbox{$X$ or $\shX$}
\]
factors through $Z_{\bar u}$ or $\shZ_{\bar u}$. We call 
$Z_{\bar u}=Z_r$, $\shZ_{\bar u}=\shZ_r$ the \emph{evaluation stratum}
for the contact order $\bar u$.

Note that we obtain an induced morphism
$r:\overline\shM_{Z_{r}}\rightarrow \NN$ as follows. Let $\bar\xi$
be the generic point of $Z_r$. Then $r\in \Hom(\overline\shM_{X,\bar\xi},
\ul{\NN})$. Thus given a section of $s$ of $\overline\shM_{Z_r}$, 
we obtain its germ $s_{\bar \xi}\in\overline{\shM}_{Z_r,\bar\xi}$ and
define $r(s)=r(s_{\bar\xi})$.
\smallskip

(2) In case that $\overline{\shM}_X^{\gp}\otimes\QQ$ is globally generated
as in Basic Setup~\ref{setup1},
the natural surjective maps 
$\Gamma(X,\overline{\shM}_X^{\gp}\otimes\QQ)\rightarrow 
\overline{\shM}_{X,x}^{\gp}\otimes\QQ$
dualize to give injective maps $\sigma_x\hookrightarrow \sigma_x^{\gp}
\rightarrow
V:=\Hom(\Gamma(X,\overline{\shM}^{\gp}_X\otimes \QQ),\RR)$. 
Compatibility with generization then yields a continuous map
$|\Sigma(X)|\rightarrow V$ 
which is not necessarily injective, but is injective on individual cones,
see \cite{ACGS18}, Prop.~3.13. By choosing any representative
$u\in N_{\sigma}$ for a global contact order $\bar u\in 
\foC(X)$, we obtain an element $u'\in V$ depending only on $\bar u$. Note that
if $u\in N_{\sigma}$ is a representative for $\bar u$,
then in the notation introduced in \S\ref{sec:mirrors},
we have
$\langle u, s\rangle=u'(s)$ for $s\in\Gamma(X,\overline{\shM}_X)$.
Further, because the map $|\Sigma(X)|\rightarrow V$ is injective
on cones, each $u'\in V$ arises from at most a finite number
of global contact orders $\bar u$.
\end{remarks}

\subsubsection{Moduli of punctured maps}
\label{subsub:moduli}
A class $\beta$ of punctured map (\cite{ACGS18}, Def.~2.44,(4)) is a choice of
curve class $A\in H_2(X)$, a genus, a number of punctured points,
and for each punctured point $x$, a choice of global contact order $\bar u_x\in
\foC(X)$.  Generally, in this paper we will specify a choice of global
contact order by specifying
an integral tangent vector $u_x$ to a cone $\sigma\in \Sigma(X)$ such that
$u_x$ is not tangent to a proper face of $\sigma$.

Most of the time, we will only need the cases considered
in Remarks~\ref{rem:contact order remarks},(1).
A point $u \in \Sigma(X)(\ZZ)$ uniquely determines a minimal cone
$\sigma\in\Sigma(X)$ containing $u_x$, and viewing $u_x\in N_{\sigma}$,
we obtain a global contact
order of a marked point. Similarly, $-u_x$ determines a contact order
of a punctured point. More general contact orders will
only appear in the arguments of \S\ref{subsec:no tails}.

Given a class $\beta$ of punctured map, 
we have a moduli space $\scrM(X,\beta)$ of stable punctured maps of
class $\beta$, with target space $X$. Here the stability condition
is just on the underlying stable map of curves.
This moduli space is a DM stack.  After forgetting the curve class
$A$ from the data of $\beta$, we have a moduli space 
$\foM(\shX,\beta)$ of pre-stable punctured maps with target space $\shX$,
the Artin fan of $X$ as described in \S\ref{subsubsec:artin fan}.
This moduli space is an algebraic stack which is
a logarithmically enhanced version of the
moduli space of pre-stable curves in ordinary Gromov-Witten theory.

Both moduli spaces $\scrM(X,\beta)$, $\foM(\shX,\beta)$ carry a canonical
choice of log structure called the \emph{basic log structure}, see
\cite{JAMS}, \S1.5, \cite{ACGS17}, \S2.3.7 and \cite{ACGS18}, 
Def.~2.31.
In particular, we call the stalk of the ghost sheaf at a point corresponding
to a given punctured map the \emph{basic monoid}.

Via composition of stable maps
with the canonical morphism $X\rightarrow \shX$, we obtain a strict morphism
\[
\scrM(X,\beta)\rightarrow \foM(\shX,\beta),
\]
and \cite{ACGS18}, \S4.2, defines a perfect
relative obstruction theory for this morphism, and hence a virtual pull-back
of cycles via \cite{Man}. In general, $\foM(\shX,\beta)$ is not
equi-dimensional, so there is not necessarily a virtual fundamental class
on $\scrM(X,\beta)$. 

In fact, \cite{ACGS18}, \S4.2, gives a somewhat more general setup.
Let ${\bf x}$ be a set of disjoint sections of the universal curve over 
$\foM(\shX,\beta)$; these may include nodal sections. Then we set
\begin{equation}
\label{eq:Mev def}
\foM^{\ev({\bf x})}(\shX,\beta)= \foM^{\ev}(\shX,\beta)
:= \foM(\shX,\beta)\times_{\prod_{x\in {\bf x}}\ul{\shX}}
\prod_{x\in {\bf x}}\ul{X}.
\end{equation}
We use the second notation when ${\bf x}$ is clear from context. The
morphism $\foM(\shX,\beta)\rightarrow \prod_{x\in {\bf x}}\ul{\shX}$
is given by evaluation at the sections in ${\bf x}$, which makes sense
at the non-logarithmic level. There is then a natural morphism
\[
\xymatrix@C=30pt
{
\scrM(X,\beta)\ar[r]^{\varepsilon} & \foM^{\ev}(\shX,\beta)\ar[r]&
\foM(\shX,\beta)
}
\]
and a perfect relative obstruction theory for $\varepsilon$.
More generally, if ${\bf x}'\subseteq {\bf x}$, we have a triple of spaces
\[
\xymatrix@C=30pt
{
\scrM(X,\beta)\ar[r]^{\varepsilon}& \foM^{\ev({\bf x})}(\shX,\beta)
\ar[r]^{\pi}&
\foM^{\ev({\bf x'})}(\shX,\beta)
}
\]
with $\varepsilon'=\pi\circ\varepsilon$. Further, $\pi$ is smooth, and hence
carries a canonical choice of perfect relative obstruction theory given by
the relative cotangent complex, and then Manolache's
virtual pull-back $\pi^!$ coincides
with flat pull-back $\pi^*$ (see \cite{Man}, Rem.~3.10). This
yields a compatible triple of obstruction theories. In particular, 
if $\foM(\shX,\beta)$
is pure-dimensional, a virtual fundamental class can be defined
using any of these obstruction theories, and they all coincide.
In this situation and in similar situations encountered later, we say
the obstruction theories for $\varepsilon$ and $\varepsilon'$ 
are \emph{compatible}.

The virtual relative dimension of $\varepsilon$ at a point of $\scrM(X,\beta)$
represented by $f:C^{\circ}\rightarrow X$ is
\begin{equation}
\label{eq:rel virtual dim}
\chi\left(f^*\Theta_{X/\kk}(-\sum_{x_i\in {\bf x}} x_i)\right),
\end{equation}
see \cite{ACGS18}, Prop.~4.5.

\begin{remark}
\label{rem:contact order determined}
One useful consequence of Remarks~\ref{rem:contact order remarks}, (2), 
crucial for various finiteness results in Lemma \ref{lem:finite}
and elsewhere, is the following. We
assume $\overline\shM_X^{\gp}\otimes\QQ$ is globally generated,
as in that remark. Suppose given a class 
of punctured map $\beta$ with punctured points $x_1,\ldots,x_n$,
with contact orders at $x_1,\ldots,x_{n-1}$ specified by global contact
orders $\bar u_{x_1},\ldots,\bar u_{x_{n-1}}$. Then $s\mapsto A
\cdot c_1(\shL_s)$ determines an element of 
$V=\Hom(\Gamma(X,\overline\shM_X^{\gp}\otimes\QQ),\RR)$,
while $\bar u_{x_1},\ldots,\bar u_{x_{n-1}}$ determine
elements $u'_{x_1},\ldots, u'_{x_{n-1}}\in V$ as
in Remarks~\ref{rem:contact order remarks}, (2). 
So by Proposition~\ref{intersectionnumbers},
if the global contact order at $x_n$ is given by $\bar u_{x_n}$,
its image $u'_{x_n}\in V$ is completely determined by
the data $A$, $u_{x_1},\ldots,u_{x_{n-1}}$, and hence again by the
discussion of Remarks~\ref{rem:contact order remarks}, (2),
$\bar u_{x_n}$ is determined up to a finite number of choices.
\end{remark}

\subsubsection{Tropical curves and maps}
\label{sec:tropical curves}
We follow \cite{ACGS17}, \S2.5 for the notion of tropical curve, with
a minor variation coming from the treatment of punctures, see 
\cite{ACGS18}, \S2.2. We consider
connected graphs $G$ with sets of vertices $V(G)$, edges $E(G)$
and legs $L(G)$. However, unlike in the marked point case, a leg
may be a compact interval or a ray. In either case, a leg has only one
endpoint in $V(G)$. A tropical curve $\Gamma=(G,{\bf g},\ell)$
of combinatorial type $(G,{\bf g})$ is the choice of a genus
function ${\bf g}:V(G)\rightarrow\NN$ and length functions
$\ell:E(G)\rightarrow \RR_{> 0}$ and $\ell:L(G)\rightarrow\RR_{>0}\cup\{
+\infty\}$. As all curves in this paper will be
of genus zero, we drop all future mention of ${\bf g}$.

We also consider \emph{families} of tropical curves parameterized
by a rational polyhedral cone $\omega$ defined in a lattice $N_{\omega}$: 
this is a choice of combinatorial
type $G$ along with a length function $\ell:E(G)\rightarrow
\Hom(\omega\cap N_{\omega},\NN)\setminus \{0\}$ as well as
a piecewise linear function $\ell(L)$ on $\omega$ for each bounded
$L\in L(G)$.
We will frequently write $\ell_E\in \omega^{\vee}$ 
instead of $\ell(E)$.

Associated to the data
$(G,\ell)$ is a generalized cone complex 
$\Gamma(G,\ell)$ along with a morphism of cone complexes
$\Gamma(G,\ell)\rightarrow \omega$ with fibre over $s\in\Int(\omega)$
being a tropical curve, i.e., a metric graph, with underlying graph
$G$ and affine edge length of $E\in E(G)$ being $\ell_E(s)$.
Associated to each vertex of $G$ is a copy $\omega_v$ of $\omega$
in $\Gamma(G,\ell)$. Associated to each edge or leg $E\in E(G) \cup L(G)$
is a cone $\omega_E \in \Gamma(G,\ell)$ with
$\omega_E\subseteq \omega\times\RR_{\ge 0}$ and the map to $\omega$
given by projection onto the first coordinate. This projection
fibres $\omega_E$ in compact intervals or rays over $\omega$
(rays only in the case of a leg representing a marked point).

A \emph{family of tropical maps}
to $\Sigma(X)$ over $\omega\in \mathbf{Cones}$ is a morphism of cone complexes
\[
h:\Gamma(G,\ell)\rightarrow \Sigma(X).
\]
If $s\in\Int(\omega)$, we may view $G$ as the fibre of
$\Gamma(G,\ell)\rightarrow \omega$ over $s$ as a metric graph, and
write
\begin{equation}
\label{eq:hs def}
h_s:G\rightarrow |\Sigma(X)|
\end{equation}
for the corresponding tropical map with domain $G$. Here $|\Sigma(X)|$
denotes the underlying topological space of the cone complex $\Sigma(X)$.
The \emph{type} of such a family consists of the data
$\tau:=(G,\bsigma,\mathbf{u})$ where
\[
\bsigma:V(G)\cup E(G)\cup L(G)\rightarrow \Sigma(X)
\]
associates to $x\in V(G)\cup E(G)\cup L(G)$ the minimal
cone of $\Sigma(X)$ containing $h(\omega_x)$. Further,
$\mathbf{u}$ associates to each (oriented) edge or leg $E\in E(G)\cup L(G)$
the corresponding \emph{contact order} $\mathbf{u}(E)\in N_{\bsigma(E)}$,
the image of the tangent vector $(0,1)\in N_{\omega_E}=
N_{\omega}\oplus\ZZ$ under the map $h$.

If $E\in E(G)$ is an edge with vertices $v, v'$, oriented from $v$
to $v'$, and $s\in \Int(\omega)$, the induced tropical map
$h_s:G\rightarrow |\Sigma(X)|$ satisfies
\begin{equation}
\label{eq:v1v2rel}
h_s(v')-h_s(v)=\ell_E(s)\cdot \mathbf{u}(E)
\end{equation}
as an equation in $N_{\bsigma(E)}\otimes_{\ZZ}\RR$.

Given a family of tropical maps $h$ as above parameterized by $\omega$, 
for $v$ any vertex of $G$, we write
\begin{align}
\label{eq:nuv}
\begin{split}
\nu_v:\omega\longrightarrow &\,\,\bsigma(v)\\
s\longmapsto &\,\,h_s(v)
\end{split}
\end{align}
Equivalently, $\nu_v$ is the restriction of $h$ to
the copy $\omega_v\in\Gamma(G,\ell)$ of $\omega$ indexed by $v$. 

\medskip

Stable log or punctured maps give families of tropical maps,
see \cite{ACGS17}, \S 2.5 and \cite{ACGS18}, \S2.2. Suppose
we have a geometric point $\bar w$ in any of the moduli spaces we consider
in this paper, for example, $\scrM(X,\beta)$ or $\foM(\shX,\beta)$ 
just mentioned. Denote the moduli space under consideration by $\scrM$. 
We obtain a diagram
\[
\xymatrix@C=30pt
{
C_{\bar w}^{\circ}\ar[r]^>>>>f\ar[d]_{\pi} & \hbox{$X$ or $\shX$}\\
\bar w
}
\]
After tropicalizing, we obtain a diagram (noting that
$\Sigma(X)=\Sigma(\shX)$ is shown in Step I of the proof
of \cite{ACGS17}, Prop.~2.10):
\[
\xymatrix@C=30pt
{
\Sigma(C_{\bar w}^{\circ})\ar[r]^>>>>>{\Sigma(f)}\ar[d]_{\Sigma(\pi)} & 
\Sigma(X)=\Sigma(\shX)\\
\Sigma(\bar w)
}
\]
Note that $\Sigma(\bar w)=\sigma_{\bar w}=\Hom(Q,\RR_{\ge 0})$ where
$Q$ is the stalk of the ghost sheaf of $\scrM$ at $\bar w$.

A fibre of $\Sigma(\pi)$ over a point in the 
interior of $\sigma_{\bar w}$
is homeomorphic to the dual intersection graph $G$ of $C^{\circ}_{\bar w}$. 
The vertices $V(G)$ are indexed by the irreducible components
of $C^{\circ}_{\bar w}$, the edges $E(G)$ are indexed by 
the nodes of $C^{\circ}_{\bar w}$, and the legs $L(G)$ are
indexed by the marked or punctured points of $C^{\circ}_{\bar w}$.
Marked points on $C^{\circ}_{\bar w}$ correspond to legs of $G$
homeomorphic to $[0,+\infty)$ and punctured points which are not
marked points correspond to legs
of $G$ homeomorphic to $[0,1]$. Further, $\Sigma(\pi)$ provides
a family of tropical curves: in particular,
length functions $\ell_E\in Q\setminus \{0\}$
for each edge $E$ of $G$. 

Further, $h:=\Sigma(f)$ now determines a tropical map to $\Sigma(X)$.
The corresponding type $\tau=(G,\bsigma,{\bf u})$ satisfies the
following. We have ${\bf u}(E)=u_q$
for each oriented edge $E$ corresponding to a node $q$, where $u_q$ is
the contact order of the node, see
\cite{ACGS17}, \S2.3.5,(iii). We have ${\bf u}(L)=u_x$  
for each leg $L$ corresponding
to a marked or punctured point $x$, where $u_x$ is the contact order of the
point $x$ as defined in \eqref{eq:ux def}.

It follows from \cite{ACGS18}, Prop.~2.23, that for $L\in L(G)$ with adjacent vertex $v\in V(G)$
giving $\omega_L,\omega_v\in \Gamma(G,\ell)$, we have
\begin{equation}
\label{eq:leg}
h(\omega_L)=(h(\omega_v)+\RR_{\ge 0}{\bf u}(L))\cap \bsigma(L)
\subseteq N_{\bsigma(L),\RR}.
\end{equation}

If $\scrM$ carries the basic log structure, then
this family of tropical maps is universal,
see \cite{ACGS18}, \S2.3. However, we may often deal with families
that are not universal in this paper.

\subsection{Some comments on stacks}

\subsubsection{Controlling the behaviour of the algebraic stacks appearing
here}
\label{sec:stacks remarks}

All stacks appearing in this paper are algebraic, but do not necessarily
satisfy the strictest notions such as those of \cite{champs}. Indeed,
if $X$ is an algebraic stack over a scheme $S$, then \cite{champs}
assumes that the diagonal morphism $\Delta:X\rightarrow X\times_S X$
is representable, separated and quasi-compact. On the other hand,
if $S$ is a log scheme, Olsson's classifying stack $\Log_S$ of log schemes over
$S$ instead satisfies the condition that $\Delta$
be representable and of finite presentation, see \cite{OlssonENS},
(1.2.4) and (1.2.5), and in particular the diagonal need not be
separated. In general, there is nothing which can be done about this:
see \cite{OlssonENS}, Rem.~3.12
for a very simple example showing $\Log_{\kk}$ has non-separated diagonal.
However, fortunately this does not cause any problems, as Kresch's
construction of Chow groups \cite{Kresch} does not require a separated
diagonal.

In general, the stacks $\foM(\shX,\beta)$ inherit all the problems
suffered by $\Log_S$. In particular, the fact that $\foM(\shX,\beta)$
is not finite type will cause persistent problems, and in particular
$\foM(\shX,\beta)$ would be expected to contain an infinite number of
logarithmic strata. We always
deal with this issue in the following way. We will frequently have
a situation where we have a strict morphism $\scrM\rightarrow \foM$
where $\scrM$ is some kind of moduli space of punctured maps to $X$,
while $\foM$ is a similar moduli space of punctured maps to $\shX$.
In every case, $\scrM$ will be of finite type over $\Spec\kk$, and hence
contain only a finite number of logarithmic strata. 
Thus we can replace $\foM$ by the complement of all logarithmic strata
whose closures do not intersect the image of $\scrM\rightarrow \foM$, and in
particular assume that $\foM$ has only a finite number of logarithmic
strata. Then $\foM$ becomes finite type, and $\Sigma(\foM)$ contains 
only a finite number of distinct
cones. This procedure will generally be applied without comment;
however, some care will be necessary at certain points. 

In order to apply results about virtual pull-backs, we will need our
stacks to be stratified by global quotients as in \cite{Man},
Construction~3.6. This condition is necessary to use Kresch's result
\cite{Kresch}, Prop.~4.3.2: if $E\rightarrow F$ is a vector bundle
stack and $F$ is stratified in global quotients, then the Chow groups
of $E$ and $F$ are isomorphic. By \cite{Kresch}, Prop.~3.5.9, 
a stack is stratified by global quotients if the stabilizer group of every
geometric point is affine. This is straightforward to check in all
cases in this paper, and we omit discussion of this point in the sequel.

\subsubsection{Log \'etale morphisms of log stacks}
\label{sec:log etale}

When a property of a (log) morphism $f:X\rightarrow Y$ is local in the
(strict) smooth topology for the domain and range, then one gets a 
corresponding notion of the property
for morphisms of (log) stacks as in \cite{stacks}, Tag 06FL. However, as we
frequently use the notion of (log) \'etale, which is not local in the
(strict) smooth topology, we recall:

\begin{definition}
A morphism $\ul{f}:\ul{X}\rightarrow \ul{Y}$ of stacks is \emph{smooth}
(resp.\ \'etale) if it is locally of finite presentation and formally
smooth (resp.\ \'etale).
\end{definition}

\begin{remark}
\label{rem:etale stack}
There are in fact a variety of definitions of \'etale morphisms of stacks
in the literature.
The definition given here is equivalent, by
\cite{rydh}, Cor.~B.9, to the definition of an \'etale morphism of
algebraic stacks given in \cite{rydh} (i.e., locally of finite presentation,
flat, and with \'etale diagonal). By \cite{stacks}, Tag 0CJ1,
Rydh's definition of \'etale is equivalent to that given in
\cite{stacks}, Tag 0CIL, and in particular an \'etale morphism
of stacks is of DM type by that definition. 
\end{remark}

Similarly, we have:

\begin{definition}
A morphism $f:X\rightarrow Y$ between log stacks
is \emph{log smooth} (resp.\ \emph{log \'etale})
if it is locally of finite presentation and formally log
smooth (resp.\ formally log \'etale), see \cite{Ogus}, IV,Def.~3.1.1.

Idealized log smooth or log \'etale morphisms are defined similarly,
as in \cite{Ogus}, IV,\S3.1.
\end{definition}

We note that while not stated in this degree of generality, the arguments
of \cite{OlssonENS}, Thm.~4.6 show that if $f:X\rightarrow Y$
is log smooth (log \'etale), then the induced morphism of ordinary
stacks $\ul{X}\rightarrow \Log_Y$ is smooth (\'etale).

\section{Construction of the invariants $N^A_{p_1p_2r}$}
\label{sec:invariant construction}
Throughout this section, we fix $X$ as in Basic Setup~\ref{setup1}.
We will also fix an element $r\in \Sigma(X)(\ZZ)$, so that $-r$ can
be interpreted as a contact order for a puncture, as discussed
in Remarks~\ref{rem:contact order remarks}, (1).
The purpose of this section is to define the structure coefficients $N^A_{p_1p_2r}$ together with some properties involved in their definition.
  
If $r=\sum_j r_jD_j^*$ is tangent to a $1$-dimensional cone $\RR_{\ge0} D_i^*$ of $\Sigma(X)$,
the invariant $N^A_{p_1p_2r}$ counts curves with negative contact $-r_i$ at a
general point of $D_i$ as follows. The curve class $A$ along with the contact
orders $p_1,p_2,-r$ defines a class $\beta$ of punctured maps
\cite{ACGS18},~Def.\,2.44,(4). The moduli space $\scrM(X,\beta)$ of punctured
maps of class $\beta$ turns out to have a virtual fundamental class of dimension
$c_1(\Theta_{X/k})\cdot A+(\dim X-1)$, and it comes with a (non-logarithmic)
evaluation map at the point of contact order $-r$
\[
\ev: \ul{\scrM(X,\beta)}\lra D_i.
\]
A first attempt that is almost correct is to define $N^A_{p_1p_2r}$ as the
degree of $[\scrM(X,\beta)]^\virt$ over $D_i$ if $c_1(\Theta_{X/k})\cdot A=0$,
and $N^A_{p_1p_2r}=0$ otherwise. With this definition, $N^A_{p_1p_2r}$ simply
counts punctured maps of class $\beta$ with one schematic point constraint at a
general point of $D_i$. Our actual definition of $N^A_{p_1 p_2 r}$ differs from
this count by logarithmic multiplicities in the computation of
$\ev_*[\ul\scrM(X,\beta)]$. At a geometric point of $\scrM(X,\beta)$ this
multiplicity equals the integral length of the image of the unit tangent vector
by the tropical evaluation map at the leg with contact order $-r$. See the
appearance of $k_\tau$ in \cite{CanScat},~Thm.~6.4. These
multiplicities arise from a logarithmic refinement of the evaluation map.

The need for a logarithmic refinement of the evaluation map becomes evident if
$r$ lies in the interior of a higher-dimensional cone of $\Sigma(X)$. In this
case $\scrM(X,\beta)$ may not be virtually pure-dimensional and may
have virtual components of strictly larger dimension than the stratum
$Z_r$. We then impose the point condition
logarithmically, akin to the definition of degrees of proper log \'etale
morphisms discussed in Remark~\ref{rem:transverse}. This is done by enhancing
$\ev$ to a logarithmic evaluation map with target a stacky enhancement
$\scrP(X,r)$ of $Z_r$ carrying logarithmic point constraints of contact order
$-r$. The stacky refinement of the evaluation map captures the map from the
$\GG_m$-torsor on $Z_r$ defined by $-r$ to the conormal bundle of the domain
curve at the punctured point. Such a refinement of evaluation maps was first
suggested in \cite{ACGM}.

The point constraint can then be imposed by a logarithmic base-change with an
appropriate map $B\GG_m^\ls\to \scrP(X,r)$. Here $B\GG_m^\ls$ is the stacky
quotient of the standard log point by $\GG_m$.\footnote{An alternative
definition of $N^A_{p_1p_2r}$ found in later work, in the case of most interest
for mirror symmetry, uses a sum over punctured invariants associated to a
realizable global type from \cite{ACGS18},~Def.~4.6. See \cite{CanScat},~\S6.}

\subsection{Evaluation spaces}
The first goal is to define a space parameterizing punctures with
contact order given by $-r$, where $r\in \Sigma(X)(\ZZ)$. This space
is analogous to the evaluation spaces constructed in \cite{ACGM}.

We define log stacks 
$\scrP(X,r),\widetilde{\scrP}(X,r)$. Here, $\scrP(X,r)$ 
is the ``moduli space of punctures in $X$ with contact order $-r$,''
and $\widetilde\scrP(X,r)$ is the universal family of such punctures.
This is made precise in
Proposition~\ref{prop:evaluation maps general}. First, we introduce
some notation which will be used throughout:

\begin{definition}
\label{def:BGm dagger}
$B\GG_m^{\dagger}$ denotes $B\GG_m$ with its universal log structure,
i.e., the pull-back to $B\GG_m$ of the divisorial log structure on 
$[\AA^1/\GG_m]$ induced by $B\GG_m\subset [\AA^1/\GG_m]$. 
Alternatively, in the notatation
of \cite{ACGS18}, (B.1), we may write this as $\shA_{\NN,\NN\setminus \{0\}}$,
and as such, can be viewed as an idealized log stack via the ideal
$\NN\setminus\{0\}\subseteq\overline\shM_{B\GG_m^{\dagger}}$.
\end{definition}

Let $Z:=Z_r$ be the closed stratum of $X$ indexed by $r$ as in Construction
\ref{const:first}, with its induced log
structure from $X$. Set
\[
\widetilde{\scrP}(X,r):= Z\times B\GG_m^{\dagger}.
\]
On the other hand, we define $\scrP(X,r)$ to have the same underlying
stack as $\widetilde{\scrP}(X,r)$, but with a sub-log structure
defined as follows.

If we abuse notation by writing $\overline{\shM}_Z$ also for the sheaf
on $\ul{\scrP(X,r)}$ given by pull-back of $\overline{\shM}_{Z}$ via
the projection $\widetilde\scrP(X,r)\rightarrow Z$, then we define
$\overline{\shM}_{\scrP(X,r)}$ as the subsheaf of
$\overline{\shM}_{\widetilde{\scrP}(X,r)}=\overline{\shM}_Z\oplus \NN$
given by 
\begin{equation}
\label{eq:MbarPdefined}
\Gamma(U,\overline{\shM}_{\scrP(X,r)})
=\{(m,r(m))\,|\, m\in \Gamma(U,\overline{\shM}_Z) \}
\end{equation}
where $r$ is viewed as an element of $\Hom(\overline{\shM}_Z,\NN)$, as
in Remarks~\ref{rem:contact order remarks}, (1). Note this gives
a canonical isomorphism between $\overline{\shM}_Z$ and 
$\overline{\shM}_{\scrP(X,r)}$.

We may now define the log structure on $\scrP(X,r)$ as
\[
\shM_{\scrP(X,r)}:=\overline{\shM}_{\scrP(X,r)}
\times_{\overline{\shM}_{\widetilde{\scrP}(X,r)}} \shM_{\widetilde{\scrP}(X,r)}.
\]

Similarly, $r$ defines a stratum $\shZ:=\shZ_r$ 
of the Artin fan $\shX$ of $X$,
and we can define $\scrP(\shX,r)$, $\widetilde{\scrP}(\shX,r)$.

\begin{remark}
\label{remark: quotient description}
The stack $\scrP(X,r)$ can also be described as a quotient 
logarithmic stack $[Z/\GG_m]$, where $\GG_m$ acts trivially on the underlying
scheme $\ul{Z}$ but acts on a torsor $\shL^{\times}_m$ over $U\subseteq Z$ open 
with weight
$r(m)$ for $m\in\Gamma(U,\overline{\shM}_Z)$.
\end{remark} 

\begin{proposition} 
\label{prop:P is punctured}
Let $\shP$ denote the log structure on $\ul{\scrP(X,r)}$
which is the pull-back of the universal log structure on $B\GG_m^{\dagger}$
under the projection
to $B\GG_m$. If $\shM$ is the log structure of $\scrP(X,r)$,
write 
\[
\shM'=\shM\oplus_{\O_{\scrP(X,r)}^{\times}} \shP.
\]
Then $\widetilde\scrP(X,r)$ is a puncturing of $\shM'$ along $\shP$,
and this puncturing is pre-stable with respect to the composed morphism
\[
f:\widetilde{\scrP}(X,r)\mapright{\pr_1} Z\hookrightarrow X
\]
in the sense of Definition~\ref{def:prestable}.
\end{proposition}

\begin{proof}
Note that $\overline\shM'=\overline\M\oplus\NN$, and there is a natural
inclusion $\overline\shM'\hookrightarrow 
\overline{\shM}_{\widetilde\scrP(X,r)}=\overline\shM_Z\oplus\NN$
given by $(m,n)\mapsto (m,r(m)+n)$. This induces an inclusion 
$\shM'\hookrightarrow \shM_{\widetilde\scrP(X,r)}$.
If $(\hat m,\hat n)\in
\shM_{\widetilde\scrP(X,r)}\setminus\shM'$ with image $(m,n)$
in the ghost sheaf, then necessarily $0\le n<r(m)$,
and $\alpha_Z(\hat m)=0$ if $r(m)>0$.
Thus $(\hat m,\hat n)$ is sent to zero by the structure map
$\alpha$ of $\widetilde{\scrP}(X,r)$. Hence we obtain a puncturing. 
This puncturing is pre-stable, as $\overline\shM_{\widetilde\P(X,r)}$
is generated by the image $\overline\shM_Z\oplus 0$ of $\bar f^{\flat}$
and $(0,1)\in \overline \shM'$.
\end{proof}

\begin{proposition}
\label{prop:evaluation maps general}
Let $f:C^{\circ}/W\rightarrow X$ be a pre-stable punctured map with a punctured point
$x:\ul{W}\rightarrow \ul{C}$ with contact order given by $-r$. 
Then there exists a canonical
morphism $\ev:W\rightarrow \scrP(X,r)$ with the property 
that
\[
W^{\circ}:=W\times^{\fine}_{\scrP(X,r)} \widetilde{\scrP}(X,r)
\]
coincides with $(\ul{W}, x^*\shM_{C^{\circ}})$.

The analogous statements for punctured maps $f:C^{\circ}/W\rightarrow\shX$
also hold.
\end{proposition}

\begin{proof}
Write $W':=(\ul{W},x^*\shM_{C^{\circ}})$ and $W^{\mathrm{m}}:=(\ul{W},
x^*\shM_C)$ (recalling that $\shM_C$ treats $x$ as a marked point). 
Then we can write $\shM_{W^{\mathrm{m}}}=\shM_W\oplus_{\O_W^{\times}}\shP_W$,
where $\shP_W$ is the pull-back of the divisorial log structure
$x(\ul{W})\subset \ul{C}$ via $x$. Thus $W'$ is a puncturing of $W^{\mathrm{m}}$
along $\shP_W$, and the induced morphism $f\circ x:W'\rightarrow X$
is pre-stable, as $f$ was pre-stable.

Next we define a morphism $W'\rightarrow \widetilde\scrP(X,r)$. First, by
Remarks~\ref{rem:contact order remarks}, (1), $Z=Z_r$ 
is the evaluation 
stratum induced by the contact order $-r$ and
the morphism $f\circ x$ factors through $Z$. 
Second, there is a canonical strict morphism $(\ul{W}, \shP_W)
\rightarrow B\GG_m^{\dagger}$. 
We note that at the level of underlying stacks, the morphism
$\ul{W}\rightarrow B\GG_m$ is the morphism induced by the torsor
on $\ul{W}$ corresponding to $1\in \overline{\shP}_W=\NN$, and this
torsor is associated to the conormal line bundle $\shN^{\vee}_{x(\ul{W})/
\ul{C}}$.
As there is a natural morphism
$W'\rightarrow (\ul{W},\shP_W)$ given by the inclusion $\shP_W\subset
\shM_{W'}$, we then obtain a morphism 
\begin{equation}
\label{eq:ev prime}
\widetilde\ev:
W'\rightarrow Z\times B\GG_m^{\dagger} =\widetilde\scrP(X,r).
\end{equation}

We now have a diagram
\begin{equation}
\label{eq:WPX}
\xymatrix@C=30pt
{
W'\ar[r]^>>>>{\widetilde\ev}\ar[d] & \widetilde\scrP(X,r)\ar[d]\ar[r]&X\\
W\ar@{.>}[r]_>>>>{\ev}&\scrP(X,r)&
}
\end{equation}
and we wish to define the dotted arrow $\ev$ to make the diagram
commute. Since the vertical arrows are the identity on underlying
stacks, it is enough to understand this diagram at the level of
ghost sheaves.

Let $\bar w\in \ul{W}$, and let
$Q=\overline{\shM}_{W,\bar w}$, $P=\overline{\shM}_{X,f(x(\bar w))}$.
Then $\overline{\shM}_{W',\bar w}\subset Q\oplus\ZZ$, and
\eqref{eq:WPX} induces
\[
\xymatrix@C=30pt
{
Q\oplus\ZZ & P\oplus\NN \ar[l]_{\ol{\widetilde\ev}^{\flat}}& P\ar[l]_>>>>>{\id\oplus 0}\\
Q\ar[u]^{\id\oplus 0} & P \ar@{.>}[l]^{\overline{\ev}^{\flat}}\ar[u]_{\id\oplus r}
}
\]
By construction
of $\widetilde\ev$, using the notation $\varphi_x:=\bar f^{\flat}:P
\rightarrow \overline{\shM}_{W',\bar w}\subseteq Q\oplus\ZZ$,
\[
\overline{\widetilde\ev}^{\flat}(m,n) = \varphi_x(m)+(0,n)=
(\pr_1\circ\varphi_x(m), -r(m)+n).
\]
Note that for $m\in P$, $\overline{\widetilde\ev}^{\flat}(m,r(m))=
(\pr_1\circ\varphi_x(m),0)$. Hence $\overline{\widetilde\ev}^{\flat}$
induces a homomorphism $\ol{\ev}^{\flat}$ making the diagram commute,
with 
\begin{equation}
\label{eq:evghost}
\ol{\ev}^{\flat}(m)=\pr_1\circ\varphi_x(m).
\end{equation}
Thus by \eqref{eq:MbarPdefined},
$\widetilde\ev$ induces the morphism $\ev$, turning \eqref{eq:WPX} into
a commutative diagram.

The commutativity of the square in \eqref{eq:WPX} 
induces a morphism $W'\rightarrow W^{\circ}$. On the other hand,
by \cite{ACGS18}, Def.~2.9,
$W^{\circ}$ is the pull-back puncturing of the
punctured log structure on $\widetilde{\scrP}(X,r)$, precisely a puncturing
of $W^{\mathrm{m}}$. Further this puncturing is pre-stable with respect
to the morphism to $X$ by Proposition~\ref{prop:P is punctured} and 
\cite{ACGS18}, Cor.~2.10. Since the same is true for $W'$, it then follows
from \cite{ACGS18}, Prop.~2.5,
that the induced morphism $W'\rightarrow W^{\circ}$ is an isomorphism,
as desired.
\end{proof}

%\begin{definition}
%\label{def: tilde W}
%In the situation of Proposition~\ref{prop:evaluation maps general},
%we write $\widetilde W\rightarrow W^{\circ}$ for the saturation
%of $W^{\circ}$. Note that
%\begin{equation}
%\label{eq: tilde W}
%\widetilde W= W\times^{\fs}_{\scrP(X,r)} \widetilde{\scrP}(X,r).
%\end{equation}
%\Mark{It is not clear we need to introduce this notation here. At the moment,
%we start using tildes again in \S7.1, but redefine them and mention this
%is compatible with this definition. But perhaps we can just delete this
%definition.}
%\end{definition}

\begin{remark}
\label{rem:P idealized}
We give $\scrP(\shX,r)$ the structure of an idealized log stack
as follows. We can view $r$ as an element
of $\Hom(\overline{\shM}_{Z_r},\NN)$ as in 
Remark~\ref{rem:contact order remarks},
(1), or equivalently as an element
of $\Hom(\overline{\shM}_{\scrP(X,r)},\NN)$. 
Let $\overline{\shI}\subset\overline{\shM}_{\scrP(\shX,r)}$ be the sheaf
of ideals given by $r^{-1}(\ZZ_{>0})$. Let $\shI$ be the
induced ideal of $\shM_{\scrP(\shX,r)}$. This defines a coherent idealized
log structure on $\scrP(\shX,r)$. 
\end{remark}

\subsection{Moduli spaces and evaluation maps}

We now fix a class of punctured map $\beta$ for $X$. We shall assume
that there is at least one punctured point $x_{\out}$ with contact order
determined by $-r$ for $r\in \Sigma(X)(\ZZ)$.

\begin{definition}
\label{def:evaluation maps}
We write
\[
\ev_X:\scrM(X,\beta)\rightarrow \scrP(X,r),\quad\quad
\ev_{\shX}:\foM(\shX,\beta)\rightarrow \scrP(\shX,r)
\]
for the evaluation maps at the point $x_{\out}$
given by Proposition~\ref{prop:evaluation maps general}
with $W=\scrM(X,\beta)$ or $\foM(\shX,\beta)$ respectively.
\end{definition}

%\begin{proposition}
%\label{prop:evaluation maps}
%With $\beta$ as given:
%\begin{enumerate}
%\item
%There are canonical evaluation maps $\ev_X$, $\ev_{\shX}$, $\widetilde{\ev}_X$,
%$\widetilde{\ev}_{\shX}$ in the following squares, each of which are
%cartesian in the fs log category:
%\[
%\xymatrix@C=15pt
%{
%\widetilde\scrM(X,\beta)\ar[r]^{\widetilde\ev_X}\ar[d] & 
%\widetilde{\scrP}(X,r)\ar[d]\\
%\scrM(X,\beta)\ar[r]_{\ev_X} & 
%\scrP(X,r)
%}
%\quad\quad\quad\quad\quad
%\xymatrix@C=15pt
%{
%\widetilde\foM(\shX,\beta)\ar[r]^{\widetilde\ev_{\shX}}\ar[d] & 
%\widetilde{\scrP}(\shX,r)\ar[d]\\
%\foM(\shX,\beta)\ar[r]_{\ev_\shX} & 
%\scrP(\shX,r)
%}
%\]
%\item
%$\widetilde\scrM(X,\beta)=\widetilde{\foM}(\shX,\beta)
%\times_{\foM(\shX,\beta)} \scrM(\shX,\beta)$.
%\end{enumerate}
%\end{proposition}
%
%\begin{proof}
%(1) follows from Proposition~\ref{prop:evaluation maps general}
%and Definition~\ref{def: tilde W}.
%For (2), note that certainly $\scrM(X,\beta)^{\circ}
%=\foM(\shX,\beta)^{\circ}\times_{\foM(\shX,\beta)}
%\scrM(X,\beta)$, as the universal curve over
%$\scrM(X,\beta)$ is the strict pull-back of the universal curve
%over $\foM(\shX,\beta)$. On the other hand, if $X\rightarrow Y$
%is a strict morphism of fine log schemes, then $X^{\sat}=X\times_Y Y^{\sat}$,
%from which the claim follows.
%\end{proof}

\begin{definition}
\label{def:Mev}
We set 
\[
\foM^{\ev}(\shX,\beta)
=\foM^{\ev(x_{\out})}(\shX,\beta):=\foM(\shX,\beta)\times_{\ul{\shX}} \ul{X},
\]
where the morphism $\foM(\shX,\beta)\rightarrow \ul{\shX}$
is (schematic) evaluation at $x_{\out}$. 
\end{definition}

\begin{lemma}
\label{lem:evaluation diagrams}
There is a cartesian diagram in all categories
\[
\xymatrix@C=30pt
{
\foM^{\ev}(\shX,\beta)\ar[d]\ar[r]^{\ev_{\shX}}&
\scrP(X,r)\ar[d]\\
\foM(\shX,\beta)\ar[r]_{\ev_{\shX}}&
\scrP(\shX,r)
}
\]
defining in particular the \emph{evaluation map} $\ev_{\shX}:
\foM^{\ev}(\shX,\beta)\rightarrow\scrP(X,r)$.
\end{lemma}

\begin{proof}
As the schematic evaluation map $\foM(\shX,\beta)\rightarrow 
\ul{\shX}$ factors through the closed substack $\shZ:=\shZ_r$
by Remarks~\ref{rem:contact order remarks}, (1), we may replace
$\ul{\shX}$ by $\ul{\shZ}$ and $\ul{X}$ by $\ul{Z}$ in the definition
of $\foM^{\ev}(\shX,\beta)$. Since both
vertical maps in the diagram are strict, it is enough
then to check the diagrams are cartesian at the level of underlying
stacks, which is immediate.
\end{proof}

%\begin{lemma}
%\label{lem:evaluation diagrams}
%There are cartesian diagrams in all categories
%\[
%\xymatrix@C=30pt
%{
%\foM^{\ev}(\shX,\beta)\ar[d]\ar[r]^{\ev_{\shX}}&
%\scrP(X,r)\ar[d]\\
%\foM(\shX,\beta)\ar[r]_{\ev_{\shX}}&
%\scrP(\shX,r)
%}
%\quad\quad\quad
%\xymatrix@C=30pt
%{
%\widetilde\foM^{\ev}(\shX,\beta)\ar[d]\ar[r]^{\widetilde\ev_{\shX}}&
%\widetilde\scrP(X,r)\ar[d]\\
%\widetilde\foM(\shX,\beta)\ar[r]_{\widetilde\ev_{\shX}}&
%\widetilde\scrP(\shX,r)
%}
%\]
%defining in particular the \emph{evaluation maps} $\ev_{\shX}$, 
%$\widetilde{\ev}_{\shX}$ in the top row of each diagram.
%In addition,\footnote{Not sure we need this last one yet.}
%there is a cartesian diagram in the fs log category
%\[
%\xymatrix@C=30pt
%{
%\widetilde\foM^{\ev}(\shX,\beta)\ar[r]^{\widetilde{\ev}_{\shX}}
%\ar[d] & \widetilde\scrP(X,r)\ar[d]\\
%\foM^{\ev}(\shX,\beta)\ar[r]_{\ev_{\shX}}& \scrP(X,r)
%}
%\]
%\end{lemma}
%
%\begin{proof}
%As the schematic evaluation map $\foM(\shX,\beta)\rightarrow 
%\ul{\shX}$ factors through the closed substack $\shZ$, we may replace
%$\ul{\shX}$ by $\ul{\shZ}$ and $\ul{X}$ by $\ul{Z}$ in the definitions
%of $\foM^{\ev}(\shX,\beta)$ and $\widetilde\foM^{\ev}(\shX,\beta)$. Since
%all vertical maps in the first two above diagrams are strict, it is enough
%then to check the diagrams are cartesian at the level of underlying
%stacks, which is then immediate. The fact the last diagram is cartesian
%in the fs log category is immediate from Proposition
%\ref{prop:evaluation maps} and the definitions.
%\end{proof}

\subsection{Imposing point constraints}
\begin{proposition}
\label{prop: BGm morphism}
Fix $r\in \Sigma(X)(\ZZ)$ and a closed point $z\in Z_r^{\circ}$,
the locally closed stratum of $X$ with closure $Z_r$.
There is a morphism of logarithmic stacks
\[
B\GG_m^{\dagger}\rightarrow \scrP(X,r)
\]
with image $z\times B\GG_m$ which, on the level of ghost sheaves, is given
by
\[
r:P:=\overline{\shM}_{X,z}\cong\overline{\shM}_{\scrP(X,r),z}\rightarrow\NN,
\]
with the isomorphism induced by \eqref{eq:MbarPdefined}. This morphism
of logarithmic stacks is also a morphism of idealized log stacks, with
domain and target carrying idealized structures as in Definition
\ref{def:BGm dagger} and Remark \ref{rem:P idealized} respectively.
\end{proposition} 

\begin{proof}
Using the definition of $\scrP(X,r)$, we can describe the restriction
of the log structure on $\scrP(X,r)$ to $z\times B\GG_m \subset \ul{Z}_r
\times B\GG_m$. The ghost sheaf of the restriction 
is the constant sheaf with stalk $P$,
and the torsor associated to $m\in P$ is $\shU^{\otimes r(m)}$, where $\shU$
is the universal torsor on $B\GG_m$. Thus we can lift the stated
map on the level of ghost sheaves to a log morphism, as the torsor
on $B\GG_m^{\dagger}$ induced by $r(m)\in \overline{\shM}_{B\GG_m^{\dagger}}$
is also $\shU^{\otimes r(m)}$. The statement that the 
morphism is idealized is immediate from the definitions of the ideals.
\end{proof}

\begin{remark}
The morphism of log stacks defined in the above proposition is
not unique. Indeed, any such morphism can be composed with a log automorphism
of an open neighbourhood of $z\times B\GG_m$.
However, the key invariants $N^{A}_{pqr}$ defined below can be seen to be 
independent of this choice by Theorem \ref{thm:NW equals N}.
Further, a priori these invariants depend on the choice of $z$, but
again Theorem \ref{thm:NW equals N} demonstrates the independence
on this choice.
\end{remark}

We are now ready to define the fundamental moduli spaces which yield
the structure constants of our ring.

\begin{definition}
\label{def:Mbetaev}
We define
\begin{align*}
\foM^{\ev}(\shX,\beta,z):= {} & \foM^{\ev}(\shX,\beta)\times_{\scrP(X,r)}
B\GG_m^{\dagger}\\
= {} & \foM(\shX,\beta)\times_{\scrP(\shX,r)} B\GG_m^{\dagger},
\end{align*}
where the morphisms $\foM^{\ev}(\shX,\beta)\rightarrow\scrP(X,r)$,
$\foM(\shX,\beta)\rightarrow\scrP(\shX,r)$ are the evaluation
maps $\ev_{\shX}$ at $x_{\out}$ defined in Lemma~\ref{lem:evaluation diagrams},
the morphism $B\GG_m^{\dagger}\rightarrow\scrP(X,r)$ in the first line is
that given in Proposition~\ref{prop: BGm morphism},
and the morphism $B\GG_m^{\dagger}\rightarrow \scrP(\shX,r)$ in the 
second line
is the composed morphism $B\GG_m^{\dagger}\rightarrow \scrP(X,r)
\rightarrow\scrP(\shX,r)$. The equality follows from 
Lemma~\ref{lem:evaluation diagrams}.
We also similarly set
\[
\scrM(X,\beta,z):=\scrM(X,\beta)\times_{\scrP(X,r)} B\GG_m^{\dagger}.
\]
\end{definition}

\begin{remark}
\label{rem:recall obstruction}
As discussed in \S\ref{subsub:moduli}, we have a perfect
relative obstruction theory for the natural morphism
$\scrM(X,\beta)\rightarrow \foM^{\ev}(\shX,\beta)$. The relative
virtual dimension at a point of $\scrM(X,\beta)$ represented by
a punctured map $f:C^{\circ}\rightarrow X$ is 
$\chi((f^*\Theta_{X/\kk})(-x_{\out}))$, see \eqref{eq:rel virtual dim}. 
\end{remark}

\begin{proposition}
\label{prop:rel virt dim}
There is a perfect relative obstruction theory
for $\scrM(X,\beta,z)\rightarrow \foM^{\ev}(\shX,\beta,z)$ of
relative virtual dimension 
\[
\chi((f^*\Theta_{X/\kk})(-x_{\out}))=A\cdot c_1(\Theta_{X/\kk})
\]
at
a point represented by a punctured map $f:C^{\circ}\rightarrow X$. 
\end{proposition}

\begin{proof}
We have a cartesian diagram (in all categories)
\[
\xymatrix@C=30pt
{
\scrM(X,\beta,z)\ar[r]\ar[d]&\scrM(X,\beta)\ar[d]\\
\foM^{\ev}(\shX,\beta,z)\ar[r]& \foM^{\ev}(\shX,\beta)
}
\]
The perfect relative obstruction theory
for the right-hand vertical arrow of Remark~\ref{rem:recall obstruction} 
pulls back to a perfect relative obstruction theory for the left-hand
vertical arrow. The dimension statement follows from Riemann-Roch.
\end{proof}

\subsection{Virtual smoothness everywhere}
Many of the arguments in this paper ultimately reduce to the idea
that while various moduli spaces of punctured maps to $X$ may be
difficult to control, moduli spaces of punctured maps to the Artin
fan $\shX$ have good local models in the smooth topology as torus invariant 
closed subschemes of toric varieties. The model for this kind of result
is \cite{ACGS18}, Thm.~3.25. We begin with a slight modification of
this result. 

We first define a space $\Mbf(\shX,r)$ which
incorporates both the moduli space of pre-stable rational curves
and the evaluation space $\scrP(\shX,r)$. This space is not important
on its own but serves as
an intermediary space for showing various morphisms are (idealized) 
log \'etale or log smooth. The point here is that $\scrP(\shX,r)$
incorporates a $B\GG_m$ factor which records the data of a $\GG_m$-torsor
on a punctured point. On the other hand, if we are given the punctured point 
as a point in a curve, then this $\GG_m$-torsor is the one corresponding
to the conormal bundle of the point in the curve: see the proof of
Proposition~\ref{prop:evaluation maps general}.
This allows us to remove the $B\GG_m$ factor via a fibre
product, as follows.

\begin{definition}
\label{def:MXr}
Fix $\beta$ a class of curve in $X$, with at least one punctured point
$x_{\out}$ with contact order $-r$, where $r\in\Sigma(X)(\ZZ)$. Assume
that $\beta$ has $n$ additional punctured points.
Let $\Mbf=\Mbf_{0,n+1}$ denote the moduli space
of pre-stable curves of genus $0$ with $n+1$ marked points, labelled
$x_1,\ldots,x_n, x_{\out}$. Give $\Mbf$ its 
basic log structure. There is a morphism $\Mbf\rightarrow
B\GG_m$ induced by the conormal bundle of $x_{\out}$ in the universal curve
over $\Mbf$. Similarly, there is a projection $\scrP(\shX,r)\rightarrow
B\GG_m$. Note here $B\GG_m$ carries the trivial log structure. We then set
\[
\Mbf(\shX,r):=\Mbf\times_{B\GG_m} \scrP(\shX,r).
\]
The pull-back of the log ideal on $\scrP(\shX,r)$ of Remark 
\ref{rem:P idealized} defines a coherent idealized log structure on
$\Mbf(\shX,r)$.
\end{definition}

\begin{proposition}
\label{prop: log etale morphism}
There is a natural idealized log \'etale morphism
\[
\Psi:\foM(\shX,\beta) \rightarrow \Mbf(\shX,r).
\]
Here the idealized log structure on $\foM(\shX,\beta)$ is given by
the puncturing log ideal $\shK$, see \cite{ACGS18}, Defs.~2.49
and 2.55 for the definition.
If the points $x_1,\ldots,x_n$ are in fact marked points, i.e., 
have non-negative contact order,\footnote{We recall the
terminology of marked versus punctured points of 
\S\ref{subsubsec:punctured maps}.} then $\Psi$ is also log \'etale.
\end{proposition}

\begin{proof}
It was shown in \cite{ACGS18}, Thm.~3.25, that 
the forgetful morphism 
$\foM(\shX,\beta)\rightarrow \Mbf$, 
taking a stable punctured map to its domain, is idealized log \'etale. 
Here $\Mbf$ is given the trivial
idealized log structure, i.e., its ideal is the empty ideal. 

Recalling that the Artin fan $\shX$ has an \'etale cover by toric
log stacks $\shA_P$, where $P$ ranges over stalks $\overline\shM_{X,x}$, 
we see that $\shZ=\shZ_r$ has a cover by 
idealized toric log stacks of the form $\shA_{P,I}$ defined in
Definition~\ref{def:APK}, where $P$ ranges
over stalks $\shM_{X,x}$ with $x\in Z$.
Here the ideal $I$ can be written as $\{m\in P\,|\, r(m)>0\}$. Thus
by \cite{ACGS18}, Lem.~B.3, $\shA_{P,I}\rightarrow \Spec\kk=
\shA_{0,\emptyset}$ is idealized log \'etale. So taking 
quotient of the domain and range of this morphism by $\GG_m$, as in 
Remark~\ref{remark: quotient description}, we see that
$\scrP(\shX,r)=[\shZ/\GG_m]\rightarrow B\GG_m$ is idealized log \'etale.
By base-change, so is the projection $\pr_1:\Mbf(\shX,r)\rightarrow \Mbf$.

To obtain the morphism $\foM(\shX,\beta)\rightarrow\Mbf(\shX,r)$,
let $\foC\rightarrow\foM(\shX,\beta)$ be the universal curve.
By the construction of the proof of
Proposition~\ref{prop:evaluation maps general}, the
composition of $\ev_{\shX}:\foM(\shX,\beta)\rightarrow\scrP(\shX,r)$ 
with the projection $\scrP(\shX,r)\rightarrow
B\GG_m$ is the morphism associated to the conormal bundle 
of $x_{\out}(\ul{\foM(\shX,\beta)})$ in $\ul{\foC}$. Since the morphism
$\Mbf\rightarrow B\GG_m$ is also given by the conormal bundle at the
section $x_{\out}$, we obtain the desired morphism $\Psi$,
with the composed morphism
\[
\xymatrix@C=30pt
{
\foM(\shX,\beta)\ar[r]^{\Psi} & \Mbf(\shX,r)\ar[r]^>>>>{\pr_1}& \Mbf
}
\]
being the forgetful morphism that only remembers the domain $\foC$. 

We need to check that $\Psi$ is an idealized morphism, i.e., that 
$\Psi^{\bullet}\pr_2^{\bullet}\shI
\subset\shK$. However, this immediately follows from the definition
of $\shK$ in \cite{ACGS18}, Def.~2.49.

Since both $\pr_1\circ\Psi$ and $\pr_1$ are idealized log \'etale
it immediately follows that $\Psi$ is idealized log \'etale by the
lifting criterion for (idealized log) \'etale morphisms.

If furthermore the point $x_1,\ldots,x_n$ are marked,
i.e., have non-negative contact order,
it follows again from the definition of $\shK$ that 
$\Psi^{\bullet}\pr_2^{\bullet}\shI=\shK$.
This is precisely the notion of ideally strict
of \cite{Ogus}, III,Def.~1.3.2. The last statement then
follows from \cite{Ogus}, IV,Variant~3.1.22.
\end{proof}

\begin{remark}
Note that in general the morphism $\Psi$ is not strict, even when it is
ideally strict. Indeed, $\Mbf(\shX,r)$ carries a product log structure coming
from the universal log structure on $\Mbf$ and the log structure on
$\scrP(\shX,r)$, while $\foM(\shX,\beta)$ carries the basic log structure
for punctured maps.
\end{remark}

\begin{theorem}
\label{thm: everything is log smooth}
The morphisms of Lemma \ref{lem:evaluation diagrams}
\[
\ev_{\shX}:\foM(\shX,\beta)\rightarrow \scrP(\shX,r),\quad
\ev_{\shX}:\foM^{\ev}(\shX,\beta)\rightarrow  \scrP(X,r)
\]
are idealized log smooth.
If the points $x_1,\ldots,x_n$ are in fact marked,
then further both these morphisms are log smooth.
\end{theorem}

\begin{proof}
From Proposition~\ref{prop: log etale morphism}, the morphism
$\foM(\shX,\beta)\rightarrow \Mbf(\shX,r)$ is idealized log \'etale, so
to show that the first morphism $\ev_{\shX}$ in the proposition is
idealized log smooth,
it is enough to check the projection $\Mbf(\shX,r)\rightarrow
\scrP(\shX,r)$ is idealized log smooth. As the morphism is ideally
strict by construction of the idealized log structure on $\Mbf(\shX,r)$,
it is sufficient to show it is log smooth by \cite{Ogus}, IV,Variant~3.1.22.
But this projection is a base-change of the morphism 
$\Mbf\rightarrow B\GG_m$. Base-changing
this latter morphism via the universal torsor $\Spec\kk\rightarrow B\GG_m$
gives the morphism $(\shN^{\vee}_{x_{\out}/\Cbf})^{\times}\rightarrow
\Spec\kk$, where the former denotes the torsor associated to the conormal
bundle of $x_{\out}$ in the universal curve $\Cbf/\Mbf$. This morphism
can be viewed as a composition of a strict smooth morphism
$(\shN^{\vee}_{x_{\out}/\Cbf})^{\times}\rightarrow \Mbf$ and the log smooth
morphism $\Mbf\rightarrow\Spec\kk$. Thus $\Mbf\rightarrow B\GG_m$, hence
$\Mbf(\shX,r)\rightarrow\scrP(\shX,r)$, is log smooth.

The morphism $\ev_{\shX}:\foM^{\ev}(\shX,\beta)\rightarrow\scrP(X,r)$
is then also idealized log smooth as it is obtained
via strict base-change from $\ev_{\shX}:\foM(\shX,\beta)\rightarrow
\scrP(\shX,r)$.

In the case that all other punctured points are in fact marked points,
then the result follows from the above discussion and the last
statement of Proposition~\ref{prop: log etale morphism}.
\end{proof}

\begin{lemma}
\label{lem: more log smoothness}
\begin{enumerate}
\item Using the definition of $\foM^{\ev}(\shX,\beta,z)$
of Definition \ref{def:Mbetaev}, the projection 
$\foM^{\ev}(\shX,\beta,z)\rightarrow B\GG_m^{\dagger}$
at the level of underlying stacks is induced by the conormal bundle of
$x_{\out}$ in the universal curve over $\foM^{\ev}(\shX,\beta,z)$.
\item
Give $\foM^{\ev}(\shX,\beta,z)$ the idealized log structure
coming from its description as a fibre product of idealized log schemes
in Definition \ref{def:Mbetaev}, that is,  the log ideal generated
by the pullbacks of the log ideals on $\foM^{\ev}(\shX,\beta)$ and
$B\GG_m^{\dagger}$. Then the
projection $\foM^{\ev}(\shX,\beta,z)\rightarrow B\GG_m^{\dagger}$
is idealized log smooth. If the points $x_1,\ldots,x_n$ are in fact marked,
then this morphism is log smooth.
\end{enumerate}
\end{lemma}

\begin{proof}
The first item follows from the construction of the proof of 
Proposition~\ref{prop:evaluation maps general}, while the
second item is
immediate by base-change from Theorem~\ref{thm: everything is log smooth}.
\end{proof}

\begin{lemma}
\label{lemma: moduli morphisms}
Let $\overline{\scrM}_{0,n+1}$ denote the Deligne-Mumford moduli space
of stable genus $0$ curves with $n+1$ marked points, with its universal
log structure, which is the divisorial log structure associated to the
boundary of $\overline{\scrM}_{0,n+1}$. Then 
there are natural morphisms
\begin{align*}
\Phi: \scrM(X,\beta,z) \rightarrow & \overline{\scrM}_{0,n+1}\times 
B\GG_m^{\dagger}=:\overline{\scrM}_{0,n+1}^{\dagger}\\
\Phi: \foM^{\ev}(\shX,\beta,z) \rightarrow & \overline{\scrM}_{0,n+1}^{\dagger}
\end{align*}
If the points $x_1,\ldots,x_n$ are in fact marked,
then the second morphism is log smooth.
\end{lemma}

\begin{proof}
There is of course a morphism $\scrM(X,\beta,z)\rightarrow 
\overline{\scrM}_{0,n+1}$ given by forgetting the morphism and
stabilizing the domain curve. Combining this with the projection
$\scrM(X,\beta,z)\rightarrow B\GG_m^{\dagger}$ gives the first desired
morphism. The second morphism is similar.

To show the second morphism $\Phi$ is log smooth if $x_1,\ldots,x_n$
are marked, note it is a composition
of morphisms 
\begin{equation}
\label{threemorphismfactorization}
\foM^{\ev}(\shX,\beta,z)=\foM(\shX,\beta)\times_{\scrP(\shX,r)}
B\GG_m^{\dagger}\rightarrow \Mbf\times_{B\GG_m} B\GG_m^{\dagger}
\rightarrow\Mbf\times B\GG_m^{\dagger}
\rightarrow \overline{\scrM}_{0,n+1}^{\dagger}.
\end{equation}
Here, the first morphism is a base-change of the morphism
$\foM(\shX,\beta)\rightarrow \Mbf(\shX,r)=\Mbf\times_{B\GG_m} \scrP(\shX,r)$
by the morphism $\Mbf\times_{B\GG_m} B\GG_m^{\dagger}\rightarrow
\Mbf(\shX,r)$ induced by $B\GG_m^{\dagger}\rightarrow \scrP(\shX,r)$. 
Since $\foM(\shX,\beta)\rightarrow \Mbf(\shX,r)$
is log \'etale by Proposition~\ref{prop: log etale morphism},
the first morphism of \eqref{threemorphismfactorization} is log \'etale.
The second morphism fits into a
cartesian diagram (in both the ordinary and fs log categories)
\begin{equation*}
\label{diagonalbasechange}
\xymatrix@C=15pt
{{\bf M}\times_{B\GG_m} B\GG_m^{\dagger}\ar[r]\ar[d] &
B\GG_m\ar[d]^{\Delta}\\
{\bf M}\times_{\Spec\kk}B\GG_m^{\dagger} \ar[r] & B\GG_m\times_{\Spec\kk}
B\GG_m
}
\end{equation*}
Since the diagonal $\Delta$ is smooth and strict (here the log structures
are trivial), the same is true of the left-hand vertical arrow.
Thus the second morphism of \eqref{threemorphismfactorization} 
is strict log 
smooth. The third morphism is induced by the stabilization morphism
$\Mbf\rightarrow \overline{\scrM}_{0,n+1}$, which is log \'etale. Hence
the composition $\Phi$ of the three morphisms is log smooth.
\end{proof}

\subsection{The invariants}
\label{subsec:the invariants}

We are now ready to define the key invariants. For the notion
of log fibre dimension, see \S\ref{sec:log fibre dim}.

\begin{proposition}
\label{prop:the invariants}
Let $\beta$ be a class of punctured map with one point $x_{\out}$
with contact order $-r$ and $n$ marked points. Then
\begin{enumerate}
\item The morphism $\Phi:\foM^{\ev}(\shX,\beta,z)\rightarrow
\overline\scrM_{0,n+1}^{\dagger}$ is of log fibre dimension $1$, and
the projection $\foM^{\ev}(\shX,\beta,z)\rightarrow B\GG_m^{\dagger}$
is log fibre dimension $n-1$.
\item
$\foM^{\ev}(\shX,\beta,z)$ is pure-dimensional of dimension $n-2$, and
hence $\scrM(X,\beta,z)$ carries a virtual fundamental class of dimension
\[
A\cdot c_1(\Theta_{X/\kk}) + n - 2.
\]
\end{enumerate}
\end{proposition}

\begin{proof}
(1) For the first morphism, we have the factorization
\eqref{threemorphismfactorization}, with the first and third morphisms
log \'etale and the second smooth of relative dimension one,
and hence by properties of log fibre dimension, Propositions
\ref{prop: flat fibre dim},(3) and \ref{prop:log etale dimensions},
we obtain $\Phi$ is log fibre dimension one.

For the second morphism, we further compose with the projection
$\overline\scrM_{0,n+1}^{\dagger}\rightarrow B\GG_m^{\dagger}$, 
which is of log fibre dimension $n-2$, giving log fibre dimension
$n-1$ for the second morphism.

(2) By Proposition~\ref{prop:rel virt dim}, it is sufficient
to prove the dimension statement for $\foM^{\ev}(\shX,\beta,z)$.
First, as $\foM^{\ev}(\shX,\beta,z)
\rightarrow B\GG_m^{\dagger}$ is log smooth by Lemma 
\ref{lem: more log smoothness}, it is log flat by Proposition
\ref{prop: flatness sorites}.
Further, 
as the ghost sheaf on $B\GG_m^{\dagger}$ is $\NN$,
this morphism is also integral by \cite{Ogus}, I,Prop.~4.6.3, 5.
So the fibre dimension and the log fibre dimension of this
morphism coincide, see Proposition~\ref{prop: flat fibre dim}, (2).
Thus the
fibre dimension of $\foM^{\ev}(\shX,\beta,z)\rightarrow B\GG_m^{\dagger}$
is $n-1$, and hence $\foM^{\ev}(\shX,\beta,z)$ is pure-dimensional of
dimension $n-2$.
\end{proof}

Recall from Basic Setup~\ref{setup1} that we do not assume $X$ is projective
over $\Spec\kk$, but only over $S$. Note that if $S$ is one-dimensional,
we do not work with the moduli space of relative punctured maps,
$\scrM(X/S,\beta)$, but with the moduli space $\scrM(X,\beta)$.
This is necessary as the contact orders we wish to make use
of do not define contact orders over $S$. In any event, $\scrM(X,\beta)$
would not be expected to be proper since $S$ is not proper. To extract
a number, we need a proper moduli space, and fortunately once we
put on a point constraint, the moduli space becomes proper:

\begin{lemma}
\label{lemma: properness}
In either the absolute or relative cases in the sense of Basic
Setup~\ref{setup1}, with $\dim \ul{S}=0$ or $1$ 
respectively, we have:
\begin{enumerate}
\item There is a morphism $\ul{\scrM(X,\beta)}\rightarrow \scrM(\ul{X}/\ul{S},
A)$ where $A$ is the underlying curve class of $\beta$. Furthermore, the
composed map
\[
\ul{\scrM(X,\beta)}\rightarrow \scrM(\ul{X}/\ul{S}, A)
\rightarrow \ul{S}
\]
is proper.
\item $\scrM(X,\beta,z)$ is proper over $\Spec \kk$.
\end{enumerate}
\end{lemma}

\begin{proof}
There is always a forgetful morphism
$\ul{\scrM(X,\beta)}\rightarrow \scrM(\ul{X},A)$, which is shown
to be proper in \cite{ACGS18}, Thm.~3.18.
Note that since $\ul{S}$ is assumed to be non-complete in the relative case,
the composition of any stable map $\ul{f}:\ul{C}/\ul{W}\rightarrow \ul{X}$ 
with $\ul{X}\rightarrow\ul{S}$ must be constant on fibres of $\ul{C}\rightarrow
\ul{W}$: see the argument of the proof of \cite{ACGS18}, Thm.~5.13. 
Thus
$\scrM(\ul{X},A)=\scrM(\ul{X}/\ul{S},A)$.
As the structure morphism
$\scrM(\ul{X}/\ul{S},A)\rightarrow\ul{S}$ is proper, 
given that $X\rightarrow S$ is projective, (1) follows.

For (2), in the absolute case, $\scrM(X,\beta)$ is proper
over $\Spec\kk$. As $B\GG_m$ is a closed substack of
$\ul{\scrP(X,r)}$, the ordinary fibre product 
$\ul{\scrM(X,\beta)}\times_{\ul{\scrP(X,r)}} B\GG_m$
is a closed substack of $\ul{\scrM(X,\beta)}$, hence proper. 
By Remark~\ref{rem:sat finite representable}, 
the fs fibre product defining $\scrM(X,\beta,z)$
is finite over this ordinary fibre product, hence is proper. 

In the relative case, one notes that
if $s=g(z)\in S$, then any family of stable maps $f:C/W\rightarrow X$ in
$\scrM(X,\beta,z)$ necessarily factors through $X_s$, and hence
$\scrM(X,\beta,z)=\scrM(X_s,\beta,z)$. This is again proper as in the absolute case.
\end{proof}

\begin{definition}
\label{def:Nbeta}
Let $p_1,p_2,r\in\Sigma(X)(\ZZ)$, and let $\beta$ be a class of
punctured map with underlying curve class $A$ and
three punctured points $x_1,x_2,x_{\out}$ with contact orders
$p_1, p_2$ and $-r$ respectively. We define
\[
N^A_{p_1p_2r} := \begin{cases}
\deg [\scrM(X,\beta,z)]^{\virt} & A\cdot c_1(\Theta_{X/\kk})=0\\
0 & \hbox{otherwise}.
\end{cases}
\] 
\end{definition}

We end this subsection by completing the calculation of these invariants
in the constant map case:

\begin{proof}[Proof of Lemma~\ref{lem:constant-maps}]
Let $\sigma_1,\sigma_2,\sigma\in\P$ be the minimal cones containing
$p_1,p_2$ and $r$ respectively. First, suppose given a 
basic stable punctured map
\[
f:(C^{\circ},x_1,x_2,x_{\out})/\Spec(Q\rightarrow\kk)\rightarrow X
\]
of class $\beta$ with $\ul{f}(x_{\out})\in Z_r^{\circ}$. Since $A=0$
the underlying morphism $\ul{f}$ is constant, and stability then implies
that $\ul{C}\cong \PP^1$. Necessarily 
$f(x_i)\in Z_{p_i}$,
$i=1,2$. Thus we have
$Z_r^{\circ}\subseteq Z_{p_1}\cap Z_{p_2}$, which is only possible if
$\sigma\supseteq \sigma_1,\sigma_2$. So $\sigma$ contains $p_1,p_2$.
It then immediately follows from the balancing condition \cite{ACGS18},
Prop.~2.27, that $p_1+p_2-r=0$, i.e., $p_1+p_2=r$. Note this balancing
can also be immediately deduced from Proposition~\ref{intersectionnumbers}. 
This proves vanishing
of $N^{0}_{p_1p_2r}$ unless this condition holds.

Now suppose $p_1,p_2,r$ lie in a common cone with $p_1+p_2=r$. We need
to show $N^0_{p_1p_2r}=1$. To do so we will show $\scrM(X,\beta,z)$
is a reduced point.

\begin{claim}
Given $z\in Z_r^{\circ}$, there
is a unique basic punctured map over $\Spec(Q\rightarrow\kk)$ (for $Q$ the
basic monoid)
of class $\beta$ with image $z$.
\end{claim}

\begin{proof}
Fix
$\ul{f}:\ul{C}=\PP^1/\Spec\kk\rightarrow \ul{X}$ to be the constant map with 
image $z$.
If $\ul{f}$ has an enhancement to a punctured map of class $\beta$,
it is immediate from the definition of the basic monoid
\cite{ACGS18}, \S2.4, that the basic monoid is $Q=P$, where
$P=\overline{\shM}_{X,z}$. Let $C$ denote the log enhancement of $\ul{C}$
to a log smooth curve over $\Spec(Q\rightarrow\kk)$ with $x_1,x_2,x_{\out}$
marked points in the log structure. We first need to determine
the only possibility for a puncturing $C^{\circ}$ of $C$. Indeed,
if there is an enhancement of $\ul{f}$ to a stable punctured map
over $\Spec(Q\rightarrow\kk)$, then in fact 
$\bar f^{\flat}:
f^{-1}\overline{\shM}_X\rightarrow \overline{\shM}_{C^{\circ}}$ is completely
determined by the contact orders at the three marked points, as follows. Write 
$\Gamma(C,\overline{\shM}_{C^{\circ}})\subset Q\oplus \NN_{x_1} \oplus \NN_{x_2}
\oplus \ZZ_{x_{\out}}$. Here the subscripts indicate the source of the
particular factor. Then by the definition of contact order \eqref{eq:ux def},
\[
\bar f^{\flat}(m)=(m, \langle p_1, m\rangle,
\langle p_2,m\rangle, - \langle r,m\rangle),
\]
for $m\in P$. 
In particular, using pre-stability of the punctured morphism
$f:C^{\circ}\rightarrow X$, we can
define the punctured log structure by choosing $\overline{\shM}_{C^{\circ}}
\subset \overline{\shM}_C^{\gp}$ so that its stalk at $x_{\out}$ is generated
by $Q\oplus\NN$ and $\{(m,-\langle r, m\rangle) \,|\, m\in P\}$, while
$\overline{\shM}_{C^{\circ}}$ and $\overline{\shM}_C$ agree away from
$x_{\out}$. We
then take $\shM_{C^{\circ}}=\overline{\shM}_{C^{\circ}}
\times_{\overline{\shM}_C^{\gp}} \shM_C^{\gp}$. Finally we define
$\alpha_{C^{\circ}}:\shM_{C^{\circ}}\rightarrow\O_C$ to agree with the structure
map $\alpha_C:\shM_C
\rightarrow \O_C$ on $\shM_C\subset \shM_C^{\circ}$ and to take the value
$0$ on $\shM_{C^{\circ}}\setminus\shM_C$. This defines the only
possible punctured curve $C^{\circ}/\Spec(Q\rightarrow\kk)$ which can
appear as a domain of a stable punctured map.

We now need to show
that up to isomorphism, there is a unique lifting of $\bar f^{\flat}$
to $f^{\flat}:f^{-1}\shM_X\rightarrow \shM_{C^{\circ}}$. Note that the torsor
$f^*\shL^{\times}_m$ (corresponding to $m\in P$) is necessarily trivial as
$\ul{f}$ is constant, and
the balancing condition is equivalent to the torsor 
$\shL^{\times}_{\bar f^{\flat}(m)}$
being trivial for all $m\in P$. Thus it is sufficient to choose a collection of
identifications $f^*\shL^{\times}_m\cong \shL^{\times}_{\bar f^{\flat}(m)}$ 
for $m\in Q^{\gp}$
compatible with multiplication in $f^*\shM_X$ and $\shM_{C^{\circ}}$. 
It is clear that such a set of identifications
can be chosen, and any two choices are related by a choice of element
in $\Hom(Q^{\gp},\kk^{\times})$. Thus in particular, any two such choices
of punctured map are isomorphic over an automorphism of 
$\Spec(Q\rightarrow\kk)$, proving the claim. 
\end{proof}

The above claim is insufficient to prove the result, as there may be
infinitesimal deformations. 
But note that as $f^*\Theta_{X/\kk}$ is trivial,
$H^1(C,f^*\Theta_{X/\kk}(-x_{\out}))=0$. This tells us that the
morphism $\scrM(X,\beta)\rightarrow \foM^{\ev}(\shX,\beta)$ is 
smooth. Moreover, the morphism $\ev_{\shX}:\foM^{\ev}(\shX,\beta)\rightarrow 
\scrP(X,r)$ is log smooth by Theorem~\ref{thm: everything is log smooth}.
Thus the composition $\scrM(X,\beta)\rightarrow\scrP(X,r)$ is log
smooth also.
From the above calculation of the basic monoid, this morphism 
is also strict over the locus $Z_r^{\circ}\times B\GG_m\subset\scrP(X,r)$. 
A strict log smooth morphism is also smooth, and hence
$\scrM(X,\beta)$ is smooth over $\scrP(X,r)$, at least over
$Z_r^{\circ}\times B\GG_m$. In turn, by base-change, $\scrM(X,\beta,z)$
is smooth over $B\GG_m$. This implies $\scrM(X,\beta,z)$ consists
of one non-singular point, with trivial stabilizer by \cite{ACGS18},
Prop.~2.41. We conclude that $N^0_{p_1p_2r}=1$.
\end{proof}

\subsection{The tropical interpretation of point constraints} 
Note that the moduli spaces $\scrM(X,\beta,z)$, $\foM^{\ev}(\shX,\beta,z)$
(Definition~\ref{def:Mbetaev})
carry a log structure coming from the fs fibre product, and thus
it is not the basic log structure associated to the family of
punctured maps. It is useful to
interpret these log structures tropically. 

In particular, let $\bar w\in \scrM(X,\beta,z)$ be a geometric
point with a corresponding
punctured map $f:C^{\circ}_{\bar w}\rightarrow X$, and images
$\bar y\in \scrM(X,\beta)$ and $\bar s\in \scrP(X,r)$. Note that 
the tropicalization of $\ev_X$ 
(see Definition~\ref{def:evaluation maps})
then induces a map $\sigma_{\bar y}\rightarrow \sigma_{\bar s}$,
with the cones $\sigma_{\bar y},\sigma_{\bar s}$ as
defined in \S\ref{sec:tropicalization}. This can be described as follows.

\begin{lemma}
\label{lem:tropical point constraint}
Let $G$ be the dual graph of $C^{\circ}_{\bar w}$ and $v_{\out}$ the
vertex of $G$ corresponding to the irreducible component of 
$C^{\circ}_{\bar w}$ containing the point $x_{\out}$. Let $h:\Gamma
\rightarrow \Sigma(X)$ be the induced universal family of tropical maps,
defined over $\sigma_{\bar y}$, i.e., defined using the basic log 
structure. Then 
$\Sigma(\ev_X):\sigma_{\bar y}\rightarrow \sigma_{\bar s}$
is given, for $m\in\sigma_{\bar y}$, by
\[
m\mapsto h_m(v_{\out}),
\]
where $h_m$ is as in \eqref{eq:hs def}.
Further, if $r\not=0$, then
$\sigma_{\bar w}$ can be identified with
\[
\{ m \in \sigma_{\bar y} \,|\, h_m(v_{\out})\in \RR_{\ge 0} r \subset 
\sigma_{\bar s}\}
\]
and the projection $\scrM(X,\beta,z)\rightarrow B\GG_m^{\dagger}$
tropicalizes to
\[
\delta:\sigma_{\bar w}\rightarrow \RR_{\ge 0}
\]
with $m\in \sigma_{\bar w}$ mapping to the number $\delta(m)$ such
that 
\[
h_m(v_{\out})=\delta(m) r.
\]
If $r=0$, then $\sigma_{\bar w}=\sigma_{\bar y}\times \RR_{\ge 0}$,
with projection to $\Sigma(B\GG_m^{\dagger})$ being projection
onto the second factor.
\end{lemma}

\begin{proof}
The description of $\Sigma(\ev_X)$ follows immediately from dualizing
\eqref{eq:evghost} and the description of tropical maps associated
to punctured maps given in \cite{ACGS18}, \S2.2.

On the other hand, the tropicalization of $B\GG_m^{\dagger}\rightarrow
\scrP(X,r)$ is clearly the map $\RR_{\ge 0}\rightarrow \sigma_{\bar s}$
given by $1\mapsto r$. The remaining statements then follow immediately
from Proposition~\ref{tropicalproduct}. 
\end{proof}

\begin{remark}
The imposition of the tropical constraint that $h_m(v_{\out})$ should
lie on the ray $\RR_{\ge 0} r$ may be understood from the point of view
of toric blowups of the target space $X$. As explained at the beginning of
the current section, if $\RR_{\ge 0}r$ is a ray of $\Sigma(X)$,
the number $N^A_{p_1p_2r}$ can be viewed as a count of curves in 
$\scrM(X,\beta)$ satisfying a schematic point constraint, up to a tropically
computed factor. For such a curve, it is automatic that $h_m(v_{\out})
\in\RR_{\ge 0}r$ since the chosen constraint point lies in $Z_r^{\circ}$.

If $\RR_{\ge 0}r$ is not a ray of $\Sigma(X)$, we can imagine performing
a toric blow-up $\widetilde X\rightarrow X$, leading to a refinement 
$\Sigma(\widetilde X)$ of $\Sigma(X)$ for which $\RR_{\ge 0}r$ is a ray.
The tropical constraint $h_m(v_{\out})\in\RR_{\ge 0}r$ then can be seen
as realising this birational modification within logarithmic geometry. 

In any event, if $r$ does not lie on a ray of $\Sigma(X)$, then
$\scrM(X,\beta)$ may not be the correct dimension or even pure-dimensional.
The imposition of the tropical constraint can also be viewed as a
transverse base-change, yielding a moduli space of the correct dimension
as in Remark~\ref{rem:transverse}.
\end{remark}

\begin{definition}
\label{def:delta}
We use the notation $\delta$ interchangeably for:
\begin{enumerate}
\item $\delta\in \overline{\shM}_{B\GG_m^{\dagger}}$ the generator;
\item $\delta\in \Gamma(\scrM(X,\beta,z),
\overline{\shM}_{\scrM(X,\beta,z)})$ the image of $\delta$ defined in (1)
under the map of ghost sheaves induced by
the projection $\scrM(X,\beta,z)\rightarrow B\GG_m^{\dagger}$;
\item The tropicalization $\delta:\Sigma(\scrM(X,\beta,z))\rightarrow \Sigma(B\GG_m^{\dagger})=
\RR_{\ge 0}$ of the projection $\scrM(X,\beta,z)\rightarrow B\GG_m^{\dagger}$.
If one views a section of the ghost sheaf on a log stack $Y$
as inducing a function $\Sigma(Y)\rightarrow\RR_{\ge 0}$, then this
is the function induced by $\delta$ as defined in (2). 
By Lemma 
\ref{lem:tropical point constraint}, if $r\not=0$, then
$\delta$ is given on a cone of $\Sigma(\scrM(X,\beta,z))$ by 
\begin{equation}
\label{eq:voutdelta}
h_m(v_{\out})=\delta(m) r.
\end{equation}
\end{enumerate}
In cases (2) and (3), we also use the same notation $\delta$
with the stack $\scrM(X,\beta,z)$ replaced by the stack
$\foM^{\ev}(\shX,\beta,z)$.
\end{definition}

More generally, the notation $\delta$ will appear in many different
places in the sequel, and will always denote either a morphism
to $B\GG_m^{\dagger}$ playing the same role as above, or denote
the element of a stalk of a ghost sheaf measuring the distance of
$v_{\out}$ from the origin, as indicated in \eqref{eq:voutdelta}.

In Figure~\ref{fig:delta}, a general
tropical map coming from tropicalizing a point in $\scrM(X,\beta)$
is depicted on the left: in this case $v_{\out}$ may fall anywhere
in a cone containing $r$, e.g., the first quadrant. On the right, 
a general tropical map coming from tropicalizing a point in
$\scrM(X,\beta,z)$ is depicted. In this case $v_{\out}$ is constrained
to lie on the dotted line through $r$, and $\delta$ measures where
$v_{\out}$ falls on this line.

\begin{figure}
\input{Figuredelta.pspdftex}
\caption{A tropical map without and with point constraints.
We note in all depictions of tropical punctured maps,
we only draw the legs corresponding to punctured
points as arrows. By pre-stability, in fact these legs should
extend as far as possible in the cone of $\Sigma(X)$ containing
the leg, but this will make the figures harder to parse.}
\label{fig:delta}
\end{figure}

\section{Sketch and first steps of the proof of the main theorem}
\label{section:sketch}

This section will break the proof of associativity and
hence the main theorems, Theorems~\ref{mainassociativity1} and
\ref{mainassociativity2}, into a number of 
steps, and illustrate why these steps are necessary, partly by way of
some simple examples. This section should be viewed as both a reader's
guide to the rest of the paper and the start of the actual proof of 
associativity.

We first divide the requirements of Theorems~\ref{mainassociativity1}
and \ref{mainassociativity2} into four pieces. 

\newpage

\subsection{Finiteness of the product rule}

We first show that \eqref{structureequation} or \eqref{structureequation2}
is always a finite sum:

\begin{lemma}
\label{lem:finite}
With the hypotheses of either Theorem~\ref{mainassociativity1} or Theorem
\ref{mainassociativity2}, with $(B,\P)=(|\Sigma(X)|,\Sigma(X))$ in the first
case, fix $p_1,p_2\in B(\ZZ)$. Then all but a finite number of
the $N^{A}_{p_1p_2r}$ for $A\in P\setminus I$, $r\in B(\ZZ)$
are zero. 
\end{lemma}

\begin{proof}
Since $I$ is assumed to be co-Artinian, there are only a finite number
of classes $A\in P\setminus I$. Thus we may fix a curve class
$A$, and only need to show there are only a finite number of choices
of $r$ such that $N^A_{p_1p_2r}\not=0$. However, this follows
immediately from Remark~\ref{rem:contact order determined}.
\end{proof}

\subsection{Commutativity}

We simply state the completely obvious lemma, which follows immediately from the
fact that the definition of the numbers $N^A_{p_1p_2r}$ is
independent of the ordering of the input points $x_1,x_2$.

\begin{lemma} 
$N^A_{p_1p_2r}=N^A_{p_2p_1r}$, hence
the multiplication operation defined using these structure constants
is commutative.
\end{lemma}

\subsection{$\vartheta_0$ is the unit element}
The fact that $\vartheta_0$ is the unit in the ring is rephrased in terms
of the numbers $N^A_{p_1p_2r}$ as follows:

\begin{theorem}
\label{thm:unit}
With the hypotheses of either Theorem~\ref{mainassociativity1} or Theorem
\ref{mainassociativity2}, with $(B,\P)=(|\Sigma(X)|,\Sigma(X))$ in the first
case, we have for $A\in P\setminus I$, $p\in B(\ZZ)$,
\[
N^A_{0pr}=\begin{cases} 0 & \hbox{$A\not=0$ or $p\not=r$,}\\
1 & A=0, r=p.
\end{cases}
\]
In particular $\vartheta_0$ is the unit in $R_I$.
\end{theorem}

The proof of this theorem is given in \S\ref{sec:forgetful}.

\subsection{Associativity: Overview and first steps}

We will rephrase the main associativity result, and then sketch
some of the main ideas of the proof while breaking the argument into
smaller pieces, which will be proved in subsequent sections.

\subsubsection{Restatement of associativity}

Associativity is implied by the following theorem,
by expanding out the products 
$(\vartheta_{p_1}\cdot\vartheta_{p_2})\cdot \vartheta_{p_3}$ 
and 
$\vartheta_{p_1}\cdot(\vartheta_{p_2}\cdot \vartheta_{p_3})$.

\begin{theorem}
\label{thm:restatement}
With the hypotheses of either Theorem~\ref{mainassociativity1} or Theorem
\ref{mainassociativity2}, with $(B,\P)=(|\Sigma(X)|,\Sigma(X))$ in the first
case, we have for each $A\in P\setminus I$, $p_1,p_2,p_3,r\in B(\ZZ)$,
\begin{equation}
\label{asseq}
\sum_{A_1,A_2,s
\atop A_1+A_2=A}
N^{A_1}_{p_1p_2s} N^{A_2}_{sp_3r}
=
\sum_{A_1,A_2,s
\atop A_1+A_2=A}
N^{A_1}_{p_1sr} N^{A_2}_{p_2p_3s}.
\end{equation}
Here we sum over $s\in B(\ZZ)$.
\end{theorem} 

So it will be sufficient to prove the above theorem to complete
the proofs of Theorems~\ref{mainassociativity1} and
\ref{mainassociativity2}.

\subsubsection{Outline}
We are proving a logarithmic analogue of associativity of ordinary quantum
cohomology. Recall that in the ordinary case associativity follows by computing
either one of $(\alpha\cdot\beta)\cdot\gamma$ and $\alpha\cdot(\beta\cdot\gamma)$ by
a sum over four-point invariants. This is done in two steps. First, one uses
virtual smoothness of the domain stabilization map $\scrM_{0,4}(X,A)\rightarrow
\overline\scrM_{0,4}\cong\PP^1$ to define a four-point invariant with fixed
stabilized domain, independently of the modulus. There are three special points
$y_i\in \ol\scrM_{0,4}$ corresponding to reducible stable domains distinguished
by the distribution of the four points into two pairs of points. 
In a second
step, one shows that the stable map moduli space over each such reducible domain
is obtained from three-point moduli spaces via matching the images in $X$ of a
pair of marked points to obtain a nodal domain:
\begin{equation}
\label{eq:ordinary gluing}
\scrM_{0,4}(X,A)\times_{\overline
\scrM_{0,4}} y_1 \cong \coprod_{A_1+A_2=A} \scrM_{0,3}(X,A_1)\times_X
\scrM_{0,3}(X,A_2).
\end{equation}
One can then use standard intersection theory techniques to relate the virtual
fundamental class of the left-hand side to the virtual fundamental classes of
the $\scrM_{0,3}(X,A_i)$. This shows that the four-point function computed by
the left-hand side of \eqref{eq:ordinary gluing} equals one of the two ways of
associating the quantum product of three classes. Performing this argument at
two of the $y_i$ then proves associativity of quantum cohomology.

We will follow the same general strategy, but everything becomes more
difficult. To prove Theorem~\ref{thm:restatement}, one starts with
a class of punctured log curve $\beta$ given by the curve class $A$
and four points with contact orders $p_1,p_2,p_3$ and $-r$, again
labelled $x_1,x_2,x_3$ and $x_{\out}$. We place a point constraint
at $x_{\out}$, obtaining moduli spaces $\scrM(X,\beta,z)
\rightarrow\foM^{\ev}(\shX,\beta,z)$ as in Definition~\ref{def:Mbetaev},
with a relative perfect obstruction theory. Statements such as virtual
smoothness of a morphism from $\scrM(X,\beta,z)$ are statements about
actual smoothness for a morphism from $\foM^{\ev}(\shX,\beta,z)$.

For example, there is a morphism $\foM^{\ev}(\shX,\beta,z)
\rightarrow\overline\scrM_{0,4}$, even at the logarithmic level, as in
the proof of Lemma~\ref{lemma: moduli morphisms}. However, this morphism
is not well-behaved. It is certainly seldom smooth, and is usually not
log smooth. Thus there is little one can do. On the other hand,
the morphism $\Phi:\foM^{\ev}(\shX,\beta,z)\rightarrow
\overline\scrM^{\dagger}_{0,4}$ of Lemma~\ref{lemma: moduli morphisms}
is log smooth. This is the utility of adding the $B\GG_m^{\dagger}$
factor. However, log smoothness does not imply flatness.
In particular, the fibres of $\Phi$ may jump in dimension over the
boundary points of $\overline\scrM_{0,4}$, leading to the
fibres of $\scrM(X,\beta,z)\rightarrow\overline\scrM^{\dagger}_{0,4}$
having a jump in virtual dimension. In \S\ref{subsubsec:extended}, we 
give an example to show this happens even for 
target spaces $X$ with very simple geometry,
and then use this example to show how our approach resolves these
difficulties.

The key idea, used over and over again, is the trick of using transverse
base-changes to obtain virtual dimension zero fibres. This is a generalization
of the idea in Remark~\ref{rem:transverse}. Also very important is the point
that the spaces $\foM^{\ev}(\shX,\beta,z)$ are controlled almost entirely
by tropical geometry: this is the benefit of log smoothness. Thus we
use tropical geometry throughout to control the pathologies of these moduli
spaces.

The second difficulty following the classical path is gluing. Gluing
stable log maps is a great deal harder than gluing ordinary stable maps.
Instead of the ordinary fibre product of \eqref{eq:ordinary gluing}, 
one instead needs to use a fibre product in the fs log category, and this
does not play well with ordinary intersection theory. Worse, again
the glued moduli space is often of the wrong virtual dimension, or even
not of pure virtual dimension. The ideas of Remark~\ref{rem:transverse}
again come to the rescue, but we need to wrestle with extreme technical
difficulties along the way.

We now go into more detail, beginning first with our promised example
demonstrating these pathologies.

\subsubsection{An extended example}
\label{subsubsec:extended}

\begin{example}
\label{ex:extended1}
We consider a sample associativity statement in the case of Examples 
\ref{runningexample1}, (1).
Consider the product
$\vartheta_{p_1}\vartheta_{p_2}\vartheta_{p_3}$ where
\[
p_1=v_1+2v_2,\quad p_2=v_2, \quad p_3=2v_3.
\]
By the observations stated in Example~\ref{runningexample2}, we
can express the corresponding theta functions as monomials in the 
$\vartheta_{v_i}$:
\[
\vartheta_{p_1}=\vartheta_{v_1}\vartheta_{v_2}^2, \quad
\vartheta_{p_2}=\vartheta_{v_2},\quad \vartheta_{p_3}=\vartheta_{v_3}^2.
\]
Thus, assuming that the ring $R_I$ is as described in \eqref{eq:first RI desc}
(in particular assuming the product is associative), we can calculate
the product of the $\vartheta_{p_i}$'s as follows, using the given
relations stated in Example~\ref{runningexample2},
\begin{equation}
\label{Eqn: triple product in first example}
\begin{aligned}
\vartheta_{p_1}\vartheta_{p_2}\vartheta_{p_3}\stackrel{\eqref{eq:monomial product}}{=} {} & 
(\vartheta_{v_1}\vartheta_{v_2} \vartheta_{v_3})\vartheta_{v_2}^2
\vartheta_{v_3}\\
\stackrel{\eqref{eq:first RI desc}}{=} {} &(t^{[L]}+t^{[L-E]}\vartheta_{v_1})\vartheta_{v_2}^2
\vartheta_{v_3}\\
\stackrel{\eqref{eq:thetav2v3 product}}{=} {} & 
t^{[L]}(\vartheta_{v_2+v_3}+t^{[L-E]})\vartheta_{v_2}
+t^{[L-E]}(\vartheta_{v_1}
\vartheta_{v_2}\vartheta_{v_3})\vartheta_{v_2}\\
\stackrel{\eqref{eq:first RI desc}}{=} {} & 
t^{[L]}(\vartheta_{v_2+v_3}+t^{[L-E]})\vartheta_{v_2}
+t^{[L-E]}(t^{[L]}+t^{[L-E]}\vartheta_{v_1})\vartheta_{v_2}\\
\stackrel{\eqref{eq:monomial product}}{=} {} & t^{[L]}\vartheta_{2v_2+v_3}+2t^{[2L-E]}\vartheta_{v_2}+t^{[2L-2E]}
\vartheta_{v_1+v_2}.
\end{aligned}
\end{equation}
We would like to study the term $2t^{[2L-E]}\vartheta_{v_2}$ and show
how the geometry of a moduli space of punctured maps exhibits the fact
that this term appears in both ways of associating the product.

To do so, set $A=2L-E$ and $r=v_2$.
This gives a class of four-pointed punctured maps $\beta$, with contact
orders $p_1, p_2, p_3$ and $-r$ at the four points. Fix a point
$z\in Z_r^{\circ}=D_2^{\circ}$. We
describe the moduli space $\scrM(X,\beta,z)$, which will in fact be
unobstructed and hence one-dimensional by Proposition~\ref{prop:rel virt dim}.
We omit the detailed analysis of this moduli space, leaving it to the
reader to use the contact orders to verify the description we give here.
Despite this moduli space being unobstructed,
it is reducible and non-reduced. Indeed, unobstructedness
means that the morphism $\scrM(X,\beta,z)\rightarrow\foM^{\ev}(\shX,\beta,z)$
is smooth, and hence $\scrM(X,\beta,z)$ inherits pathologies from
$\foM^{\ev}(\shX,\beta,z)$. Shortly, we will examine the structure of
the latter moduli space. However, one finds that
the underlying reduced stack $\scrM(X,\beta,z)_{\red}$
of $\scrM(X,\beta,z)$ has two 
irreducible components isomorphic
to $\PP^1$ attached at a single point. Call these 
reduced irreducible components $\scrM_1$
and $\scrM_2$. 

The general curve parameterized by $\scrM_1$ is as depicted on the
left in Figure~\ref{Figure1}, with its tropicalization  depicted on
the right. The image of the map is the union of $D_2$ and the strict transform
of the line through the center of the blowup and $D_2\cap D_3$.
The $\PP^1$ worth of moduli is given by the location of the
point $x_2$. There are three ways this curve can degenerate: $x_2$
can either move towards the node of the domain contained in the same component
as $x_2$; $x_2$ can move to $x_{\out}$; and $x_2$ can move to $x_1$. In
Figures~\ref{Figure2} and \ref{Figure3} we draw the first and third of
these possibilities; the second possibility will not play an important role
in this discussion. 

The moduli space $\scrM(X,\beta,z)$ is in fact generically non-reduced along the
irreducible component $\scrM_1$,
a thickening of multiplicity $2$. As explained above, this
arises from the structure of $\foM^{\ev}(\shX,\beta,z)$. 
Note that $\foM^{\ev}(\shX,\beta,z)\rightarrow 
B\GG_m^{\dagger}$ is log smooth by Lemma~\ref{lem: more log smoothness}.
By Lemma~\ref{lem:tropical point constraint}, if $Q=\NN$ is the ghost sheaf
at a general point of $\scrM_1$, then $\delta:Q^{\vee}\rightarrow\NN$
is given by \eqref{eq:voutdelta}.
Now $Q$ is generated by $\ell\in Q$, which viewed as a function on the
moduli space $Q^{\vee}_{\RR}$ of tropical maps of the given type
depicted in Figure~\ref{Figure1},
yields the affine length of either of
the two edges (as opposed to legs) of the tropical curve.
In particular, the location of the vertex adjacent to $x_{\out}$
is at distance $2\ell$ from the origin. Thus $\delta=2$. From this, it follows
that \'etale locally at a general point of $\scrM_1$, $\scrM(X,\beta,z)$
is of the form $\AA^1\times\Spec\kk[x]/(x^2)$, see 
\cite{ACGS18},~Prop.~B.4.

\begin{figure}
\input{Figure1.pspdftex}
\caption{Punctured map at a general point of $\scrM_1$}
\label{Figure1}
\end{figure}

\begin{figure}
\input{Figure2.pspdftex}
\caption{The limit $x_2\to\text{node}$ yields a punctured map in $\scrM_1\cap\scrM_2$}
\label{Figure2}
\end{figure}

\begin{figure}
\input{Figure3.pspdftex}
\caption{The limit $x_2\to x_1$ in $\scrM_1$}
\label{Figure3}
\end{figure}

The curve of Figure~\ref{Figure2} lies at the intersection of the 
two components $\scrM_1$ and $\scrM_2$, and we obtain $\scrM_2$ by
smoothing the node at the intersection of the components containing
$x_2$ and $x_3$. The generic behaviour is then illustrated in
Figure~\ref{Figure4}, with two further degenerations shown in
Figures~\ref{Figure5} and \ref{Figure6}. The notation $q$, $E_q$ will
be explained later.

\begin{figure}
\input{Figure4.pspdftex}
\caption{Punctured map at a general point of $\scrM_2$}
\label{Figure4}
\end{figure}

\begin{figure}
\input{Figure5.pspdftex}
\caption{Another special point in $\scrM_2$ with $x_2$ and $x_3$ in Fig.~\ref{Figure2} swapped}
\label{Figure5}
\end{figure}

\begin{figure}
\input{Figure6.pspdftex}
\caption{The limit $x_2\to x_3$ in $\scrM_2$}
\label{Figure6}
\end{figure}

The obvious forgetful and stabilization map $\scrM(X,\beta,z)\rightarrow
\overline{\scrM}_{0,4}$ induces an isomorphism $\scrM_1\rightarrow
\overline{\scrM}_{0,4}$ and maps $\scrM_2$ to the divisor
$D(x_2,x_3\,|\,x_1,x_{\out})$. Again we use standard notation for the irreducible
divisor in the moduli space of (pre-) stable curves whose general point
represents a curve with two irreducible components with marked points
distributed as indicated.

Thus the coefficient $2$ of the term $t^{[2L-E]}\vartheta_{v_2}$ in
the computation \eqref{Eqn: triple product in first example} of
$\vartheta_{p_1}\vartheta_{p_2}\vartheta_{p_3}$ is the degree of the map
$\scrM(X,\beta,z)\to \ol\scrM_{0,4}$, or the length of the fiber over a general
point. The proof of associativity in ordinary Gromov-Witten theory now works by
comparing the (virtual) length of the fibers over two of the three special points of $\ol\scrM_{0,4}$. Specifically, equality of this virtual length over $D(x_2,x_3\,|\,x_1,x_{\out})$ and over $D(x_1,x_2\,|\,x_3,x_{\out})$ would show
equality of the corresponding terms in $\vartheta_{p_1}\cdot
(\vartheta_{p_2}\vartheta_{p_3})$ and in $(\vartheta_{p_1}\vartheta_{p_2})
\cdot\vartheta_{p_3}$.

Indeed, over the divisor $D(x_1,x_2\,|\,x_3,x_{\out})$ sits the curve of
Figure~\ref{Figure3}. Normalizing at the node
$q$, we obtain two curves. Roughly speaking,
the component containing $x_1$ and $x_2$ contributes $t^0\vartheta_{p_1+p_2}
=\vartheta_{v_1+3v_2}$, and the other connected component calculates the
contribution $t^{A}\vartheta_r=t^{[2L-E]}\vartheta_{v_2}$ of the product
$\vartheta_{p_1+p_2}\cdot \vartheta_{p_3}$. In the tropical picture the
normalization at the node amounts to dividing the tropical curve into two
connected components by replacing the edge $E_q$ by a pair of legs. The component containing the outgoing leg is then of the form on the right of 
Figure~\ref{fig:delta}, and can be shown to readily compute the invariant
$N^A_{sp_3r}$ with $s=p_1+p_2$. The other component obtained from splitting,
however, has outgoing leg translated away from the origin, thus is of the
form on the left of Figure~\ref{fig:delta}.
To show that this moduli space computes
$N^0_{p_1 p_2s}$ requires the consideration of a two-parameter family of
tropical curves that tropically involves taking the length of $E_q$ to $\infty$
while keeping the location of the vertex attached to $x_{\out}$ fixed.

On the other hand, the fibre of the forgetful morphism over $D(x_2,x_3\,|\,x_1,
x_{\out})$ is $\scrM_2$, which is of course one-dimensional, despite
the moduli space being unobstructed, and there is no way to assign
a number to this fibre. Rather, as we shall see, this fibre is genuinely
virtually one-dimensional as well as actually one-dimensional. How, then,
do we select out a contribution to $\vartheta_{p_1}\cdot(\vartheta_{p_2}
\cdot \vartheta_{p_3})$ which will also yield the same contribution 
$t^{A}\vartheta_r$?

The answer lies in 
using the power of log geometry to select out the
part of the fibre over boundary points of $\overline{\scrM}_{0,4}^{\dagger}$
with the correct tropical behaviour. 
Morally, this proves to be similar to the situation described in 
Remark~\ref{rem:transverse}.

Thus we turn to a more detailed analysis
of the tropical maps appearing here.
\end{example}

\subsubsection{The relevant tropical maps}
\label{sec:tropical curves relevant}
We will consider families of tropical
maps to $\Sigma(X)$ as summarized in \S\ref{sec:tropical curves}. At
this stage, we will not assume they arise as tropicalization of logarithmic
maps, but the reader should keep in mind that these tropical maps
should arise from a logarithmic curve with three marked points $x_1,x_2,x_3$
and a puncture $x_{\out}$ as previously, with a point constraint on 
$x_{\out}$.

We fix $p_1,p_2,p_3,r\in \Sigma(X)(\ZZ)$ as usual and a domain graph $G$
of genus $0$ with precisely four legs,
labelled $E_{x_1},E_{x_2},E_{x_3},E_{\out}$
respectively, with the first three legs being rays and the fourth
being a ray only if $r=0$; otherwise $E_{\out}$ is a line segment.

\begin{definition}
The \emph{spine} $G'$ of $G$ is the minimal connected subgraph
of $G$ containing all the legs of $G$.
\end{definition}

See Figure~\ref{Figurefourpointed} for an example of a spine.

\begin{figure}
\input{fourpointedtropical.pspdftex}
\caption{The spine of a four-pointed tropical curve. The solid lines are
the spine, while the dotted lines indicate additional pieces of the 
tropical curve which are not contained in the spine. This tropical
curve has a terminal tail, drawn as a fat red path.}
\label{Figurefourpointed}
\end{figure}

Note that the underlying topological space of the spine $G'$ of $G$ is 
homeomorphic to the underlying topological space of
either: (1) a graph $G''$ with one quadrivalent vertex, with
four legs attached to this vertex; (2)
a graph $G''$ with two trivalent vertices, $v$ and $w$,
one edge connecting $v$ and $w$, and four legs, two adjacent to $v$
and two adjacent to $w$. 
Note in case (2) that as $G'$ is a subgraph of $G$, we can view
$v$, $w$ as vertices of $G$. 

\begin{definition}
\label{def:leg adjacent}
In the above situation, suppose $G''$ has two trivalent
vertices, $v$ and $w$. If a leg $E_{x_i}$ or $E_{\out}$ of $G$
is contained in the leg of the underlying topological space of
$G''$ adjacent to $v$ (resp.\ $w$), we say that the leg $E_{x_i}$
or $E_{\out}$ \emph{is adjacent to $v$} (resp. $w$), even if such
a leg is not adjacent to $v$ (resp.\ $w$) in $G$.
\end{definition}

For the remainder of this subsection, we assume given a family of tropical 
maps $h:\Gamma(G,\ell)\rightarrow \Sigma(X)$
over a cone $\sigma$ and with domain $G$
in the sense of \S\ref{sec:tropical curves}, with type given
by $(G,\bsigma,\mathbf{u})$.
We write as usual for $m\in \sigma$ the map
\[
h_m:G\rightarrow|\Sigma(X)|
\]
parameterized by $m$. We also assume given a function $\delta:\sigma
\rightarrow\RR_{\ge 0}$ which will play the same role as the function
$\delta$ in Lemma~\ref{lem:tropical point constraint}.

\begin{definition}
\label{def:tropical constraint}
We say $h$, along with the function $\delta$, is a 
family of tropical maps \emph{with contact
orders $p_1,p_2,p_3$ and $-r$ and constraint $\delta$}
if the following conditions hold.
\begin{enumerate}
\item
Let $\sigma_1,\sigma_2,\sigma_3$ and $\sigma_{\out}$ be the minimal
cones of $\Sigma(X)$ containing $p_1,p_2,p_3$ and $r$ respectively.
Then $\sigma_i\subseteq \bsigma(E_{x_i})$ and $\mathbf{u}(E_{x_i})=p_i$ 
for each $i$.
Further $\sigma_{\out}\subseteq \bsigma(E_{\out})$ and $\mathbf{u}(E_{\out})
=-r$.
\item Let $v_{\out}$ be the endpoint of $E_{\out}$. Then
\begin{equation}
\label{eq:deltanu}
h_m(v_{\out})=\delta(m) r.
\end{equation}
Note that $\delta$ is determined by this condition if $r\not=0$, but
if $r=0$, $\delta$ is genuinely an extra piece of data.
\end{enumerate}
\end{definition}

\begin{definition}
\label{def:miniversal}
We say the family $h$ of tropical maps parameterized by $\sigma$
with contact orders
$p_1,p_2,p_3$ and $-r$ and constraint $\delta$ is \emph{universal} if for
every family of tropical maps $h'$ of the same type as $h$ 
parameterized by $\sigma'$, with constraint $\delta'$, there is a unique map
$\sigma'\rightarrow\sigma$ such that $h'$ is the pull-back of
$h$ and $\delta'$ is the pull-back of $\delta$. 
We write $h_{\bas}$, $\sigma_{\bas}$, $\delta_{\bas}$ for the universal
family.

We say a family $h$ parameterized by $\sigma$ is \emph{miniversal}
if the tautological morphism $\sigma\rightarrow\sigma_{\bas}$ induces
an isomorphism of tangent spaces of cones (not necessarily respecting
integral structures).
\end{definition}

\begin{remark}
\label{rem:delta basic}
Note that if the family $h_{\bas},\delta_{\bas}$ is universal and 
$r\not=0$, then as before $\delta_{\bas}$ is not additional data,
being determined by \eqref{eq:deltanu}. However, if $r=0$, then
$\sigma_{\bas}=\sigma'_{\bas}\times\RR_{\ge 0}$, where $\delta_{\bas}:
\sigma_{\bas}\rightarrow\RR_{\ge 0}$ is the second projection, and
$\sigma'_{\bas}$ itself is the universal family of tropical maps constrained
by $r$, i.e., $h_m(v_{\out})=0$ for $m\in \sigma'_{\bas}$,
without the additional data of $\delta_{\bas}$.
\end{remark}

\begin{definition}
\label{def:common trop}
Suppose the spine $G'$ of
$G$ has two distinct trivalent vertices $v$ and $w$, with 
$E_{\out}$ adjacent to $v$.
We say that $G$ \emph{lies in $D(x_1,x_2\,|\,x_3,x_{\out})$} if 
$E_{x_3}$ is adjacent to $v$, so that $E_{x_1},E_{x_2}$ are 
adjacent to $w$. Here, all adjacency is in the sense of
Definition \ref{def:leg adjacent}. In this case,
let $v=v_1, v_2,\ldots, v_n=w$ be the sequence of vertices
traversed in passing from $v$ to $w$ in $G$, and let 
$E_i$ be the edge connecting $v_i$ to $v_{i+1}$. 

Let $\ell_i=\ell(E_i):
\sigma\rightarrow
\RR_{\ge 0}$ be the length function associated to the edge $E_i$, and set
\[
\ell=\sum_{i=1}^{n-1} \ell_i.
\]

See Figure~\ref{Figurefourpointed} for a picture of a graph $G$ lying
in $D(x_1,x_2\,|\,x_3,x_{\out})$.

We remark we can similarly say that $G$ lies in 
$D(x_2,x_3\,|\,x_1,x_{\out})$
if $E_{x_2}, E_{x_3}$ are adjacent to $w$ instead, and everything said
below for the case of $D(x_1,x_2\,|\,x_3,x_{\out})$ also goes for
$D(x_2,x_3\,|\,x_1,x_{\out})$ with the obvious modifications.
\end{definition}

\begin{definition}
In the situation of Definition~\ref{def:common trop}, 
we say that $G$ has a \emph{terminal tail} if
$v_{\out}\not=v$.
\end{definition}

We will see that terminal tails cause trouble at several points in the
proof. Another type of tail, discussed later,
called an \emph{internal tail} (see Definition \ref{def:internal tails}), 
will also cause problems. Let us observe
for now that families of curves with terminal tails can indeed arise:

\begin{example}
\label{ex:terminal tail}
In Example 
\ref{runningexample1}, (1), consider $A=D_1+E$, where $E$
is the exceptional curve of the blow-up, and take 
\[
p_1=v_1+v_2,\quad p_2=v_1,\quad p_3=v_3,\quad r=v_1.
\]
Fix $z\in Z_r=D_1$ to be the point of intersection of $E$ with $D_1$.
Then in the one-dimensional moduli space $\scrM(X,\beta,z)$, there
is a curve as depicted in Figure~\ref{Figure7}. There the component
containing $x_3$ maps isomorphically to $D_1$, the component containing
$x_1$ and $x_2$ is contracted to $D_1\cap D_2$, the component containing
$x_{\out}$ is contracted to $D_1\cap E$, and the remaining unmarked component
is mapped isomorphically to $E$.

The tropicalization of
this curve has a terminal tail, as is visible on the right in Figure
\ref{Figure7}.
Here it is clear that if we had chosen $z$ more generally, this tail wouldn't
have appeared, perhaps supplying some intuition as to why such curves
will not cause trouble at a virtual level. Indeed, this hints
at the proof of what we call the no-tail lemma, Lemma~\ref{thm:no-tail-lemma}, 
which concludes that the contributions from tails vanish
because they involve gluing a moduli space with negative virtual dimension.

\begin{figure}
\input{Figure7.pspdftex}
\caption{}
\label{Figure7}
\end{figure}
\end{example}

\subsubsection{\textsc{Step I}: Independence of modulus}
\label{sec:ind of modulus}

We are now in a situation where we have a moduli space $\scrM(X,\beta,z)$
of virtual dimension one
along with the morphism $\Phi:\scrM(X,\beta,z)\rightarrow
\overline\scrM_{0,4}^{\dagger}$ of Lemma \ref{lemma: moduli morphisms}.
As we have seen, $\Phi$ is virtually log smooth, in the sense that the
corresponding morphism $\foM^{\ev}(\shX,\beta,z)\rightarrow
\overline\scrM_{0,4}^{\dagger}$ is log smooth, but may not be virtually flat.
In particular, the fibre of $\Phi:\scrM(X,\beta,z)\rightarrow
\overline\scrM_{0,4}^{\dagger}$ may be virtual dimension one over the
boundary points of $\overline\scrM_{0,4}^{\dagger}$. 

To rectify this,
for each $y\in \overline\scrM_{0,4}$, 
we will define below a log stack $B\GG_m^{\ddagger}$ with a log morphism
$\psi_y:B\GG_m^{\ddagger}\rightarrow \overline\scrM_{0,4}^{\dagger}$.
The underlying morphism of stacks of $\psi_y$ is a closed embedding, identifying
$B\GG_m$ with $y\times B\GG_m$. Further, we require $\psi_y$ be
transverse to $\Phi:\foM^{\ev}(\shX,\beta,z)\rightarrow \overline
\scrM_{0,4}^{\dagger}$ in the sense of Definition~\ref{def: transverse}.
This will guarantee that the base-change $\scrM_y^{\ddagger}$
of $\scrM(X,\beta,z)$ by $\psi_y$
is virtual dimension zero. We will ultimately see that the degree of the
virtual fundamental class of this base-change is independent of the 
specific parameters which go into the choice of $B\GG_m^{\ddagger}$
and $\psi_y$. This replaces the virtual flatness argument in the usual
proof of associativity of quantum cohomology.

Second, to calculate the virtual degree of $\scrM^{\ddagger}_y$ when
$y$ is a boundary point of $\overline\scrM_{0,4}$, we need to be able
to split each domain curve in this moduli space into two pieces, contributing
to terms on the (say) left-hand side of \eqref{asseq}. To do this, we
need to have a natural choice of node at which to split. We shall see
in \S\S\ref{sec:key tropical analysis} and \ref{sec:second step} how
this node is selected. However, the key point is that we wish to
select a node whose corresponding edge in the tropicalized curve has a
very specific behaviour. In particular, given a generic
point of $\foM^{\ev}(\shX,\beta,z)$ lying over a boundary point of
$\overline\scrM_{0,4}^{\dagger}$,
we obtain a family of tropical maps parameterized by
a cone $\sigma$ as in Definition \ref{def:tropical constraint}. 
We wish to select an edge of $G$ whose length is unbounded in $\sigma$
given a fixed value of $\delta$. Note it is insufficient to just ask
that the length of the edge be unbounded without fixing the value of 
$\delta$, as all edges have unbounded length via rescaling the tropical
map. 

We note that this idea of considering four-pointed tropical curves and
looking at families where the tropical modulus goes to infinity is a
familiar one
in tropical geometry: the first occurence of this that we are aware of is
in Gathmann and Markwig's proof of WDVV for tropical curves in the
projective plane \cite{GM}.

Since we cannot fix the value of $\delta$ in logarithmic geometry,
instead we demand that the ratio of the length of the chosen edge
to the value of $\delta$ be unbounded. If this quantity is unbounded
for an edge contained in the spine of $G$, then the quantity $\ell/\delta$,
where $\ell$ is as in Definition \ref{def:common trop}, is also unbounded. This
motivates the following definition.

\begin{definition}
\label{def: BGmdagger}
Let $y\in \overline{\scrM}_{0,4}$. We define a log stack
$B\GG_m^{\ddagger}$ and a morphism
\[
\psi_y:B\GG_m^{\ddagger}\rightarrow\overline{\scrM}_{0,4}^{\dagger},
\]
depending on $y$, where 
$\overline{\scrM}_{0,4}^{\dagger}=\overline{\scrM}_{0,4}\times 
B\GG_m^{\dagger}$
as defined in Lemma~\ref{lemma: moduli morphisms}.

If $y$ is not one of the three boundary points of $\overline{\scrM}_{0,4}$,
then we set $B\GG_m^{\ddagger}=B\GG_m^{\dagger}$, and take
$\psi_y$ to be an isomorphism with the closed substack 
$y\times B\GG_m^{\dagger}$ of $\overline{\scrM}_{0,4}^{\dagger}$.

If $y$ is a boundary point, we select an auxiliary positive integer
$\lambda$, which in general will be taken to be suitably large.
Consider the group $\ZZ^2$ with basis $\ell,\delta$ and define submonoids
of $\ZZ^2$ and their duals given by
\begin{align}
\label{eq:R Rlambda def}
\begin{split}
R:= \NN\ell\oplus\NN\delta,& \quad R_{\lambda}:=\NN(\ell-\lambda\delta)
+\NN\delta\subset R^{\gp},\\
R^{\vee}= \NN\ell^*\oplus\NN\delta^*& \quad
R_{\lambda}^{\vee}=\NN(\lambda\ell^*+\delta^*)+\NN\ell^*\subset R^*.
\end{split}
\end{align}
Let $B\GG_m^{\ddagger}$ denote $B\GG_m$ with log structure
having ghost sheaf $R_{\lambda}$, with torsor associated to $a\ell+b\delta\in 
R_{\lambda}$
given by $\shU^{\otimes b}$, where $\shU$ is the universal torsor on
$B\GG_m$.
\end{definition}

Note that if $y$ is a boundary point of $\ol\scrM_{0,4}$ then 
the log structure on 
$\overline{\scrM}_{0,4}^{\dagger}$ pulled back to $y\times B\GG_m$
has ghost sheaf $R=\NN\ell\oplus\NN\delta$ as in \eqref{eq:R Rlambda def}, where the factor $\NN\ell$ 
comes from the
divisorial log structure induced by $y\in \overline{\scrM}_{0,4}$
and the factor $\NN\delta$ comes from $B\GG_m^{\dagger}$. Thus the
torsor associated to $a\ell+b\delta$ is again $\shU^{\otimes b}$, and 
there is a canonical log
morphism $\psi_y:B\GG_m^{\ddagger}\rightarrow \overline{\scrM}_{0,4}^{\dagger}$
with image $y\times B\GG_m$ and
$\bar\psi_y^{\flat}:R\rightarrow R_{\lambda}$ the inclusion.
The tropicalization
\begin{equation}
\label{eq:psiy tropical}
\Sigma(\psi_y):\Sigma(B\GG_m^{\ddagger}) \rightarrow 
\Sigma(\overline{\scrM}_{0,4}^{\dagger})
\end{equation}
has image the subcone of $\RR_{\ge 0}^2$ generated by $\lambda\ell^*
+\delta^*$ and $\ell^*$. See Figure~\ref{sigmam} below for some illustration.

\begin{definition}
\label{def: Mydef}
Let $y\in\overline\scrM_{0,4}$. Assume further fixed 
$p_1,p_2,p_3,r\in B(\ZZ)$ and a class of
curve $\beta$ with four punctured points of contact orders 
$p_1,p_2,p_3$ and $-r$.
Assume $z\in Z_r^{\circ}$ given. Then set
\begin{align*}
\scrM_y^{\ddagger}=\scrM_y^{\ddagger}(X,\beta,z):= {} & \scrM(X,\beta,z) \times_{\overline\scrM_{0,4}^{\dagger}}
B\GG_m^{\ddagger},\\
\foM^{\ddagger,\ev}_y=\foM^{\ddagger,\ev}_y(\shX,\beta,z):= {} & 
\foM^{\ev}(\shX,\beta,z) 
\times_{\overline\scrM_{0,4}^{\dagger}}
B\GG_m^{\ddagger}.
\end{align*}
Here the morphisms $\scrM(X,\beta,z), \foM^{\ev}(\shX,\beta,z)\rightarrow
\overline\scrM_{0,4}^{\dagger}$ are the morphisms $\Phi$ given by Lemma
\ref{lemma: moduli morphisms}. We use the short-hand notation on the left
when clear from context, which will essentially always be the case
as we never consider more than one four-pointed curve class at a time.

As $\scrM_y^{\ddagger}=\scrM(X,\beta,z)\times_{\foM^{\ev}(\shX,\beta,z)} 
\foM_y^{\ddagger,\ev}$,
there is a perfect relative obstruction theory for $\scrM_y^{\ddagger}
\rightarrow
\foM_y^{\ddagger,\ev}$ pulled back from that of $\scrM(X,\beta,z)\rightarrow
\foM^{\ev}(\shX,\beta,z)$.
\end{definition}

The first of the key results, to be proved in \S\ref{sec: modulus invariance}, 
is the following.

\begin{theorem}
\label{theorem: independence of modulus}
Suppose that the virtual dimension of $\scrM(X,\beta,z)$ is $1$
and either $\pm c_1(\Theta_{X/D})$ is nef or $(X,D)$ is log Calabi-Yau. 
Let $y\in \overline{\scrM}_{0,4}$. If either (1) $y$ is a boundary point
of $\overline{\scrM}_{0,4}$ and $\lambda$ is chosen 
to be sufficiently large in Definition 
\ref{def: BGmdagger} or (2) $y$ is not a boundary point, then 
$\scrM_y^{\ddagger}$ carries
a virtual fundamental class of dimension $0$. Furthermore,
$\deg [\scrM^{\ddagger}_y]^{\virt}$ is then independent of $y$, 
so that in particular $\deg [\scrM^{\ddagger}_y]^{\virt}=\deg
[\scrM^{\ddagger}_{y'}]^{\virt}$ for $y$ a boundary point and
$y'$ a non-boundary point.
\end{theorem}

\begin{remark} 
\label{rem: tropicalization of m}
It is useful to keep in mind the nature of the tropicalization
\[
\Sigma(\Phi):\Sigma(\scrM(X,\beta,z))\rightarrow
\Sigma(\overline\scrM_{0,4}^{\dagger}).
\]
The tropicalization of the map $\scrM(X,\beta,z)\rightarrow B\GG_m^{\dagger}$
has already been discussed in Definition~\ref{def:delta}, and is
written as $\delta:\Sigma(\scrM(X,\beta,z))\rightarrow
\Sigma(B\GG_m^{\dagger})=\RR_{\ge 0}$. On the
other hand, the morphism $\scrM(X,\beta,z)\rightarrow \overline{\scrM}_{0,4}$
factors through $\Mbf=\Mbf_{0,4}$, the moduli space of pre-stable 
curves. 

We first review a standard description
of the tropicalization of the stabilization map
$\Sigma(\Mbf_{0,4})\rightarrow
\Sigma(\overline{\scrM}_{0,4})$, see e.g., \cite{ACP15} or 
\cite{CCUW},~Lem.~8.8. Given
a geometric point $\bar x$ of $\Mbf_{0,4}$, let $C$ be the corresponding basic
pre-stable curve. Suppose the stabilization $C'$ of $C$ is reducible, 
necessarily
having two components.
If the dual graph of $C$ is $G$, then the dual graph of $C'$
is $G'$, the spine of $G$, with two-valent vertices removed.
Further, $G'$ has two distinct trivalent
vertices, $v$ and $w$,
and let $v=v_1,v_2,\ldots,v_n=w$, $E_i$ be as in Definition~\ref{def:common trop}.
Recall that $\overline{\shM}_{\Mbf,\bar x}$ is the
free monoid with one generator for each node of $C$, with $\ell(E)$
the generator corresponding to an edge $E$ of $G$. Set
$\ell_i=\ell(E_i)$.
Meanwhile, the ghost sheaf of $\overline{\scrM}_{0,4}$ at the image of 
$\bar x$ is $\NN$, generated by an element $\ell$.
Then the stabilization map at the level of ghost sheaves
sends $\ell$ to $\ell_1+\cdots+\ell_{n-1}$.
Put another way, the tropicalization of the stabilization map at 
$\bar x$ is given by
\begin{align*}
\Hom(\overline{\shM}_{\Mbf,\bar x},\RR_{\ge 0})
\rightarrow & \RR_{\ge 0}\ell^*,\\
m \mapsto  &\left( \sum_i \ell_i(m)\right)\ell^*.
\end{align*}
So $\ell$ is the sum of the lengths of the edges $E_1,\ldots,E_{n-1}$,
and may be viewed as the tropical modulus of the tropical stabilization
of the tropical curve 
parameterized by a point of $\Hom(\overline{\shM}_{\Mbf,\bar x},\RR_{\ge 0})$.

Under the forgetful map $\scrM(X,\beta,z)\rightarrow \Mbf_{0,4}$,
the elements $\ell_1,\ldots,\ell_{n-1}$ are taken to the elements
of the basic monoid corresponding to the nodes indexed by
$E_1,\ldots,E_{n-1}$, see \cite{ACGS18},~\S2.3.
So if instead $\bar x$ is a geometric point of $\scrM(X,\beta,z)$,
$C\rightarrow \bar x$ the corresponding curve, then we can write
$\ell_1,\ldots,\ell_{n-1}$ for the elements of the basic monoid indexed
by the nodes corresponding to $E_1,\ldots,E_{n-1}$. We continue to
use the notation 
\[
\ell:=\ell_1+\cdots+\ell_{n-1}.
\]

Putting this together, we see that given $\bar x\in \scrM(X,\beta,z)$
mapping to a boundary point of 
$\overline\scrM_{0,4}^{\dagger}$ with notation as above, 
the tropicalization is given by
\begin{equation}
\label{eq:trop stuff}
\sigma_{\bar x}\rightarrow \RR_{\ge 0}\ell^*\oplus\RR_{\ge 0}\delta^*,\quad
m\mapsto  \ell(m)
\ell^*+\delta(m)\delta^*.
\end{equation}

Given this description, along with Proposition~\ref{tropicalproduct}, we
can now tropically interpret the fibre product
$\scrM(X,\beta,z)\times^{\fs}_{\overline\scrM_{0,4}^{\dagger}} 
B\GG_m^{\ddagger}$ defining $\scrM^{\ddagger}_y$ for $y$ a boundary point.
From \eqref{eq:psiy tropical}, we see that the fibre product
defining $\scrM^{\ddagger}_y$ selects out that part of the moduli space 
$\scrM(X,\beta,z)$ where the value of $\ell$ (which from the above
discussion can be viewed as the tropical modulus of the $4$-pointed
domain) is very large compared to $\delta$, the distance of the image of
the vertex $v_{\out}$ to the origin in $\Sigma(X)$.
As long as the slope $\lambda$ is
taken sufficiently large, we will be able to guarantee that for any
tropical map appearing in the family parameterized by 
$\Sigma(\scrM^{\ddagger}_y)$,
the tropical modulus can be taken to infinity while keeping the distance of
$v_{\out}$ from the origin fixed.
\end{remark}

\begin{remark}
\label{rmk:first place tails}
The possible existence of tails (in this case terminal tails) is important,
and a certain amount of work is necessary in the course of the proofs
here to eliminate their influence. The first time they cause potential
problems is in the proof of Theorem~\ref{theorem: independence of modulus}.
Indeed, the existence of curves in $\scrM(X,\beta,z)$ whose tropicalizations
have terminal tails may result in correction terms, leading to the possibility
that $\deg [\scrM^{\ddagger}_y]$ 
as in Theorem~\ref{theorem: independence of modulus} 
in fact depends on $y$. See Step 3 of the proof of Lemma 
\ref{lemma: invarianceI}, where we need to show that the contribution
to the first Chern class of the conormal bundle at $x_{\out}$ coming
from the divisor $D(x_1,x_2,x_3\,|\,x_{\out})\subseteq \Mbf_{0,4}$
vanishes.

This vanishing is proved in in the no-tail lemma, 
Lemma~\ref{thm:no-tail-lemma},
Essentially, we are able to prove there is no contribution because
any contribution would arise from gluing two moduli spaces of punctured maps,
one of which is of negative virtual dimension.

However, it is expected that tails should play
a role in defining a differential in a definition of an all-degree
analogue of symplectic cohomology in our setting. 
In the course of the proofs of theorems in this paper comparing various
invariants, the vanishing of
potential corrections caused by tails might then be reinterpreted
in terms of vanishing of the differential from degree zero to degree
one, at least in the log Calabi-Yau case.
\end{remark}

\begin{example}
\label{ex:extended2}
We continue with Example~\ref{ex:extended1}. Recall that
the component $\scrM_2$ sits over the boundary divisor
$D(x_2,x_3\,|\,x_1,x_{\out})$ of $\overline{\scrM}_{0,4}$. Let us
consider $\Sigma(\Phi):\Sigma(\scrM(X,\beta,z))\rightarrow\Sigma(\overline
\scrM_{0,4}^{\dagger})$ when restricted to $\Sigma(\scrM_2)$. 
Figure~\ref{sigmam} shows the quadrant $\RR_{\ge 0}\ell^*\oplus
\RR_{\ge 0}\delta^*$ with coordinates
given by $\ell$ and $\delta$, and the tropical maps corresponding
to points in the inverse image of $\Sigma(\Phi)$. Here the red chain
of edges is the chain $E_1,\ldots,E_{n-1}$. The one curve
depicted in the cone where $\delta\ge \ell$ is the curve shown in Figure
\ref{Figure2}, the curve over the diagonal where $\delta=\ell$
is shown in Figure~\ref{Figure4}, and the two curves lying
over the cone where $\ell\ge\delta$ are the curves of Figures~\ref{Figure5}
and \ref{Figure6}.

As long as the image of $\Sigma(\psi_y)$, 
as in \eqref{eq:psiy tropical},
is contained in the cone
where $\ell\ge\delta$ 
(i.e., as long as we have taken $\lambda\ge 1$)
$\scrM^{\ddagger}_y$ will consist of two points projecting to the two
boundary points of $\scrM(X,\beta,z)$ corresponding to Figures
\ref{Figure5} and \ref{Figure6}. In this case $\scrM^{\ddagger}_y$ is 
non-singular
and length two. 

Taking $y$ instead to be the boundary point
in the divisor $D(x_1,x_2\,|\,x_3,x_{\out})$, we may take
$\lambda\ge 0$, and then $\scrM^{\ddagger}_y$ consists of one non-reduced point
of length $2$, with tropicalization given by Figure~\ref{Figure3}.
So in both cases, $\deg [\scrM^{\ddagger}_y]^{\virt}=2$, as desired.

\begin{figure}
\input{sigmam.pspdftex}
\caption{The $(\ell,\delta)$-quadrant for Example~\ref{ex:extended1}. The
integral length of the fat red sequence of edges is $\ell$. Base change with
$B\GG_m^{\ddagger}$ restricts to a narrow cone with $\ell\gg\delta$ as
indicated by the shaded region.}
\label{sigmam}
\end{figure}
\end{example}

\subsubsection{\textsc{Step II}: The tropical analysis
of splitting edges} 
\label{sec:key tropical analysis}

Here we show that given a certain type of family of four-pointed
tropical maps to $\Sigma(X)$, there is a canonical choice of edge 
along which to split. A bit more explicitly,
such families of tropical maps can be viewed as arising by tropicalizing
a punctured map lying over a generic point of an irreducible component
of $\foM^{\ddagger,\ev}_y$ from Definition~\ref{def: Mydef}. We will eventually see that such
a family of tropical maps is two-dimensional,
parameterized by $\ell$ and $\delta$, parameters described in
\S\ref{sec:tropical curves relevant}. What we will show is that with
certain assumptions on this family, the length $\ell_i$ of at least
one of the edges $E_i$ described
in Definition~\ref{def:common trop} will be unbounded even when $\delta$
is fixed. Further, in nice cases, there will be a unique
choice of such an edge.
Splitting the domain of the tropical map at $E_i$ into two curves
allows us to decompose the tropical map as a gluing of two other
tropical maps. Ultimately, we will then need to do the same splitting
at the level of punctured maps, but here we focus only on the requisite
tropical geometry.

\begin{assumptions}
\label{ass:tropical}
For the remainder of this subsection, we assume given a family of tropical maps 
$h$ to $\Sigma(X)$ with contact orders $p_1,p_2,p_3$ and $-r$ and
constraint $\delta$ as in Definition~\ref{def:tropical constraint}
parameterized by a two-dimensional cone $\sigma$
and with underlying domain $G$. We write as usual for $m\in \sigma$ the map
\[
h_m: G \rightarrow |\Sigma(X)|
\]
parameterized by $m$. We assume further that this family satisfies the 
following properties, referring to Definition \ref{def:common trop} 
for notation:
\begin{enumerate}
\item $G$ lies in $D(x_1,x_2\,|\,x_3,x_{\out})$.
\item The family $\Sigma$ is miniversal in the sense of Definition
\ref{def:miniversal}.
\item Using the notation of Definition~\ref{def:common trop},
we assume that the image of the map
\[
\sigma\rightarrow \RR_{\ge 0}\ell^*\oplus \RR_{\ge 0}\delta^*,\quad
m\mapsto  \ell(m)\ell^*+\delta(m)\delta^*
\]
as in \eqref{eq:trop stuff}
has image a two-dimensional cone containing $\ell^*$.
\end{enumerate}
Note (3) implies that the image of $(\ell,\delta)$ may be written as
the cone generated by $\lambda\ell^*+\delta^*$ and $\ell^*$ for some
$\lambda \ge 0$, $\lambda\in\QQ$.

Given such a family, we recall the notation of Definition \ref{def:common trop}
for $v_1,\ldots,v_n$ and $E_1,\ldots,E_{n-1}$.
For each $1\le i\le n$, let 
\[
\tau_i:=\nu_{v_i}(\sigma)\subseteq \bsigma(v_i),
\]
where $\nu_{v_i}$ is defined in \eqref{eq:nuv}.
Let $u_i=\mathbf{u}(E_i)\in \bsigma(E_i)^{\gp}$ be the 
contact order of the edge $E_i$, oriented from $v_i$ to $v_{i+1}$.
\end{assumptions}

The following is perhaps the single most important
fact in the proof 
of associativity. It will tell us that provided the tropical modulus
of a tropical curve is sufficiently large, there will be a good
choice of edge along which to split.

\begin{proposition}
\label{claim:unique}
Given Assumptions~\ref{ass:tropical}, assume in addition that $G$
does not have a terminal tail. Then
there exists a unique $i$ such that $u_i$ is
tangent to $\tau_{i+1}$. Further, for this $i$,
\begin{enumerate}
\item $\ell_1,\ldots,\ell_{i-1}$ are proportional to $\delta$.
\item $\ell_i$ and $\delta$ are linearly independent.
\item $u_i\in \bsigma(E_i)$ (as opposed to in $\bsigma(E_i)^{\gp}$).
\end{enumerate}
\end{proposition}

\begin{proof}
We first show existence.  
Write $\nu_i:=\nu_{v_i}$ for brevity. If $r=0$, then by the no-tail assumption,
$\nu_1(m)=0$ for any $m$ by \eqref{eq:deltanu}. Hence $u_1$ is proportional
to $\nu_2(m)$ for any $m$ by \eqref{eq:v1v2rel}. Since $\nu_2$ is linear, it follows in
particular that $u_1$ is tangent to $\tau_2$, and we can take $i=1$.

If $r\not=0$, then $\nu_1(m)$ is determined by $\delta$ via
\eqref{eq:deltanu}. By \eqref{eq:v1v2rel}, we have
\begin{equation}
\label{eq:i to iplus1}
\nu_{i+1}(m)=\nu_i(m)+\ell_i(m) u_i.
\end{equation}
If there exists
$m, m'\in \sigma$ with $\nu_i(m)=\nu_i(m')$ but $\nu_{i+1}(m)
\not=\nu_{i+1}(m')$, then
$\ell_i(m)\not=\ell_i(m')$. In this case,
$\left(\ell_i(m')-\ell_i(m)\right) u_i$ is a non-zero vector in
$\tau_{i+1}^{\gp}$, and
hence $u_i$ is tangent to $\tau_{i+1}$. 
Conversely, if there is no edge
$E_i$ with $u_i$ tangent to $\tau_{i+1}$, then inductively
$\nu_1(m),\ldots,\nu_n(m)$
are determined entirely by $\nu_1(m)$. As
$\nu_1(m)$ is determined by $\delta$, \eqref{eq:i to iplus1} tells us
inductively that $\ell_i$ is a function of $\delta$, i.e., is
constant on the sets $\delta=C$ in $\sigma$ for any constant real number $C$. 
However, $\ell_i$
is a linear function on $\sigma$, and hence must be rationally
proportional to $\delta$. Thus
$\ell$ is also rationally proportional to $\delta$, violating
condition (3) of Assumptions~\ref{ass:tropical}. Note this also
shows (1).

To show uniqueness, suppose there are two edges $E_i$, $E_j$, $i<j$,
with $u_i,u_j$ tangent to $\tau_{i+1}, \tau_{j+1}$ respectively. We will
argue that we can deform $\Gamma\rightarrow\Sigma(X)$ in a 
two-dimensional family if $r=0$ and a three-dimensional family if $r\not=0$,
so that in both cases $\sigma_{\bas}$ is at least three-dimensional.
Since $\dim\sigma=2$ by assumption, this
contradicts the miniversality assumption (2) of
Assumptions~\ref{ass:tropical}. We will
do this by constructing a two-dimensional family on which $\nu_1$ is
constant regardless of the value of $r$. Once this is accomplished, 
if $r=0$, we are done. Indeed, in this case $\nu_1$ is identically
zero on the entire moduli space in any event, and $\delta$ is an independent
parameter, so the full moduli space is at least three-dimensional.
If $r\not=0$, we again obtain a three-dimensional
family as $\delta$, hence $\nu_1$, 
is not constant on $\sigma$ by \eqref{eq:deltanu}. 

Take $m\in\Int(\sigma)$.
Since $u_i$ is tangent to $\tau_{i+1}$ and $\nu_{i+1}(m)\in\Int(\tau_{i+1})$,
we have $\nu_{i+1}(m)+\epsilon u_i\in \tau_{i+1}$ for $\epsilon$ sufficiently
close to $0$. Thus for each $\epsilon$, there is an $m'\in\sigma$
with $\nu_{i+1}(m')=\nu_{i+1}(m)+\epsilon u_i$. This allows us to
create a new tropical map $h_{\epsilon}:\Gamma\rightarrow
\Sigma(X)$ as follows. Removing the edge $E_i$ splits $\Gamma$
into two connected components, $\Gamma_1$ and $\Gamma_2$, with
$\Gamma_2$ containing the leg $E_{\out}$. Then we define 
$h_{\epsilon}|_{\Gamma_2}=h_m|_{\Gamma_2}$, and 
$h_{\epsilon}|_{\Gamma_1}=h_{m'}|_{\Gamma_1}$. Finally
$h_{\epsilon}$ maps $E_i$ to the segment joining 
$\nu_i(m)$ and $\nu_{i+1}(m')$. This gives a tropical map of the same
type as $h_m$, since $\nu_{i+1}(m')-\nu_i(m)$ is a positive
multiple of $u_i$ for sufficiently small $\epsilon$. In addition,
this tropical map has the same value as $h_m$ on $v_1$.

The same procedure can be carried out at the edge $E_j$, completely 
independently of the adjustment to the length of $E_i$. This shows
that we do have a two-parameter family as claimed, with $\nu_1$
constant in this family.

To prove (2), note that if $\ell_i$
were proportional to $\delta$, then \eqref{eq:i to iplus1}
tells us that $\nu_{i+1}(m)$ also only depends on $\delta(m)$.
Proceeding inductively as before, using the fact that $u_j$ is not tangent to
$\tau_{j+1}$ for any $j>i$, \eqref{eq:i to iplus1}
shows that $\ell_{i+1},\ldots,\ell_{n-1}$ only are functions of
$\delta$.
Thus all $\ell_1,\ldots,\ell_{n-1}$ are proportional to $\delta$, violating
Assumptions~\ref{ass:tropical}.

For (3), note that as $\ell_i$ and $\delta$ are linearly
independent and $\sigma$ is two-dimensional,
we can write $\ell=a\ell_i+b\delta$ for some rational
numbers $a$ and $b$. Further, from Assumptions~\ref{ass:tropical},
(3), $a\not=0$ and the image of the map $(\ell,\delta)$ 
is generated by $\lambda\ell^*+\delta^*$ and $\ell^*$ for some
$\lambda\ge 0$. Thus in particular, this image
contains
the half-line $L=\{\lambda'\ell^*+\delta^*\,|\,\lambda'\ge \lambda\}$,
on which $\delta$ takes the constant value $1$
and $\ell$ is positive and unbounded. Identify $\sigma$ with
its image under $(\ell,\delta)$, so that we may view $\ell_i$
as a linear function on this image. Since $\ell_i$ is non-negative on
$\sigma$ and $\ell_i=(\ell-b\delta)/a$,
we see that necessarily
$a>0$ and $\ell_i$ is unbounded on $L$. 
For $m\in L$, 
$\nu_i(m)$ is determined as before and hence independent of the choice
of $m\in L$ since
$\delta$ is constant on $L$, but $\nu_{i+1}(m)=\nu_i(m)+
\ell_i(m) u_i\in \bsigma(E_i)$. Since $\ell_i(m)$ is unbounded,
necessarily $u_i\in \bsigma(E_i)$, as claimed.
\end{proof}

\begin{example}
\label{ex:extended3}
Continuing on from Example~\ref{ex:extended2}, one sees that none
of the tropical maps appearing there have terminal tails. In this
case, Proposition~\ref{claim:unique} can be applied for
the two types of tropical maps parameterized by the cone
$\sigma = \RR_{\ge 0} \ell^* + \RR_{\ge 0} (\delta^*+\ell^*)$. These curves
are those appearing in Figures~\ref{Figure5} and \ref{Figure6}, also
depicted on the right-hand side of Figure~\ref{sigmam}. 
In the case of these two families of curves, the unique $i$ given by
Proposition~\ref{claim:unique} is $i=1$ and $2$ for the first
and second respectively. In Figures~\ref{Figure5},\ref{Figure6}, 
we label this edge
as $E_q$, and also indicate the corresponding node $q$ of the domain
curve.
\end{example}

\begin{remark}
The first place terminal tails caused problems was described in 
Remark~\ref{rmk:first place tails}, where they interfere with independence
of modulus in $\overline\shM_{0,4}$.
We now see a second place where terminal tails may cause problems.
If $G$ has
a terminal tail in Proposition~\ref{claim:unique}, 
there may not be an edge along which to split. This will cause
problems in the final stage of the proof, namely in the proof 
of Theorem~\ref{thm:main comparison}, as not all curves in the moduli
space $\scrM^{\ddagger}_y$ will have a choice of splitting. Again, we need to know
contributions from such curves vanish. This requires another
application of the no-tail lemma, Lemma~\ref{thm:no-tail-lemma}.
\end{remark}

As mentioned earlier, there is another type of tail 
which we will also have to deal
with in the last steps of the proof:

\begin{definition}
\label{def:internal tails}
In the situation of Assumptions~\ref{ass:tropical}, suppose $G$ does not
have a terminal tail. We say an edge $E_i$ is \emph{of splitting type}
if $\ell_i$ and $\delta$ are linearly independent as linear
functions on $\sigma^{\gp}$.
We say that the family of tropical
maps parameterized by $\sigma$ \emph{has an internal tail} if there is 
more than one edge $E_i$ of splitting type.
If either $G$ has a terminal
tail or the family has an internal tail, then we say the family of 
tropical maps parameterized by $\sigma$ \emph{has a tail}. Otherwise,
we say the family is \emph{tail-free}.\footnote{The reader may 
view the terminal tail as the union of edges
contained in the spine of $G$ connecting the vertex $v$ with the leg
$E_{\out}$, as depicted in Figure~\ref{Figurefourpointed} in red. 
The reader may view an internal tail as the collection of edges
from the second edge of splitting type to the first edge of splitting
type, including the second edge but not the first. However, we will not
need this language.}
\end{definition}

\begin{example}
We have already seen an example of a terminal tail in
Example~\ref{ex:terminal tail} and in fact
internal tails can also happen. Take $\overline X=\PP^1\times\PP^1$
and $\overline D$ the toric boundary, $\overline D=\overline D_1+\cdots+
\overline D_4$ in
cyclic order. Blow up a point $x$ in the interior of $\overline D_1$ to
get $X$, and take $D$ to be the inverse image of $\overline D$ (rather
than the strict transform!) See the left of Figure~\ref{Figure8}.
So this is a non-minimal log Calabi-Yau
manifold. We can write $D=D_1+\cdots+D_4+E$, where $E$ is the exceptional
divisor, and $K_X+D=E$. In this case the Kontsevich-Soibelman
skeleton $B$ is naturally identified with $\Sigma(\overline X)$. On the
other hand, $\Sigma(X)$ is obtained, abstractly as a cone complex, 
from $\Sigma(\overline X)$ by adding an additional two-dimensional standard
cone with rays generated by $D_1^*$ and $E^*$, so that this cone intersects
$B$ in the ray generated by $D_1^*$. Note $|\Sigma(X)|$ is not a manifold.

Consider 
\[
p_1=v_1+v_2, \quad p_2=v_4, \quad p_3=v_3, \quad r=0,
\]
where $v_i$ is the primitive generator of the ray of $\Sigma(X)$
corresponding to $D_i$. In addition, take
the curve class 
\[
A\sim D_3+D_4\sim D_1+(D_4-E)+2E.
\] 
The resulting curve class $\beta$
will potentially contribute to the $\vartheta_0$
term in the product $\vartheta_{p_1}\vartheta_{p_2}\vartheta_{p_3}$.
If we place $z\in Z_r^{\circ}=X\setminus D$ to lie on the strict transform
$L$
of the line of class $\overline D_4$ passing through the point $x$,
then Figure~\ref{Figure8} depicts a general curve in the moduli space
$\scrM^{\ddagger}_y$ as defined in Definition~\ref{def: Mydef}, $y\in
D(x_1,x_2\,|\, x_3,x_{\out})\subset \overline{\scrM}_{0,4}$. 
The picture is slightly misleading. The component
drawn as an ellipse is a degree $2$ cover of $E$. This cover must be
ramified over $E\cap D_1$ but need not be ramified where $E$
meets $L$. In particular, this moduli space is one-dimensional.

The essential part of the 
tropicalization of this curve is depicted on the right-hand
side of Figure~\ref{Figure8}, with the quadrant corresponding to the
double point $E\cap D_1$ depicted. All vertices of the tropical curve
are mapped into this quadrant, whereas some of the legs are not, and
we draw those legs as arrows. The reader should keep in mind these 
correspond to half-lines in other cones of $\Sigma(X)$. Both of the edges of the
curve are of splitting type. Indeed, as $r=0$, the value of $\delta$ has
no impact on the tropical map, and $\ell_1=\ell_2$.

The claim, shown in \S\ref{section:key comparison},
using the no-tail lemma, Lemma~\ref{thm:no-tail-lemma},
is that curves with
tails will not contribute to any virtual calculaton. This may seem a bit 
puzzling,
as in fact we expect a constant term in $\vartheta_{p_1}\vartheta_{p_2}
\vartheta_{p_3}$, as is immediately visible if we had chosen $z$
somewhat more generally. The answer lies in the magic of virtual irreducible
decompositions, which we now explain. 

If we degenerate the double cover of $E$ to a reducible
double cover, we obtain a particularly degenerate member of the family
of curves $\scrM^{\ddagger}_y$,
which tropicalizes to the left-hand picture in Figure~\ref{Figure9}.
Along with the parameter $\delta$, which is not visible in this picture, this
yields a three-dimensional family of tropical maps.
This degenerate curve lies over the intersection of two irreducible components
of the moduli space $\foM_y^{\ddagger,\ev}$ defined in Definition~\ref{def: Mydef}.
One of these irreducible components contains the image of the
family of curves described in Figure~\ref{Figure8}. Because $\foM^{\ddagger,\ev}_y$ is insensitive
to the precise location of $z\in Z_r^{\circ}=X\setminus D$, the other
irreducible component can be seen as the one containing the image of what would
have been
$\scrM^{\ddagger}_y$ if $z$ had been chosen generally. In that case, the domain consists
of three irreducible components, mapping to $D_1$, $E$ and a curve
linearly equivalent to $D_4$ passing through $z$. The tropicalization
of this curve is shown on the right in Figure~\ref{Figure9}.
It then makes sense that it is only this second irreducible component
which contributes to $\vartheta_{p_1}\vartheta_{p_2}\vartheta_{p_3}$,
as the first one would not contribute if $z$ had been chosen generally.

\begin{figure}
\input{Figure8.pspdftex}
\caption{}
\label{Figure8}
\end{figure}

\begin{figure}
\input{Figure9.pspdftex}
\caption{}
\label{Figure9}
\end{figure}
\end{example}

\subsubsection{\textsc{Step III}: Gluing}
\label{sec:second step}

The basic idea is now as follows. Let $y$ be a boundary point of 
$\overline{\scrM}_{0,4}$,
either the divisor $D(x_1,x_2\,|\,x_3,x_{\out})$ or
$D(x_2,x_3\,|\,x_1,x_{\out})$. Then, up to
the issue of terminal tails, Proposition~\ref{claim:unique}
essentially tells us that the moduli space
$\scrM^{\ddagger}_y$ described in \S\ref{sec:ind of modulus} should decompose
into (virtual) irreducible components over which the universal curve
splits into two curves, at least if $\lambda$ is large enough. 
If the original curve class was $\beta$, with
underlying curve class $A$ and with four marked points of
contact order $p_1,p_2,p_3$ and $-r$, then we are interested in situations
where the curve splits into two curves of classes
$\beta_1$ and $\beta_2$.
These classes
have underlying curve classes $A_1,A_2$ satisfying
$A_1+A_2=A$, and
$\beta_1$ having punctured points $x_1,x_2,x_s$ of contact orders
$p_1,p_2$ and $-s$ for some $s\in \Sigma(X)(\ZZ)$ and $\beta_2$
having punctured points $x_s'$, $x_3$, $x_{\out}$ of contact orders
$s, p_3$ and $-r$.

For example, if we continue with Example~\ref{ex:extended1}, we might see
the curve of Figure~\ref{Figure5} as obtained by gluing two pieces,
a curve of class $\beta_1$ with $A_1=L-E$ and $s=v_3$,
and a curve of class $\beta_2$ with $A_2=L$. The curve
of class $\beta_1$ is just the 
connected component of the partial normalization
of the domain at the labelled node $q$ containing the points $x_2,x_3$.

In particular, we have moduli spaces $\scrM(X,\beta_1)$ and $\scrM(X,\beta_2,z)$,
with $z\in Z_r^{\circ}$. Note that we don't put a point constraint in the
first moduli space. We can see that this is reasonable at a tropical level.
In Figure~\ref{Figure5}, if one restricts the tropical map to the sub-graph
of $G$ corresponding to the curve of class
$\beta_1$, one sees that the punctured leg (labelled $E_q$) 
with tangent vector $-s$
is not constrained to lie in the ray $\RR_{\ge 0}s$.

Reviewing the methods of \cite{ACGS18}, we construct
in \S\ref{sec:gluing review} a moduli space of glued curves 
$\scrM^{\gl}(X,\beta_1,\beta_2,z)$, parameterizing
curves in $\scrM(X,\beta,z)$ with a splitting into curves of class 
$\beta_1$ and $\beta_2$. Doing the same thing at the level of
maps to the Artin fan gives moduli spaces $\foM(\shX,\beta_1)$, 
$\foM^{\ev}(\shX,\beta_2,z)$ and a glued moduli space
$\foM^{\gl,\ev}(\shX,\beta_1,\beta_2,z)$, along with a morphism
$\scrM^{\gl}(X,\beta_1,\beta_2,z)\rightarrow \foM^{\gl,\ev}
(\shX,\beta_1,\beta_2,z)$.
There is a commutative diagram
\[
\xymatrix@C=30pt
{
\scrM^{\gl}(X,\beta_1,\beta_2,z)\ar[r]\ar[d]&\scrM(X,\beta,z)\ar[d]\\
\foM^{\gl,\ev}(X,\beta_1,\beta_2,z)\ar[r]&\foM^{\ev}(\shX,\beta,z)
}
\]
which is not cartesian; rather, it induces an open and closed
embedding of $\scrM^{\gl}(X,\beta_1,\beta_2,z)$ into the
fibre product. This is because as classes of maps to the Artin fan
$\shX$, $\beta_1$ and $\beta_2$ do not remember the curve classes
in $X$, and hence we need to pick out those curves in the
fibre product which split into curves of classes $A_1$ and
$A_2$.
Now ideally we would like $\scrM^{\gl}(X,\beta_1,\beta_2,z)$
to be of virtual dimension zero, and then the point of gluing would
be to express the degree of this virtual fundamental class as
$N^{A_1}_{p_1p_2s}N^{A_2}_{sp_3r}$. Unfortunately, this is not
always the case. 

This is precisely the same problem as we had in 
\S\ref{sec:ind of modulus}, where fibres of the morphism $\Phi$
were of too high virtual dimension.
Again, we rectify this problem by using a careful choice of transverse
map to select a part of $\scrM^{\gl}(X,\beta_1,\beta_2,z)$ which is
virtual dimension zero. We will give some more details below, but at
the tropical level, the idea is as follows.

For the tropicalization of any punctured map in 
$\scrM^{\gl}(X,\beta_1,\beta_2,z)$
there is a tropical parameter $\ell_q$ measuring the
length of the glued edge $E_q$. We wish to impose tropically the condition
that the ratio $\ell_q/\delta$ is unbounded, and this is done via a 
base-change which selects the part of $\scrM^{\gl}(X,\beta_1,\beta_2,z)$
for which $\ell_q/\delta$ is bounded below by an additional parameter $\mu$
which needs to be taken to be sufficiently large.

\begin{example}
We see that this is already an issue in
Example~\ref{ex:extended1}.
Indeed, the curve of Figure~\ref{Figure5} appears as a degeneration of
a family of curves shown in Figure~\ref{Figure4}. The node $q$
in Figure~\ref{Figure5} is a limit of the node $q$ indicated
in Figure~\ref{Figure4}. As a consequence, $\scrM^{\gl}(X,\beta_1,\beta_2,z)$,
in this instance, is one-dimensional and can be checked to be unobstructed. 

Note in the tropical map in Figure~\ref{Figure4}, the edge
$E_q$ corresponding to the node $q$
has length equal to the distance $\delta$ of the vertex adjacent to $x_{\out}$
to the origin. However, if we fix $\mu>1$ and only consider punctured maps 
where the
ratio $\ell_q/\delta$ is always at least $\mu$,
then we select out 
only the curve depicted in Figure~\ref{Figure5}.
\end{example}

We impose this tropical condition at the moduli space level as follows. In
\S\ref{sec:key gluing} (see \eqref{eq:Phidef}), we construct a morphism
\[
\Psi:\foM^{\gl,\ev}(\shX,\beta_1,\beta_2,z)\rightarrow
B\GG_m^{\dagger}\times B\GG_m^{\dagger}.
\]
The composition of $\Psi$ with projection onto the first factor is just
the composition $\foM^{\gl,\ev}(\shX,\beta_1,\beta_2,z)
\rightarrow\foM^{\ev}(\shX,\beta,z)\rightarrow B\GG_m^{\dagger}$
with the second morphism the projection from the fibre product
definition of $\foM^{\ev}(\shX,\beta,z)$ in Definition \ref{def:Mbetaev}.
The composition of $\Psi$ with projection onto the second factor is new, and
the induced tropical map $\Sigma(\foM^{\gl,\ev}(\shX,\beta_1,\beta_2,z))
\rightarrow \Sigma(B\GG_m^{\dagger})=\RR_{\ge 0}$ takes a tropical map
to the length $\ell_q$ of the edge $E_q$ corresponding to the glued node $q$.

Analogously to \eqref{eq:R Rlambda def}, we define
\begin{align}
\label{eq:T Tmu}
\begin{split}
T:=\NN\delta\oplus \NN\ell_q,&\quad T_{\mu}:=\NN\delta+\NN(\ell_q-\mu\delta)
\subseteq T^{\gp},\\
T^{\vee}=\NN\delta^*\oplus\NN\ell_q^*&\quad
T_{\mu}^{\vee}=\NN \ell_q^*+\NN (\delta^*+\mu\ell_q^*)\subset T^*,
\end{split}
\end{align}
where of course $\delta^*,\ell_q^*$ is the dual basis to $\delta,\ell_q$.
Here $\mu$ is a fixed positive integer, which we will eventually take
to be large.
We view $T$ as the stalk of the ghost sheaf of
$B\GG_m^{\dagger}\times B\GG_m^{\dagger}$, with $\delta$ and $\ell_q$
coming from the first and second factors respectively. 
We can then
define a new log structure on $B\GG_m^2$, which we write as $B\GG_m^{2,\mu}$,
with $\shM^{\gp}_{B\GG_m^{2,\mu}}=\shM^{\gp}_{B\GG_m^{\dagger}
\times B\GG_m^{\dagger}}$ but
$\overline{\shM}_{B\GG_m^{2,\mu}}=T_{\mu}$.  Thus there is a morphism
\[
B\GG_m^{2,\mu}\rightarrow B\GG_m^{\dagger}\times B\GG_m^{\dagger}.
\]
This tropicalizes to an inclusion of cones
\[
T_{\mu}^{\vee})_{\RR}=\RR_{\ge 0}\ell_q^*+\RR_{\ge 0}(\delta^*+\mu\ell_q^*)\hookrightarrow
\RR_{\ge 0}\delta^*\oplus \RR_{\ge 0}\ell_q^*.
\]

We can then define modified moduli spaces via the fs log cartesian
diagrams
\[
\xymatrix@C=30pt
{
\scrM^{\mu,\gl}(X,\beta_1,\beta_2,z)\ar[r]\ar[d]&
\scrM^{\gl}(X,\beta_1,\beta_2,z)\ar[d]\\
B\GG_m^{2,\mu}\ar[r]&B\GG_m^{\dagger}\times B\GG_m^{\dagger}
}
\]
and
\[
\xymatrix@C=30pt
{
\foM^{\mu,\gl,\ev}(\shX,\beta_1,\beta_2,z)\ar[r]\ar[d]&
\foM^{\gl,\ev}(\shX,\beta_1,\beta_2,z)\ar[d]\\
B\GG_m^{2,\mu}\ar[r]&B\GG_m^{\dagger}\times B\GG_m^{\dagger}
}
\]

We can now state the key gluing result:

\begin{theorem}
\label{thm: main gluing theorem}
Suppose $\pm c_1(\Theta_{X/\kk})$ is nef or $(X,D)$ is log Calabi-Yau. Then
we have a diagram cartesian in all categories
\[
\xymatrix@C=30pt
{
\scrM^{\mu,\gl}(X,\beta_1,\beta_2,z)\ar[r]\ar[d] &
\scrM^{\gl}(X,\beta_1,\beta_2,z)\ar[d]\\
\foM^{\mu,\gl,\ev}(\shX,\beta_1,\beta_2,z)\ar[r]&  \foM^{\gl,\ev}
(\shX,\beta_1,\beta_2,z)
}
\]
The relative obstruction theory for the right-hand vertical arrow pulls
back to a relative obstruction theory for the left-hand vertical arrow.
If $\mu\gg 0$, then
$\scrM^{\mu,\gl}(X,\beta_1,\beta_2,z)$ is proper of virtual dimension
$A\cdot c_1(\Theta_{X/\kk})$, and assuming this virtual
dimension is zero, then
\[
\deg [\scrM^{\mu,\gl}(X,\beta_1,\beta_2,z)]^{\virt}=
\begin{cases}
N^{A_1}_{p_1p_2s}N^{A_2}_{sp_3r} & s\in B(\ZZ)\\
0 & s\not\in B(\ZZ).
\end{cases}
\]
\end{theorem}

The proof of this theorem is given in \S\ref{sec:proof of main gluing}.
This is the first significant application of the logarithmic gluing
technology set up in \cite{ACGS18}.

\subsubsection{\textsc{Step IV}: Comparison of moduli spaces}

Now let $y\in \overline{\scrM}_{0,4}$ be either the
boundary point
$D(x_1,x_2\,|\,x_3,x_{\out})$ or $D(x_2,x_3\,|\, x_1,x_{\out})$.
Without loss of generality, we shall assume $y$ is always the former.
We fix $p_1,p_2,p_3,r\in
B(\ZZ)$ and an underlying curve
class $A$, giving a curve class $\beta$. We fix $z\in Z_r^{\circ}$
as usual. After fixing a $\lambda\gg 0$, we obtain,
from this data, moduli spaces
$\scrM^{\ddagger}_y$ and $\foM_y^{\ddagger,\ev}$, as defined
in Definition~\ref{def: Mydef}.

On the other hand,
given $s\in \Sigma(X)(\ZZ)$ and a decomposition $A=A_1+A_2$,
we obtain two classes of punctured map $\beta_1$, $\beta_2$
with the given underlying curve classes, each with three punctures of
contact orders $p_1,p_2,-s$ and $s,p_3,-r$ respectively. Fixing $\mu\gg 0$,
this gives rise to the glued moduli spaces
$\scrM^{\mu,\gl}(X,\beta_1,\beta_2,z)$, $\foM^{\mu,\gl,\ev}(\shX,
\beta_1,\beta_2,z)$, with the virtual degree of the former calculated
by Theorem~\ref{thm: main gluing theorem}. 

Note that there is a relationship between the 
parameters $\mu$ and $\lambda$.
The parameter $\lambda$ is used to bound
below the ratio $\ell/\delta$, where $\ell=\sum_i \ell_i$ as in 
Definition~\ref{def:common trop}. In our situation, 
$E_q$ will be an edge of splitting
type and hence $\ell_q$ agrees with $\ell_i$ 
for some $i$, so that $\ell_q$ is a summand
of $\ell$. Since the $\ell_i$ are all non-negative, we see that 
$\ell/\delta\ge \ell_q/\delta$ in our situation. As a consequence, it
will be natural to take $\lambda\ge\mu$. 
We typically first fix $\mu$ large
enough to obtain a virtual zero-dimensional moduli space, as we 
explained in \S\ref{sec:second step},
and then usually require
$\lambda\gg \mu$ in what follows to avoid having to worry about the
values of the $\ell_i$ for $\ell_i\not=\ell_q$. 

We wish to compare the virtual degree of $\scrM^{\ddagger}_y$ with the 
virtual degree of $\coprod_{A_1,A_2,s} \scrM^{\mu,\gl}
(X,\beta_1,\beta_2,z)$.
Here, when we index a union or sum over $A_1,A_2,s$,
we mean we take the union or sum over all choices of $s\in \Sigma(X)(\ZZ)$
 and all decompositions
$A=A_1+A_2$, with the classes
$\beta_1,\beta_2$ then being determined by
this data as in the previous paragraph. The main remaining theorem is:

\begin{theorem}
\label{thm:main comparison}
Suppose $\pm c_1(\Theta_{X/\kk})$ is nef or $(X,D)$ is log
Calabi-Yau. 
Suppose that $\mu$ is chosen sufficiently large, and that 
$\lambda\gg\mu$ is chosen sufficiently large compared to $\mu$.
Then, assuming $A\cdot c_1(\Theta_{X/\kk})=0$, we have
\[
\sum_{A_1,A_2,s}\deg [\scrM^{\mu,\gl}(X,\beta_1,\beta_2,z)]^{\virt}
=\deg [\scrM^{\ddagger}_y]^{\virt}.
\]
\end{theorem}

To prove this, we will need another intermediate moduli space.
There are canonical morphisms
\begin{equation}
\label{eq:modular maps}
\scrM^{\mu,\gl}(X,\beta_1,\beta_2,z)\rightarrow \scrM(X,\beta,z),
\quad
\foM^{\mu,\gl,\ev}(\shX,\beta_1,\beta_2,z)\rightarrow 
\foM^{\ev}(\shX,\beta,z).
\end{equation}
We may compose each of these with the canonical morphism $\Phi$ of
Lemma~\ref{lemma: moduli morphisms} to
$\overline\scrM_{0,4}^{\dagger}$.

Fix $\lambda,\mu> 0$. We obtain from $\lambda$
a morphism $\psi_y:B\GG_m^{\ddagger}
\rightarrow \overline{\scrM}_{0,4}^{\dagger}$ via Definition
\ref{def: BGmdagger}, and hence moduli spaces $\scrM^{\ddagger}_y$,
$\foM_y^{\ddagger,\ev}$ by Definition~\ref{def: Mydef}. We then obtain
a diagram
\begin{equation}
\label{eq:comparison diagram}
\xymatrix@C=30pt
{
\coprod_{A_1,A_2,s} \scrM^{\mu,\gl}(X,\beta_1,\beta_2,z)\ar[d]_{k_1} &
\coprod_{A_1,A_2,s} \scrM^{\mu,\gl}(X,\beta_1,\beta_2,z)
\times^{\fs}_{\overline\scrM_{0,4}^{\dagger}} B\GG_m^{\ddagger}\ar[l]_>>>>>>>{i'}
\ar[d]_{k_2} \ar[r]^>>>>>{j'} &
\scrM^{\ddagger}_y\ar[d]^{k_3}\\
\coprod_{s} \foM^{\mu,\gl,\ev}(\shX,\beta_1,\beta_2,z) &
\coprod_{s} \foM^{\mu,\gl,\ev}(\shX,\beta_1,\beta_2,z) 
\times^{\fs}_{\overline\scrM_{0,4}^{\dagger}} B\GG_m^{\ddagger}\ar[l]^>>>>>{i}
\ar[r]_>>>>>j &
\foM^{\ddagger,\ev}_y
}
\end{equation}
All vertical arrows are strict. For $k_1$ and $k_2$, note that
as punctured maps to the Artin fan do not remember the curve class,
specifying $s$ alone is enough to specify the classes $\beta_1$,
$\beta_2$ of maps to $\shX$.
Each square is cartesian in all
categories, as will be shown in Lemma~\ref{lem:finite representable}. 
Here $i'$ and $i$ are projections onto the first factor
and $j'$, $j$ are induced by the morphisms \eqref{eq:modular maps}.
The morphisms $k_1$ and $k_3$ carry relative obstruction theories by
Theorem~\ref{thm: main gluing theorem} and Lemma~\ref{lem:Mtau proper}
respectively. We shall see that these obstruction theories pull back
to the same relative obstruction theory for $k_2$.

It may be illustrative to interpret the horizontal arrows in
Diagram~\eqref{eq:comparison diagram} tropically. Assuming the outgoing contact
order to be non-zero for simplicity of the discussion, fix the distance
$\delta$ of $v_{\out}$ from the origin to be $1$.
Then $\foM_y^{\ddagger,\ev}$ on
the lower right restricts to tropical curves with sufficiently large tropical
modulus $\ell=\sum_i \ell_i$ (Definition~\ref{def:common trop}). The lower
left spaces $\foM^{\mu,\gl,\ev}$ on the other hand are moduli spaces of glued
punctured maps with length $\ell_q$ of the tropical gluing edge at least
$\mu$. The fiber products $\foM^{\mu,\gl,\ev}
\times^\fs_{\ol\scrM_{0,4}^\dagger} B\GG_m^\ddagger$ in the middle records
both $\ell_q,\ell$ and enforces $\ell_q\ge \mu$, $\ell\ge \lambda$.

In an ideal world, we would show that $i$ and $j$
are degree $1$ morphisms, so that we may apply Costello's result
\cite{Costello},~Thm.~5.0.1 or \cite{Man}. 
Unfortunately, this is not true in general as there may be irreducible
components of the various moduli spaces involving tails, and these components
cause problems for $i$ and $j$. However, we will show that such 
components contribute
$0$ to the degree of the virtual fundamental classes, again
using Lemma~\ref{thm:no-tail-lemma}.
It is only at this step that internal tails cause trouble.
For the details, see \S\ref{section:key comparison}.

\subsubsection{The proof of Theorem~\ref{thm:restatement}}

Theorem~\ref{thm:restatement}, hence the main results of the paper,
Theorems~\ref{mainassociativity1} and
\ref{mainassociativity2}, now follows immediately. Indeed, fix
$A$, $p_1,p_2,p_3,r\in B(\ZZ)$. If $A\cdot c_1(\Theta_{X/\kk})
\not=0$, then for each splitting $A=A_1+A_2$,
we have $A_i\cdot c_1(\Theta_{X/\kk})\not=0$ for some $i$, and hence
one of the two factors in the product
 $N^{A_1}_{p_1p_2s}N^{A_2}_{sp_3r}$
is defined to be $0$. Hence both sides
of \eqref{asseq} are zero. Thus we may assume that $A\cdot
c_1(\Theta_{X/\kk})=0$, so that the virtual dimension of
$\scrM(X,\beta,z)$ is $1$.

Fixing $\lambda\gg \mu \gg 0$ and taking
$y\in\overline{\scrM}_{0,4}$ to 
be the boundary point $D(x_1,x_2\,|\,x_3,x_{\out})$
and $y'$ to be the boundary point $D(x_2,x_3\,|\,x_1,x_{\out})$, 
Theorem~\ref{theorem: independence of modulus} now tells us
that $\deg [\scrM^{\ddagger}_y]^{\virt}=\deg [\scrM^{\ddagger}_{y'}]^{\virt}$.
However, combining Theorems~\ref{thm: main gluing theorem}
and~\ref{thm:main comparison}, we see that these two
degrees coincide with the left and right-hand sides of 
\eqref{asseq}, hence the desired equality.

\section{Forgetful maps: $\vartheta_0$ is the unit}
\label{sec:forgetful}
We give here the proof of Theorem~\ref{thm:unit}. The proof is
morally the same as the proof that the fundamental class of a smooth
variety is the unit in quantum cohomology. However, as usual, life is
made harder by difficulties in forgetting marked points in logarithmic
geometry.

\subsection{The case of curves}
We discuss forgetting marked points of punctured
log maps, following some of the ideas of \cite{AMW} and \cite{AW},
but generalising to the punctured case.

\begin{lemma}
\label{lem:forgetful punctured}
Let $C^\circ\rightarrow W$ be a punctured curve with punctures 
$x,x_1,\ldots,x_n$. Let $C^{\circ}_{\mathrm{for}}\rightarrow W$
be the curve obtained by forgetting the puncture $x$, and let
$\ul{\tau}:\ul{C}\rightarrow \ul{C}'$ be a partial stabilization
such that for each geometric point $\bar w\in W$, the induced map of pointed
curves
\[
\ul{\tau}:(\ul{C}_{\bar w},x,x_1,\ldots,x_n)\rightarrow
(\ul{C}'_{\bar w}, x_1,\ldots,x_n)
\]
is stable, i.e., the only component possibly contracted is the component
containing $x$. Then $\tau_*\shM_{C^{\circ}_{\mathrm{for}}}$ defines
a log structure $C'^{\circ}$ on $\ul{C}'$ making $C'^{\circ}\rightarrow W$
a punctured curve. 
%Further, let $Q:=\overline{\shM}_{W,\bar w}$ and
%$Q^{\circ}:=\overline{\shM}_{C^{\circ}_{\bar w},x_i}\subseteq Q\oplus\ZZ$.
%If $x_i\in \ul{C}_{\bar w}$ is contained
%in a component contracted by $\ul{\tau}$, then there is necessarily
%a unique node $q\in \ul{C}_{\bar w}$ contained in this component, and
%we may write $\overline{\shM}_{C^{\circ},q} = Q\oplus_{\NN}\NN^2$ with
%the maps $\NN\rightarrow Q,\NN^2$ given by $1\mapsto \rho_q\in Q$, $1\mapsto
%(1,1)$ respectively. Then 
%\begin{equation}
%\label{eq:pushdown monoid}
%\overline{\shM}_{C'^{\circ}_{\bar w},\tau(x_i)}=\left\{(m,r)\in Q\oplus\ZZ\,|\,
%(m+r\rho_q,r)\in Q^{\circ}\right\}.
%\end{equation}
%If $x_i$ is not contained in a contracted component, then
%$\overline{\shM}_{C'^{\circ}_{\bar w},\tau(x_i)}=Q^{\circ}$.
\end{lemma}

\begin{proof}
This is a slight variation of ideas in \cite{AMW},~App.~B. By 
\cite{AMW}, Lem.~B.3, both $\ul{\tau}_*\shM_{C_{\mathrm{for}}}$ and 
$\ul{\tau}_*\shM_{C^{\circ}_{\mathrm{for}}}$ define
log structures on $\ul{C}'$. (Here we continue to use the convention
that $C_{\mathrm{for}}$ is the log curve with marked points $x_1,\ldots,x_n$,
rather than punctured points.) Write these two log schemes as $C'$ and
$C'^{\circ}$, respectively; both are log schemes over $W$.

By \cite{AMW},~Lem.~B.5, $C'$ is log smooth over $W$. We claim 
that $C'$ has marked points $\ul\tau\circ x_i$, $1\le i\le n$, and that
$C'^{\circ}$ is a puncturing of $C'$.

We will first compute the stalks of the ghost sheaves of 
$\ul{\tau}_*\shM_{C_{\mathrm{for}}}$ and 
$\ul{\tau}_*\shM_{C^{\circ}_{\mathrm{for}}}$ at points in the image
of $\ul{\tau}\circ x_i$, $1\le i\le n$. 
As formation of $\tau_*\shM_{C^{\circ}_{\mathrm{for}}}$ commutes with strict base
change in $W$ by \cite{AMW},~Lem.~B.4, we may assume that $W$ is
a geometric point. Further, if $\ul{\tau}$ is an isomorphism
in a neighbourhood of $x_i$, then of course there is
nothing to compute.

Otherwise, $\ul{C}\rightarrow \ul{C}'$ contracts
an irreducible component $D$ containing $x$ and $x_i$, and necessarily
$D$ contains precisely one node of $\ul{C}_{\mathrm{for}}$.
Let $q$ be the node
of $\ul{C}$ contained in $D$. Let $y\in \ul{C}'$ be the image
of $D$ under $\ul{\tau}$. 

Note that $\overline{\shM}_{C'^{\circ},y}\subseteq (\ul{\tau}_*
\overline\shM_{C^{\circ}_{\mathrm{for}}})_y=\Gamma(D,
\overline\shM_{C^{\circ}_{\mathrm{for}}}|_D)$, and we have
\[
\overline{\shM}_{C'^{\circ},y}=\{\bar s\in
\Gamma(D,\overline{\shM}_{C^{\circ}_{\mathrm{for}}}|_{D})\,|\,
\shL^{\times}_{\bar s}\cong \O_D^{\times}\}.
\]
Let $\eta$ be the generic point of $D$. The generization map
$\chi_{x_i\eta}:Q^{\circ}\rightarrow Q$ is given by projection onto
the first factor, and the generization map
$\chi_{q\eta}:Q\oplus_{\NN}\NN^2\rightarrow Q$ is given by
$\left(m,(a,b)\right) \mapsto m + a\rho_q$, see \cite{JAMS},~Rmk.~1.2. 
Then
\begin{align*}
&\Gamma(D,\overline{\shM}_{C^{\circ}_{\mathrm{for}}}|_D) \\
= {} &
\left\{\big((m',r), (m,(a,b))\big)\in Q^{\circ}\oplus (Q\oplus_{\NN}\NN^2)
\,|\,\chi_{x_i\eta}(m',r)=
\chi_{q\eta}(m,(a,b))\right\}\\
= {} & 
\left\{\big((m+a\rho_q,r), (m,(a,b))\big)\in Q^{\circ}\oplus 
(Q\oplus_{\NN}\NN^2)\right\}.
\end{align*}

On the other hand, if $\bar s\in 
\Gamma(D,\overline{\shM}_{C^{\circ}_{\mathrm{for}}}|_D)$ is given
as $\big((m+a\rho_q,r),(m,(a,b))\big)$, then as in 
\cite{JAMS},~Lem.~1.14,
the degree of the corresponding line bundle $\shL_{\bar s}$ is
$-r+a-b$. Thus if this degree is $0$, we may write $\bar s$
as $\big((m+a\rho_q,r), (m,(a,a-r))
\big)$, necessarily with $(m+a\rho_q,r)\in Q^{\circ}$,
$m\in Q$ and $a,a-r\ge 0$. 
Thus we may write
\[
\overline{\shM}_{C'^{\circ},y} 
=\left\{ \big((m+a\rho_q,r), (m,(a,a-r))\,|\, (m+a\rho_q,r)\in Q^{\circ},
m\in Q, a,a-r\ge 0
\right\}.
\]
We may map $\overline{\shM}_{C'^{\circ},y}$
injectively into $Q\oplus \ZZ$ via
\begin{equation}
\label{eq:map without a name}
\big((m+a\rho_q,r), (m,(a,a-r))\big) \mapsto (m+(a-r)\rho_q, r).
\end{equation}
Identifying $\overline{\shM}_{C'^{\circ},y}$ with its image, we obtain
\begin{equation}
\label{eq:set without a name}
\overline{\shM}_{C'^{\circ},y} =\left\{
(m,r)\in Q\oplus\ZZ\,|\, (m+r\rho_q,r)\in Q^{\circ}\right\}.
\end{equation}
Indeed, any element of the image of \eqref{eq:map without a name}
lies in the right-hand side of \eqref{eq:set without a name}. Conversely,
if $(m,r)$ lies in the right-hand side of \eqref{eq:set without a name}, and
if $r\ge 0$, then $(m,r)$ is the image of $\big((m+r\rho_q,r),(m,(r,0))\big)$,
while if $r<0$, then $(m,r)$ is the image of
$\big((m+r\rho_q,r),(m+r\rho_q,(0,-r))\big)$.
%This gives the desired description in the statement of the lemma.

Note under this identification, the map $Q=\overline{\shM}_W
\rightarrow \overline\shM_{C'^{\circ},y}$ is given by $m\mapsto (m,0)$,
and that if $x_i$ were itself a marked point, so that $Q^{\circ}=Q\oplus
\NN$, then also $\overline{\shM}_{C'^{\circ},y}=Q\oplus\NN$. 

Returning to general $W$, applying the above discussion to
$\ul{\tau}_*\shM_{C_{\mathrm{for}}}$, since $C'\rightarrow W$ is log smooth, 
necessarily
$\ul{\tau}\circ x_i$ is a marked point in the log structure of $C'$. 
Further, from the description of $\overline{\shM}_{C'^{\circ},y}$,
one sees that $C'^{\circ}$ is a puncturing of $C'$.
\end{proof}

In the following lemma,
we reverse this procedure, starting with a punctured curve
$\pi':C'^{\circ}\rightarrow W$ and a choice of logarithmic section $x'$.
Stabilization of this family modifies the curve when a point in the image
of $x'$ is a node or lies in the image of $x_i$ for some $i$. In the
latter case, a $\PP^1$ is bubbled off, producing a new node. At the
logarithmic level, the smoothing parameter is determined by
$\rho_i$ defined in (4) of the following lemma.

\begin{lemma}
\label{lem:punctured stabilization}
Let $\pi':C'^{\circ}\rightarrow W$ be a punctured curve with
punctures $x_1,\ldots,x_n$
and let $x':W\rightarrow C'^{\circ}$ be a section of $\pi'$
as a morphism of log schemes. Then there exists a
punctured curve $C^{\circ}\rightarrow W$ with punctures
$x_1,\ldots,x_n$ and marked point $x$, unique up to choice of puncturing,
along with a diagram
\begin{equation}
\label{eq:forgetting diagram}
\xymatrix@C=30pt
{
C^\circ\ar[r]\ar@/^2pc/[rr]^{\tau'}
&C^{\circ}_{\mathrm{for}}\ar[r]_{\tau}& C'^{\circ}
}
\end{equation}
where the first morphism forgets the marked point $x$, satisfying
the following properties:
\begin{enumerate}
\item 
$\ul{\tau}:(\ul{C},x,x_1,\ldots,x_n)\rightarrow (\ul{C}',x_1,\ldots,x_n)$
is a partial stabilization.
\item 
For each geometric point $\bar w\in W$, the map of marked
curves induced by $C^{\circ}_{\bar w}\rightarrow (C')^{\circ}_{\bar w}$
\[
\ul{\tau}:(\ul{C}_{\bar w},x,x_1,\ldots,x_n)\rightarrow 
(\ul{C}'_{\bar w},x_1,\ldots,x_n)
\]
is stable.
\item If $x:W\rightarrow C^{\circ}_{\mathrm{for}}$ is the 
logarithmic section corresponding to the marked point $x$, then
$x'=\tau\circ x$.
\item Let $\bar s_i$ be the global
section of $\overline{\shM}_{C'}$ which is $0$ away
from the image of $x_i$ and has germ $(0,1)\in 
\overline{\shM}_{C',x_i(\bar w)} =\overline{\shM}_{W,\bar w}\oplus \NN$
for any $\bar w$ a geometric point of $W$. 
Define $\rho_i \in \Gamma(W,\overline{\shM}_W)$ to be 
\[
\rho_i = \overline{(x')}^{\flat}(\bar s_i).
\]
Then identifying $\ul{W}$ with
its image in $\ul{C}$ and $\ul{C}'$ under $x_i$, we have an isomorphism
of line bundles
\[
\shN^{\vee}_{x_i/C} \otimes \shL_{\rho_i}\cong \shN^{\vee}_{x_i/C'}.
\]
Further, if $W_i^{\circ}$, $W_i'^{\circ}$ denote the log schemes
obtained by restricting the log structures of $C^{\circ}$ and $C'^{\circ}$ 
to the images of $x_i$, the restricted morphism $\tau':W_i^{\circ}
\rightarrow W_i'^{\circ}$ induces, on the level of ghost sheaves,
for $(m,r) \in \overline{\shM}_{W_i'^{\circ},\bar w}\subseteq
\overline{\shM}_{W,\bar w}\oplus\ZZ$, $\bar w$ a geometric point of $W$,
\[
(\bar\tau')^{\flat}(m, r) = (m+r\rho_i, r)
\in
\overline{\shM}_{W_i^{\circ},\bar w}\subseteq
\overline{\shM}_{W,\bar w}\oplus\ZZ.
\]
\end{enumerate}
\end{lemma}

\begin{proof}
{\bf The diagram \eqref{eq:forgetting diagram} at the level of marked
log curves.}
The existence of the diagram \eqref{eq:forgetting diagram} at the level of
underlying schemes is classical, see \cite{Knud},~Thm.~2.4,
or \cite{ACG}, Chap.~10,\S8.
Further, if $C'$ is the corresponding log smooth (marked) curve over $W$,
then $x'$ also induces a section $x':W\rightarrow C'$.
Then the diagram \eqref{eq:forgetting diagram} exists with
$C^{\circ}_{\mathrm{for}}$, $C^{\circ}$ replaced with the corresponding
marked log smooth curves $C_{\mathrm{for}}, C$. This is 
standard, the morphism $C_{\mathrm{for}}\rightarrow C'$ being seen as a
logarithmic modification. As we need a more detailed analysis 
at marked points, we
explain how this may be described as a logarithmic blow-up when the images 
of $x'$ and $x_i$ coincide.

Let $\bar s_i$, $\rho_i$ be as in item (4) of the lemma.
Let $w\in W$ be a scheme-theoretic point with $x_i(w)=x'(w)$,
and replace $W$ with any Zariski neighbourhood of $w$ so that
the image of $x'$ is disjoint from the image of $x_j$ for any 
$j\not=i$. 
Further shrinking $W$, we may assume the image of $x'$ 
is disjoint from any node of any fibre of $\pi'$. 
Then $\overline{(\pi')}^{\flat}(\rho_i), \bar s_i$
generate a sheaf of monoid ideals in $\overline{\shM}_{C'}$. The pull-back
of this sheaf of ideals under the quotient map $\shM_{C'}\rightarrow
\overline\shM_{C'}$ is then a coherent sheaf of monoid ideals $\shK$,
and we take $C_{\mathrm{for}}\rightarrow C'$ to be the logarithmic
blow-up of $\shK$ (see \cite{FKato},~\S3 or \cite{Ogus},~III,\S2.6). Note that 
the pull-back ideal $(x')^{\bullet}\shK$ is locally principal, so by
the universal property of log blow-up (see \cite{FKato},~Def.~3.8), 
$x'$ factors as
\[
\xymatrix@C=20pt
{
W\ar[r]^{x}&C_{\mathrm{for}}\ar[r]&C'.
}
\]
As a log blow-up is log \'etale, the composed map 
$C_{\mathrm{for}}\rightarrow W$ is still log smooth. To verify that
this is a logarithmic curve, a local calculation, which we now
carry out, suffices to check that all fibres are reduced and of dimension one.

\'Etale locally around $w$ we may write $W$ as $\Spec A$ with chart 
$\alpha:Q\rightarrow A$, and then \'etale locally at $x_i(w)$, 
we may write $C'=\Spec A[u]$, with chart $Q\oplus\NN\rightarrow A[u]$
given by $(m,r)\mapsto \alpha(m)u^r$. Further, the section $x'$
induces a monoid homomorphism $Q\oplus \NN\rightarrow Q$ given by
$(m,r)\mapsto m+r \rho_i$, and necessarily $\ul{x}'$ is then induced
by the $A$-algebra homomorphism
$A[u]\rightarrow A$ defined by $u\mapsto \varphi\cdot \alpha(\rho_i)$
for some invertible function $\varphi$. The logarithmic ideal $\shK$
corresponds to the monoid ideal of $Q\oplus\NN$ generated by
$(\rho_i,0)$ and $(0,1)$, and the underlying scheme of the
logarithmic blow-up is then 
\[
\Proj A[u][v,w]/(vu-w\alpha(\rho_i)),
\]
with $v,w$ of degree $1$,
covered by the open affines 
\begin{align*}
U_1:={} & \Spec A[u,v]/(vu-\alpha(\rho_i))\\
U_2:={} &  \Spec A[u,w]/(u-w\alpha(\rho_i))=\Spec A[w],
\end{align*}
with logarithmic structures given by charts
\begin{align*}
Q_1:=(Q\oplus\NN) + \NN(\rho_i,-1)  \rightarrow & A[u,v]/(vu-\alpha(\rho_i))\\
Q_2:=(Q\oplus\NN) + \NN(-\rho_i,1)  \rightarrow & A[u,w]/(u-w\alpha(\rho_i)).
\end{align*}
Here the monoids $Q_1,Q_2$ are submonoids of $Q^{\gp}\oplus\ZZ$. In each of
the two cases, the chart takes $(m,r)\in Q\oplus\NN$ to $\alpha(m)u^r$. 
In the first case, $(\rho_i,-1)\mapsto v$ and in the second case
$(-\rho_i,1)\mapsto w$.
Note further that there is an isomorphism
$Q\oplus_{\NN}\NN^2 \cong Q_1$ (with the push-out defined using
the map $\NN\rightarrow Q$ given by $1\mapsto \rho_i$). This isomorphism
is given by the map $\big(m,(a,b)\big)\mapsto (m+b\rho_i,a-b)$.
There is also an isomorphism 
\begin{align}
\label{eq:Q2 isomorphism}
\begin{split}
Q\oplus\NN \rightarrow {}  &  Q_2\\
(m,r)\mapsto  {} & (m-r\rho_i,r).
\end{split}
\end{align}
Finally, note that the lift $x$ of the section 
$x'$ has image in $U_2$, and $x:\Spec A\rightarrow U_2$ is defined by
$w\mapsto \varphi$.

From the above explicit descriptions we see that indeed 
$C_{\mathrm{for}}\rightarrow W$ is a log smooth curve. Marking
the image of the section of $x$ gives $C$. This now
gives the diagram \eqref{eq:forgetting diagram} at the level of log
smooth curves, and now we need to construct the diagram for punctured
curves.

\medskip

{\bf The diagram \eqref{eq:forgetting diagram} at the level of punctured
curves.}
We define
\[
C^{\circ}_{\mathrm{for}}:=C_{\mathrm{for}}\times^{\fine}_{C'} C'^{\circ}.
\]
We need to check this is a puncturing of
$C_{\mathrm{for}}$. Note this just needs
to be checked \'etale locally in neighbouhoods of points $\bar p$ of
$\ul{C'}$ which lie in the image of both $x'$ and $x_i$, and thus
we may use the previous local description for $C'$ and $C_{\mathrm{for}}$
to calculate the fibre product. Write the stalk of the ghost
sheaf of $C'$ at $\bar p$ as $Q\oplus\NN$ and the stalk of the ghost
sheaf of $C'^{\circ}$ at $\bar p$ as $Q^{\circ}\subseteq Q\oplus\ZZ$.
Noting the inclusion $Q\oplus\NN \hookrightarrow Q^{\circ}$
yields an isomorphism at the level of groups, we may use
Proposition \ref{prop:fibre product properties}, (3), to identify
$Q_i\oplus^{\fine}_{Q\oplus\NN} Q^{\circ}$ with a submonoid of
$Q_i^{\gp}$. 

Explicitly, for $Q_1$, we obtain
\[
Q_1\oplus^{\fine}_{Q\oplus\NN} Q^{\circ} = Q_1.
\]
Indeed, this follows if $Q^{\circ}\subseteq Q\oplus\ZZ$ is
contained in $Q_1=(Q\oplus\NN)+
\NN(\rho_i,-1)\subseteq Q\oplus\ZZ$. But the section $x'$ induces
a map $Q^{\circ}\rightarrow Q$ given by $(m,r)\mapsto m+r \rho_i\in Q$, 
so if $(m,r)\in Q^{\circ}$ with $r<0$, we may write 
$(m,r)= (m+r\rho_i,0)-r(\rho_i,-1)\in Q_1$.

Using the identification of $Q_2$ with $Q\oplus\NN$ of 
\eqref{eq:Q2 isomorphism}, we have
\begin{equation}
\label{eq:blown up puncture}
Q_2\oplus^{\fine}_{Q\oplus \NN} Q^{\circ} = (Q\oplus \NN)+
\{(m,r)\in Q\oplus\ZZ\,|\, (m-r\rho_i,r)\in Q^{\circ}\} \subseteq Q\oplus\ZZ.
\end{equation}

Now the fine fibre product is constructed by integralizing the
log fibre product, which on the two charts $U_i:=\Spec B_i$ 
has underlying scheme $\Spec B_i\times_{A[u]} A[u] =\Spec B_i$, and 
log structure given by strict morphisms
$\psi_i:\Spec B_i\times_{A[u]} A[u] \rightarrow
A_{Q_i\oplus_{Q\oplus\NN} Q^{\circ}}$. The integralization is
obtained by base-change via the closed embedding
$A_{Q_i\oplus^{\fine}_{Q\oplus\NN} Q^{\circ}} \rightarrow
A_{Q_i\oplus_{Q\oplus\NN} Q^{\circ}}$. However it is easy to check that
in fact the morphisms $\psi_i$ factor through these closed embeddings
by showing the charts $Q_i\oplus_{Q\oplus \NN}Q^{\circ}\rightarrow
B_i\otimes_{A[u]} A[u]$ for the log fibre product factor through the
integralization morphism $Q_i\oplus_{Q\oplus \NN}Q^{\circ}\rightarrow
Q_i\oplus^{\fine}_{Q\oplus \NN}Q^{\circ}$.
So integralization is an isomorphism. Hence we see
that $\ul{C}^{\circ}_{\mathrm{for}}=\ul{C}_{\mathrm{for}}$,
and the form of the push-outs $Q_i\oplus^{\fine}_{Q\oplus\NN} Q^{\circ}$
shows that $C^{\circ}_{\mathrm{for}}$ is indeed a punctured
curve, with ghost sheaf at the image of $x_i$ being given by
\eqref{eq:blown up puncture}. 

This construction now provides diagram \eqref{eq:forgetting diagram}
at the level of punctured curves, satisfying conditions (1)--(3).

\medskip

{\bf Proof of (4).} The given form for $(\bar\tau')^{\flat}$ follows
from the description of the local models above, with 
$(\bar\tau')^{\flat}$
described at the level of the marked curves as the inclusion $Q\oplus\NN
\hookrightarrow Q_2=(Q\oplus \NN)+\NN(-\rho_i,1) \cong Q\oplus\NN$ 
where the latter isomorphism is given by \eqref{eq:Q2 isomorphism}.
Now writing $\overline{\shM}_{W_i^{\circ}},
\overline{\shM}_{W_i'^{\circ}}
\subseteq \overline\shM_W\oplus \ZZ$, the line bundle corresponding
to the section $(0,1)$ of $\overline{\shM}_{W_i'^{\circ}}$ is
$\shN^{\vee}_{x_i/C'}$, which is then isomorphic to the line bundle
corresponding to the section $(\bar\tau')^{\flat}(0,1)=(\rho_i,1)$
of $\overline\shM_{W_i^{\circ}}$,
i.e., to $\shN^{\vee}_{x_i/C}\otimes \shL_{\rho_i}$, as claimed.
\end{proof}

\begin{remark}
\label{rem:forgetful tropical 1}
The sections $\rho_i$ of $\overline{\shM}_W$ can be interpreted 
tropically as follows. Assume $W$ is a geometric log point, and let
$G$, $G'$ be the dual intersection graphs of $C^{\circ}$, $C'^{\circ}$
respectively.
If the section $x'$ is disjoint from nodes of $C'^{\circ}$
or punctures $x_1,\ldots,x_n$, then $G$ is obtained from $G'$ 
simply by adding an additional unbounded leg attached to some
vertex of $G'$. If, on the
other hand, $x$ coincides with $x_i$, then $G$ is obtained from $G'$
by subdividing the leg of $G$ corresponding to $x_i$ by adding a new
vertex to this leg, and attaching an additional leg to this vertex.
In this case, $\rho_i$ can be viewed as the length of the new edge
created by this process.
\end{remark}

\subsection{The case of punctured maps}

We now turn to forgetful maps for punctured maps. In this case,
we may only forget marked points with contact order $0$. Thus we consider
for the remainder of this subsection the following fixed situation.
Let $\beta$ be a class of punctured curve for a target space $X$
with punctures $x,x_1,
\ldots,x_n$ with
$x$ having contact order given by $0\in \Sigma(X)(\ZZ)$. 
Let $\beta'$ be the class 
of punctured curve with the same underlying curve class and
genus, but with only punctured points $x_1,\ldots,x_n$. Let ${\bf x}\subseteq
\{x_1,\ldots,x_n\}$ be a subset of these marked points.

To understand the forgetful map at the virtual level,
we need an intermediate space of maps to the Artin fan $\shX$. 
This is a minor variant of
the stack $\foM'(\shY\rightarrow\shX)$ of \cite{AW},~\S3.

\begin{definition}
We define $\foM(\shX,\beta/\beta')$ to be the log stack as follows.
The objects of this stack over a log scheme $W$ are
pairs of pre-stable punctured maps 
\[
\hbox{$f:C^{\circ}/W\rightarrow \shX$ and
$f':C'^{\circ}/W\rightarrow\shX$}
\]
of class $\beta$ and $\beta'$
respectively along with a commutative diagram
\[
\xymatrix@C=30pt
{
C^{\circ}\ar[rd]^f\ar[d]\ar@/_2pc/[dd]_{\tau'}&\\
C^{\circ}_{\mathrm{for}}\ar[d]^{\tau}&\shX\\
C'^{\circ}\ar[ur]_{f'}
}
\]
As in Lemmas 
\ref{lem:forgetful punctured},\ref{lem:punctured stabilization},
$C^{\circ}_{\mathrm{for}}$ is the same curve as $C^{\circ}$ without
the marked point $x$, $C^{\circ}\rightarrow C^{\circ}_{\mathrm{for}}$
is the canonical morphism induced by the obvious inclusion of
log structures, and $\ul{\tau}$ satisfies conditions (1) and (2)
of Lemma \ref{lem:punctured stabilization}.

As usual, set 
\[
\foM^{\ev({\bf x})}(\shX,\beta/\beta')=
\foM(\shX,\beta/\beta')\times_{\prod_{x_i \in {\bf x}} \ul{\shX}}
\prod_{x_i\in {\bf x}}\ul{X}.
\]
\end{definition}

The stack $\foM(\shX,\beta/\beta')$ is an algebraic log stack via precisely
the same argument given in \cite{AW},~\S3 for $\foM'(\shY\rightarrow
\shX)$. Alternatively, the following lemma gives an explicit
description of $\foM^{\ev({\bf x})}(\shX,\beta/\beta')$.
Also, $\foM(\shX,\beta/\beta')$ carries the basic log
structure inherited from
the log morphism $f:C^{\circ}\rightarrow\shX$, as this morphism determines
$f'$: see the argument of \cite{AW},~Lem.~4.1.

For brevity of notation, write
\[
\foM:=\foM^{\ev({\bf x})}(\shX,\beta/\beta'), \quad
\foM':=\foM^{\ev({\bf x})}(\shX,\beta').
\]

\begin{theorem}
\label{thm:forgetful}
In the above situation, assume the underlying curve class
$A$ of $\beta$ is non-zero. Then there is a diagram
\[
\xymatrix@C=30pt
{
\scrM(X,\beta)\ar[r]^{\phi'}\ar[d]_{\varepsilon}&\scrM(X,\beta')
\ar[d]^{\varepsilon'}\\
\foM\ar[d]_{\psi}\ar[r]^{\phi} & 
\foM'\\
\foM^{\ev({\bf x})}(\shX,\beta)
}
\]
such that:
\begin{enumerate}
\item $\psi$ is strict \'etale.
\item Let $\pi':C'^{\circ}_{\foM'} \rightarrow \foM'$
be the universal prestable punctured curve over
$\foM'$. Then there is a factorization
\[
\xymatrix@C=30pt
{
\foM
\ar[r]^>>>>{\mathrm{sat}}&
C'^{\circ}_{\foM'}
\ar[r]^>>>>>{\pi'}& \foM'
}
\]
of $\phi$ where
$\mathrm{sat}$ is the saturation of 
$C'^{\circ}_{\foM'}$.
In particular, $\phi$ is proper and representable.
\item The square is cartesian in all categories, and $\varepsilon'$
and $\psi\circ\varepsilon$ are the standard morphisms 
as in \S\ref{subsub:moduli}. 
The relative obstruction theory for $\varepsilon'$ pulls back to give
a relative obstruction theory for $\varepsilon$ which is compatible with
the standard relative obstruction theory for $\psi\circ\varepsilon$.
\end{enumerate}
\end{theorem}

\begin{proof}
{\bf Construction of $\psi$ and (1).} 
Given a log morphism $W\rightarrow \foM$,
corresponding to the data of $f, f'$ and $\tau$ over
$W$, we may forget the data of $f'$ and $\tau$ to get
$W\rightarrow \foM^{\ev({\bf x})}(\shX,\beta)$.
This defines $\psi$ as a morphism of log stacks.
Because the basic log structure on $\foM$ is basic for
$f$, this makes $\psi$ strict. It
is proved to be log \'etale precisely as in \cite{AW},~\S5.1. This
gives (1).

\medskip

{\bf Construction of $\phi'$.}
As the contact order of $x$ is $0$, we see from \eqref{eq:ux def}
that the image of $\bar f^{\flat}$ at $x$ is contained in
$\overline\shM_{W,\bar w}$, for any geometric point $\bar w$ of $W$.
Thus we may forget the marked point without affecting the existence of the
morphism $f$. So we may replace $C^{\circ}$ with $C^{\circ}_{\mathrm{for}}$.
One then has a partial stabilization map $\ul{\tau}:
\ul{C}\rightarrow
\ul{C}'$ defined over $\ul{W}$ so that
$\ul{f}$ factors as $\ul{f}'\circ \ul{\tau}$ for an ordinary stable map
$\ul{f}':\ul{C}'
\rightarrow X$ with marked points $x_1,\ldots,x_n$. Necessarily,
over each geometric point $\bar w$ of $W$, $\ul{\tau}$ contracts
at most the irreducible component of $\ul{C}_{\bar w}$ containing $x$.
Note that here we use $A\not=0$ as otherwise $\ul{C}$ may
be contracted to a point.

Thus we are in the situation of Lemma \ref{lem:forgetful punctured},
obtaining a log structure $C'^{\circ}$ on $\ul{C}'$ via
$\shM_{C'^{\circ}}=\ul{\tau}_*\shM_{C^{\circ}_{\mathrm{for}}}$,
making $C'^{\circ}\rightarrow W$ into a punctured curve.
Further, by adjunction, the morphism $f:C^{\circ}_{\mathrm{for}}
\rightarrow X$ induces a morphism $f':C'^{\circ}\rightarrow X$.
The punctured map $f'$ need not be pre-stable, but
after shrinking the puncturing $C'^{\circ}$ via
\cite{ACGS18}, Prop.~2.5, we may assume also that $f'$ is pre-stable.

This gives the morphism $\phi'$.

\medskip

{\bf The morphism $\phi$.}
The morphism $\phi$ takes the data of $f, f'$ and $\tau$ 
over $W$ to the punctured map $f'$. Note this morphism need not be strict,
as the basic log structure on $\foM$ need not be basic for $f'$.

\medskip

{\bf Proof of (2).}
We first construct the factorization of $\phi$ 
through $\pi'$. Let 
\[
\tau':C^{\circ}_{\foM}\mapright{}
C^{\circ}_{\foM,\mathrm{for}}\mapright{\tau} C'^{\circ}_{\foM}
\]
be the universal data for $\foM$. Note that 
\[
C'^{\circ}_{\foM}=\foM
\times^{\fine}_{\foM'} C'^{\circ}_{\foM'}
\]
by \cite{ACGS18},~Prop.~2.16, (2). 
We can view the section $x:\ul{\foM}\rightarrow 
\ul{C}^{\circ}_{\foM}$ as inducing a strict log morphism
$x:\foM\rightarrow C^{\circ}_{\foM,\mathrm{for}}$, as $x$ is not marked
on $C^{\circ}_{\foM,\mathrm{for}}$. Hence we obtain 
\[
\tau\circ x:\foM\rightarrow C'^{\circ}_{\foM}.
\]
Projection to $C'^{\circ}_{\foM'}$ then gives the desired morphism
$\sat:\foM\rightarrow C'^{\circ}_{\foM'}$ such that
$\phi=\pi'\circ\sat$.

As $\foM$ is saturated, we thus obtain a factorization of $\tau\circ x$ as
\[
\foM\rightarrow C'^{\circ,\sat}_{\foM'}\rightarrow C'^{\circ}_{\foM'},
\]
and we need to construct an inverse to this first morphism to show
it is an isomorphism. To do so, we consider a test fs
log scheme $T$ along with a morphism $g:T\rightarrow C'^{\circ}_{\foM'}$: 
by the universal 
property of saturation, $g$ must factor uniquely through the saturation
of $C'^{\circ}_{\foM'}$. We need to construct a morphism
$T\rightarrow \foM$, which we do by constructing the relevant data over
$T$.

By composition of $g$ with $\pi'$, we may pull back 
$C'^{\circ}_{\foM'}$ to $T$ to obtain 
\[
C'^{\circ}:=T\times^{\fine}_{\foM'}C'^{\circ}_{\foM'},
\]
along with a pre-stable punctured map $f':C'^{\circ}\rightarrow\shX$.
Note there is a section $x':T\rightarrow C'^{\circ}$
induced by $\id:T\rightarrow T$ and 
$g:T\rightarrow C'^{\circ}_{\foM'}$. 
Thus, by Lemma \ref{lem:punctured stabilization}, we obtain
the data of \eqref{eq:forgetting diagram}. We may then
take $f=f'\circ\tau'$; it follows from the explicit description
\eqref{eq:blown up puncture} of the puncturing monoids on $C^{\circ}$
that pre-stability of $f'$ implies pre-stability of $f$.
Thus we obtain an object of $\foM$, defining the desired inverse
to the morphism $\mathrm{sat}$. That this morphism is indeed the inverse to
the morphism $\mathrm{sat}$ follows immediately by construction.
Certainly, saturation is always 
finite and representable, and $\pi'$ is proper and representable (by
algebraic spaces),
being a family of curves. Thus $\phi$ is proper and representable (by
algebraic spaces).
This completes the proof of (2).

\medskip

{\bf Proof of (3).}
The morphism $\varepsilon$ takes a stable punctured map $f:C^{\circ}\rightarrow
X$ with $\phi(f)$ given by $f':C'^{\circ}\rightarrow X$ to the pair of the 
compositions $C^{\circ}\rightarrow X\rightarrow \shX$ and
$C'^{\circ}\rightarrow X\rightarrow \shX$, with $\tau:C^{\circ}
\rightarrow C'^{\circ}$ the stabilization map.
In particular, $\varepsilon$ is
strict. The square of the theorem is then cartesian exactly as in the
proof of \cite{AW},~Prop.~1.6.3. 

The obstruction theory for $\varepsilon'$ is given, for
$\pi':\scrC'\rightarrow \scrM(X,\beta')$ the universal curve
and $f':\scrC'\rightarrow X$ the universal map,
by a morphism 
\[
\left[R\pi'_*(f')^*\Theta_{X/\kk}(-\sum_{x_i\in {\bf x}} x_i)\right]^{\vee}
\rightarrow \LL_{\scrM(X,\beta')/\foM^{\ev({\bf x})}(\shX,\beta')},
\]
so, using the projection formula and $R\tau'_*\O_C=\O_{C'}$, the pull-back 
obstruction theory is similarly given by
\[
\left[R\pi_*f^*\Theta_{X/\kk}(-\sum_{x_i\in {\bf x}} x_i)\right]^{\vee}
\rightarrow \LL_{\scrM(X,\beta)/\foM^{\ev({\bf x})}(\shX,
\beta/\beta')}\cong \LL_{\scrM(X,\beta)/\foM^{\ev({\bf x})}(\shX,\beta)},
\]
the latter
isomorphism as $\psi$ is strict log \'etale, hence \'etale. One may
check
that the resulting obstruction theory for $\psi\circ\varepsilon$ then agrees
with the relative obstruction theory for
$\scrM(X,\beta)\rightarrow\foM^{\ev({\bf x})}(\shX,\beta)$ precisely
as in \cite{AW},~\S6.2
\end{proof}

\subsection{Forgetful maps with point constraints and the proof of
Theorem \ref{thm:unit}}

We now consider a curve class $\beta$ as in Theorem~\ref{thm:unit}:
assume $A\not=0$ and $\beta$ has three punctured points
$x_1,x_2,x_{\out}$ with contact orders $0,p$ and $-r$ respectively,
with $p,r\in B(\ZZ)$.

We can then apply Theorem~\ref{thm:forgetful} to get the diagram
given there, taking ${\bf x}=\{x_{\out}\}$, with $\phi$ forgetting
the marked point $x_1$. We write $\ev$
instead of $\ev({\bf x})$ as is usual.
 
We wish to obtain an analogous diagram for the corresponding moduli
spaces with point constraints. For this purpose, note we have
evaluation morphisms 
\begin{align}
\label{eq:two evs}
\begin{split}
\ev_{\shX}:\foM^{\ev}(\shX, \beta/\beta')\rightarrow &\scrP(X,r)\\
\ev'_{\shX}:\foM^{\ev}(\shX, \beta')\rightarrow &\scrP(X,r),
\end{split}
\end{align}
with $\ev'_{\shX}$ the usual evaluation morphism at $x_{\out}$ given by 
Definition~\ref{def:evaluation maps}, while $\ev_{\shX}$
is given by evaluation at $x_{\out}
\in C^{\circ}$ (rather than $x_{\out}\in C'^{\circ}$), using 
Proposition~\ref{prop:evaluation maps general}. Unfortunately,
we do not have $\ev_{\shX}'\circ\phi=\ev_{\shX}$.
This arises 
due to the change of normal bundle at $x_{\out}$ under stabilization,
see Lemma~\ref{lem:punctured stabilization}, (4). Thus, the argument is a 
bit more subtle.
We first examine the difference between $\ev_{\shX}'$ and $\ev_{\shX}$
at the ghost sheaf level.

A key ingredient for this is the section $\rho_{\out}$ of the ghost sheaf of
$\foM^{\ev}(\shX,\beta/\beta')$ obtained from 
Lemma~\ref{lem:punctured stabilization}, (4), with $\tau\circ x_1$ playing
the role of the section $x'$ and $x_{\out}$ playing the role of $x_i$.

\begin{lemma}
\label{lem:ev comparison}
In the setup of this subsection, 
we may identify the ghost sheaf of $\scrP(X,r)$ with $\overline{\shM}_{Z_r}$
and view $r$ as defining a homomorphism $r:\overline{\shM}_{Z_r}\rightarrow
\NN$ as in Remark~\ref{rem:contact order remarks},(1). Then in the notation of \eqref{eq:two evs} and Theorem 
\ref{thm:forgetful},
for $\bar s$ a section of $\overline\shM_{\scrP(X,r)}$, we have
\[
\overline{\ev'_{\shX}\circ\phi}^{\flat}(\bar s)=\overline{\ev}_{\shX}^{\flat}
(\bar s) + r(\bar s)\rho_{\out}.
\]
\end{lemma}

\begin{proof}
Fix a morphism $W\rightarrow \foM^{\ev}(\shX,\beta/\beta')$.
This means we are given the
data of $C^{\circ}\rightarrow C'^{\circ}$, compatible maps 
$f:C^{\circ}\rightarrow \shX$, $f':C'^{\circ}\rightarrow \shX$,
and a factorization of $\ul{f}\circ x_{\out}$ through 
$\ul{X}\rightarrow\ul{\shX}$.
In particular, if $W^{\circ}$, $W'^{\circ}$ denote the
log structures on $\ul{W}$ obtained by pull-back via $x_{\out}$
from $C^{\circ}$ and $C'^{\circ}$ respectively, we obtain a commutative
diagram
\begin{equation}
\label{eq:TTprime diag}
\xymatrix@C=30pt
{
W^{\circ}\ar[r]^{f_W}\ar[d]_{\tau}&Z\\
W'^{\circ}\ar[ur]_{f'_W}&
}
\end{equation}
with $Z=Z_r$.

The morphism $f_W$, $f'_W$ yield evaluation maps
$\ev,\ev':W\rightarrow \scrP(X,r)$ via 
Proposition~\ref{prop:evaluation maps general}. Note that then
$\ev$ and $\ev'$ coincide with the composition of $W\rightarrow 
\foM^{\ev}(\shX,\beta/\beta')$ with $\ev_{\shX}$ and $\ev'_{\shX}\circ\phi$,
respectively. Thus it is sufficient to understand the relationship
between $\ev$ and $\ev'$, which we do by following the
construction of Proposition~\ref{prop:evaluation maps general}.

We first construct $\widetilde\ev:W^{\circ}\rightarrow \widetilde\scrP(X,r)
=Z\times B\GG_m^{\dagger}$, $\widetilde\ev':W'^{\circ}\rightarrow 
\widetilde\scrP(X,r)$
as $\widetilde\ev=f_W\times \beta$, $\widetilde\ev'=f'_W\times \beta'$,
where the morphisms $\beta$ and $\beta'$, with target $B\GG_m^{\dagger}$,
are determined  by
sections $(0,1)\in \overline{\shM}_{W^{\circ}}\subseteq
\overline{\shM}_W\oplus\ZZ$, $(0,1)\in \overline{\shM}_{W'^{\circ}}
\subseteq \overline\shM_W\oplus\ZZ$ respectively. We may write,
for $\bar s$ a section of  $(f'_W)^{-1}\overline{\shM}_Z$,
\[
\overline{f'_W}^{\flat}(\bar s)= (\varphi(\bar s),-r(\bar s))
\]
for some homomorphism $\varphi:(f_W')^{-1}\overline\shM_Z\rightarrow 
\overline{\shM}_W$.
Then for $(\bar s,n)\in (f'_W)^{-1}\overline\shM_Z\oplus\NN$, we have
\[
\overline{(\widetilde\ev')}^{\flat}(\bar s,n)=(\varphi(\bar s),n-r(\bar s))
\]

On the other hand, 
\begin{align*}
\overline f_W^{\flat}(\bar s) = {} & \overline\tau^{\flat}\circ 
\overline f'^{\flat}_W(\bar s)\\
= {} & \overline\tau^{\flat}(\varphi(\bar s),-r(\bar s))\\
= {} & (\varphi(\bar s)-r(\bar s)\rho_{\out},-r(\bar s)),
\end{align*}
the last equality by
Lemma \ref{lem:punctured stabilization}, (4).
Thus
\[
\overline{\widetilde\ev}^{\flat}(\bar s,n)=
(\varphi(\bar s)-r(\bar s)\rho_{\out}, n-r(\bar s)).
\]
Restricting to $\overline\shM_{\scrP(X,r)}\subseteq
\overline\shM_{\widetilde\scrP(X,r)}$, we consider sections
$(\bar s,r(\bar s))$ instead of the more general $(\bar s,n)$. This immediately
gives the claimed result.
\end{proof}

\begin{remark}
Continuing with the tropical interpretation of 
Remark~\ref{rem:forgetful tropical 1}, the above comparison can
be interpreted as follows. Again assume $W$ is a geometric log point
with a morphism $W\rightarrow \foM^{\ev}(\shX,\beta/\beta')$
giving data of $C^{\circ},C'^{\circ}$, etc.
Suppose further that $\tau':C^{\circ}\rightarrow
C'^{\circ}$ is such that $x_1$ and $x_{\out}$ lie 
in the same irreducible
component of $C^{\circ}$, and this component is contracted by $\tau'$.
Tropically, the maps $f:C^{\circ}\rightarrow \shX$, $f':C'^{\circ}\rightarrow
\shX$ then give rise to two families of tropical maps
$h_s:G\rightarrow \Sigma(X)$, $h'_s:G'\rightarrow \Sigma(X)$. 
Let $v'_{\out}$ be the vertex of $G'$  adjacent to the leg $L'_{\out}$
corresponding to $x_{\out}\in C'^{\circ}$, and similarly let $v_{\out}$ be the
vertex of $G$ adjacent to the leg $L_{\out}$ corresponding to 
$x_{\out} \in C^{\circ}$. Then as in Remark~\ref{rem:forgetful tropical 1},
$v_{\out}$ can be viewed as lying in the interior of
$L'_{\out}$, and $\rho_{\out}$ is the length of the segment of
$L'_{\out}$ joining $v'_{\out}$ and $v_{\out}$. Now by 
Lemma~\ref{lem:tropical point constraint}, $\Sigma(\ev_{\shX})(s)
=h_s(v_{\out})$, but by the above discussion $h_s(v_{\out})=
h_s(v'_{\out})-r \rho_{\out}(s)= \Sigma(\ev'_{\shX})(s)-r\rho_{\out}(s)$.
This is the transpose of the relation of Lemma~\ref{lem:ev comparison}.
\end{remark}

\begin{theorem}
\label{thm:forget with point constraints}
Define
\[
\foM^{\ev}(\shX,\beta/\beta',z)
:= 
\foM^{\ev}(\shX,\beta/\beta')\times_{\scrP(X,r)}^{\fs}
B\GG_m^{\dagger},
\]
where the morphism $\foM^{\ev}(\shX,\beta/\beta')\rightarrow\scrP(X,r)$
is $\ev_{\shX}$. There is a diagram
\[
\xymatrix@C=30pt
{
\scrM(X,\beta,z)\ar[r]^{\phi'_z}\ar[d]_{\varepsilon_z} & \scrM(X,\beta',z)\ar[d]^{\varepsilon_z'}\\
\foM^{\ev}(\shX,\beta/\beta',z)
\ar[d]_{\psi_z}\ar[r]_>>>>>{\phi_z} & \foM^{\ev}(\shX,\beta',z)\\
\foM^{\ev}(\shX,\beta,z)
}
\]
with the square cartesian in all categories, $\phi_z$ proper and
representable, $\psi_z$ strict \'etale, and with 
the pull-back relative obstruction 
theory for $\varepsilon_z'$ constructed in Proposition~\ref{prop:rel virt dim} 
being compatible with that for 
$\psi_z\circ\varepsilon_z$ also constructed in Proposition~\ref{prop:rel virt dim}.
\end{theorem}

\begin{proof}
{\bf The construction of $\psi_z$ and $\phi_z$.}
\medskip

Here $\psi_z$ is just a base-change of $\psi$ in 
Theorem~\ref{thm:forgetful}, hence is strict
\'etale.
To construct $\phi_z$, it is convenient to view
\begin{align*}
\foM^{\ev}(\shX,\beta/\beta',z)= {} & \foM(\shX,\beta/\beta')
\times^{\fs}_{\scrP(\shX,r)}B\GG_m^{\dagger}\\
\foM^{\ev}(\shX,\beta',z)= {} & \foM(\shX,\beta')
\times^{\fs}_{\scrP(\shX,r)}B\GG_m^{\dagger}
\end{align*}
as in Definition~\ref{def:Mbetaev}. It is then sufficient to
construct a ($2$-)commutative diagram
\begin{equation}
\label{eq:forget comm}
\xymatrix@C=30pt
{
\foM^{\ev}(\shX,\beta/\beta',z)\ar[r]^>>>>>{\phi\circ \pr_1}\ar[d]_{\gamma}&
\foM(\shX,\beta')
\ar[d]^{\ev'_{\shX}}
\\
B\GG_m^{\dagger}\ar[r]_{\zeta}
&\scrP(\shX,r)
}
\end{equation}
We write $\pr_1,\pr_2$ for the projections of
$\foM^{\ev}(\shX,\beta/\beta',z)$ onto $\foM(\shX,\beta/\beta')$ and
$B\GG_m^{\dagger}$ respectively, and $\zeta$ is the map of 
Proposition~\ref{prop: BGm morphism}.
We first define $\gamma$. Let $C^{\circ}\rightarrow C^{\circ}_{\mathrm{for}}
\stackrel{\tau}{\rightarrow} C'^{\circ}$ be the universal data
over $\foM(\shX,\beta/\beta')$. Recall the section $\rho_{\out}$ of
the ghost sheaf of $\foM(\shX,\beta/\beta')$. Similarly, denoting
as usual by $\delta$ the generator of the ghost sheaf $\NN$ of 
$B\GG_m^{\dagger}$, we take $\gamma$ to be the log morphism 
which, at the level of ghost sheaves, is given by
$\bar\gamma^{\flat}(\delta)
= 
\overline\pr_1^{\flat}(\rho_{\out})+\overline\pr_2^{\flat}(\delta)$; this
completely determines the morphism $\gamma$.

We now argue that commutativity of \eqref{eq:forget comm} can
essentially be checked by checking at the level of ghost sheaves.
Indeed, necessarily
$\ev'_{\shX}\circ\phi\circ\pr_1$ factors through the open stratum
of $\scrP(\shX,r)$, i.e., $[\shZ^{\circ}_r/\GG_m]$ as per the description
of Remark~\ref{remark: quotient description}. Note that
$\shZ^{\circ}_r\cong \shA_{P,P\setminus\{0\}}$ in the notation of
Definition~\ref{def:APK}, where $P$ is the stalk of $\overline{\shM}_X$
at the generic point of the stratum $Z_r$. Now $\scrP(\shX,r)$ was constructed
as a sub-log structure of $\shZ_r\times B\GG_m^{\dagger}$,
so $[\shZ^{\circ}_r/\GG_m]$ may be viewed as a sub-log structure
of $\shA_{P,P\setminus\{0\}}\times B\GG_m^{\dagger}\cong
\shA_{P\oplus\NN, (P\oplus\NN)\setminus\{(0,0)\}}$.
This sub-log structure is induced by the sub-monoid $P'\subseteq P\oplus\NN$,
$P'=\{(p,r(p))\,|\,p\in P\}$. We then have a canonical morphism
$[\shZ_r^{\circ}/\GG_m]\rightarrow \shA_{P'}$ via 
Proposition~\ref{prop:APrep}, which necessarily factors through
$\shA_{P',P'\setminus \{0\}}$. Further we have the  projection 
$[\shZ_r^{\circ}/\GG_m]
\rightarrow B\GG_m$, the latter with no log structure. These two maps
in fact induce an isomorphism
\[
[\shZ_r^{\circ}/\GG_m]\cong \shA_{P',P'\setminus\{0\}} \times B\GG_m,
\]
corresponding to a splitting of the algebraic torus
$\Spec \kk[(P\oplus\NN)^{\gp}]$ as 
$\Spec \kk[(P')^{\gp}]\times \Spec \kk[\ZZ(0,1)]$.

Identifying $P'$ with $P$, we have thus reduced 
the task of constructing $\phi_z$ to 
determining commutativity of the diagram
\begin{equation}
\label{eq:final comm diagram}
\xymatrix@C=30pt
{
\foM^{\ev}(\shX,\beta/\beta',z)\ar[r]^>>>>>>>{\gamma}
\ar[rd]_{\ev'_{\shX}\circ\phi\circ\pr_1}
& B\GG_m^{\dagger}\ar[d]^{\zeta}\\
&\shA_{P,P\setminus\{0\}}\times B\GG_m
}
\end{equation}
This commutativity can be checked by composing with the
first and second projections $\pr_1'$ and $\pr_2'$
to $\shA_{P,P\setminus\{0\}}$ and
$B\GG_m$. After composing with the first projection, it is enough
to check commutativity at the ghost sheaf level by 
Proposition~\ref{prop:APrep}, considering the maps
$\overline{\ev'_{\shX}\circ\phi\circ\pr_1}^{\flat}:P\rightarrow
\overline\shM$ and
$\overline{\zeta\circ\gamma}^{\flat}:P\rightarrow \overline\shM$,
where $\overline\shM$ is the ghost sheaf of $\foM^{\ev}(\shX,
\beta/\beta',z)$. However, by Lemma \ref{lem:ev comparison},
\[
\overline{\ev'_{\shX}\circ\phi\circ\pr_1}^{\flat}(p)
=\overline\pr_1^{\flat}(\overline{\ev}_{\shX}^{\flat}(p)+r(p)\rho_{\out}),
\]
while
\[
\overline{\zeta\circ\gamma}^{\flat}(p)=r(p)\big(\overline\pr_1^{\flat}(\rho_{\out})
+ \overline\pr_2^{\flat}(\delta)\big),
\]
the value of $\overline{\zeta}$ given by Proposition~\ref{prop: BGm morphism}.
Note however that because the fibre product defining
$\foM^{\ev}(\shX,\beta/\beta',z)$ is given by the maps
$\ev_{\shX}:\foM(\shX,\beta/\beta')\rightarrow \scrP(\shX,r)$
and $\zeta:B\GG_m^{\dagger}\rightarrow \scrP(\shX,r)$,
we in fact have
\[
\overline\pr_1^\flat(\overline{\ev}_{\shX}^{\flat}(p))=
\overline\pr_2^{\flat}(r(p)\delta),
\]
which shows that in fact 
$\overline{\ev'_{\shX}\circ\phi\circ\pr_1}^{\flat}=
\overline{\zeta\circ\gamma}^{\flat}$, hence commutativity of
\eqref{eq:final comm diagram} after composing with the projection 
$\pr_1'$.

On the other hand, 
$\pr'_2\circ\ev'_{\shX}\circ\phi\circ\pr_1$ is determined
by the conormal bundle $\shN^{\vee}_{x_{\out}/C'}$. But
$\pr'_2\circ\zeta$ is the identity on underlying stacks, 
and thus $\pr'_2\circ\zeta\circ\gamma$
is determined by the line bundle $\shL_{\overline\pr_1^{\flat}(\rho_{\out})
+\overline{\pr}_2^{\flat}(\delta)}
\cong \pr_1^*\shL_{\rho_{\out}}\otimes \shN^{\vee}_{x_{\out}/C}$.
By Lemma \ref{lem:punctured stabilization}, (4), we then see
that $\pr'_2\circ\ev'_{\shX}\circ\phi\circ\pr_1$ and 
$\pr'_2\circ\zeta\circ\gamma$ agree (or rather are $2$-isomorphic). 

This completes the construction of $\phi_z$. 

\medskip
{\bf Completion of the proof.}
\medskip

The cartesianity of the square in the statement of the theorem now
follows from Theorem~\ref{thm:forgetful},(3).
Now note that from the construction of $\phi_z$ there is a commutative 
(but not cartesian) diagram
\begin{equation}
\label{eq:foMs commute}
\xymatrix@C=30pt
{
\foM^{\ev}(\shX,\beta/\beta',z)\ar[r]^{\phi_z}\ar[d]_{\pr_1}
& \foM^{\ev}(\shX,\beta',z)\ar[d]^{\pr_1}\\
\foM^{\ev}(\shX,\beta/\beta')\ar[r]_{\phi}&\foM^{\ev}(\shX,\beta')
}
\end{equation}
Both vertical arrows in this diagram are
finite and representable by Remark~\ref{rem:sat finite representable},
being a base-change of a (non-strict) closed embedding. Further,
$\phi$ is proper and representable by Theorem~\ref{thm:forgetful}, (2).
Thus $\phi\circ\pr_1=\pr_1\circ\phi_z$ is proper and representable, hence
$\phi_z$ is proper and representable.

The final statement involves compatibility of obstruction theories.
There is a similar commutative square as \eqref{eq:foMs commute}
involving the moduli spaces
$\scrM(X,\beta)$ etc., and a cube whose top face is this latter
square and whose bottom face is \eqref{eq:foMs commute}, with all
side faces cartesian. The statement on compatibility of obstruction
theories then follows immediately from the construction
of the obstruction theories for $\scrM(X,\beta,z)\rightarrow
\foM^{\ev}(\shX,\beta,z)$, $\scrM(X,\beta',z)\rightarrow
\foM^{\ev}(\shX,\beta',z)$ in Proposition~\ref{prop:rel virt dim} 
and the compatibility statement of Theorem~\ref{thm:forgetful}.
\end{proof}

\begin{proof}[Proof of Theorem~\ref{thm:unit}]
If $A=0$, the claim follows immediately from Lemma 
\ref{lem:constant-maps}. Otherwise, we may apply
Theorem~\ref{thm:forget with point constraints}.

If the virtual dimension of $\scrM(X,\beta,z)$ is non-zero, then
$N^{\beta}_{0pr}=0$ regardless by definition. Otherwise, 
by push-pull, (\cite{Man},~Thm.~4.1,~(i)), writing
$\varepsilon_z^!$ and $(\varepsilon_z')^!$ for Manolache's virtual
pull-back, we have
\begin{align*}
N^{\beta}_{0pr}= \deg [\scrM(X,\beta,z)]^{\virt} = {} & 
\deg (\phi'_{z})_*[\scrM(X,\beta,z)]^{\virt}\\
= {} 
& \deg
(\phi'_z)_*\varepsilon_z^! [\foM^{\ev}(\shX,\beta/\beta',z)]\\
={} & \deg(\varepsilon_z')^!\phi_{z*}[\foM^{\ev}(\shX,\beta/\beta',z)].
\end{align*}
However, $\dim\foM^{\ev}(\shX,\beta/\beta',z)=\dim\foM^{\ev}(\shX,\beta,z)
=0$ by the fact that $\psi_z$ is \'etale and Proposition
\ref{prop:the invariants}, (2), while $\dim\foM^{\ev}(\shX,\beta',z)=-1$
by the same proposition. Thus $\phi_{z*}[\foM^{\ev}(\shX,\beta/\beta',z)]=0$,
and hence $N^{\beta}_{0pr}=0$.
\end{proof}

\section{Independence of modulus}
\label{sec: modulus invariance}

\subsection{First steps: definitions of invariants}

Our goal is to prove Theorem~\ref{theorem: independence of modulus}.  
In this subsection, we will define some invariants in
a slightly more general situation than
that of this theorem.

Fix a boundary point $y\in \overline{\scrM}_{0,4}$, and let
$\sigma_y=\RR_{\ge 0}\ell^*\oplus\RR_{\ge 0}\delta^*$ denote the cone 
$\Hom(R,\RR_{\ge 0})$, where $R$ is the stalk of the ghost sheaf of
$\overline{\scrM}_{0,4}^{\dagger}=\overline{\scrM}_{0,4}\times 
B\GG_m^{\dagger}$ at $y$, 
with generators $\ell, \delta$
as in \eqref{eq:R Rlambda def}.
Let $\tau\subset \sigma_y$ be any subcone isomorphic to
the standard two-dimensional cone. Dually, this defines an inclusion
of monoids 
\[
R\subset R_{\tau} := \tau^{\vee}\cap R^{\gp}.
\]
We may then define a log stack $B\GG_m^{\tau}$ whose ghost sheaf
is the constant sheaf with stalk $R_{\tau}$, 
with the torsor associated to $a\ell+b\delta\in R_{\tau}$
being $\shU^{\otimes b}$, where $\shU$ is the
universal torsor on $B\GG_m$. Thus we obtain a morphism
\begin{equation}
\label{eq:psitau def}
\psi_{\tau}:B\GG_m^{\tau}\rightarrow \overline{\scrM}_{0,4}^{\dagger}
\end{equation}
with image $y\times B\GG_m$ and $\overline{\psi}^{\flat}_{\tau}:
R\rightarrow R_{\tau}$ the canonical inclusion.

\begin{definition} 
\label{def:Mtaudef}
Let $\beta$ be a curve class as in Definition~\ref{def: Mydef}.
We set
\begin{align*}
\scrM_\tau:= {} & \scrM(X,\beta,z) \times^{\fs}_{\overline\scrM_{0,4}^{\dagger}}
B\GG_m^{\tau}\\
\foM^{\ev}_\tau:= {} & \foM^{\ev}(\shX,\beta,z) 
\times^{\fs}_{\overline\scrM_{0,4}^{\dagger}}
B\GG_m^{\tau}
\end{align*}
\end{definition}

\begin{remark}
The spaces
$\scrM_y^{\ddagger}$, $\foM^{\ddagger,\ev}_y$ of Definition~\ref{def: Mydef} 
are special cases  of the above definition, with 
$\tau$ the cone generated by $\lambda\ell^*+\delta^*$ and $\ell^*$.
\end{remark}

\begin{lemma}
\label{lem:Mtau proper}
$\scrM_{\tau}$ is a proper Deligne-Mumford stack over $\Spec\kk$
and has a perfect relative
obstruction theory over $\foM_{\tau}^{\ev}$, with relative virtual
dimension at a point of $\scrM_\tau$ represented by a punctured
map $f:C^{\circ}\rightarrow X$ being
\[
\chi(f^*\Theta_{X/\kk}(-x_{\out})).
\]
\end{lemma}

\begin{proof}
Write $\overline\scrM_{0,4}\times B\GG_m$ for the underlying stack
of $\overline{\scrM}_{0,4}^{\dagger}$.
First note that $\ul{\psi_{\tau}}:B\GG_m \rightarrow\overline\scrM_{0,4}\times
B\GG_m$
is representable, hence of DM type, and being of DM type is stable
under base-change, so 
\[
\scrM(X,\beta,z) \times_{\overline\scrM_{0,4}\times B\GG_m} B\GG_m
\rightarrow\scrM(X,\beta,z)
\]
is of DM type. Since $\scrM(X,\beta,z)$ itself is Deligne-Mumford,
i.e., of DM type over $\Spec\ZZ$, and a composition of DM type
morphisms is DM type, it follows that 
$\scrM(X,\beta,z) \times_{\overline{\scrM}_{0,4}
\times B\GG_m} B\GG_m $ is Deligne-Mumford. 
From Remark~\ref{rem:sat finite representable}, the morphism
\[
\scrM_{\tau}=\scrM(X,\beta,z)\times^{\fs}_{\overline{\scrM}_{0,4}^{\dagger}}
B\GG_m^{\tau}
\rightarrow
\scrM(X,\beta,z)\times_{\overline{\scrM}_{0,4}^{\dagger}}
B\GG_m^{\tau}
\]
is finite and representable. Thus $\scrM_{\tau}$
is Deligne-Mumford, as promised. Furthermore, since $\ul{\psi_{\tau}}$
is also proper, being a closed immersion, and $\scrM(X,\beta,z)$ is
proper over $\Spec \kk$ by Lemma~\ref{lemma: properness}, it follows
that $\scrM_{\tau}$ is proper over $\Spec \kk$.

From the diagram cartesian in all categories
\[
\xymatrix@C=15pt
{
\scrM_{\tau}\ar[r]\ar[d]&\scrM(X,\beta,z)\ar[d]\\
\foM_{\tau}^{\ev}\ar[r]&\foM^{\ev}(\shX,\beta,z)
}
\]
the perfect relative obstruction theory for the right-hand vertical
arrow pulls back to a perfect relative obstruction theory for
the left-hand vertical arrow. In particular, the relative virtual
dimension remains unchanged, given in Proposition~\ref{prop:rel virt dim}.
\end{proof}

Here we will prove the part of Theorem~\ref{theorem: independence of modulus}
involving virtual dimensions.

\begin{lemma}
\label{lemma: tau VFC}
\begin{enumerate}
\item
Suppose $y\in\overline\scrM_{0,4}$ is a boundary point, $\tau\subseteq
\sigma_y$ a cone determining a morphism
$\psi_{\tau}:B\GG_m^{\tau}\rightarrow \overline\scrM_{0,4}^{\dagger}$.
Suppose $\psi_{\tau}$ is transverse to 
\[
\Phi:\foM^{\ev}(\shX,\beta,z)\rightarrow
\overline\scrM_{0,4}^{\dagger}
\]
in the sense of Definition 
\ref{def: transverse}. Then $\foM_{\tau}^{\ev}$ (Definition 
\ref{def:Mtaudef}) is of pure dimension zero,
and $\scrM_{\tau}$ carries a virtual fundamental class
of dimension $\chi((f^*\Theta_{X/\kk})(-x_{\out}))$. In particular,
as a special case, if $B\GG_m^{\ddagger}\rightarrow \overline{\scrM}_{0,4}^{\dagger}$ is transverse to $\Phi$, then the same dimension statements
hold for $\foM^{\ddagger,\ev}_y$
and $\scrM^{\ddagger}_y$ (Definition~\ref{def: Mydef}).
\item Suppose $y\in\overline\scrM_{0,4}$ is not a boundary point.
Then $\foM_y^{\ddagger,\ev}$ is of pure dimension zero, and
$\scrM^{\ddagger}_y$ carries a virtual fundamental class of dimension 
$\chi((f^*\Theta_{X/\kk})(-x_{\out}))$.
\end{enumerate}
\end{lemma}

\begin{proof}
By Lemma~\ref{lemma: moduli morphisms}, the base-change 
$\widetilde \Phi:\foM^{\ev}_{\tau}\rightarrow B\GG_m^{\tau}$ of
$\Phi$ is log smooth, and the
transversality condition implies that it is integral. Now for
$x\in \foM^{\ev}_{\tau}$, we have, using Proposition~\ref{prop: flat fibre dim},
(2) for the second equality,
\begin{align*}
\dim\foM^{\ev}_{\tau} = {} & \dim_{x}\ul{\widetilde \Phi}^{-1}(\ul{\widetilde 
\Phi}(x))
+\dim B\GG_m^{\tau}\\
= {} & \dim^{\log}_x \widetilde \Phi^{-1}(\widetilde \Phi(x)) - 1.
\end{align*}
However, by the base-change result of Proposition 
\ref{prop: flat fibre dim}, (5) and
Proposition~\ref{prop:the invariants}, (1), the log fibre
dimension of $\widetilde \Phi$ is $1$.
This yields $\dim \foM^{\ev}_{\tau}=0$. 

In (2), the fact that the morphism $\psi_y:B\GG_m^{\ddagger}\rightarrow
\overline\scrM_{0,4}^{\dagger}$ of Definition~\ref{def: BGmdagger}
is transverse to $\Phi$ is automatic, and a 
similar argument gives $\dim \foM^{\ddagger,\ev}_y=0$. The dimension of
the virtual fundamental class in both cases
then immediately follows from the form of
the obstruction theory given in Proposition~\ref{prop:rel virt dim}
or Lemma~\ref{lem:Mtau proper}.
\end{proof}

\subsection{Comparing invariants for adjacent cones}
The first step in proving
Theorem~\ref{theorem: independence of modulus} is
to compare enumerative invariants obtained from Lemma
\ref{lemma: tau VFC}, (1) for different choices of $\tau$. 
An important ingredient is the following lemma:

\begin{lemma}
\label{lemma: PoneGmChow}
Let $\GG_m$ act with weight $w$ on $\PP^1$ for some positive integer
$w$, i.e.,  by $(x_0:x_1)\mapsto  (\lambda^w x_0: x_1)$
for $\lambda\in \GG_m$. Write $0=(0:1)$, $\infty = (1:0)$. Let
$\shL$ be the universal line bundle on $[\PP^1/\GG_m]$, i.e., the pull-back
of the universal line bundle under the canonical morphism $[\PP^1/\GG_m]
\rightarrow B\GG_m$. Then the relation
\[
[0/\GG_m]-[\infty/\GG_m]+wc_1(\shL)\cap [\PP^1/\GG_m]=0
\]
holds in $A_{-1}([\PP^1/\GG_m])$, where $A_*$ denotes Kresch's Chow group
of Artin stacks, see \cite{Kresch}.
\end{lemma}

\begin{proof}
As $\PP^1 \rightarrow [\PP^1/\GG_m]$ is the universal torsor
over the latter stack, the total space of the associated universal line bundle 
is
$L=\PP^1\times^{\GG_m} \AA^1$, the quotient of $\PP^1\times\AA^1$
by the action 
$\big((x_0:x_1),y\big)\mapsto \big((\lambda^w x_0:x_1),\lambda^{-1}y\big)$
for $\lambda\in\GG_m$. Here we are using the convention that the
universal line bundle over $B\GG_m$ is the quotient of $\AA^1$ by
$\GG_m$ acting with weight $-1$. This choice of convention will have
no effect on the sequel.
As in \cite{Kresch}, we write $A^{\circ}_*(L)$
for the naive Chow group of $L$. By \cite{Kresch},~Thm.~2.1.12, 
there is then a natural homomorphism $A^{\circ}_0(L)\rightarrow A_0(L)$ and 
an isomorphism $A_{-1}([\PP^1/\GG_m]) \cong A_0(L)$ given by
pull-back of Chow classes, hence a homomorphism
\[
s:A^{\circ}_0(L)\rightarrow A_{-1}([\PP^1/\GG_m]).
\]
Note the field of rational functions of $L$ is
\[
k(L)=k(\PP^1\times\AA^1)^{\GG_m} = k(x=x_0/x_1,y)^{\GG_m} = k(xy^w).
\]
Thus there is a relation in $A_0^{\circ}(L)$ given by
\[
0 = (xy^w) = [(0\times \AA^1)/\GG_m]-[(\infty\times\AA^1)/\GG_m)]
+w[(\PP^1\times 0)/\GG_m].
\]
On the other hand, by \cite{Kresch},~Def.~2.4.2, $c_1(\shL)\cap
[\PP^1/\GG_m]$ is the image of the zero section of $L\rightarrow [\PP^1/\GG_m]$
under $s$. This zero section is $[(\PP^1\times 0)/\GG_m]$, hence the result.
\end{proof}

The procedure for comparing $\scrM^{\ddagger}_y$ with 
$\scrM^{\ddagger}_{y'}$ for $y'$ not
a boundary point is as follows.
Morally, we would like to perform a logarithmic blow-up
$\widetilde{\scrM}_{0,4}^{\dagger}\rightarrow
\overline{\scrM}_{0,4}^{\dagger}$ by refining $\sigma_y$ so that
$\widetilde{\scrM}_{0,4}^{\dagger}\rightarrow \overline{\scrM}_{0,4}^{\dagger}$
is transverse to $\Phi$, at least locally near $y$,
and so that the fan $\Sigma$ refining $\sigma_y$ contains
the cone $\RR_{\ge 0}(\lambda\ell^*+\delta^*)+\RR_{\ge 0}\ell^*$. Note
that if $\tau$ is taken to be this cone, then $\scrM^{\ddagger}_y
=\scrM_{\tau}$ by definition.
This blow-up, up to the quotient by $\GG_m$,
would produce a chain of $\PP^1$'s over $y$. Each double point
in this chain of $\PP^1$'s then corresponds logarithmically to
a point $B\GG_m^{\tau}$ in the log blow-up 
for the corresponding two-dimensional cone $\tau$, and then we obtain
a moduli space $\scrM_{\tau}$. Integrality implies this moduli space
is of virtual dimension zero, and we then just need to show that
the virtual degree of $\scrM_{\tau}$ is independent of the two-dimensional
cone $\tau\in \Sigma$. We do this by comparing these virtual degrees for
adjacent cones. Unfortunately, the fact that we don't have a chain of
$\PP^1$'s but a $\GG_m$-quotient thereof makes the argument a bit more
difficult, causing potential problems with terminal tails: these 
problems are caused by the previous lemma.

\begin{lemma} 
\label{lemma: invarianceI}
Assume we are in the situation of Theorems~\ref{mainassociativity1} or 
\ref{mainassociativity2}, and $\beta$ is a curve class 
as in Definition \ref{def: Mydef}.
Assume further that $A\cdot c_1(\Theta_{X/\kk})=0$.
Fix a boundary point $y\in\overline{\scrM}_{0,4}$ and standard two-dimensional
cones $\tau,\tau'
\subset \sigma_y$ such that $\tau\cap \tau'=\omega$ is a one-dimensional
ray in $\sigma_y$. Suppose further that the maps
\[
\psi_{\tau}:B\GG_m^{\tau}\rightarrow \overline{\scrM}_{0,4}^{\dagger},\quad
\psi_{\tau'}:B\GG_m^{\tau'}\rightarrow \overline{\scrM}_{0,4}^{\dagger}
\]
are both transverse to $\Phi:\foM^{\ev}(\shX,\beta,z)\rightarrow
\overline{\scrM}_{0,4}^{\dagger}$. So $\scrM_{\tau}$, $\scrM_{\tau'}$
carry a virtual fundamental class of dimension zero by
Lemma~\ref{lemma: tau VFC}, (1). Then
\[
\deg [\scrM_{\tau}]^{\virt}=\deg [\scrM_{\tau'}]^{\virt}.
\]
\end{lemma}

\begin{proof}
{\bf Step 1. The interpolating moduli space $\scrM_{\omega}$.}

Let $\Sigma$ be the fan in $\RR\sigma_y=\Hom(R,\RR)=\RR\ell^*\oplus
\RR\delta^*$ consisting of $\tau$,$\tau'$
and their faces, and let $X_{\Sigma}$ be the corresponding toric variety.
This contains a one-dimensional closed stratum $Z_{\omega}$ 
corresponding to the common ray $\omega$ of $\tau,\tau'$, with
$Z_{\omega}\cong \PP^1$. Recalling the generators $\ell,\delta$ of 
$R$, we define the slope of a ray
$\omega'$
contained in $\sigma_y$ to be $\delta(p)/\ell(p)$ for any $p\in\omega'\setminus
\{0\}$. Order $\tau,\tau'$ so that the boundary ray of $\tau$ which is
not $\omega$ has slope greater than the slope of $\omega$.

Giving $X_{\Sigma}$ and $X_{\sigma_y}=\AA^2$ their toric log structures,
the canonical morphism of fans from $\Sigma$ to $\sigma_y$ gives a log
morphism $\tilde\psi_{\omega}:X_{\Sigma}\rightarrow X_{\sigma_y}$,
with $\tilde\psi_{\omega}(Z_{\omega})=\{O\}$, where $O$ is the origin
in $X_{\sigma_y}$. 

Let $\GG_m$ act on $X_{\Sigma}$ via the one-parameter
subgroup $\ZZ\delta^*\otimes \GG_m\subset R^*\otimes\GG_m$. 
Note that if $\omega=\RR_{\ge 0} (a\ell^*+b\delta^*)$ with $a,b$
relatively prime, then $z^{-b\ell
+a\delta}$ is a coordinate on $Z_{\omega}$ having a simple zero at $Z_{\tau}$
and a simple pole at $Z_{\tau'}$. From this, we see that if we identify
$Z_{\tau}$ with $0\in\PP^1=Z_{\omega}$ and $Z_{\tau'}$ with $\infty\in\PP^1$,
then the weight of this $\GG_m$ action on $Z_{\omega}$ is $w(\omega):=a$,
in the sense of Lemma~\ref{lemma: PoneGmChow}.

Now $\tilde\psi_{\omega}$ descends to a morphism of log stacks
\[
\psi_{\omega}:[X_{\Sigma}/\GG_m]\rightarrow [\AA^2/\GG_m] = \AA^1\times
[\AA^1/\GG_m]
\]
which induces by restriction and inclusion a composed morphism
\[
\psi_{\omega}:[Z_{\omega}/\GG_m]\rightarrow [O/\GG_m] =
\Spec\kk^{\dagger}\times B\GG_m^{\dagger}\hookrightarrow
\overline{\scrM}_{0,4}^{\dagger},
\]
the latter inclusion being strict with image $y\times B\GG_m^{\dagger}$.
Note that the zero-dimensional strata $Z_{\tau},Z_{\tau'}\subset Z_{\omega}$
are invariant under the torus action, and the compositions of the inclusions
\[
B\GG_m^{\tau}=[Z_{\tau}/\GG_m]\hookrightarrow [Z_{\omega}/\GG_m],
\quad
B\GG_m^{\tau'}=[Z_{\tau'}/\GG_m]\hookrightarrow [Z_{\omega}/\GG_m]
\]
with $\psi_{\omega}$ coincide with $\psi_{\tau}$ and $\psi_{\tau'}$
respectively. 

We now define
\begin{align*}
\scrM_{\omega}:= {} &[Z_{\omega}/\GG_m]
\times_{\overline{\scrM}_{0,4}^{\dagger}}^{\fs}\scrM(X,\beta,z),\\
\foM^{\ev}_{\omega}:= {} &[Z_{\omega}/\GG_m]
\times_{\overline{\scrM}_{0,4}^{\dagger}}^{\fs}\foM^{\ev}(\shX,\beta,z).
\end{align*}
Note that since $[Z_{\omega}/\GG_m]\rightarrow\overline\scrM_{0,4}^{\dagger}$
is proper, the projection $\scrM_{\omega}\rightarrow\scrM(X,\beta,z)$,
as in the proof of Lemma~\ref{lem:Mtau proper},
is proper. Hence, by Lemma~\ref{lemma: properness},(2), $\scrM_{\omega}$ is proper
over $\Spec\kk$.

We now obtain a diagram
\[
\xymatrix@C=30pt
{
\scrM_{\tau}\ar[r]^q\ar[d]_h&\scrM_{\omega}\ar[d]_{\bar h}&
\scrM_{\tau'}\ar[l]_{q'}\ar[d]^{h'}\\
\foM^{\ev}_\tau
\ar[r]^p\ar[d]_{\widetilde \Phi}&\foM^{\ev}_{\omega}\ar[d]_{\overline\Phi}&
\foM_{\tau'}^{\ev}
\ar[l]_{p'}\ar[d]^{\widetilde \Phi'}\\
B\GG_m^{\tau}\ar[r]&[Z_{\omega}/\GG_m]&B\GG_m^{\tau'}\ar[l]
}
\]
Here $h$, $\bar h$ and $h'$, as well as all horizontal arrows,
are strict. All squares are
cartesian in both the category of underlying stacks and the fs log category.
By the transversality assumption, $\widetilde\Phi$ and $\widetilde\Phi'$ are 
integral. On the other
hand, as any morphism $\NN\rightarrow P$ is integral if $P$ is an integral
monoid, it then follows that
$\overline\Phi$ is integral as the ghost sheaf of 
the log structure on $[Z_{\omega}/\GG_m]$
is $\NN$ except at $B\GG_m^{\tau}$ and $B\GG_m^{\tau'}$. By Lemma
\ref{lemma: moduli morphisms}, $\overline{\Phi}$ is also log smooth, hence flat
by Proposition~\ref{prop: flatness sorites}.

\medskip

{\bf Step 2. Pull-back of cycles.}

As in the proof of Lemma~\ref{lemma: tau VFC},
the relative dimension of $\overline{\Phi}$ 
is $1$, so we obtain a flat pull-back of cycles 
\[
\overline{\Phi}^*:A_*([Z_{\omega}/\GG_m])\rightarrow A_{*+1}(\foM_{\omega}^{\ev}).
\]

Note also that $\bar h$ carries a relative obstruction theory pulled
back from that of $\scrM(X,\beta,z)\rightarrow \foM^{\ev}(\shX,\beta,z)$,
which further pulls back to the relative obstruction theories for $h$
and $h'$. Denote by $h^!$, $(h')^!$ and $\bar h^!$ the virtual
pull-backs of \cite{Man} defined by these obstruction theories.
It then follows from \cite{Man},~Thm.~4.1,(i) and (iii) that
\[
q_*[\scrM_\tau]^{\virt}= \bar h^!\overline{\Phi}^* [0/\GG_m],\quad 
q_*[\scrM_{\tau'}]^{\virt}=  \bar h^!\overline{\Phi}^* [\infty/\GG_m],
\]
so it will be sufficient to show the equality
\[
\deg \bar h^!\overline{\Phi}^* [0/\GG_m] = \deg \bar h^!\overline{\Phi}^* 
[\infty/\GG_m].
\]
Now by Lemma~\ref{lemma: PoneGmChow}, we have in $A_{-1}([Z_{\omega}/\GG_m])$,
\[
[0/\GG_m]-[\infty/\GG_m]+w(\omega)c_1(\shL)\cap[Z_{\omega}/\GG_m]=0.
\]
We thus must show that
\[
\deg \bar h^!\overline{\Phi}^*(c_1(\shL)\cap[Z_{\omega}/\GG_m]) =
\deg \bar h^!(c_1(\overline{\Phi}^*\shL)\cap[\foM_{\omega}^{\ev}])=0
\]
to complete the proof.
\medskip

\pagebreak

{\bf Step 3. No tails: the vanishing of 
$\deg \bar h^!(c_1(\overline{\Phi}^*\shL)\cap[\foM_{\omega}^{\ev}])$.}

\medskip

What is $\overline{\Phi}^*\shL$? Recall from Lemma 
\ref{lem: more log smoothness}, (1)
that the projection
$\foM^{\ev}(\shX,\beta,z) \rightarrow B\GG_m^{\dagger}$
is given by the conormal bundle of $x_{\out}$ in the universal curve
over $\foM^{\ev}(\shX,\beta,z)$.
Thus $\overline{\Phi}^*\shL$ is this conormal bundle.
Equivalently, if $\phi:\foM^{\ev}_{\omega}\rightarrow\Mbf=\Mbf_{0,4}$ is the
forgetful map just remembering the domain, then $\overline{\Phi}^*\shL=
\phi^*\shN^{\vee}_{x_{\out}/{\mathbf C}}$, where ${\mathbf C}\rightarrow
\Mbf$ is the universal curve. It is standard that
$c_1(\shN^{\vee}_{x_{\out}/{\mathbf C}})$ can be written
as the divisor $D'+
D(x_1,x_2,x_3\,|\,x_{\out})$, where $D'$ is the pull-back of a point
under the stabilization morphism $\Mbf\rightarrow\overline{\scrM}_{0,4}$.
Indeed, if $(\ul{C},x_1,x_2,x_3,x_{\out})$ is a four-pointed rational
curve, then the map to its stabilization is an isomorphism in a neighbourhood of
$x_{\out}$ unless the curve lies in the irreducible
divisor $D(x_1,x_2,x_3\,|\,x_{\out})$. 
Thus $c_1(\shN^{\vee}_{x_{\out}/{\mathbf C}})$ can be written
as $n'D'+n D(x_1,x_2,x_3\,|\,x_{\out})$ for some integers $n,n'$.
It is standard that the conormal bundle of
any marked point of the universal curve over $\overline\scrM_{0,4}=\PP^1$
is $\O_{\PP^1}(1)$, so $n'=1$. On the other hand, $n$ 
is easily seen to be $1$ by considering a one-parameter family of curves
in $\Mbf$ transverse to $D(x_1,x_2,x_3\,|\,x_{\out})$.

Let $\foM$ be an irreducible component of $\foM_{\omega}^{\ev}$
with the reduced induced stack structure, $f:\foC^{\circ}/\foM\rightarrow\shX$
the induced family of punctured maps. Let $\phi:\foM\rightarrow \Mbf$ be
the forgetful morphism as before, and 
let $\scrM=\scrM_{\omega}\times_{\foM^{\ev}_{\omega}} \foM$. It
will be sufficient to show that $\deg \bar h^!(c_1(\overline{\Phi}^*\shL)
\cap[\foM])=0$.

As the image of $\foM_\omega^{\ev}$ under the forgetful and stabilization
map $\foM_\omega^{\ev}\rightarrow \overline{\scrM}_{0,4}$ is the point
$y$, we may choose $D'$ so that $\phi^{-1}(D')=
\emptyset$. Thus we have two possibilities. Either
(1) $\phi^{-1}(D(x_1,x_2,x_3\,|\,x_{\out}))=\foM$, or (2)
$c_1(\bar m^*\shL)\cap [\foM]$ is supported
on the codimension one closed substack $\phi^{-1}(D(x_1,x_2,x_3\,|\,x_{\out}))$.
As whether a given point of $\foM$ maps into
$D(x_1,x_2,x_3\,|\,x_{\out})$ can be determined from the tropicalization
of the corresponding punctured map, necessarily
in this case $\phi^{-1}(D(x_1,x_2,x_3\,|\,x_{\out}))$ is a union of
strata of $\foM$.

Now set $\foM'=\foM$ in case (1) of the previous paragraph or
set $\foM'$ to be an irreducible component of
$\phi^{-1}(D(x_1,x_2,x_3\,|\,x_{\out}))$ in case (2). Let $\scrM'=
\foM'\times_{\foM}\scrM$. We are now in a position to apply
the no-tail lemma, Lemma~\ref{thm:no-tail-lemma}. In that theorem,
take $F=\foM'$, $S^{\circ}$ to be the 
stratum of $\foM^{\ev}(\shX,\beta)$
containing the image of the generic point of $F$ under the composed morphism
\[
\foM'\hookrightarrow \foM\hookrightarrow\foM_{\omega}^{\ev}\rightarrow
\foM^{\ev}(\shX,\beta,z)\rightarrow\foM^{\ev}(\shX,\beta),
\]
where the last two maps are the obvious projections,
and take $W=B\GG_m^{\dagger}$ with $g:W\rightarrow \scrP(X,r)$ the
standard morphism determined by $z$ by Proposition~\ref{prop: BGm morphism}.
Then the hypotheses
of Lemma~\ref{thm:no-tail-lemma} are satisfied; in particular,
all transversality statements in hypothesis (5) of that theorem
are trivial since $W$ has ghost sheaf $\NN$. We thus conclude
that $\bar h^![\foM']=0$ in rational Chow.

Thus in case (1),
\[
\deg \bar h^!(c_1(\overline{\Phi}^*\shL)\cap [\foM])=
\deg c_1(\bar h^*\overline{\Phi}^*\shL)\cap \bar h^![\foM]=0.
\]
Here the second identity follows from $\foM=\foM'$ in this case.

In case (2), we have
\[
\deg \bar h^!(c_1(\overline{\Phi}^*\shL)\cap [\foM])
=\deg \sum_i a_i\bar h^!([\foM_i']),
\]
where $\foM_i'$ ranges over codimension one strata of $\foM$ in
$\phi^{-1}(D(x_1,x_2,x_3\,|\,x_{\out}))$ and $a_i$ are some rational 
numbers. The vanishing then follows as before from $\bar h^![\foM_i']=0$ in
rational Chow.
\end{proof}

\subsection{Comparing the invariant at different points of
$\overline{\scrM}_{0,4}$}
The second step of the proof is to compare the moduli spaces
$\scrM^{\ddagger}_{y'}$ for $y'\in \overline{\scrM}_{0,4}$ for general $y'$ and
$\scrM_{\tau}$, for $\tau\subseteq\sigma_y$ a cone with two boundary
rays: $\RR_{\ge 0} \delta^*$ and another ray of very large slope.
The following lemma makes this comparison.

\begin{lemma}
\label{lemma: invarianceII}
Assume we are in the situation of Theorems~\ref{mainassociativity1} or 
\ref{mainassociativity2}, and $\beta$ is a curve class as in Definition 
\ref{def: Mydef}.
Assume further that $A\cdot c_1(\Theta_{X/\kk})=0$.
Suppose
that for each boundary point $y\in\overline{\scrM}_{0,4}$, we are
given a standard
two-dimensional cone $\tau_y\subset \sigma_y$ such that (1) one boundary
ray of $\tau_y$ is $\RR_{\ge 0}\delta^*$ and (2)
$\psi_{\tau_y}:B\GG_m^{\tau_y}\rightarrow\overline\scrM_{0,4}^{\dagger}$
is transverse to $\Phi:\foM^{\ev}(\shX,\beta,z)\rightarrow
\overline{\scrM}_{0,4}^{\dagger}$. Let $y'\in \overline{\scrM}_{0,4}$ be a
non-boundary point. Note $\scrM_{\tau_y}$, $\scrM^{\ddagger}_{y'}$ each carry
a virtual fundamental class of dimension zero by Lemma~\ref{lemma: tau VFC}. 
Then for each boundary point $y\in \overline{\scrM}_{0,4}$,
\[
\deg [\scrM_{\tau_y}]^{\virt} = \deg [\scrM^{\ddagger}_{y'}]^{\virt}.
\]
\end{lemma}

\begin{proof}
The argument is very similar to that of Lemma~\ref{lemma: invarianceI},
but in fact simpler because tails are not an issue here. Define
a new log structure on $\overline{\scrM}_{0,4}\times B\GG_m$ as follows.
If $\shM$ and $\overline{\shM}$ denote the log structure and its ghost sheaf on 
$\overline{\scrM}_{0,4}^{\dagger}$, we let $\overline{\shM}'\subset
\overline{\shM}^{\gp}$ be a subsheaf defined as follows.
First, $\overline{\shM}'$ agrees with
$\overline\shM$ away from the boundary points.
At the boundary point $y$, the stalk $\overline\shM'_y$
is the inverse image of $\overline{\shM}_{B\GG_m^{\tau_y}}$
under the homomorphism $\overline{\psi}_{\tau_y}^{\flat}:
\overline{\shM}_y^{\gp} \rightarrow \overline{\shM}_{B\GG_m^{\tau_y}}^{\gp}$,
where $\psi_{\tau_y}$ is as defined in \eqref{eq:psitau def}.
We define a sheaf of monoids
\[
\shM'=\shM^{\gp}\times_{\overline\shM^{\gp}} \overline\shM' \subset \shM^{\gp}.
\]
It is easy to check that the structure map $\alpha:\shM\rightarrow
\O_{\overline{\scrM}_{0,4}^{\dagger}}$ extends to $\shM'$ by sending
all elements of $\shM'\setminus\shM$ to zero. We remark that this new
log structure is a puncturing of the log structure on
$\overline{\scrM}_{0,4}^{\dagger}$. Write $\overline{\scrM}'_{0,4}$
for $\overline{\scrM}_{0,4}\times B\GG_m$ with this new log structure.
We then have a factorization of $\psi_{\tau_y}$
\[
B\GG_m^{\tau_y}\rightarrow\overline{\scrM}_{0,4}' \rightarrow
\overline{\scrM}_{0,4}^{\dagger}
\]
with the first morphism strict.

Defining 
\begin{align*}
\scrM':= {} &\overline{\scrM}_{0,4}'
\times_{\overline{\scrM}_{0,4}^{\dagger}}^{\fs}\scrM(X,\beta,z)\\
\foM':= {} &\overline{\scrM}_{0,4}'
\times_{\overline{\scrM}_{0,4}^{\dagger}}^{\fs}\foM^{\ev}(\shX,\beta,z)
\end{align*}
we see that the morphism
\[
\foM'\rightarrow \overline\scrM_{0,4}'
\]
is integral. Indeed, this is the case at each boundary point $y$ 
because of the transversality
assumption and elsewhere because the stalk of the ghost sheaf of
$\overline{\scrM}_{0,4}'$ is $\NN$. In particular, we have a diagram
with $\overline{\Phi}$ flat and all squares cartesian in both the ordinary and
fs log categories:
\[
\xymatrix@C=30pt
{
\scrM_{\tau_y}\ar[r]^q\ar[d]_h&\scrM'\ar[d]_{\bar h}&\scrM^{\ddagger}_{y'}\ar[l]_{q'}\ar[d]^{h'}\\
\foM_{\tau_y}^{\ev}\ar[r]^p\ar[d]& \foM' \ar[d]_{\overline{\Phi}} & 
\foM_{y'}^{\ddagger,\ev}\ar[l]_{p'}
\ar[d]\\
B\GG_m^{\tau_y}\ar[r] & \overline{\scrM}_{0,4}' & y'\times B\GG_m^{\dagger}\ar[l]
}
\]
Now $[B\GG_m^{\tau_y}]$ and $[y'\times B\GG_m]$ define the same element
in $A_{-1}(\overline{\scrM}_{0,4}')$. Thus we have
\[
q_*[\scrM_{\tau_y}]^{\virt} = \bar h^! \overline{\Phi}^* [B\GG_m^{\tau_y}]=
\bar h^! \overline{\Phi}^* [y'\times B\GG_m] = q'_* [\scrM^{\ddagger}_{y'}]^{\virt},
\]
hence the result.
\end{proof}

\begin{proof}[Proof of Theorem~\ref{theorem: independence of modulus}]
Fix a boundary point $y\in\overline\scrM_{0,4}^{\dagger}$.
Consider the map $\Sigma(\Phi):\Sigma(\foM^{\ev}(\shX,\beta,z))
\rightarrow\Sigma(\overline{\scrM}_{0,4}^{\dagger})$. We can find a refinement
$\Sigma$ of $\sigma_y=\RR_{\ge 0}\ell^*\oplus\RR_{\ge 0}\delta^*$ 
into standard cones with the property that
for any $\sigma \in \Sigma(\foM^{\ev}(\shX,\beta,z))$ with
$\Sigma(\Phi)(\sigma)\subset\sigma_y$, we can write $\Sigma(\Phi)(\sigma)$ as
a union of cones in $\Sigma$. We can further assume that one cone
in $\Sigma$ is $\RR_{\ge 0}(\lambda \ell^*+\delta^*)+\RR_{\ge 0} \ell^*$
for sufficiently large $\lambda$. Note to do this we need to shrink
$\foM^{\ev}(\shX,\beta,z)$ as in \S\ref{sec:stacks remarks}
to ensure $\Sigma(\foM^{\ev}(\shX,\beta,z))$ only contains
a finite number of cones.

It now follows from Theorem~\ref{theorem: transversality} that
for any $\tau\in \Sigma$ of dimension two, $B\GG_m^{\tau}\rightarrow 
\overline{\scrM}_{0,4}^{\dagger}$ is transverse to $\Phi$. Thus we can
apply Lemma~\ref{lemma: invarianceI} repeatedly to show that 
$\deg[\scrM_{\tau}]^{\virt}$ is independent of $\tau\in\Sigma$,
and then apply Lemma~\ref{lemma: invarianceII} to show 
Theorem~\ref{theorem: independence of modulus}.
\end{proof}

\section{Gluing}
\label{sec:key gluing}

\subsection{Review of gluing}
\label{sec:gluing review}
We will begin by reviewing certain constructions in \cite{ACGS18}
for gluing punctured maps, and adapt those constructions
to our needs. 
That reference considers quite general gluing situations, and
as here we only ever have to glue two punctured maps to obtain
an additional node in the domain, we give a brief summary of 
the results of \cite{ACGS18} in our context, using slightly
different notation so as to avoid introducing the general notation
of \cite{ACGS18}.

We assume given two classes $\beta_1,\beta_2$ of punctured map to $X$,
with underlying curve classes $A_1,A_2$ and
with $\beta_i$ having a puncture $x_i$ and an additional set of
punctured points ${\bf x}_i$. We assume the global contact orders
$\bar u_{x_i}$ are related by $u_{x_1}=-u_{x_2}$ for some
choice of representatives $u_{x_i}\in N_{\sigma}$ of $\bar u_{x_i}$
and $\sigma\in\Sigma(X)$.
Furthermore, we may take $\beta$
to be the class of punctured map with underlying curve class
$A=A_1+A_2$ and a set of punctured points ${\bf x}={\bf x}_1\cup
{\bf x}_2$ with the punctures in ${\bf x}_i$ having the same 
global contact orders as the corresponding points of the class 
$\beta_i$. As usual, all genera are assumed to be $0$.

We may then define the stack $\foM^{\gl}(\shX,\beta_1,\beta_2)$
as follows. This is the stack of punctured maps $f:C^\circ/W\rightarrow
\shX$ along with data of pre-stable curves $\ul{C}_i\rightarrow \ul{W}$ 
with marked points indexed by the set ${\bf x}_i\cup \{x_i\}$ and
an isomorphism of marked pre-stable curves $\ul{C}\cong 
\ul{C}_1\amalg_{x_1=x_2} \ul{C}_2$. Here, the latter notation denotes
the gluing of $\ul{C}_1$ and $\ul{C}_2$ along the sections $x_1,x_2$
to produce a node, producing a nodal section $q:\ul{W}\rightarrow \ul{C}$. 
This data should satisfy the condition that the splitting of
$f$ at $q$ (see \cite{ACGS18},~Prop.~5.2) yields punctured maps
$f_i:C_i^{\circ}/W\rightarrow \shX$ of class $\beta_i$.

Note that there is a canonical morphism
\[
\foM^{\gl}(\shX,\beta_1,\beta_2)\rightarrow \foM(\shX,\beta).
\]
We may then define 
\[
\scrM^{\gl}(X,\beta_1,\beta_2)\subseteq
\foM^{\gl}(\shX,\beta_1,\beta_2)
\times_{\foM(\shX,\beta)} \scrM(X,\beta)
\]
to be the open and closed substack of the right-hand side
defined as follows.
Any object in the fibre product consists of a stable punctured map
$f:C^{\circ}/W\rightarrow X$ along with a decomposition
$\ul{C}=\ul{C}_1 \amalg_{x_1=x_2} \ul{C}_2$, and we require
that $\ul{f}|_{\ul{C}_i}$ be a stable map of curve class $A_i$.
Thus we have a commutative diagram
\[
\xymatrix@C=30pt
{
\scrM^{\gl}(X,\beta_1,\beta_2)\ar[d]_{\varepsilon'}\ar[r]&
\scrM(X,\beta)\ar[d]^{\varepsilon''}\\
\foM^{\gl}(\shX,\beta_1,\beta_2)\ar[r]&\foM(\shX,\beta)
}
\]
and the relative obstruction theory for $\varepsilon''$ pulls
back to a relative obstruction theory for $\varepsilon'$,
see \cite{ACGS18},~Prop.~5.21. Note this is possible since 
$\scrM^{\gl}(X,\beta_1,\beta_2)$ is a union of connected components of the
fibre product of the other three vertices of the square (since fixing
curve classes is an open and closed condition).

\begin{remark}
In the language of \cite{ACGS18}, the moduli spaces
$\foM^{\gl}(\shX,\beta_1,\beta_2)$ and $\scrM^{\gl}(X,\beta_1,\beta_2)$ are
written as $\foM'(\shX,\tau)$ and $\scrM'(X,\boldsymbol{\tau})$, where
$\tau=(G,{\bf g},\bsigma,\bar{\bf u})$ with $G$ a graph
with two vertices, one edge connecting the two vertices,
and a number of legs corresponding to ${\bf x}_1$ and ${\bf x}_2$.
The genus function ${\bf g}$ takes value $0$ on all vertices,
the map $\bsigma$ assigns $0\in \Sigma(X)$ to every leg, edge or
vertex of $G$. Next $\bar{\bf u}$ assigns to each leg of $G$
the corresponding global contact order $\bar{\bf u}_x$ for
$x\in {\bf x}_1\cup {\bf x}_2$ and to the edge of $G$
the global contact order represented by $\pm u_{x_1}$, depending on
orientation. Further, the assignment of curve
classes $A_1$ and $A_2$ to the two vertices of $G$ determines
$\boldsymbol{\tau}$.
\end{remark}

We need to understand gluing at a virtual level. To do so,
after choosing (possibly empty) subsets ${\bf x}_1'\subseteq {\bf x}_1$
and ${\bf x}_2'\subseteq {\bf x}_2$, we define as in \eqref{eq:Mev def}:
\begin{align*}
\foM^{\ev(x_i,{\bf x}'_i)}(\shX,\beta_i):= {}  &
 \foM(\shX,\beta_i) \times_{\ul{\shX}
\times \prod_{x\in {\bf x}'_i} \ul{\shX}}\left( \ul{X}
\times \prod_{x\in {\bf x}'_i} \ul{X}\right),\\
\foM^{\gl,\ev(q,{\bf x}_1'\cup {\bf x}_2')}(\shX,\beta_1,\beta_2):= {}
& \foM^{\gl}(\shX,\beta_1,\beta_2)
\times_{\ul{\shX}
\times \prod_{x\in {\bf x}'_1\cup {\bf x}'_2} \ul{\shX}}\left( \ul{X}
\times \prod_{x\in {\bf x}'_1\cup {\bf x}'_2} \ul{X}\right).
\end{align*}
Here, the map $\foM(\shX,\beta_i)\rightarrow\ul{\shX}
\times \prod_{x\in {\bf x}_i'}\ul{\shX}$ is given by
evaluation at $x_i$ and $x\in {\bf x}_i'$, 
while the map $\foM^{\gl}(\shX,\beta_1,\beta_2)
\rightarrow\ul{\shX}\times \prod_{x\in {\bf x}'_1\cup {\bf x}'_2}\ul{\shX}$ 
is given by evaluation at the node $q$ obtained
by gluing and the points $x\in {\bf x}_1'\cup {\bf x}_2'$.
We omit the decoration on $\ev$ when clear from context.
We then obtain a natural factorization of $\varepsilon'$:
\[
\xymatrix@C=30pt
{
\scrM^{\gl}(X,\beta_1,\beta_2)
\ar[r]^>>>>>{\varepsilon}& 
\foM^{\gl,\ev(q,{\bf x}_1'\cup {\bf x}_2')}(\shX,\beta_1,\beta_2)
\ar[r]^>>>>>{\pi} & 
\foM^{\gl}(\shX,\beta_1,\beta_2).
}
\]
As in \S\ref{subsub:moduli} above, 
\cite{ACGS18},~\S\S4.2,4.3 constructs a relative
obstruction theory for $\varepsilon$. As $\pi$ is smooth, there is a canonical
choice of relative obstruction theory for $\pi$ such that virtual
pull-back coincides with flat pull-back of cycles. Together
with the relative obstruction theories for $\varepsilon$ and $\varepsilon'$,
we obtain a compatible
triple of obstruction theories,
and we say the obstruction theories for $\varepsilon$ and $\varepsilon'$ are compatible.

The main gluing result of \cite{ACGS18},~Thm.~C,
then states in this case:

\begin{theorem}
\label{thm:gluing recall}
Given the above setup, there is a diagram cartesian in all categories
\begin{equation}
\label{eq:first gluing}
\xymatrix@C=40pt
{
\scrM^{\gl}(X,\beta_1,\beta_2)\ar[r]^>>>>>>>>>>>>>{\psi'}\ar[d]_\varepsilon&
\scrM(X,\beta_1)\times\scrM(X,\beta_2)\ar[d]^{\varepsilon_1\times \varepsilon_2}\\
\foM^{\gl,\ev(q,{\bf x}_1'\cup {\bf x}_2')}(\shX,\beta_1,\beta_2)
\ar[r]_>>>>>>{\psi}&
\foM^{\ev(x_1,{\bf x}_1')}(\shX,\beta_1)
\times\foM^{\ev(x_2,{\bf x}_2')}(\shX,\beta_2)
}
\end{equation}
where $\psi$ is representable and finite, defined by splitting
a punctured map at the glued node.
Further, the product of the perfect relative obstruction theories for 
$\varepsilon_1$ and
$\varepsilon_2$ pulls back to the perfect relative obstruction theory for
$\varepsilon$.
\end{theorem}

As we will often need to glue when one of the moduli spaces has a point
constraint, we need the following minor modification of the above
result:

\begin{corollary}
\label{cor:gluing recall}
Given the above setup, suppose in addition that one of the punctures
$x_{\out}\in {\bf x}_2'$ has contact order specified by $-r$, for
$r\in \Sigma(X)(\ZZ)$, and let $z\in Z_r^{\circ}$. Set
\begin{align*}
\scrM^{\gl}(X,\beta_1,\beta_2,z):= {} & \scrM^{\gl}(X,\beta_1,\beta_2)
\times^{\fs}_{\scrP(X,r)} B\GG_m^{\dagger}\\
\foM^{\gl,\ev(q,{\bf x}_1'\cup{\bf x}_2')}(\shX,\beta_1,\beta_2,z):= {} & 
\foM^{\gl,\ev(q,{\bf x}_1'\cup{\bf x}_2')}(\shX,\beta_1,\beta_2)
\times^{\fs}_{\scrP(X,r)} B\GG_m^{\dagger}
\end{align*}
where the morphisms $\scrM^{\gl}(X,\beta_1,\beta_2),
\foM^{\gl,\ev(q,{\bf x}_1'\cup{\bf x}_2')}(\shX,\beta_1,\beta_2)\rightarrow\scrP(X,r)$ are
the evaluation maps given by Proposition~\ref{prop:evaluation maps general}
at $x_{\out}$
(as in Lemma~\ref{lem:evaluation diagrams} in the case
of $\foM^{\gl,\ev(q,{\bf x}_1'\cup{\bf x}_2')}(\shX,\beta_1,\beta_2)
\rightarrow\scrP(X,r)$),
and the morphism $B\GG_m^{\dagger}\rightarrow\scrP(X,r)$ is given
by Proposition~\ref{prop: BGm morphism}.
Then there is a diagram cartesian in all categories
\[
\xymatrix@C=40pt
{
\scrM^{\gl}(X,\beta_1,\beta_2,z)\ar[r]^>>>>>>>>>>>>>>{\psi'}\ar[d]_\varepsilon&
\scrM(X,\beta_1)\times\scrM(X,\beta_2,z)\ar[d]^{\varepsilon_1\times \varepsilon_2}\\
\foM^{\gl,\ev(q,{\bf x}_1'\cup{\bf x}_2')}(\shX,\beta_1,\beta_2,z)\ar[r]_>>>>>>{\psi}&
\foM^{\ev(x_1,{\bf x}_1')}(\shX,\beta_1)\times
\foM^{\ev(x_2,{\bf x}_2')}(\shX,\beta_2,z)
}
\]
with $\scrM(X,\beta_2,z)$, $\foM^{\ev(x_2,{\bf x}_2')}(\shX,\beta_2,z)$
as in Definition~\ref{def:Mbetaev}.
Further, $\varepsilon_2$ carries a perfect relative obstruction theory
pulled back from the right-hand morphism in 
the cartesian diagram
\[
\xymatrix@C=20pt
{
\scrM(X,\beta_2,z)\ar[r]\ar[d]_{\varepsilon_2}&\scrM(X,\beta_2)\ar[d]\\
\foM^{\ev(x_2,{\bf x}_2')}(\shX,\beta_2,z)\ar[r]&\foM^{\ev(x_2,{\bf x}_2')}(\shX,\beta_2)
}
\]
that is compatible with the obstruction theory defined in
Proposition~\ref{prop:rel virt dim}. The product of the
relative obstruction theories for $\varepsilon_1$ and $\varepsilon_2$
pulls back to a relative obstruction theory for
$\varepsilon$, which is compatible 
with the pull-back perfect relative obstruction theory
from the right-hand arrow in the diagram
\[
\xymatrix@C=30pt
{
\scrM^{\gl}(X,\beta_1,\beta_2,z)\ar[d]_{\varepsilon}\ar[r]&
\scrM^{\gl}(X,\beta_1,\beta_2)\ar[d]\\
\foM^{\gl,\ev(q,{\bf x}_1'\cup {\bf x}_2')}(\shX,\beta_1,\beta_2,z)\ar[r]&
\foM^{\gl,\ev(q,{\bf x}_1'\cup {\bf x}_2')}(\shX,\beta_1,\beta_2)
}
\]
\end{corollary}

\begin{proof}
The diagram \eqref{eq:first gluing} is
a diagram of log stacks over $\scrP(X,r)$ cartesian in all categories. 
We have already explained how to view $\scrM^{\gl}(X,\beta_1,\beta_2)$,
$\foM^{\gl,\ev(q,{\bf x}_1'\cup {\bf x}_2')}(\shX,\beta_1,\beta_2)$ 
as stacks over $\scrP(X,r)$,
while we have a composed morphism 
\[
\scrM(X,\beta_1)
\times \scrM(X,\beta_2)\rightarrow \scrM(X,\beta_2)\rightarrow
\scrP(X,r)
\]
with first morphism the projection and second morphism
the similarly defined evaluation map. One defines a morphism
\[
\foM^{\ev(x_1,{\bf x}_1')}(\shX,\beta_1)\times 
\foM^{\ev(x_2,{\bf x}_2')}(\shX,\beta_2)
\rightarrow\scrP(X,r)
\]
similarly. Thus, applying the fibre product
$\times^{\fs}_{\scrP(X,r)} B\GG_m^{\dagger}$ to \eqref{eq:first gluing},
we obtain the claimed cartesian diagram of the corollary. Note it is
cartesian in all categories as the vertical arrows are still strict.
The statements on obstruction theories then follows, as the obstruction
theories involved are all defined as pull-backs of the obstruction
theories appearing in Theorem~\ref{thm:gluing recall}.
\end{proof}

\medskip

We next review the fibre product description of the glued moduli spaces
given in \S5.2 of \cite{ACGS18}.
If $C^{\circ}\rightarrow \foM^{\gl}(\shX,\beta_1,\beta_2)$ is the universal
curve, then the glued node provides a nodal section
$q:\ul{\foM^{\gl}(\shX,\beta_1,\beta_2)}\rightarrow \ul{C}$, and we
define 
\[
\widetilde\foM^{\gl}(\shX,\beta_1,\beta_2):=
\big(\ul{\foM^{\gl}(\shX,\beta_1,\beta_2)},q^*\shM_{C^{\circ}}\big),
\]
the log structure on the moduli space pulled back from the node.
Similarly, if $C_i^{\circ}\rightarrow \foM(\shX,\beta_i)$ is the
universal curve, then we define
\[
\widetilde\foM(\shX,\beta_i):=
\big(\ul{\foM(\shX,\beta_i)},x_i^*\shM_{C^{\circ}_i}\big)^{\sat},
\]
recalling that the log structure at a puncture need only be fine, and thus
we saturate this pull-back log structure. This may change the underlying
stack. Note we then have evaluation morphisms
$\ev_i:\widetilde\foM(\shX,\beta_i)\rightarrow \shX$
given by composing the saturation morphism with $f_i\circ x_i$.

Similarly, we may define $\widetilde{\foM}^{\gl,\ev(q,{\bf x}_1'
\cup{\bf x}_2')}(\shX,\beta_1,\beta_2)$ and
$\widetilde{\foM}^{\ev(x_i,{\bf x}_i')}(\shX,\beta_i)$.
Note because $\ev_i$ is defined at the logarithmic level, we 
in particular have
\[
\widetilde{\foM}^{\ev(x_i,{\bf x}_i')}(\shX,\beta_i)
=
\widetilde{\foM}(\shX,\beta_i)\times_{\shX\times \prod_{x\in {\bf x}_i'}
\ul{\shX}} X\times \prod_{x\in {\bf x}_i'}\ul{X}.
\]
Thus we obtain similarly an evaluation morphism
$\ev_i:\widetilde{\foM}^{\ev(x_i,{\bf x}_i')}(\shX,\beta_i)
\rightarrow X$.

We have the following description,
a special case of \cite{ACGS18},~Thm.~5.8 and Cor.~5.15.

\begin{theorem}
\label{thm:gluing review}
Suppose that the global contact orders $\bar u_{x_i}$, $i=1,2$, are
monodromy free. Then we have
\[
\widetilde\foM^{\gl}(\shX,\beta_1,\beta_2)
\cong
\widetilde\foM(\shX,\beta_1)\times^{\fs}_{\shX}
\widetilde\foM(\shX,\beta_2),
\]
and
\[
\widetilde\foM^{\gl,\ev(q,{\bf x}_1'\cup{\bf x}_2')}(\shX,\beta_1,\beta_2)
\cong
\widetilde\foM^{\ev(x_1,{\bf x}_1')}(\shX,\beta_1)\times^{\fs}_{X}
\widetilde\foM^{\ev(x_2,{\bf x}_2')}(\shX,\beta_2),
\]
with the fibre product defined using the morphisms
$\ev_i$ defined above.
\end{theorem}

\begin{remark}
\label{rem:glued log structure}
The log stack $\foM^{\gl}(\shX,\beta_1,\beta_2)$ has the same
underlying stack as $\widetilde \foM^{\gl}(\shX,\beta_1,\beta_2)$
by the construction of the latter, and the log structure of the
former is a sub-log structure of the latter. 
This may be seen explicitly as follows; for a deeper understanding
of the gluing process, see
the proof of \cite{ACGS18},~Thm.~5.8.

If locally near the
node the universal curve over $\foM^{\gl}(\shX,\beta_1,\beta_2)$
is given by the equation $xy=0$, then the coordinates $x$ and $y$
lift to sections $s_x$ and $s_y$ of the log structure on the universal
curve with a relation $s_xs_y=s_{t}$ for some section $s_t$ of the
log structure pulled back from $\foM^{\gl}(\shX,\beta_1,\beta_2)$.
Thus the log structure on $\widetilde \foM^{\gl}(\shX,\beta_1,\beta_2)$
may be seen as being generated by the log structure on
$\foM^{\gl}(\shX,\beta_1,\beta_2)$ along with sections $s_x,s_y$ subject
to the relation $s_xs_y=s_t$. 

On the other hand, there are morphisms
\[
\xymatrix@C=20pt
{\widetilde\foM^{\gl}(\shX,\beta_1,\beta_2)\ar[r]^>>>>{\pr_i}&
\widetilde\foM(\shX,\beta_i)\ar[r]^{s_i}&\foM(\shX,\beta_i)
}
\]
with $\pr_i$ the projection and $s_i$ the morphism obtained by
composing the saturation morphism with the canonical map
$\big(\ul{\foM(\shX,\beta_i)},x_i^*\shM_{C^{\circ}_i}\big)
\rightarrow \foM(\shX,\beta_i)$. Then
$s_i^{-1}\overline{\shM}_{\foM(\shX,\beta_i)}\oplus\NN \subseteq
\overline{\shM}_{\widetilde\foM(\shX,\beta_i)}\subseteq
s_i^{-1}\overline{\shM}_{\foM(\shX,\beta_i)}\oplus\ZZ$ by
definition of a punctured log structure. In particular, there
is a sub-log structure $\shP_i$ of the log structure on
$\widetilde\foM(\shX,\beta_i)$
with ghost sheaf given by the factor $\NN$, generated by sections $s_x$
or $s_y$ lifting coordinates $x$ and $y$ as before. These sections
map via $\pr_i^{\flat}$ to the corresponding sections of the log structure on
$\widetilde\foM^{\gl}(\shX,\beta_1,\beta_2)$.

There is then an inclusion of  
log structures $\shN\subset
\pr_1^*\shP_1\oplus_{\O^{\times}}\pr_2^*\shP_2$ 
on the underlying space of $\widetilde\foM(\shX,\beta_1)\times 
\widetilde\foM(\shX,\beta_2)$ corresponding to the inclusion of ghost
sheaves the diagonal $\NN\subseteq \NN\oplus\NN$. 
Then the log structure of $\foM^{\gl}(\shX,\beta_1,\beta_2)$
is the fine saturated sub-log structure
of the log structure on $\widetilde\foM^{\gl}(\shX,\beta_1,\beta_2)$ generated 
by the images under $(s_i\circ \pr_i)^{\flat}$ of $\shM_{\foM(\shX,\beta_i)}$, 
$i=1,2$, and the image under $(\pr_1\times \pr_2)^{\flat}$ of 
$\shN$.
\end{remark}

It will be useful to recast this fibre product description
in the most important case we will need here.
 
\begin{theorem}
\label{thm:new gluing}
Suppose that the global contact orders $\bar u_{x_i}$, $i=1,2$,
are monodromy-free. Suppose further that $\bar u_{x_1}$ is 
specified by $-s$ with $s\in \Sigma(X)(\ZZ)$, and $\bar u_{x_2}$
is specified by $s$. Then there exist morphisms
\begin{equation}
\label{eq:special eval}
\foM(\shX,\beta_2)\times B\GG_m^{\dagger}\rightarrow\scrP(\shX,s),
\quad
\scrM(X,\beta_2)\times B\GG_m^{\dagger}\rightarrow\scrP(X,s)
\end{equation}
yielding
diagrams cartesian in the fs log category
\[
\xymatrix@C=30pt
{
\foM^{\gl}(\shX,\beta_1,\beta_2)\ar[r]\ar[d] &
\foM(\shX,\beta_1)\ar[d]^{\ev}\\
\foM(\shX,\beta_2)\times B\GG_m^{\dagger}\ar[r]&\scrP(\shX,s)
}
\]
and
\[
\xymatrix@C=30pt
{
\scrM(X,\beta_1,\beta_2)\ar[r]\ar[d]&
\scrM(X,\beta_1)\ar[d]^{\ev}\\
\scrM(X,\beta_2)\times B\GG_m^{\dagger}\ar[r]&\scrP(X,s)
}
\]
where $\ev$ is as given in Proposition~\ref{prop:evaluation maps general}
using the puncture $x_1$.
\end{theorem}

\begin{proof}
We consider the first diagram.
For notational convenience, we write $W_i=\foM(\shX,\beta_i)$,
$\widetilde W_i=\widetilde\foM(\shX,\beta_i)$, $W=\foM^{\gl}(\shX,\beta_1,
\beta_2)$ and $\widetilde W=\widetilde{\foM}^{\gl}(\shX,\beta_1,\beta_2)$.

First consider the evaluation morphisms $\ev_i:\widetilde W_i\rightarrow
\shX$. By Remark~\ref{rem:contact order remarks}, (1), $\ev_i$
factors through the stratum $\shZ:=\shZ_s$ of $\shX$. Second,
note that since $x_1$ has contact order given by $-s$, it follows
from Proposition~\ref{prop:evaluation maps general} and the construction
of $\widetilde W_1$ that $\widetilde W_1$ is the saturation of
$W_1\times^{\fine}_{\scrP(\shX,s)} \widetilde{\scrP}(\shX,s)$, i.e.,
\begin{equation}
\label{eq: tilde W}
\widetilde W_1 \cong W_1\times^{\fs}_{\scrP(X,s)} \widetilde{\scrP}(X,s).
\end{equation}

On the other hand, note that the underlying stack of $\widetilde W_2$ coincides
with $\ul{W}_2$; indeed, 
the contact order $s$ of $x_2$ lies in $\Sigma(X)$,
so that $x_2$ is a marked point. In particular,
$\overline{\shM}_{\widetilde W_2}=\overline{\shM}_{W_2}\oplus\NN$.

We then have by Proposition~\ref{thm:gluing review} that $\widetilde W$
coincides with
\begin{align*}
\widetilde W_1\times^{\fs}_{\shX} \widetilde W_2\cong {} & 
\widetilde W_1\times^{\fs}_{\shZ} \widetilde W_2\\
\cong {} & \widetilde W_1\times^{\fs}_{\shZ\times B\GG_m^{\dagger}} 
(\widetilde W_2\times B\GG_m^{\dagger})\\
\cong {} & 
\widetilde W_1\times^{\fs}_{\widetilde\scrP(\shX,s)} 
(\widetilde W_2\times B\GG_m^{\dagger})\\
\cong {} & 
(W_1\times^{\fs}_{\scrP(\shX,s)}\widetilde{\scrP}(\shX,s))
\times^{\fs}_{\widetilde\scrP(\shX,s)} 
(\widetilde W_2\times B\GG_m^{\dagger})\\
\cong {} & W_1\times^{\fs}_{\scrP(\shX,s)} 
(\widetilde W_2\times B\GG_m^{\dagger}).
\end{align*}
Here for the second isomorphism we use the morphism
$\widetilde\ev:\widetilde W_1\rightarrow \shZ\times B\GG_m^{\dagger}
=\widetilde\scrP(\shX,s)$
of \eqref{eq:ev prime}
through which $\widetilde W_1\rightarrow \shZ$ factors.
The morphism $\widetilde W_2\times B\GG_m^{\dagger}
\rightarrow \shZ\times B\GG_m^{\dagger}$ is just given by the identity
on the $B\GG_m^{\dagger}$ factor. For the fourth isomorphism, we use
\eqref{eq: tilde W}.

We next consider the sub-log structure on 
$\widetilde W_2\times B\GG_m^{\dagger}$ defined as follows.
The ghost sheaf of $\widetilde W_2\times B\GG_m^{\dagger}$ 
is $\overline{\shM}_{W_2}\oplus
\NN\oplus\NN$, and the torsors corresponding to
$(0,1,0)$ and $(0,0,1)$ are $\shN_{x_2/C_2}^{\vee}$ and $\shU$
respectively, where $C_2\rightarrow W_2$ is the universal domain over
$W_2$ and as usual $\shU$ is the universal torsor on
$B\GG_m$. Consider the subsheaf of this ghost sheaf
generated by $\overline{\shM}_{W_2}$ and $(0,1,1)$. This defines
the sub-log structure, which we call $W'_2$. Note there is an isomorphism
\begin{equation}
\label{eq:W2 prime} W'_2\stackrel{\cong}{\longrightarrow} W_2\times B\GG_m^{\dagger}.
\end{equation}
This is defined by specifying
morphisms $W'_2\rightarrow W_2$ and $W'_2\rightarrow B\GG_m$.
The first morphism is defined by noting that
the composition $\widetilde W_2\times B\GG_m^{\dagger}\rightarrow
\widetilde W_2\rightarrow W_2$ factors through $W_2'$. On the other
hand, the morphism $\ul{W}'_2\rightarrow B\GG_m$ is given by
the line bundle $\shU\otimes\shN^{\vee}_{x_2/C_2}$, and this immediately
lifts to a morphism $W'_2\rightarrow B\GG_m^{\dagger}$ which on the level
of ghost sheaves takes $1\in \NN$ to $(0,1,1)\in \overline{\shM}_{W_2'}$.

Note the morphism $\widetilde W_2\times B\GG_m^{\dagger}\rightarrow
\scrP(\shX,s)$ factors through $W_2'$. Indeed, by construction
we have a factorization
\begin{equation}
\label{eq:factorization}
\widetilde W_2\times B\GG_m^{\dagger}
\rightarrow \shZ\times B\GG_m^{\dagger} \rightarrow \scrP(\shX,s),
\end{equation}
which translates at the level of stalks of ghost sheaves as
\begin{equation}
\label{eq:ghost sheaf map}
Q_2\oplus\NN\oplus\NN \leftarrow P\oplus\NN \leftarrow P
\end{equation}
with the first map given by $(p,n)\mapsto (q,s(p),n)$ for some $q$
depending on $p$
and the second map given by $p\mapsto (p,s(p))$. 
Thus the composition is $p\mapsto
(q,s(p),s(p))$, which lies in the ghost sheaf of $W_2'$ by construction.

We now define
$W'$ to make the right-hand square of the following diagram cartesian in the
fs log category:
\[
\xymatrix@C=30pt
{
\widetilde W_1\times_{\shZ}^{\fs}\widetilde W_2\ar[d]\ar[r]
& W' \ar[d]\ar[r]& W_1\ar[d]\\
\widetilde W_2\times B\GG_m^{\dagger}\ar[r]&W_2'\ar[r]&\scrP(\shX,s)
}
\]
where $W_2'$ is as in \eqref{eq:W2 prime}.
Since the big rectangle is also cartesian, so is the left-hand rectangle.
If we show that $W'\cong W$, then the right-hand cartesian square
is the desired diagram of the theorem.

First note that by construction,
$\widetilde W_2\times B\GG_m^{\dagger}\rightarrow W_2'$ is
an integral and saturated morphism, and the identity on underlying stacks.
Thus the stack underlying the fs fibre product 
of the left square agrees with the ordinary
fibre product, so $\widetilde W_1\times_{\shZ}^{\fs}\widetilde W_2$
and $W'$ have the same underlying stack, and $W'$ is a sub-log
structure of $\widetilde W_1\times_{\shZ}^{\fs}\widetilde W_2$,
as can easily be checked on the level of stalks of ghost sheaves.

Thus it is enough to check that $W'$ is the same sub-log structure
as that of $W$ as described in Remark~\ref{rem:glued log structure}.
Indeed, the ghost sheaf of $W'$ is the smallest fine saturated 
subsheaf of
the ghost sheaf of $\widetilde W_1\times_{\shZ}^{\fs}\widetilde W_2$
containing the images of
$\overline{\shM}_{W_1}$ and $\overline{\shM}_{W_2'}$. The latter
is generated by $\overline{\shM}_{W_2}$ and an $\NN$ factor.
By tracing through the definition of the composed morphism
$\widetilde W_1\times_{\shZ}^{\fs}\widetilde W_2\rightarrow
\widetilde W_2\times B\GG_m^{\dagger}\rightarrow W_2'$,
one sees that this $\NN$ factor generates $\overline{\shN}$, 
where $\shN$ is as described in Remark~\ref{rem:glued log structure}.
Hence the result for the first diagram of the lemma.

For the second diagram, we proceed exactly as above, but replace the
application of Proposition \ref{thm:gluing review} with the analogous
statement of \cite{ACGS18},~Cor.~5.15, and use Theorem~\ref{thm:gluing recall}
to pass from moduli spaces of maps to $\shX$ to moduli spaces of maps to
$X$. 
\end{proof}

\begin{remark}
\label{rem:trop of gluing}
From the explicit description of the morphisms \eqref{eq:special eval}
given in the above proof, especially \eqref{eq:ghost sheaf map}, we
can understand the cartesian diagrams of Theorem \ref{thm:new gluing}
at a tropical level. Let $\bar w_i$ be a geometric point of 
$\foM(\shX,\beta_i)$ or $\scrM(X,\beta_i)$, $i=1,2$, with corresponding cones
$\sigma_{\bar w_i}$ parameterizing 
families of tropical maps $h_i$ to $\Sigma(X)$
with domain $G_i$. Let $v_{i,\out}$ be the vertex of $G_i$ adjacent to the
leg corresponding to $x_i$. 
Then the morphisms of \eqref{eq:special eval} tropicalize to
the map
\[
\sigma_{\bar w_2}\times\RR_{\ge 0}
\rightarrow \Sigma(X), \quad (m,t) \mapsto (h_2)_m(v_{2,\out})+t s.
\]
Thus we see via the tropical interpretation of $\ev$ in 
Lemma~\ref{lem:tropical point constraint} that if there is a geometric point
$\bar w$ of $\foM^{\gl}(\shX,\beta_1,\beta_2)$ mapping to $\bar w_1$ and
$\bar w_2$ respectively, then the cone $\sigma_{\bar w}$ associated
to this point satisfies
\[
\sigma_{\bar w}=\{(m_1,(m_2,t))\in \sigma_{\bar w_1}\times
(\sigma_{\bar w_2}\times\RR_{\ge 0})\,|\,
(h_1)_{m_1}(v_{1,\out})=(h_2)_{m_2}(v_{2,\out})+t s\}.
\]
It is clear this parameterizes a family of tropical maps with domain
obtained by gluing $G_1,G_2$ along the legs corresponding to $x_1,x_2$.
\end{remark}

\begin{corollary}
\label{cor:new gluing}
Suppose we are in the setup of Corollary~\ref{cor:gluing recall}. Suppose 
further that the contact
order of $x_1$ is specified by $-s$, with $s\in \Sigma(X)(\ZZ)$, and
the contact order of $x_2$ is specified by $s$. Then
there are fs cartesian diagrams
\[
\xymatrix@C=30pt
{
\foM^{\gl,\ev(q,{\bf x}_1'\cup {\bf x}_2')}(\shX,\beta_1,\beta_2,z)\ar[r]\ar[d]&
\foM^{\ev(x_1,{\bf x}_1')}(\shX,\beta_1)\ar[d]\\
\foM^{\ev(x_2,{\bf x}_2')}(\shX,\beta_2,z)\times B\GG_m^{\dagger}\ar[r]&\scrP(X,s)
}
\]
and
\[
\xymatrix@C=30pt
{
\foM^{\gl,\ev({\bf x}_1'\cup {\bf x}_2')}(\shX,\beta_1,\beta_2,z)\ar[r]\ar[d]&
\foM^{\ev({\bf x}_1')}(\shX,\beta_1)\ar[d]\\
\foM^{\ev({\bf x}_2')}(\shX,\beta_2,z)\times B\GG_m^{\dagger}\ar[r]&\scrP(\shX,s)
}
\]
\end{corollary}

\begin{proof}
We do the first case, the second being slightly easier.
First, note that each stack in the first diagram of
Theorem~\ref{thm:new gluing} comes with a morphism to $\ul{\shX}$,
given either by evaluation at the glued node, at 
$x_1$ or $x_2$, or the canonical
projection $\scrP(\shX,s)\rightarrow \ul{\shZ}_s$ followed by the closed
embedding $\ul{\shZ}_s\hookrightarrow \ul{\shX}$. Further, all these morphisms
are compatible. Noting that $\scrP(X,s)=\scrP(\shX,s)\times_{\ul{\shX}}\ul{X}$,
we then have from Theorem~\ref{thm:new gluing} an fs cartesian diagram
\[
\xymatrix@C=30pt
{
\foM^{\gl}(\shX,\beta_1,\beta_2)\times_{\ul{\shX}} \ul{X}\ar[r]\ar[d] &
\foM(\shX,\beta_1)\times_{\ul{\shX}} \ul{X}\ar[d]\\
(\foM(\shX,\beta_2)\times_{\ul{\shX}} \ul{X})\times B\GG_m^{\dagger}\ar[r]
&\scrP(X,s)
}
\]
It then follows from general properties of fibre products that
\[
\xymatrix@C=30pt
{
\foM^{\gl,\ev(q,{\bf x}_1'\cup {\bf x}_2')}(\shX,\beta_1,\beta_2)\ar[r]\ar[d] &
\foM^{\ev(x_1,{\bf x}_1')}(\shX,\beta_1)\ar[d]\\
\foM^{\ev(x_2,{\bf x}_2')}(\shX,\beta_2)\times B\GG_m^{\dagger}\ar[r]&\scrP(X,s)
}
\] 
is fs cartesian, and a further base-change on the left by the projection
morphism
$\foM^{\ev(x_2,{\bf x}_2')}(\shX,\beta_2,z)\times B\GG_m^{\dagger}\rightarrow
\foM^{\ev(x_2,{\bf x}_2')}(\shX,\beta_2)\times B\GG_m^{\dagger}$ then gives the desired
fs cartesian diagram.
\end{proof}

\subsection{Alternative calculations for $N^{A_1}_{p_1p_2s}$}

We fix here a punctured map class $\beta_1$, of genus zero with three marked
points $x_1,x_2,x_{\out}$, with contact orders specified by $p_1,p_2$ and $-s$
for $p_1,p_2,s\in \Sigma(X)(\ZZ)$, and curve class $A_1$. 
We return to the conventions of
Definitions~\ref{def:Mev} and \ref{def:Mbetaev} for the notation
$\foM^{\ev}(\shX,\beta_1)$, $\foM^{\ev}(\shX,\beta_1,z)$.

In the proof of the key gluing result, Theorem~\ref{thm: main gluing theorem},
we glue the moduli space $\scrM(X,\beta_2,z)$ to a moduli space
$\scrM(X,\beta_1)$ without a point constraint. Suppose given a punctured map 
$f_2:(C_2,x_s',x_3,x_{\out})\rightarrow X$ in $\scrM(X,\beta_2,z)$ 
which we wish to glue to punctured maps of the form
$f_1:(C_1,x_1,x_2,x_s)\rightarrow X$ in $\scrM_{\beta_1}(X,\beta_1)$, gluing
$x'_s$ to $x_s$. In order to glue, at the very least one needs
$f_1(x_s)=f_2(x_s')$, and this effectively imposes a point constraint
on the maps of class $\beta_1$. In particular, if $z'=f_2(x_s')$, 
we would expect we would need to look at punctured maps of class
$\beta_1$ with a constraint given by $z'$.
Unfortunately, there is no reason why the point $z'$ should lie in
$Z_s^{\circ}$, rather than just $Z_s$,
and as a result it does not necessarily impose a 
point constraint of the type we have previously considered.
So here we find a broader description of the numbers $N^{A_1}_{p_1p_2s}$
which allows more general point constraints. Effectively, this number
should be thought of as a kind of ``virtual log degree'' of the 
morphism $\ev_X:\scrM(X,\beta_1)\rightarrow\scrP(X,s)$. This
point of view extends the philosophy of Remark~\ref{rem:transverse}.
However, $\ev_X$
has virtual log fibre dimension $1$, so some care is needed in
defining this degree, with the definition given in 
Definition~\ref{def: alternateN} below.

In Theorem \ref{thm:NW equals N}, we show this gives a more flexible way
of defining the numbers $N^{A}_{pqr}$. This definition has been
used to obtain a deeper understanding of our construction in 
\cite{CanScat}, \cite{J22a}, \cite{J22b} and \cite{J24}.

\begin{lemma}
\label{lem: ev transversal}
Let $g:W\rightarrow \scrP(X,s)$ be a morphism from an fs log stack
which is transverse to $\ev_{\shX}:\foM^{\ev}(\shX,\beta_1)\rightarrow
\scrP(X,s)$ in the sense of Definition~\ref{def: transverse}. Then the
projection
\begin{equation}
\label{eq:proj to W}
W\times_{\scrP(X,s)}^{\fs}\foM^{\ev}(\shX,\beta_1)\rightarrow W
\end{equation}
is flat of fibre dimension $1$. 
\end{lemma}

\begin{proof}
As $\ev_{\shX}:\foM^{\ev}(\shX,\beta_1)\rightarrow\scrP(X,s)$ is log
smooth by Theorem~\ref{thm: everything is log smooth} (since $x_1$, $x_2$ are
marked), $\ev_{\shX}$ is also log flat. 
By transversality the projection \eqref{eq:proj to W}
is integral, so \eqref{eq:proj to W} is in fact flat by Proposition
\ref{prop: flatness sorites}, (2).

We compute the log fibre dimension of $\ev_{\shX}$.
By Lemma~\ref{lem:evaluation diagrams} and Proposition~\ref{prop: flat fibre dim}, (5), this log fibre dimension is the same as that of
$\ev_{\shX}:\foM(\shX,\beta_1)\rightarrow\scrP(\shX,s)$.
Recall $\mathbf{M}(\shX,s)$ from Definition~\ref{def:MXr}.
The morphism $\ev_{\shX}$ factors as
$\foM(\shX,\beta_1)\rightarrow \mathbf{M}(\shX,s)\rightarrow \scrP(X,s)$
as in the proof of Theorem~\ref{thm: everything is log smooth},
with the first morphism being log \'etale and the second morphism,
a base-change of $\mathbf{M}\rightarrow B\GG_m$, being of fibre and
log fibre dimension $1$.  Thus $\ev_{\shX}$ has log fibre dimension $1$ by
Proposition~\ref{prop: flat fibre dim}, (3), and the same holds for
\eqref{eq:proj to W}, by Proposition~\ref{prop: flat fibre dim}, (5).
Thus the result follows from Proposition~\ref{prop: flat fibre dim}, (2).
\end{proof}

\begin{definition}
\label{def: alternateN}
Suppose $g:W\rightarrow \scrP(X,s)$ is a 
morphism with $W$ of pure dimension $-1$. 
Suppose further that $g$ is transverse
to $\ev_{\shX}:\foM^{\ev}(\shX,\beta_1)\rightarrow \scrP(X,s)$.
Define a relative obstruction theory for 
\[
W\times^{\fs}_{\scrP(X,s)}\scrM(X,\beta_1)
\rightarrow
W\times^{\fs}_{\scrP(X,s)}\foM^{\ev}(\shX,\beta_1)
\]
by pull-back of that for $\scrM(X,\beta_1)\rightarrow\foM^{\ev}(\shX,\beta_1)$.
By Lemma~\ref{lem: ev transversal}, $W\times^{\fs}_{\scrP(X,s)}
\foM^{\ev}(\shX,\beta_1)$ is pure dimension zero, and hence it follows
by Riemann-Roch as in the proof of Proposition~\ref{prop:rel virt dim}
that $W\times^{\fs}_{\scrP(X,s)}\scrM(X,\beta_1)$
is virtual dimension $A\cdot c_1(\Theta_{X/k})$.
If $W\times^{\fs}_{\scrP(X,s)}\scrM(X,\beta_1)$ is proper over $\kk$,
we define
\[
N^{A_1,W}_{p_1p_2s}=\begin{cases}
\deg [W\times^{\fs}_{\scrP(X,s)}
\scrM(X,\beta_1)]^{\virt}& \virt.\dim W\times^{\fs}_{\scrP(X,s)}
\scrM(X,\beta_1)=0,\\
0 & \hbox{otherwise.}
\end{cases} 
\]
\end{definition}

\begin{theorem}
\label{thm:NW equals N}
Let $\ul{W}=B\GG_m$, with log structure $\shM_W$ so that $\overline\shM_W=Q$
for an fs monoid $Q$, and suppose there exists
$\ell \in Q^{\vee}$ such that for
$q\in Q$, the torsor contained in $\shM_W$ corresponding to $q$
is $\shL_q\cong \shU^{\otimes\langle \ell,q\rangle}$. Suppose further
that $\ell$ is a generator of a one-dimensional face of $Q^{\vee}$, and
that $g:W\rightarrow \scrP(X,s)$ is given so that (1) $g$ is transverse
to $\ev_{\shX}$; (2) $\Sigma(g)(\ell)=s$ and (3) $\ul{g}$ is a closed
embedding identifying $\ul{W}$ with $z\times B\GG_m$ for some
$z\in \ul{Z}$. Then
\[
N^{A_1,W}_{p_1p_2s}=N^{A_1}_{p_1p_2s}.
\]
\end{theorem}

\begin{proof}
We proceed in three steps. In the first step, we show that the
relevant moduli spaces are in fact proper, so that in particular
the number $N_{p_1p_2s}^{A_1,W}$ makes sense. Our goal is then to
deform the point constraint given by the morphism $g$ to a point constraint
of the type used to define $N_{p_1p_2s}^{A_1}$. This is done by finding
a discrete valuation ring $R$, a log structure $T$ on $\ul{T}=\Spec R$,
and a moprhism $[T/\GG_m]\rightarrow \scrP(X,s)$ whose restriction to
the closed point of $T$ gives the morphism $g$, and whose restriction to
the generic point is (close to) a morphism of the type constructed in
Proposition~\ref{prop: BGm morphism}. Here, assumption (2) of the statement
of the theorem is crucial. This allows us to deform
the moduli space defining $N^{A_1,W}_{p_1p_2s}$ to the moduli space
defining $N^{A_1}_{p_1p_2s}$. In Step 3, we then use standard
intersection theory techniques to show equality of these two numbers.

\medskip

{\bf Step 1. Properness.}

\medskip

By the transversality of the assumption (1) of the
statement of the theorem
and the fact that $\ul{W}=B\GG_m$, we see from 
Definition \ref{def: alternateN} that $N^{A_1,W}_{p_1p_2s}$
is defined provided that $W\times^{\fs}_{\scrP(X,s)}\scrM(X,\beta_1)$
is proper over $\kk$. For use later in the proof, we show a slightly
more general properness, assuming that $\ul{W}$ is a stack $\ul{T}\times
B\GG_m$ carrying some log structure and the morphism $\ul{W}\rightarrow
\ul{\scrP(X,s)}=\ul{Z}\times B\GG_m$ is of the form $\ul{g}'\times \id$
for some morphism $\ul{g}':\ul{T}\rightarrow\ul{Z}$.
We then show that $W\times^{\fs}_{\scrP(X,s)}\scrM(X,\beta_1)$ is
proper over $\ul{T}$.

As $W\times^{\fs}_{\scrP(X,s)}\scrM(X,\beta_1)$ is representable and
finite over 
\[
\ul{W}\times_{\ul{\scrP(X,s)}} \ul{\scrM(X,\beta_1)}
=\ul{T}\times_{\ul{Z}} \ul{\scrM(X,\beta_1)},
\]
it is enough to show that the latter is proper over $\ul{T}$.
But the evaluation map
$\ul{\scrM(X,\beta_1)}\rightarrow \ul{Z}$ is proper as the composition
with the proper map $\ul{Z}\rightarrow\ul{S}$ is proper by Lemma
\ref{lemma: properness}, (1). Thus $\ul{T}\times_{\ul{Z}}
\ul{\scrM(X,\beta_1)}$ is proper over $\ul{T}$. In particular, with
the choice of $W$ given in the statement of the theorem, we
see that $N^{A_1,W}_{p_1p_2s}$ is defined.

Note that $W\times^{\fs}_{\scrP(X,s)}\scrM(X,\beta_1)$ is of virtual
dimension zero precisely when $A_1\cdot c_1(\Theta_{X/\kk})=0$,
which is also the case when $\scrM(X,\beta_1,z)$ is virtual dimension zero.
Since both $N^{A_1,W}_{p_1p_2s}$ and $N^{A_1}_{p_1p_2s}$ are defined
to be zero if the virtual dimension of the relevant moduli space is
non-zero, we restrict now to the case that
$A_1\cdot c_1(\Theta_{X/\kk})=0$.

\medskip
{\bf Step 2. Construction of $T$.}
\medskip

To make the comparison, we will choose a discrete valuation ring $R$,
a log structure $T$ on $\ul{T}=\Spec R$, and a careful choice of
morphism $T\rightarrow Z=Z_s$. This allows us to deform the
point constraint determined by $g$ to a point constraint determined
by a morphism of the form given by Proposition~\ref{prop: BGm morphism}.
With the image of $g$ being $z\times B\GG_m$,
let $P=\overline{\shM}_{Z,\bar z}$ and
let $P_K$ be the stalk of $\overline{\shM}_Z$ at the generic
point of $Z$. These coincide with the stalks of the ghost sheaf of
$\scrP(X,s)$ at the corresponding points.
We have the generization map $\chi_P:P\rightarrow P_K$
and face $F=\chi_P^{-1}(0)$ of $P$. Note that
\[
F^{\vee}=(P^{\vee}+P_K^*)/P_K^*.
\]

As $\bar g^{\flat}:P\rightarrow Q$ is a local homomorphism between
sharp monoids, i.e., $(\bar g^{\flat})^{-1}(0)=\{0\}$,
$(\bar g^{\flat})^t(Q^{\vee})$ is not contained in a proper
face of $P^{\vee}$. Choose $u\in\Int(Q^{\vee})$ and
set
\begin{equation}
\label{eq:uprime def}
u'=(\bar g^{\flat})^t(u).
\end{equation}
Note $u'$ is not contained in a proper face of $P^{\vee}$, and thus
the image of $u'$ in $F^{\vee}$ under the projection map $P^{\vee}
\rightarrow F^{\vee}$ lies in the interior of $F^{\vee}$. We will denote
this element also by $u'$.

\'Etale locally at $\bar z$ we can write
$\ul{Z}$ as $\Spec\kk[F]\times \ul{\AA}_{\kk}^d$, by 
\cite[IV.3.3.1,3]{Ogus} and \cite[Tag~039P]{stacks}
for some $d\ge 0$,
with $\bar z$ corresponding to
the point $(z,0)$ where $z$ is the torus fixed-point
of $\Spec\kk[F]$. Write $\ul{\AA}_{\kk}^d=\Spec \kk[x_1,\ldots,x_d]$. 
Note that the log structure on $\Spec\kk[F] \times \ul{\AA}^d_{\kk}$
arises from the natural chart $P\rightarrow \kk[F] = \kk[P]/(P\setminus F)$,
$p\mapsto z^p$, and the factor $\ul{\AA}^d_{\kk}$ does not carry a log
structure.

We now choose $\ul{T}$ and a morphism $\ul{T}\rightarrow\ul{Z}$ as
follows. We give a morphism $\AA^1=\Spec\kk[t]
\rightarrow \Spec\kk[F]\times \ul{\AA}^d_{\kk}$ given at the ring level
by $z^f\mapsto t^{u'(f)}$, $x_i\mapsto 0$. After passing to an 
\'etale neighbourhood $\ul{C}\rightarrow \AA^1$ of the origin in $\AA^1$,
with $x\in \ul{C}$ mapping to the origin in $\AA^1$, we may then
assume that there is a morphism 
\begin{equation}
\label{eq:C to Z}
\ul{C} \rightarrow \ul{Z}
\end{equation}
with $x$ mapping to $\bar z$ and such that the induced homomorphism
on strict Henselian local rings $\O_{\ul{Z},\bar z}\rightarrow
\O_{\ul{C},\bar x}$ concides with the corresponding homomorphism
of strict Henselian rings given by the constructed morphism
$\AA^1\rightarrow \Spec\kk[F]\times\AA^d_k$.
Set $R=\O_{\ul{C},x}$,
a discrete valuation ring.

We now have a morphism $\ul{T}=\Spec R\rightarrow \ul{Z}$.
If $T'$ denotes the pull-back log structure on $\ul{T}$ from
$X$, we note that this log structure is determined by the data
$P, F$ and $u'$ of Lemma~\ref{lem:log dash}. In what
follows, $\xi$ denotes the generic point of $\ul{T}$, and we write
$\xi^{\dagger}$ for the standard log point structure on $\xi$,
which is, by construction, the restriction of the log structure on
$T$ to $\xi$.

On the other hand, let $\bar u\in 
\Int((\ell^{\perp}\cap Q)^{\vee})$ denote the restriction of
$u\in\Int(Q^{\vee})$ to the face $\ell^{\perp}\cap Q$ of $Q$.
Consider a log structure $T$ on $\ul{T}$ given
from Lemma~\ref{lem:log dash} by
the data of the monoid $Q$, its face $\ell^{\perp}\cap Q$, and
the element $\bar u$.
Note that $Q_K=\NN$ as
$\ell$ generates a one-dimensional face of $Q^{\vee}$ by the
assumption stated in the theorem.

We can now define a morphism $T\rightarrow T'$, hence a morphism $T\rightarrow
Z$, using Lemma~\ref{lem:log dash morphism}. Such a morphism is
determined, in this lemma, by maps $\varphi,\varphi_K$. We take
$\varphi=\bar g^{\flat}:P\rightarrow Q$ and
$\varphi_K=s:P_K\rightarrow Q_K=\NN$. We  verify the required
properties of these maps as stated in the lemma. First, we need
to check the commutativity of the lemma's diagram of monoids.
Writing $\chi_P,\chi_Q$ for the two generization maps, we need to check
that $\chi_Q\circ\varphi=\varphi_K\circ \chi_P$. Instead we prove equality
of the transpose maps, noting that $\chi_Q^{t}$ and $\chi_P^{t}$ are
just inclusions of faces. Thus we need only to show that
$\varphi^t|_{Q_K^{\vee}}=\varphi_K^t$. However, $\ell$ is the generator
of $Q_K^{\vee}$, again by the stated assumptions of the theorem,
and $\varphi^t(\ell)=
(\bar g^{\flat})^t(\ell)=\Sigma(g)(\ell)=s=\varphi_K^t(\ell)$ by the
assumption (2) of the theorem. Hence we have the desired equality.
Second, we need to check that $\bar u(\varphi(p))=u'(p)$ for
$p\in F$. However, this is immediate from \eqref{eq:uprime def}.
Thus Lemma~\ref{lem:log dash morphism} now gives us the desired morphism
$T\rightarrow T'$. 

\medskip
{\bf Step 3. The comparison.}
\medskip

We can now act on $T$ and $Z$ with $\GG_m$. The action is taken to be
trivial on the underlying schemes, but acts on the log structure
by acting on a torsor $\shL_{q}$ over $T$ with weight $\ell(q)$ for
$q\in Q$
and acting on a torsor $\shL_p$ over $U\subseteq Z$ for
$p\in\Gamma(U,\overline{\shM}_Z)$ with weight $\langle s,p\rangle$.
By Remark~\ref{remark: quotient description},
with this action we have $[Z/\GG_m]\cong \scrP(X,s)$.

Because of compatibility of the weights of the $\GG_m$ action,
the morphism $T\rightarrow Z$ descends to a morphism $[T/\GG_m]
\rightarrow \scrP(X,s)$ which is in fact transverse to 
$\ev_{\shX}:\foM^{\ev}(\shX,\beta_1)\rightarrow \scrP(X,s)$.
To check this, we can use Theorem~\ref{theorem: transversality},
and it is then sufficient to check transversality of the 
morphisms obtained by restricting $[T/\GG_m]\rightarrow
\scrP(X,s)$ to each of the two points of $[T/\GG_m]$.
However, there is a closed embedding $W\hookrightarrow [T/\GG_m]$ identifying
$W$ with the closed point of $[T/\GG_m]$, and the composition 
$W\rightarrow \scrP(X,s)$ agrees at the level of ghost sheaves
with the morphism $g$, which we have assumed is transverse.
Transversality is obvious at the
generic point $[\xi^{\dagger}/\GG_m]$ of 
$[T/\GG_m]$, as the monoid at that generic point is
$\NN$. 

Let us first complete the proof under the assumption that
the composition $W\hookrightarrow [T/\GG_m]\rightarrow
\scrP(X,s)$ agrees precisely with $g$.
We then have a diagram
\[
\xymatrix@C=30pt
{
[\xi^{\dagger}/\GG_m] \times^{\fs}_{\scrP(X,s)} \scrM(X,\beta_1)
\ar[r]^j\ar[d]_{p_{\xi}}&
[T/\GG_m]\ar[d]_p \times^{\fs}_{\scrP(X,s)} \scrM(X,\beta_1)
&
W \times^{\fs}_{\scrP(X,s)}\scrM(X,\beta_1)
\ar[d]^{p_0}\ar[l]_i\\
[\xi^{\dagger}/\GG_m] \times^{\fs}_{\scrP(X,s)} \foM^{\ev}(\shX,\beta_1)
\ar[r]^j
\ar[d]_{q_{\xi}}&
[T/\GG_m] \times^{\fs}_{\scrP(X,s)} \foM^{\ev}(\shX,\beta_1)
\ar[d]_q&
W \times^{\fs}_{\scrP(X,s)} \foM^{\ev}(\shX,\beta_1)
\ar[d]^{q_0}\ar[l]_i\\
\xi\ar[r]_j&\ul{T}&\Spec\kk\ar[l]^i
}
\]
Here the morphism $q$ is the composition of the projection to
$[T/\GG_m]$ whose underlying stack is $\ul{T}\times B\GG_m$,
followed by the further projection to $\ul{T}$. The morphisms
$q_0$, $q_{\xi}$ are defined similarly. All squares are cartesian in
all categories, with all horizontal arrows strict. The morphisms
$q_{\xi}$, $q$, $q_0$
are flat of fibre dimension zero by Lemma~\ref{lem: ev transversal}.

The relative obstruction theory for $p$ given by Definition
\ref{def: alternateN} pulls back via the closed
embedding $i$ and the open embedding $j$ to the corresponding
relative obstruction theories for $p_{\xi}$ and $p_0$. These give
virtual fundamental classes on the three moduli spaces of the top
row, which we write as $\alpha_{\xi}$, $\alpha$ and $\alpha_0$ respectively.
By properties of virtual pull-back \cite{Man},~Thm.~4.1,
it follows that $\alpha_\xi=j^*\alpha$, $\alpha_0=i^!\alpha$. 

From this
it follows by a standard argument that $\deg\alpha_{\xi}=\deg \alpha_0$.
Explicitly, note that by Step 1,
$q\circ p$, $q_{\xi}\circ
p_{\xi}$ and $q_0\circ p_0$ are all proper. We can define
$\deg \alpha\in \QQ$ by the formula
$(q\circ p)_*\alpha = (\deg\alpha) [\ul{T}]$, as $\alpha$ is a one-dimensional
class.\footnote{Here
we work in the rational Chow groups of these Deligne-Mumford stacks
so that we can work with proper, rather than projective, push-forward.}
Further, we have 
\[
(q_0\circ p_0)_*\alpha_0= (\deg\alpha_0)[\Spec\kk],\quad 
(q_{\xi}\circ p_{\xi})_*\alpha_{\xi}=(\deg\alpha_{\xi})
[\xi].
\]
Then $\deg\alpha=\deg\alpha_{\xi}$ because
\[
(\deg \alpha)[\xi]=j^*(q\circ p)_*\alpha
=(q_{\xi}\circ p_{\xi})_* j^* \alpha = (\deg \alpha_{\xi})[\xi]
\]
by compatibility of flat pull-back and proper push-forward.
Also, $\deg\alpha=\deg\alpha_0$, as
\[
(\deg\alpha)[\Spec\kk]=i^!(q\circ p)_*\alpha=(q_0\circ p_0)_* i^!\alpha=(\deg\alpha_0)
[\Spec\kk],
\]
by compatibility of the Gysin map with proper push-forward,
see \cite{Kresch}, Thm.~2.1.12,(xi).
Thus we conclude that $N^{A_1,W}_{p_1p_2s}=
N^{A_1,[\xi^{\dagger}/\GG_m]}_{p_1p_2s}$ under the assumption
that the composition
$W\hookrightarrow [T/\GG_m]\rightarrow\scrP(X,s)$ agrees with $g$.

If the composition $W\hookrightarrow [T/\GG_m]\rightarrow\scrP(X,s)$
does not agree with $g$, we note that in fact the number 
$N^{A_1,W}_{p_1p_2r}$ only depends on the morphism $g$ at the
level of ghost sheaves. Indeed, in the above argument, take instead 
$\ul{T}=\Spec \kk[P^{\gp}]$,
with constant log structure given by $T=\Spec(Q\rightarrow\kk)\times \ul{T}$.
Then $\Gamma(T,\O_T^{\times})=P^{\gp}\times\kk^{\times}$, and there
is a universal morphism $T\rightarrow \Spec(P\rightarrow\kk)$
given at the logarithmic level by $P\rightarrow \Gamma(T,\shM_T)=Q\times 
\Gamma(T,\O_T^{\times})$
defined by $p\mapsto (\bar g^{\flat}(p),(p,1))$. This then gives rise 
to a morphism
$[T/\GG_m]\rightarrow \scrP(X,s)$ as before. Any morphism $g':W\rightarrow
\scrP(X,s)$ agreeing with $g$ at the level of underlying stacks and
at the level of ghost sheaves then factors through some strict
embedding $W\hookrightarrow [T/\GG_m]$.  The same argument as in the
previous paragraphs then
shows that $N^{A_1,W}_{p_1p_2s}$ defined using
the composition $W\hookrightarrow [T/\GG_m]\rightarrow \scrP(X,s)$
is independent of the choice of strict embedding $W\hookrightarrow [T/\GG_m]$.

A similar but simpler argument then allows comparison of
$N^{A_1,[\xi^{\dagger}/\GG_m]}_{p_1p_2s}$ with 
$N^{A_1,B\GG_m^{\dagger}}_{p_1p_2s}$
for a morphism $B\GG_m^{\dagger}\rightarrow  \scrP(X,s)$
given by Proposition~\ref{prop: BGm morphism}. 
Indeed, using the morphism $\ul{C}\rightarrow \ul{Z}$
given in \eqref{eq:C to Z}, we may take a closed point $x' \in \ul{C}$
such that its image in $\ul{Z}$ lies in $\ul{Z}^{\circ}$. Indeed, such
a point exists by the construction of this morphism. We now take
$R=\O_{\ul{C},x'}$,
set $\ul{T}=\Spec R'$ and equip $\ul{T}$
with two log structures. The first log structure, $T$,
is constructed from Lemma~\ref{lem:log dash} from the data $Q=\NN$, $F=0$, 
$u=0$. The second log structure $T'$ is pulled back from 
$Z$, and arises from the data $Q=P_K$, $F=0$, and $u=0$. Thus by
Lemma~\ref{lem:log dash morphism}, there is a morphism $T\rightarrow T'$
induced by $s:P_K\rightarrow \NN$. We may then repeat the argument above
to show that $N^{A_1,[\xi^{\dagger}/\GG_m]}_{p_1p_2s}= 
N^{A_1}_{p_1p_2s}$, which completes the proof.
\end{proof}

\subsection{The proof of Theorem~\ref{thm: main gluing theorem}}
\label{sec:proof of main gluing}

Fix $p_1,p_2,p_3,r\in B(\ZZ)$ and $s\in \Sigma(X)(\ZZ)$.
We further fix punctured map classes $\beta_1,\beta_2$, each of genus zero with 
three marked points and curve classes $A_1,A_2$. 
We label the punctured points of $\beta_1$ as
$x_1,x_2$ and $x_s$ with contact orders $p_1,p_2,-s$, and
label the punctured points of $\beta_2$ as $x_3,x_s'$ and $x_{\out}$
with contact orders $p_3, s$ and $-r$. We also fix $z\in Z^{\circ}_r$.
In this section, when we write $\foM^{\ev}(\shX,\beta_1)$, 
$\foM^{\ev}(\shX,\beta_2,z)$, we mean the moduli spaces defined
in Definitions~\ref{def:Mev} and \ref{def:Mbetaev}; otherwise, we will
be careful in decorating the $\ev$ superscripts.

This data gives rise to moduli spaces
\[
\scrM_1:= \scrM(X,\beta_1), \quad
\scrM_2:= \scrM(X,\beta_2,z), \quad 
\scrM^{\gl}:=\scrM^{\gl}(X,\beta_1,\beta_2,z)
\]
and
\[
\foM_1:=  \foM^{\ev(x_s)}(\shX,\beta_1),\quad
\foM_2:=  \foM^{\ev(x_s',x_{\out})}(\shX,\beta_2,z),
\quad
\foM^{\gl}:=\foM^{\gl,\ev(q,x_{\out})}(\shX,\beta_1,\beta_2,z).
\]
Note we have a diagram cartesian in the fs log category
\begin{equation}
\label{eq:brexit sucks}
\xymatrix@C=30pt
{
\foM^{\gl}\ar[d]_{\Pi}\ar[r] & \foM_1\ar[d]^{\ev_{\shX}}\\
\foM_2\times B\GG_m^{\dagger}\ar[r] & \scrP(X,s)
}
\end{equation}
by Corollary~\ref{cor:new gluing}. 

Write $\delta:\foM_2\rightarrow B\GG_m^{\dagger}$
for the composition
\begin{equation}
\label{eq:delta def}
\foM_2\rightarrow \foM^{\ev}(\shX,\beta_2,z)=\foM^{\ev}(\shX,\beta_2)
\times^{\fs}_{\scrP(X,r)} B\GG_m^{\dagger}\rightarrow B\GG_m^{\dagger},
\end{equation}
where the first morphism forgets evaluation at $x_s'$ and
is strict smooth and the last morphism is the second projection.
We then have a composed morphism
\begin{equation}
\label{eq:Phidef}
\xymatrix@C=30pt
{
\Psi:\foM^{\gl}\ar[r]^>>>>{\Pi} & \foM_2\times B\GG_m^{\dagger}\ar[r]^{\delta\times\id}&
B\GG_m^{\dagger}
\times B\GG_m^{\dagger}.
}
\end{equation}

Fix a non-negative integer $\mu$ and recall the definition
of $B\GG_m^{2,\mu}$ from \S\ref{sec:second step}
and the morphism
\[
B\GG_m^{2,\mu} \rightarrow B\GG_m^{\dagger}\times B\GG_m^{\dagger}.
\]
We use $\delta$, $\ell_q$ for the generators of the ghost sheaves
of the first and second factors of $B\GG_m^{\dagger}\times B\GG_m^{\dagger}$.
We may now define 
\begin{align}
\label{def:Mmu moduli}
\begin{split}
\foM^{\mu,\gl}= \foM^{\mu,\gl,\ev(q,x_{\out})}(\shX,\beta_1,\beta_2,z)
:= {} & \foM^{\gl}\times^{\fs}_{B\GG_m^{\dagger}\times
B\GG_m^{\dagger}} B\GG_m^{2,\mu}\\
\foM_2^{\mu}=\foM^{\mu,\ev(x_s',x_{\out})}(\shX,\beta_2,z):= {} & (\foM_2\times B\GG_m^{\dagger})
\times^{\fs}_{B\GG_m^{\dagger}\times
B\GG_m^{\dagger}} B\GG_m^{2,\mu}\\
\scrM^{\mu,\gl}=\scrM^{\mu,\gl}(X,\beta_1,\beta_2,z):= {} &
\scrM^{\gl}(X,\beta_1,\beta_2,z)\times_{\foM^{\gl}} \foM^{\mu,\gl}.
\end{split}
\end{align}
In the remainder of this section, we use the short-hand notations,
spelling out these moduli spaces more fully in subsequent sections.
We also write
\begin{equation}
\label{eq:def Pimu}
\Pi^{\mu}:\foM^{\mu,\gl}\rightarrow \foM_2^{\mu}
\end{equation}
for the base-change of $\Pi$ from \eqref{eq:brexit sucks}.

At times, we will also need to remove the evaluation at the
node $q$, so we also define
\begin{align}
\label{eq:Mmu no ev}
\begin{split}
\foM^{\mu,\gl,\ev}(\shX,\beta_1,\beta_2,z)
:=  {} &
\foM^{\gl,\ev(x_{\out})}(\shX,\beta_1,\beta_2,z)\times^{\fs}_{B\GG_m^{\dagger}\times
B\GG_m^{\dagger}} B\GG_m^{2,\mu}\\
\foM^{\mu,\ev}(\shX,\beta_2,z)
:= {} & (\foM^{\ev(x_{\out})}(\shX,\beta_2,z)\times B\GG_m^{\dagger})
\times^{\fs}_{B\GG_m^{\dagger}\times B\GG_m^{\dagger}} B\GG_m^{2,\mu}.
\end{split}
\end{align}
For the structure of $\foM_2^{\mu}$, we note:

\begin{lemma}
\label{lem: underlying}
Given an fs log stack $W$ equipped with a morphism $\delta:W\rightarrow
B\GG_m^{\dagger}$, the projection
\[
W^{\mu}:=
(W\times B\GG_m^{\dagger})\times^{\fs}_{B\GG_m^{\dagger}\times B\GG_m^{\dagger}}
B\GG_m^{2,\mu}\rightarrow W\times B\GG_m^{\dagger}
\]
is an isomorphism on underlying stacks.
\end{lemma}

\begin{proof}
First observe the morphism $\delta:W\rightarrow B\GG_m^{\dagger}$ is
integral, as the ghost sheaf monoid on the target is $\NN \delta$.
As usual, we use the same notation $\delta$ for the morphism,
the generator of the ghost sheaf on the target, and the pull-back
of this generator to a section of $\overline\shM_W$. So
$W\times B\GG_m^{\dagger}\rightarrow B\GG_m^{\dagger}\times
B\GG_m^{\dagger}$ is an integral morphism. Therefore if $W^{\mu}$
were defined using the fibre product in the category of fine log structures,
the statement would be true. Thus it is sufficient to check
that this fibre product is already saturated. The
construction of the fine fibre product (see \cite{Ogus},~III,\S2.1)
implies the following. In a smooth neighbourhood
of a geometric point $\bar x$ of $W$, there is a chart for
the log structure of the fine fibre product with monoid the fine push-out
\[
(\overline{\shM}_{W,\bar x}\oplus\NN\ell_q)\oplus^{\fine}_{\NN\delta
\oplus\NN\ell_q}
(\NN\delta+\NN(\ell_q-\mu\delta)).
\]
By Lemma~\ref{prop:fibre product properties}, this 
is the submonoid of $(\overline{\shM}_{W,\bar x}\oplus\NN\ell_q)^{\gp}$
generated by $\overline{\shM}_{W,\bar x}\oplus\NN\ell_q$ and
the image of $\ell_q-\mu\delta$ under the map $\ZZ\delta\oplus
\ZZ\ell_q\rightarrow
(\overline{\shM}_{W,\bar x}\oplus\NN\ell_q)^{\gp}$. This is precisely
the submonoid generated by $(-\mu\delta, \ell_q)
\in \overline{\shM}^{\gp}_{W,\bar x}\oplus\ZZ \ell_q$ and $\overline{\shM}_{W,\bar
x}\oplus 0$. This monoid is clearly saturated.
\end{proof}

We turn to the structure of $\foM^{\mu,\gl}$.

\begin{proposition}
\label{prop:transverse M2}
\begin{enumerate}
\item 
If $B\GG_m^{2,\mu}\rightarrow B\GG_m^\dagger\times B\GG_m^{\dagger}$
is transverse to $\Psi$ defined in \eqref{eq:Phidef}, then $\foM^{\mu,\gl}$ is
pure-dimensional of dimension $\dim X$. Further, the stalk of the ghost
sheaf of $\foM^{\mu,\gl}$ at the generic point of any irreducible component
of $\foM^{\mu,\gl}$ is rank two.
\item If $B\GG_m^{2,\mu}\rightarrow B\GG_m^\dagger\times B\GG_m^{\dagger}$
is transverse to $\Psi$, and
$S_2^{\mu}\subset \foM_2^{\mu}$ is the union of open strata of
$\foM_2^{\mu}$, then 
the composition
\[
S_2^{\mu}\hookrightarrow \foM_2^{\mu}\stackrel{\pr_1}{\longrightarrow}
\foM_2\times
B\GG_m^{\dagger}
\]
is transverse to $\Pi$.
\item For $\mu$ sufficiently large, 
$B\GG_m^{2,\mu}\rightarrow B\GG_m^\dagger\times B\GG_m^{\dagger}$
is transverse to $\Psi$.
\end{enumerate}
\end{proposition}

\begin{proof}
For (1), note
that as $\ev_{\shX}:\foM_1\rightarrow\scrP(X,s)$ is 
log smooth by Theorem~\ref{thm: everything is log smooth}, 
so is $\Pi$ by base-change. Further, the morphism $\delta$ of
\eqref{eq:delta def} is log smooth
by Lemma~\ref{lem: more log smoothness}. Thus $\Psi$
is log smooth, and we will first calcuate log fibre dimensions.
The log fibre dimension of $\Pi$ is the same as that of
$\ev_{\shX}$ by Proposition~\ref{prop: flat fibre dim}, (5).
The log fibre dimension of $\ev_{\shX}$ was
already calculated in the proof of Lemma~\ref{lem: ev transversal} to be $1$.
Thus $\Pi$ has log fibre dimension $1$.
On the other hand, in the definition of $\delta$,
\eqref{eq:delta def}, the first morphism is strict and smooth
of relative dimension $\dim X$ and the second is of log fibre dimension $1$,
by Proposition~\ref{prop:the invariants}, (1). 
Putting this all together, we see that $\Psi$
is log smooth of log fibre dimension $1+\dim X + 1=\dim X+2$. The same is then
true of the base-change $\foM^{\mu,\gl}\rightarrow B\GG_m^{2,\mu}$.

Now let $\bar\xi$ be a generic point of an irreducible component 
$\foM'$ of $\foM^{\mu,\gl}$, and set $Q:=\overline\shM_{\foM^{\mu,\gl},
\bar\xi }$.
We now apply the definition \eqref{eq:log fibre dim def} of log fibre
dimension for $\Psi$, taking $x=\bar\xi$ in that formula. However,
\eqref{eq:log fibre dim def} 
needs to be appropriately
interpreted for stacks, which one can do by passing to smooth charts
and keeping track of the change of dimension.
Since $\bar\xi$ is a generic point of an irreducible component, 
the stalk of the structure sheaf
appearing in \eqref{eq:log fibre dim def} should be interpreted as
zero-dimensional, whereas the term $\mathrm{tr.deg.}\,\kappa(x)/\kappa(y)$
should be interpreted in this case as $\dim \foM'-\dim B\GG_m^{2,\mu}
=\dim \foM'+2$. With this interpretation, \eqref{eq:log fibre dim def}
then gives
\[
\dim X+ 2= \dim \foM'+2+\rank Q-2=\dim\foM'+\rank Q,
\]
with the two terms in the middle expression arising from the monoids.
So it suffices to show that the assumption that $B\GG_m^{2,\mu}
\rightarrow B\GG_m^{\dagger}\times B\GG_m^{\dagger}$ is
transverse to $\Psi$ implies $\rank Q=2$.

For any choice of $\mu$, if $\rank Q>2$, consider the induced map
\begin{align*}
\sigma_{\bar\xi}=\Hom(Q,\RR_{\ge 0})\rightarrow &
\tau_{\mu}:=\Sigma(B\GG_m^{2,\mu})
=\RR_{\ge 0} (\delta^*+\mu\ell_q^*)+\RR_{\ge 0}\ell_q^*\\
\subseteq&\tau:=
\Sigma(B\GG_m^{\dagger}\times B\GG_m^{\dagger})=\RR_{\ge 0}^2.
\end{align*}
The image of $\sigma_{\bar\xi}$ must intersect the interior of $\tau_{\mu}$,
as the induced morphism of stalks of ghost sheaves is local. Since
$\rank Q>2$, it 
follows that there is a proper face $\sigma'$ of $\sigma_{\bar\xi}$ whose image
intersects the interior of $\tau_{\mu}$. By Proposition 
\ref{prop: realisability} applied to the map
$\foM^{\mu,\gl}\rightarrow B\GG_m^{2,\mu}$, it then follows that 
there exists a point $\bar\xi'\in |\foM^{\mu,\gl}|$ with $\sigma_{\bar\xi'}=\sigma'$ and
with $\bar\xi'$ specializing to $\bar\xi$. However, since $\bar\xi$
is the generic point of an irreducible component, there is no such
point $\bar\xi'$.

If in fact $B\GG_m^{2,\mu}\rightarrow B\GG_m^{\dagger}\times B\GG_m^{\dagger}$
is transverse to $\Psi$, then by Proposition
\ref{prop:integral morphism}, $\sigma_{\bar\xi}$ must surject onto a face
of $\tau_{\mu}$. Again, by the local property, the image of
$\sigma_{\bar\xi}$ intersects the interior of $\tau_{\mu}$, and hence the image
is $\tau_{\mu}$. This shows that $\rank Q\ge 2$. Thus $\rank Q=2$
and we obtain the desired dimension for $\foM^{\mu,\gl}$.

(2) Take $S_2$ to be the union of open strata of $\foM_2$. As
$\delta$ is log smooth, the ghost sheaf on $S_2$ has stalk
$\NN$, and every two-dimensional cone of $\Sigma(S_2\times B\GG_m^{\dagger})$ 
surjects onto $\tau$ under the map $\Sigma(\delta\times\id)$. This
surjection is in 
particular an isomorphism of rational polyhedral cones (not necessarily
preserving the integral structure). If $S_2^{\mu}\subseteq \foM_2^{\mu}$
is the corresponding open subset of $\foM_2^{\mu}$, then it follows that
$\Sigma(S_2^{\mu})\rightarrow \Sigma(B\GG_m^{2,\mu})$ maps each two-dimensional
cone isomorphically to $\tau_{\mu}$, by the tropical description of the
fibre product, Proposition~\ref{tropicalproduct}. If 
$B\GG_m^{2,\mu}\rightarrow B\GG_m^\dagger\times B\GG_m^{\dagger}$
is transverse to $\Psi$, then 
every cone of $\Sigma(\foM^{\mu,\gl})$ surjects onto a face of $\tau_{\mu}$,
and via the factorization 
$\Sigma((\Pi^{\mu})^{-1}(S_2^{\mu}))\rightarrow
\Sigma(S_2^{\mu})\rightarrow\Sigma(B\GG_m^{2,\mu})$ with
$\Pi^{\mu}$ given in \eqref{eq:def Pimu}, we see that
every cone of $\Sigma((\Pi^{\mu})^{-1}(S_2^{\mu}))$ surjects onto a cone
of $\Sigma(S_2^{\mu})$, hence the desired transversality.

(3) To obtain transversality, we choose $\mu$ as follows.
Consider the map
\[
\Sigma\left(\Psi\right):
\Sigma(\foM^{\gl})\rightarrow \Sigma(B\GG_m^{\dagger}\times B\GG_m^{\dagger})=\tau.
\]
By \S\ref{sec:stacks remarks}, we can assume $\Sigma(\foM^{\gl})$ contains
only a finite number of one-dimensional cones. We can thus choose
$\mu$ so that for any one-dimensional cone $\rho\in\Sigma(\foM^{\gl})$, 
the intersection of $\Sigma(\Psi)(\rho)$ with $\tau_{\mu}$
is either $\{0\}$ or the boundary ray $\RR_{\ge 0}\ell_q^*$. Transversality
then follows from Theorem~\ref{theorem: transversality}.
\end{proof}

\begin{proof}[The proof of Theorem~\ref{thm: main gluing theorem}.]

{\bf Step 1. The relative obstruction theories and calculation of
the virtual dimension.}

\medskip

Note the diagram given in the theorem is obviously cartesian in all
categories. However, $\foM^{\mu,\gl,\ev}(\shX,\beta_1,\beta_2,z)$
as defined in \eqref{eq:Mmu no ev} does not involve evaluation at $q$. But we 
can replace the diagram of the theorem with the diagram
\[
\xymatrix@C=30pt
{
\scrM^{\mu,\gl}(X,\beta_1,\beta_2,z)\ar[r]\ar[d] &
\scrM^{\gl}(X,\beta_1,\beta_2,z)\ar[d]\\
\foM^{\mu,\gl,\ev(q,x_{\out})}(\shX,\beta_1,\beta_2,z)\ar[r]&  
\foM^{\gl,\ev(q,x_{\out})}(\shX,\beta_1,\beta_2,z)
}
\]
without any harm to the statements. Note in the notational convention
in this subsection, this is
\[
\xymatrix@C=15pt
{
\scrM^{\mu,\gl}\ar[r]\ar[d]&\scrM^{\gl}\ar[d]\\
\foM^{\mu,\gl}\ar[r]&\foM^{\gl}
}
\]

First we spell out carefully the spaces involved in the relative obstruction
theory for $\scrM^{\mu,\gl}\rightarrow\foM^{\mu,\gl}$.
We have a diagram
\begin{equation}
\label{eq:Mgl obstruction diagram}
\xymatrix@C=30pt
{
\scrM^{\gl}(X,\beta_1,\beta_2,z)\ar[r]\ar[d] & \scrM(X,\beta,z)\ar[d]\\
\foM^{\gl,\ev(x_{\out})}(\shX,\beta_1,\beta_2,z)\ar[r]&\foM^{\ev(x_{\out})}
(\shX,\beta,z)
}
\end{equation}
which is not cartesian but for which the upper left corner is an
open and closed substack of the fibre product consisting of
those punctured maps which split into maps of degree classes $A_1,A_2$. 
The relative
obstruction theory on the right pulls back to one on the left, of
relative dimension $\chi((f^*\Theta_{X/\kk})(-x_{\out}))$. Now we
have a factorization
\[
\scrM^{\gl}(X,\beta_1,\beta_2,z)\rightarrow \foM^{\gl,\ev(q,x_{\out})}
(\shX,\beta_1,\beta_2,z)
\rightarrow \foM^{\gl,\ev(x_{\out})}(\shX,\beta_1,\beta_2,z),
\]
with a compatible obstruction theory for the first morphism of
relative virtual dimension
$\chi((f^*\Theta_{X/\kk})(-x_{\out}-q))$, where $q$ is the glued node.
This relative obstruction theory then pulls back to give the desired
relative obstruction theory for $\scrM^{\mu,\gl}\rightarrow \foM^{\mu,\gl}$. 
In particular, the relative virtual dimension is still
\[
\chi((f^*\Theta_{X/\kk})(-x_{\out}-q))=A\cdot c_1(\Theta_{X/\kk})
-\dim X,
\]
see \eqref{eq:rel virtual dim}.

Now choose $\mu$ as in Proposition~\ref{prop:transverse M2}, (3). 
Then by (1) of that proposition, $\foM^{\mu,\gl}$ is pure-dimensional
of dimension $\dim X$, and hence the virtual dimension of $\scrM^{\mu,\gl}$
is $A\cdot c_1(\Theta_{X/\kk})$, as claimed.

\medskip

\pagebreak

{\bf Step 2. The setup for gluing.}
\medskip

We have a diagram
\begin{equation}
\label{eq:big diagram}
\xymatrix@C=30pt
{
\scrM^{\mu,\gl}\ar[r]\ar[d]\ar@{}[dr] | I&\scrM'_1\ar[d]
\ar[rr]\ar@{}[drr] | {II}&&\scrM_1\ar[d]\\
\scrM_2\times_{\foM_2}\foM^{\mu,\gl}\ar[d]
\ar[r]\ar@{}[dr] | {III} &\foM^{\mu,\gl}\ar[d]\ar[r]
\ar@{}[dr] | {IV} &\foM^{\gl}\ar[d]\ar[r]\ar@{}[dr] | {V} & \foM_1\ar[d]\\
\scrM_2\times_{\foM_2}\foM_2^{\mu}\ar[r]\ar[d]
\ar@{}[dr] | {VI}&\foM_2^{\mu}\ar[r]\ar[d]&\foM_2\times B\GG_m^{\dagger}
\ar[r] & \scrP(X,s)\\
\scrM_2\ar[r]&\foM_2
}
\end{equation}
Here squares $IV$ and $V$ are cartesian in the fs log category.
On the other hand, 
$\scrM'_1$ is defined so that the rectangle $II$ is
fs log cartesian, but as the vertical arrows are strict, it is
cartesian in all categories. Similarly, squares $III$ and $VI$
are cartesian in all categories, with horizontal arrows being strict.
Moreover, the large rectangle with vertices $\scrM^{\mu,\gl}$, $\scrM_2
\times_{\foM_2} \foM_2^{\mu}$, $\scrM_1$ and $\scrP(X,s)$
is fs log cartesian. Indeed, using the second diagram of 
Theorem~\ref{thm:new gluing} and the definitions of the various stacks, we have
\begin{align*}
(\scrM_2\times_{\foM_2}\foM^{\mu}_2)\times^{\fs}_{\scrP(X,s)} \scrM_1
\cong {} & 
\big(\scrM_2\times_{\foM_2}((\foM_2\times B\GG_m^{\dagger})
\times^{\fs}_{B\GG_m^{\dagger}
\times B\GG_m^{\dagger}}B\GG_m^{2,\mu})\big)\times^{\fs}_{\scrP(X,s)} \scrM_1\\
\cong {} & ((\scrM_2\times B\GG_m^{\dagger})\times^{\fs}_{B\GG_m^{\dagger}
\times B\GG_m^{\dagger}} B\GG_m^{2,\mu})\times^{\fs}_{\scrP(X,s)}\scrM_1\\
\cong {} & \big((\scrM_2\times B\GG_m^{\dagger})\times^{\fs}_{\scrP(X,s)}
\scrM_1\big)\times^{\fs}_{B\GG_m^{\dagger}\times B\GG_m^{\dagger}}
B\GG_m^{2,\mu}
\\
\cong {} & \scrM^{\gl}\times^{\fs}_{B\GG_m^{\dagger}\times B\GG_m^{\dagger}}
B\GG_m^{2,\mu}\\
\cong {} & \scrM^{\mu,\gl}.
\end{align*}
Thus it follows
that square $I$ is also cartesian in the fs log category, and in this
square all arrows are strict, hence it is cartesian in all
categories.

We now extract the sub-diagram 
\begin{equation}
\label{eq:sub-diagram}
\xymatrix@C=30pt
{
\scrM^{\mu,\gl}\ar[r]^{\alpha''}\ar[d]_{p''}&\scrM'_1\ar[d]^{p'}
\ar[r]^{\beta''}&\scrM_1\ar[d]^p\\
\scrM_2\times_{\foM_2}\foM^{\mu,\gl}\ar[r]_>>>>>>>{\alpha'}
\ar[d]_{q''}&\foM^{\mu,\gl}\ar[r]_{\beta'}\ar[d]^{q'}&\foM_1\\
\scrM_2\ar[r]_{\alpha}&\foM_2
}
\end{equation}
Here all squares are cartesian in all categories. This diagram also
induces a cartesian diagram in all categories
\[
\xymatrix@C=50pt
{
\scrM^{\mu,\gl}\ar[d]_{\alpha'\circ p''=p'\circ\alpha''}
\ar[rr]^{(\beta''\circ\alpha'')\times (q''\circ p'')}&&
\scrM_1\times \scrM_2\ar[d]^{p\times\alpha}\\
\foM^{\mu,\gl}\ar[rr]_{\beta'\times q'} &&\foM_1\times\foM_2
}
\]
with the obstruction theory for $\alpha'\circ p''$ the pull-back of
that for $p\times\alpha$ by Corollary~\ref{cor:gluing recall}.
\medskip

We now proceed with the necessary virtual intersection calculus.
By \cite{Man},~Thms.~4.8 and 4.3, 
\begin{equation}
\label{eq:Muvirt}
[\scrM^{\mu,\gl}]^{\virt}=(\alpha'\circ p'')^![\foM^{\mu,\gl}]=
(p\times \alpha)^![\foM^{\mu,\gl}]= p^!\alpha^![\foM^{\mu,\gl}]=\alpha^!p^![\foM^{\mu,\gl}].
\end{equation}

In Step 3, we will show that $q'\circ p'$ and $q''\circ p''$ are
proper of DM type. This allows us to push-forward rational Chow classes.
This is a technical issue.
In \cite{Kresch}, push-forward of Chow classes on Artin stacks
were constructed only for
projective morphisms. However, in \cite{sko}, \cite{BS22}, 
proper push-forward for
DM type morphisms is constructed, assuming one works with rational Chow
groups. As we are interested only in the degree of cycles, this is sufficient
for our purpose. The compatibility of proper push-forward with 
virtual pull-back of \cite{Man},~Thm.~4.1,(i), can then be checked
to hold for this push-forward also. 

We may thus compute
\begin{align}
\label{eq:push-pull}
\begin{split}
\deg[\scrM^{\mu,\gl}]^{\virt}=\deg \alpha^!p^![\foM^{\mu,\gl}] = {} & \deg (q''\circ p'')_* 
\alpha^!p^![\foM^{\mu,\gl}]\\
= {} &\deg \alpha^! (q'\circ p')_* p^![\foM^{\mu,\gl}]\\
= {} &\deg \alpha^! (q'\circ p')_* (p')^![\foM^{\mu,\gl}].
\end{split}
\end{align}
Here, the last equality follows since the virtual pull-back
$(p')^!$ is defined using the pull-back of the relative obstruction
theory defining $p^!$.

As Theorem~\ref{thm: main gluing theorem} claims a value for
$\deg [\scrM^{\mu,\gl}]^{\virt}$ when the virtual dimension of $\scrM^{\mu,\gl}$
is zero, by the calculation of this virtual dimension in Step 1,
we may now assume that $A\cdot c_1(\Theta_{X/\kk})=0$. 
With this assumption, we will show 
\begin{claim}
\label{claim:claim}
\begin{enumerate}
\item If $s\in B(\ZZ)$, then both $\scrM_1(X,\beta_1,z')$ and $\scrM_2$
are virtual dimension zero for $z'\in Z_s^{\circ}$.
\item If $s\not\in B(\ZZ)$, then $\scrM_2$ is of negative virtual dimension.
\item $(q'\circ p')_*(p')^![\foM^{\mu,\gl}]= N^{A_1}_{p_1p_2s} [\foM_2].$
\end{enumerate}
\end{claim}

The theorem will then follow from \eqref{eq:push-pull} and these claims, 
as $\alpha^![\foM_2]= [\scrM_2]^{\virt}$ and $\deg[\scrM_2]^{\virt}
=N^{A_2}_{sp_3r}$.
We will show (1) and (2) of the claim in Step 4 and (3) of the claim in
Step 5.

\medskip

{\bf Step 3. $q'\circ p'$ and $q''\circ p''$ are proper of DM type.} 

There is a diagram of ordinary stacks
\[
\xymatrix@C=30pt
{
(\ul{\foM}_2\times B\GG_m)\times_{\ul{\scrP(X,s)}} \ul{\scrM}_1
\ar[r]\ar[d]&\ul{\scrM}_1\ar[d]\\
\ul{\foM}_2\times B\GG_m\ar[d]\ar[r]&\ul{\scrP(X,s)}\ar[d]\\
\ul{\foM}_2\ar[r]&\ul{Z}_s
}
\]
Here the top square is obviously cartesian and the bottom square,
with vertical maps being projections using
$\ul{\scrP(X,s)}=\ul{Z}_s\times B\GG_m$, is also cartesian.
But $\scrM'_1=\foM_2^{\mu}\times^{\fs}_{\scrP(X,s)}
\scrM_1$ by definition, and $\ul{\foM}_2^{\mu}=\ul{\foM}_2\times
B\GG_m$ by Lemma~\ref{lem: underlying}. Hence there is a finite
representable morphism
\[
\ul{\scrM}'_1\rightarrow 
(\ul{\foM}_2\times B\GG_m)\times_{\ul{\scrP(X,s)}} \ul{\scrM}_1.
\]
Thus, since $\ul{\scrM}_1\rightarrow \ul{Z}_s$
is a DM type morphism, as $\ul{\scrM}_1$ is already DM,
we see that $q'\circ p'$ is DM type. Similarly, as $\ul{\scrM}_1
\rightarrow \ul{Z}_s$ is proper by Lemma~\ref{lemma: properness},
it follows that $q'\circ p'$ is proper. By base-change, $q''\circ p''$
enjoys the same properties.

Note that as $\scrM_2$ is proper over $\Spec\kk$, this also shows that
$\scrM^{\mu,\gl}$ is proper over $\Spec\kk$, so that $\deg[\scrM^{\mu,\gl}]^{\virt}$
makes sense.

\medskip

{\bf Step 4. Calculation of virtual dimensions.}

We now show (1) and (2) of Claim~\ref{claim:claim}.
Recalling that we are assuming that $A \cdot c_1(\Theta_{X/\kk})=0$,
we see that if $\pm c_1(\Theta_{X/\kk})$ is nef, and both moduli spaces
$\scrM_1$ and $\scrM_2$ are non-empty, then necessarily $A_i
\cdot c_1(\Theta_{X/\kk})=0$ for $i=1,2$, and so $\scrM_1(X,\beta_1,z')$ and
$\scrM_2$ are both virtual dimension zero, for
$z'\in Z_s^{\circ}$. 

On the other hand, in the log Calabi-Yau case, where $c_1(\Theta_{X/\kk})
\equiv_{\QQ}-\sum_i a_i D_i$ for $a_i\ge 0$, we need to show that either
(1) $s\in B(\ZZ)$ and $A_i\cdot 
c_1(\Theta_{X/\kk})=0$ for $i=1,2$ as above,
or (2) $s\not\in B(\ZZ)$ and the virtual dimension of $\scrM_2$ is negative.

Indeed, by Corollary~\ref{cor:intersectionnumbers},
$A_1\cdot D_i= \langle p_1, D_i\rangle + \langle p_2, D_i\rangle
-\langle s, D_i\rangle$ for any $i$. Now $p_1,p_2\in B(\ZZ)$ by assumption, so
if $a_i\not=0$, the right-hand side is $-\langle s,D_i\rangle$.
If in fact $s\in B(\ZZ)$, then further the right-hand side is $0$,
so $A_1\cdot D_i=0$. Thus $A_1\cdot c_1(\Theta_{X/\kk})=0$.
So necessarily $A_2\cdot c_1(\Theta_{X/\kk})=0$, showing (1)
of Claim~\ref{claim:claim}.

On the other hand, if $s\in \Sigma(X)(\ZZ)\setminus 
B(\ZZ)$, we similarly have $A_2
\cdot D_i=\langle s, D_i\rangle+\langle p_3,D_i\rangle - \langle r,D_i\rangle$.
So if $a_i>0$, the last two terms on the right-hand side are zero
and $A_2\cdot D_i\ge 0$.
But there is at least one $i$ with $a_i\not=0$ and $\langle s,D_i\rangle
\not=0$ since $s\not\in B(\ZZ)$. Thus we see that $A_2\cdot
c_1(\Theta_{X/\kk})<0$, and hence the virtual dimension of $\scrM_2$ is
negative, showing (2).

\medskip

{\bf Step 5. Proof that $(q'\circ p')_*(p')^![\foM^{\mu,\gl}] =
N^{A_1}_{p_1p_2s} [\foM_2]$.}

We now turn to proving (3) of Claim~\ref{claim:claim}. The idea is 
fairly simple: morally, we have the cycle $(p')^![\foM^{\mu,\gl}]$
sitting over $\foM_2$, and we need to calculate the degree of this
cycle over $\foM_2$. To do this, it will be sufficient to restrict to a general
point of $\foM_2$ and calculate the degree of $(p')^![\foM^{\mu,\gl}]$
over such a point. We will do this in two steps, first by restricting
to a connected open stratum of $\foM_2$, and then passing to a point
in this stratum.

Note $\foM_2$ is of pure dimension $\dim X$, by Proposition 
\ref{prop:the invariants},(2). Thus $A_{\dim X}(\foM_2)\otimes_{\ZZ} \QQ$ is 
generated by the classes of the irreducible components of
$\foM_2$. Indeed, as $\foM_2$ is stratified
by global quotients (see \S\ref{sec:stacks remarks}), we can use excision 
(\cite{Kresch},~Thm.\ 2.1.12,(iv),(v)) to replace $\foM_2$
with a dense open substack which is a disjoint union of 
irreducible global quotients
without changing the group $A_{\dim X}$. However, for a global quotient
stack, $A_*$ coincides with the Edidin-Graham-Totaro Chow groups
by \cite{Kresch},~Thm.~2.1.12,(iii), and it is clear from the definition
of the latter that the top-dimensional Edidin-Graham-Totaro Chow group
of an irreducible global quotient
is rationally generated by the fundamental class of the global quotient.

Let $S_2$ be the union of open strata of $\foM_2$, so
$\dim \foM_2\setminus S_2< \dim X$. Hence, again by excision,
$A_{\dim X}(\foM_2)=A_{\dim X}(S_2)$. Let $S_2^{\mu}$ be the corresponding
open substack of $\foM_2^{\mu}$, and 
$S^{\mu}:=(\Pi^{\mu})^{-1}(S_2^{\mu})$
with $\Pi^{\mu}:\foM^{\mu,\gl}\rightarrow \foM_2^{\mu}$ 
as in \eqref{eq:def Pimu}.
We then have a diagram
\begin{equation}
\label{eq:M2 to S2}
\xymatrix@C=30pt
{
\scrM'_1\times_{\foM_2}S_2\ar[d]_{j_4}\ar[r]^>>>>>{p'}&S^{\mu}\ar[r]
\ar@/^2pc/[rr]^{q'}
\ar[d]_{j_3}&
S_2^{\mu}\ar[r]\ar[d]_{j_2}&S_2\ar[d]_{j_1}\\
\scrM'_1\ar[r]_{p'}& \foM^{\mu,\gl}\ar[r]\ar@/_2pc/[rr]_{q'}&\foM_2^{\mu}\ar[r]&\foM_2
}
\end{equation}
with all squares cartesian in all categories. Then
\begin{align}
\label{eq:bunch of equalities}
\begin{split}
j_1^*(q'\circ p')_*(p')^![\foM^{\mu,\gl}]= {} & (q'\circ p')_* j_4^*
(p')^![\foM^{\mu,\gl}]\\
= {} & (q'\circ p')_* (p')^!j_3^*[\foM^{\mu,\gl}]\\
= {} & (q'\circ p')_* (p')^! [S^{\mu}],
\end{split}
\end{align}
the first equality by standard compatibility of flat pull-back and
proper push-forward, the second by \cite{Man},~Thm.~4.1,(ii). Thus,
it is sufficient to show that 
\[
(q'\circ p')_*(p')^![S^{\mu}]=N^{A_1}_{p_1p_2s} [S_2].
\]

Without loss of generality, we may replace $S_2$ with a connected component 
of $S_2$. Next, we show that $S_2$ can be replaced by its reduction.
By the transversality of Proposition~\ref{prop:transverse M2}, (2) and
log smoothness of $\Pi$, we see that
$S^{\mu}\rightarrow S_2^{\mu}$ is flat. As the underlying stack
of $S_2^{\mu}$ is $\ul{S}_2\times B\GG_m$, the projection $S_2^{\mu}\rightarrow
S_2$ is flat. Thus $q':S^{\mu}\rightarrow S_2$ is flat.
By the definition of flat pull-back of cycles, 
$(q')^*[S_2]=[S^{\mu}]$. Writing $[S_2]=\nu [(S_2)_{\red}]$, for some
$\nu\in\QQ$, we see we can write $[S^{\mu}]=\nu \cdot (q')^*[(S_2)_{\red}]$.
Thus, we may replace $S_2$ by $(S_2)_{\red}$, $S_2^{\mu}$ by
$(S_2^{\mu})_{\red}$, and $S^{\mu}$ by $S^{\mu}\times_{S_2^{\mu}}
(S_2^{\mu})_{\red}$. Note now in particular that $S_2$ is smooth
over $\Spec \kk$.

Now let $s_2 \rightarrow S_2$ be a choice of strict morphism from
a standard log point with underlying scheme $\Spec\kk$.
Because $\ul{S}_2$
is smooth, the underlying stack morphism is regular of codimension $\dim X$.
We then have a diagram
\begin{equation}
\label{eq:S2 to p2}
\xymatrix@C=30pt
{
\scrM'_1\times_{\foM_2}s_2\ar[d]_{i_4}\ar[r]^>>>>>{p'}&
S^{\mu}\times_{S_2} s_2 \ar[r]
\ar@/^2pc/[rr]^{q'}
\ar[d]_{i_3}&
S_2^{\mu}\times_{S_2} s_2\ar[r]\ar[d]_{i_2}&s_2\ar[d]_{i_1}\\
\scrM'_1\times_{\foM_2} S_2\ar[r]_{p'}& S^{\mu}\ar[r]\ar@/_2pc/[rr]_{q'}&S_2^{\mu}\ar[r]&S_2
}
\end{equation}
with all squares cartesian in all categories.
Note that, with $i_1^!, i_3^!$ the usual Gysin pull-backs,
\begin{equation}
\label{eq:bernd requested}
i_1^!(p')^![S^{\mu}]=i_3^!(p')^![S^{\mu}]=(p')^!i_3^![S^{\mu}].
\end{equation}
Indeed, the first equality
follows from \cite{Fulton}, Thm.~6.2,(3), as flatness of $q'$ 
implies that $i_3$ is a regular morphism
of the same codimension as $i_1$. The second equality follows
from \cite{Man},~Thm.~4.3, as Gysin pull-back is a special case
of virtual pull-back, see \cite{Man},~Rmk.~3.9. 

Note that 
$(p')^![S^{\mu}]\in A_*(\scrM_1'\times_{\foM_2} S_2)$ (see the lower 
left-hand corner of \eqref{eq:S2 to p2}).
As $i_1^![S_2]=[s_2]$,
we wish to calculate the coefficient of $s_2$ in
\begin{align}
\label{eq:a whole bunch of equalities}
\begin{split}
i_1^!(q'\circ p')_*(p')^![S^{\mu}] = {} & (q'\circ p')_* i_1^! (p')^![S^{\mu}]\\
= {} &(q'\circ p')_* (p')^! i_3^![S^{\mu}]\\
= {} &(q'\circ p')_* (p')^![S^{\mu}\times_{S_2}s_2].
\end{split}
\end{align}
The first equality follows from the standard compatibility of
Gysin pull-back with push-forward, the second from \eqref{eq:bernd requested}
and the third from flatness and $i_1^![S_2]=[s_2]$.

Thus, it is sufficient to show that 
\[
(q'\circ p')_* (p')^![S^{\mu}\times_{S_2} s_2]=
N^{A_1}_{p_1p_2s} [s_2].
\]

Now take $W=S_2^{\mu}\times_{S_2} s_2$ in Definition~\ref{def: alternateN}.
Then by \eqref{eq:big diagram}, 
\begin{align}
\label{eq:bunch of isomorphisms}
\begin{split}
W\times_{\scrP(X,s)}^{\fs} \foM_1\cong  {} &
W\times_{\foM_2^{\mu}} (\foM^{\mu}_2\times^{\fs}_{\scrP(X,s)}\foM_1)\\
\cong {} & (S_2^{\mu}\times_{S_2} s_2)\times_{\foM_2^{\mu}} \foM^{\mu,\gl}\\
\cong {} & S^{\mu}\times_{S_2} s_2.
\end{split}
\end{align}
In particular, $\deg (p')^![S^{\mu}\times_{S_2} s_2]=
\deg [W\times^{\fs}_{\scrP(X,s)}\scrM(X,\beta_1)]^{\virt}=
N^{A_1,W}_{p_1p_2s}$ by Definition~\ref{def: alternateN}, and thus 
\[
(q'\circ p')_* (p')^![S^{\mu}\times_{S_2} s_2]=N^{A_1,W}_{p_1p_2s}[s_2].
\]
To conclude, it is sufficient to check that 
the induced morphism $g:W\rightarrow \scrP(X,s)$ satisfies
the hypotheses of Theorem~\ref{thm:NW equals N}, as then
$N^{A_1,W}_{p_1p_2s}=N^{A_1}_{p_1p_2s}$, proving the
desired result.

To check these hypotheses, first note that as $\ul{S}_2^{\mu}=\ul{S}_2
\times B\GG_m$ and $s_2\rightarrow S_2$ is strict, $\ul{W}=B\GG_m$.
Further, the stalk of the ghost sheaf of
$S_2\times B\GG_m^{\dagger}$ at any point is $\NN\delta'\oplus\NN\ell_q$,
with 
\begin{equation}
\label{eq:def nu}
\delta=\nu\delta'
\end{equation}
the image of 
$1$ under $\bar\delta^{\flat}$
with $\delta:S_2\rightarrow B\GG_m^{\dagger}$, as in the proof
of Lemma~\ref{lem: underlying}. Then the stalk $Q$ of
the ghost sheaf of $S_2^{\mu}$ is generated by 
$-\mu\delta+\ell_q$ and $\delta'$ in $\ZZ\delta'\oplus\ZZ\ell_q$,
see the proof of Lemma~\ref{lem: underlying}. Then $Q^{\vee}$
is generated by $\ell_q^*$ and $(\delta')^*+\mu\nu\ell_q^*$.
We see from the description of the map \eqref{eq:ghost sheaf map}
that $\Sigma(g)(\ell_q^*)=s$. The required transversality
follows since, by Proposition~\ref{prop:transverse M2}, (2),
$S_2^{\mu}\rightarrow \scrP(X,s)$ is transverse to $\ev_{\shX}$.
The remaining hypotheses
are clear.
\end{proof}

\subsection{The no-tail lemma}
\label{subsec:no tails}

Here we prove the crucial no-tail lemma mentioned a number of times
in \S\ref{section:sketch} and already applied in the proof of
Step 3 of the proof of Lemma~\ref{lemma: invarianceI}. Since the basic
idea is needed in a number of different contexts throughout the paper,
we prove here the various cases needed simultaneously. 

We first introduce some additional terminology.
\begin{definition}
Let, as usual, $\Mbf_{0,n+1}$ denote the moduli stack of pre-stable 
genus zero curves with 
$n+1$ marked points, $x_1,\ldots,x_n,x_{\out}$. We denote by
$D(x_1,\ldots,x_n\,|\,x_{\out})$ the reduced divisor on $\Mbf_{0,n+1}$ which is
the closure of the (locally closed) stratum parameterizing curves with two irreducible
components, with $x_1,\ldots,x_n$ contained in one irreducible component
and $x_{\out}$ contained in the other. 
\end{definition}

\begin{definition}
\label{def:strongly transverse}
Let $\beta$ be a class of punctured map of curve class $A$,
with punctured points $x_1,\ldots,x_n,x_{\out}$ having contact orders
$p_1,\ldots,p_n,-r$ respectively, with $p_1,\ldots,p_n,r\in B(\ZZ)$.
Suppose given a splitting $A=A_1+A_2$ with $A_1$ and $A_2$ effective curve
classes on $X$. This defines a class of punctured map 
$\beta_1$ with underlying curve
class $A_1$ and two punctured points, labelled $x_{\out}$
and $x$ respectively, with $x_{\out}$ having contact order specified
by $-r$. By Remark~\ref{rem:contact order determined}, the contact order
at $x$ is then specified up to a finite number of choices. Each choice
gives a moduli space and evaluation morphism at $x_{\out}$, 
\[
\ev_{\shX}: \foM^{\ev}(\shX,\beta_1)\rightarrow \scrP(X,r).
\]

Let $g:W\rightarrow\P(X,r)$ be a log morphism.
We say that $g$ is \emph{strongly $ev_{\shX}$-transverse} if $g$
is transverse to $\ev_{\shX}$ for each choice of splitting $A=A_1+A_2$
and choice of contact order for $x$.
\end{definition}

\begin{lemma}
\label{thm:no-tail-lemma}
Let $\beta$ be a class of punctured map with $A
\cdot c_1(\Theta_{X/\kk})=0$ and 
$\pm c_1(\Theta_{X/\kk})$ nef or $X$ log Calabi-Yau. Suppose further that $\beta$
has punctured points $x_1,\ldots,x_n,x_{\out}$ with contact orders
$p_1,\ldots,p_n,-r$, respectively, with $p_1,\ldots,p_n,r\in B(\ZZ)$. 
Set $Z=Z_r$. Suppose also given a log stack $W$ with $\ul{W}=B\GG_m$
and a morphism $g:W\rightarrow \scrP(X,r)$.
Let $F$ be an fs log Artin stack with $\ul{F}$ integral and
stratified by global quotients,
equipped with a morphism
$\psi:F\rightarrow \foM^{\ev}(\shX,\beta)$ such that there exists
a (locally closed) stratum $S^{\circ}\subset \foM^{\ev}(\shX,\beta)$
with closure $S$ such that: 
\begin{enumerate}
\item $\psi(F)\subset S$;
\item $\ev_{\shX}(S^{\circ})\subset \ul{Z}^{\circ}\times B\GG_m
\subset \ul{\scrP(X,r)}$.
\item The forgetful morphism $\phi:\foM^{\ev}(\shX,\beta)\rightarrow
\Mbf_{0,n+1}$ satisfies 
\[
\phi(S)\subset D(x_1,\ldots,x_n\,|\, x_{\out}).
\]
\item There is a morphism given by the dotted arrow
in the diagram
\[
\xymatrix@C=30pt
{
F\ar[r]^>>>>{\psi}\ar@{-->}[d]& \foM^{\ev}(\shX,\beta)\ar[d]^{\ev_{\shX}}\\
W\ar[r]_>>>>>g&\scrP(X,r)
}
\]
making the square commute.
\item
$W$ and $g$ satisfy the following
two conditions. First, they satisfy all hypotheses of Theorem
\ref{thm:NW equals N}; in particular, the required transversality is
with respect to the morphism $\ev_{\shX}:\foM^{\ev}(\shX,\beta)
\rightarrow\scrP(X,r)$. 
Second, $g$ is strongly $\ev_{\shX}$-transverse in the sense of 
Definition~\ref{def:strongly transverse}.
\end{enumerate}
Let $\varepsilon:M:=F\times_{\foM^{\ev}(\shX,\beta)} \scrM(X,\beta)\rightarrow
F$ be the base-change of $\scrM(X,\beta)\rightarrow\foM^{\ev}(\shX,\beta)$,
and let $\varepsilon$ carry the pull-back relative obstruction theory.
Then
\[
\varepsilon^![F] = 0
\]
in the rational Chow group of $M$. 
\end{lemma}

\begin{proof}
{\bf Step 1. Splitting the maps.}
\medskip

By Condition (1),
we have a factorization 
$M\rightarrow F\rightarrow S\rightarrow \foM^{\ev}(\shX,\beta)$,
and have universal punctured maps
$f_M:C_M^{\circ}/M\rightarrow X$, 
$f_F:C_F^{\circ}/F\rightarrow \shX$ and $f_S:C_S^{\circ}/S\rightarrow \shX$. 
We choose a node $q$ of $C_S$
as follows. 
Let
$\eta$ be the generic point of $S$, and let $C_{\bar\eta}\rightarrow\bar\eta$
be the generic fibre. Let $C'_{\bar\eta}\subset C_{\bar\eta}$
be the irreducible component of $C_{\bar\eta}$ containing $x_{\out}$.
By Condition (3), there is a 
unique node of $C_{\bar\eta}$ contained in
$C'_{\bar\eta}$ whose removal disconnects $C_{\bar\eta}$ into
two connected components, one containing $x_{\out}$ and the other containing
$x_1,\ldots,x_n$. Taking the closure of this node, we obtain a node of
the universal curve $C_S/S$. Since $C_F$ and $C_M$ are pull-backs of this
universal curve, we obtain a choice of node on these curves. We call
this node $q$.

As the base $F$ is integral, we can thus assume that $C_F/F$ can be
split at $q$ using \cite{ACGS18},~Prop.~5.2, giving morphisms
$f_{F,i}:C_{F,i}^{\circ}/F\rightarrow \shX$ with punctures $x\in C_{F,1}^{\circ}$,
$x'\in C_{F,2}^{\circ}$ mapping to the node $q$, and $C_{F,1}^{\circ}$ 
containing $x_{\out}$. We then also obtain a splitting of
$C_M$, getting $f_{M,i}:C_{M,i}^{\circ}/M\rightarrow X$.

Without loss of generality, 
we can fix a splitting $A=A_1+A_2$
and replace $M$ with the union
of connected components of $M$ such that 
the maps $f_{M,1}, f_{M,2}$ represent $A_1$, $A_2$.
This is necessary as a priori the
degree data of these two punctured maps is only locally constant on $M$.

Let $\beta_i$ be the class of the punctured map 
$f_{M,i}:C_{M,i}^{\circ}\rightarrow X$.
\medskip

{\bf Step 2. $A_1\cdot c_1(\Theta_{X/\kk})\le 0$.}
\medskip

In the case that $c_1(\Theta_{X/\kk})$ is nef or anti-nef, the fact
that $A\cdot c_1(\Theta_X)=0$ implies that $A_1\cdot
c_1(\Theta_X)=0$, provided that $\scrM(X,\beta_i)$ is non-empty
for $i=1,2$. (If one of these moduli spaces is empty, then $M$ is empty.)

In the log Calabi-Yau case (where we assume $D$ is normal crossings), 
we have $-c_1(\Theta_{X/\kk})\equiv_{\QQ}
\sum_i a_i D_i$, $a_i\ge 0$. Suppose that
$A_1\cdot c_1(\Theta_{X/\kk})>0$, i.e., 
\[
A_1\cdot \sum_i a_i D_i < 0.
\]
Then there exists an $i$ such that $a_i>0$ but $A_1\cdot D_i<0$.
However, by Corollary~\ref{cor:intersectionnumbers}, 
\[
A_1\cdot D_i = \langle u_{x_{\out}},D_i\rangle + \langle u_x, D_i
\rangle.
\]
Now $u_{x_{\out}}=-r$, and since $r$ is assumed to lie in $B(\ZZ)$,
rather than more broadly in $\Sigma(X)(\ZZ)$, we have $\langle r, D_i\rangle=0$.
Thus $\langle u_x, D_i\rangle <0$. By \cite{ACGS18},~Rmk.~2.20,
this can only happen if the irreducible
component $C'_{\bar\eta}$ is mapped by $f_S$ into $\shD_i$, the divisor
of $\shX$ corresponding to $D_i$.
But by Condition (2), $f_{\bar\eta}(x_{\out})\in \shZ^{\circ}$
(where $\shZ$ is the stratum of $\shX$ corresponding to $Z$). Hence 
$C'_{\bar\eta}$ cannot have image contained in $\shD_i$, as $\shD_i$ does not
contain $\shZ=\shZ_r$ as $\langle r,D_i\rangle=0$.
This is a contradiction, so in fact $A_1\cdot
c_1(\Theta_{X/\kk})\le 0$ in this case also.

\medskip

{\bf Step 3. The gluing diagram.}

\medskip
Note that the choice of node $q$ induces an evaluation morphism
$\ev_q:F\rightarrow\ul{\shX}$ at $q$, and we can then define
\[
F^{\ev(q)}:=F\times_{\ul{\shX}} \ul{X}.
\]
As usual, we have a factorization 
\[
\xymatrix@C=30pt
{
M\ar[r]_{\varepsilon^{\ev}} \ar@/^1pc/[rr]^{\varepsilon} & F^{\ev(q)}\ar[r]_{p_1}& F
}
\]
with the second morphism $p_1$ the projection, smooth since
it is the base-change of the smooth morphism $\ul{X}\rightarrow\ul{\shX}$.
We now have a diagram
\begin{equation}
\label{eq:splitting diagram}
\xymatrix@C=30pt
{
M\ar[r]\ar[d]_{\varepsilon^{\ev}} & \scrM^{\gl}(X,\beta_1,\beta_2)\ar[r]\ar[d]&
\scrM(X,\beta_1)\times \scrM(X,\beta_2)\ar[d]\\
F^{\ev(q)}\ar[r]&\foM^{\gl,\ev(q,x_{\out})}(\shX,\beta_1,\beta_2)
\ar[r]&\foM^{\ev(x_{\out},x)}(\shX,\beta_1)\times 
\foM^{\ev(x')}(\shX,\beta_2)
}
\end{equation}
with each square cartesian, and the relative obstruction theory
for the right vertical arrow pulling back to the natural relative
obstruction theories for the other two vertical arrows, as follows
from the discussion of \S\ref{sec:gluing review}. In particular,
the virtual fundamental class of $M$ may be calculated as either
$\varepsilon^![F]$ or $(\varepsilon^{\ev})^![F^{\ev(q)}]$.

By condition (4), the morphism
$F^{\ev(q)}\rightarrow \foM^{\ev(x_{\out},x)}(\shX,\beta_1)$ factors through
the composition
\[
\foM_1^{\ev}:=
(\foM^{\ev}(\shX,\beta_1)\times_{\scrP(X,r)}^{\fs} W)\times_{\ul{\shX}}
\ul{X}
\cong
\foM^{\ev(x_{\out},x)}(\shX,\beta_1)\times_{\scrP(X,r)}^{\fs} W
\rightarrow \foM^{\ev(x_{\out},x)}(\shX,\beta_1).
\]
Set
\[
\scrM_1:=\scrM(X,\beta_1)\times_{\foM^{\ev(x_{\out},x)}(\shX,\beta_1)} 
\foM^{\ev}_1.
\]
Then the outer cartesian rectangle of \eqref{eq:splitting diagram}
factors into two cartesian diagrams
\[
\xymatrix@C=30pt
{
M\ar[r]\ar[d]_{\varepsilon^{\ev}} & \scrM_1\times \scrM(X,\beta_2)\ar[d]\ar[r]&
\scrM(X,\beta_1)\times \scrM(X,\beta_2)\ar[d]\\
F^{\ev(q)}\ar[r]&
\foM^{\ev}_1\times \foM^{\ev(x')}(\shX,\beta_2)\ar[r] & 
\foM^{\ev(x_{\out},x)}(\shX,\beta_1)\times \foM^{\ev(x')}(\shX,\beta_2)
}
\]
The relative obstruction theory for the right vertical arrow
pulls back to a relative
obstruction theory for the middle vertical arrow also. 
Now note that the relative virtual dimension of $\scrM(X,\beta_1)$
over $\foM^{\ev(x_{\out},x)}(\shX,\beta_1)$ is 
\[
\chi\big((f_{M,1}^*\Theta_{X/\kk})(-x_{\out}-x)\big)
= A_1\cdot c_1(\Theta_{X/\kk})-\dim X\le -\dim X,
\]
the equality by \eqref{eq:rel virtual dim} and the inequality
by Step 2.
This is also the relative virtual dimension of $\scrM_1$ over
$\foM_1^{\ev}$. 
Thus, if we can show that 
\begin{equation}
\label{eq: necessary dim}
\dim \foM_1^{\ev} -\dim X <0,
\end{equation}
then the virtual dimension of $\scrM_1$ is negative and
Theorem~\ref{thm:virtual vanishing} 
shows that $\varepsilon^![F]=(\varepsilon^{\ev})^!([F^{\ev(q)}])=0$.

\medskip
{\bf Step 4. The dimension of $\foM_1^{\ev}$.}
\medskip

Note that
\[
\dim \foM_1^{\ev}= \dim(\foM^{\ev}(\shX,\beta_1)\times^{\fs}_{\scrP(X,r)} W)
\times_{\ul{\shX}}\ul{X}
=\dim (\foM^{\ev}(\shX,\beta_1)\times^{\fs}_{\scrP(X,r)} W)+\dim X.
\]
On the other hand, by Proposition~\ref{prop: log etale morphism},
$\foM(\shX,\beta_1)$ is idealized log \'etale over 
$\Mbf(\shX,r)=\Mbf\times_{B\GG_m}
\scrP(\shX,r)$, and hence 
\[
\foM(\shX,\beta_1)\times^{\fs}_{\scrP(\shX,r)} W = 
\foM^{\ev}(\shX,\beta_1)\times^{\fs}_{\scrP(X,r)} W
\]
is idealized
log \'etale over $\Mbf\times_{B\GG_m}W$.
Since the punctured map of class $\beta_1$ has two marked points
and the underlying curve is genus $0$, $\Mbf=\Mbf_{0,2}$ is of dimension $-1$,
so $\dim \Mbf\times_{B\GG_m}W$ is of dimension $-1$. Further,
$\Mbf\rightarrow B\GG_m$ is log \'etale (being log smooth
as in the proof of Theorem~\ref{thm: everything is log smooth}, and
an isomorphism on a dense open subset of $\Mbf$).
Thus the projection $\Mbf\times_{B\GG_m}W\rightarrow W$
is log \'etale. Putting this together, we see that 
$\foM^{\ev}(\shX,\beta_1)
\times^{\fs}_{\scrP(X,r)}W \rightarrow W$
is idealized log \'etale and by the transversality hypothesis of 
(5), is integral.
Thus $\dim\foM^{\ev}(\shX,\beta_1)\times^{\fs}_{\scrP(X,r)}W\le \dim W=-1$ 
by Proposition
\ref{prop:log etale dimensions}, (3), and so 
$\dim \foM_1^{\ev}<\dim X$, showing \eqref{eq: necessary dim} as desired.
\end{proof}

\section{Comparison of moduli spaces}
\label{section:key comparison}

We shall prove Theorem~\ref{thm:main comparison}, so we fix in this
section $p_1,p_2,p_3,r\in
B(\ZZ)$ and an underlying curve
class $A$, giving a punctured map class $\beta$. 
We fix $z\in Z_r^{\circ}$ 
as usual, and assume $c_1(\Theta_{X/\kk})\cdot A=0$.

The main point will be to study the diagram \eqref{eq:comparison diagram},
with a detailed analysis of the morphisms $i$ and $j$.

\subsection{First steps}

After fixing the data $A$, $p_1,p_2,p_3,r\in B(\ZZ)$,
we also fix a choice of splitting
$A=A_1+A_2$, and punctured map classes
$\beta_1$ with three points of contact orders determined by
$p_1,p_2,-s$ with $s\in \Sigma(X)(\ZZ)$, and $\beta_2$ with three points
of contact orders $s,p_3,-r$. In general, we continue to use the notation
of \S\ref{sec:key gluing}.

We will use short-hand notation for the stacks appearing in 
\eqref{eq:comparison diagram},
\begin{align*}
\scrM^{\ddagger,\mu,\gl}(X,\beta_1,\beta_2,z):= {} &
\scrM^{\mu,\gl}(X,\beta_1,\beta_2,z)\times^{\fs}_{\overline{\scrM}_{0,4}^{\dagger}}
B\GG_m^{\ddagger}\\
\foM^{\ddagger,\mu,\gl,\ev}(\shX,\beta_1,\beta_2,z):= {} &
\foM^{\mu,\gl,\ev}(\shX,\beta_1,\beta_2,z)\times^{\fs}_{\overline{\scrM}_{0,4}^{\dagger}}
B\GG_m^{\ddagger}
\end{align*}
Thus \eqref{eq:comparison diagram}
becomes
\begin{equation}
\label{eq:comparison diagram2}
\xymatrix@C=30pt
{
\coprod_{A_1,A_2,s} \scrM^{\mu,\gl}(X,\beta_1,\beta_2,z)\ar[d]_{k_1} &
\coprod_{A_1,A_2,s} \scrM^{\ddagger,\mu,\gl}(X,\beta_1,\beta_2,z)
\ar[l]_>>>>>>>{i'}
\ar[d]_{k_2} \ar[r]^>>>>>{j'} &
\scrM^{\ddagger}_y\ar[d]^{k_3}\\
\coprod_{s} \foM^{\mu,\gl,\ev}(\shX,\beta_1,\beta_2,z) &
\coprod_{s} \foM^{\ddagger,\mu,\gl,\ev}(\shX,\beta_1,\beta_2,z) 
\ar[l]^>>>>>{i}
\ar[r]_>>>>>j &
\foM^{\ddagger,\ev}_y
}
\end{equation}

Recall from \eqref{eq:T Tmu} the ghost sheaves $T=\NN\delta\oplus
\NN\ell_q$, $T_{\mu}\subseteq T^{\gp}$ 
of $B\GG_m^{\dagger} \times B\GG_m^{\dagger}$ 
and $B\GG_m^{2,\mu}$ respectively, with $T_{\mu}$
generated by $\delta$ and
$\ell_q-\mu\delta$. Here $\ell_q$ tropically measures the length
of the glued edge while $\mu$ is the parameter we use to ensure that this
length is sufficiently large, see the discussion of \S\ref{sec:second step}.

We also recall the notation from \eqref{eq:R Rlambda def} of the
monoid $R=\NN\ell\oplus\NN\delta$ and $R_{\lambda}\subseteq R^{\gp}$ 
generated by $\ell-\lambda\delta$ and $\delta$. Here $R$ is the stalk
of the ghost sheaf of $y\times B\GG_m^{\dagger}\subseteq 
\overline{\scrM}_{0,4}^{\dagger}$ and $R_{\lambda}$ is the stalk of
the ghost sheaf of $B\GG_m^{\ddagger}$. Now $\ell$ measures the tropical 
modulus of the stabilizaton of a four-pointed curve, 
see Remark~\ref{rem: tropicalization of m},
while $\lambda$ is the parameter we use to ensure this tropical modulus
is sufficiently large, see the discussion of \S\ref{sec:ind of modulus}.

\begin{remark}
We recall here that in a number of points of the arguments throughout
the paper we need to
replace the various moduli spaces of punctured maps to
$\shX$ with finite type approximations. This does not always need to
be done consistently in general: for example, the virtual fundamental
class of $\scrM^{\mu,\gl}(X,\beta_1,\beta_2,z)$ is impervious to which
particular finite type open substack
of $\foM^{\mu,\gl,\ev}(\shX,\beta_1,
\beta_2,z)$ one uses, provided it contains the image of
$\scrM^{\mu,\gl}(X,\beta_1,\beta_2,z)$. Thus, in the proof of
Theorem~\ref{thm: main gluing theorem}, one needs to choose
an open set of $\foM^{\mu,\gl,\ev}(\shX,\beta_1,\beta_2,z)$
by choosing open subsets of the moduli spaces
$\foM^{\ev(x_s)}(\shX,\beta_1)$ and $\foM^{\ev(x_s',x_{\out})}(\shX,
\beta_2,z)$ from which the glued moduli space is constructed.
However, starting in Lemma~\ref{lem:finite representable}, we need
to be more careful. In particular, to guarantee the claimed finiteness
of the maps $i,j$ in that lemma, we will need to make more careful
choices of the three
moduli spaces in the bottom row of \eqref{eq:comparison diagram2}.

To do this, note all three moduli spaces come with a canonical
morphism to $\foM^{\ev}(\shX,\beta,z)$ and $i$ and $j$ are compatible
with these canonical morphisms. We choose a suitable finite
type open subset of $\foM^{\ev}(\shX,\beta,z)$ and pull these back
to the three moduli spaces to obtain the desired finite type approximations.
This choice will be made implicitly from now on.

Note that one consequence is that at most a finite number of 
$s\in\Sigma(X)(\ZZ)$ now appear in the disjoint union in 
\eqref{eq:comparison diagram2}. Indeed, once $\foM^{\ev}(\shX,\beta,z)$
is taken to be finite type, it has a finite number of strata, and hence
only a finite number of types of punctured map occur. Each type has
a finite number of nodes at which to split, giving at most a finite number
of contact orders for these nodes.
\end{remark}

In order to use push-pull formulas 
for the various squares of \eqref{eq:comparison diagram2}, we need:

\begin{lemma}
\label{lem:finite representable}
All squares in \eqref{eq:comparison diagram2} are cartesian in all categories,
and the morphisms $i$ and $j$ are finite and representable.
\end{lemma}

\begin{proof}
The left-hand square in \eqref{eq:comparison diagram2} is obviously
fs log cartesian, and the vertical arrows are strict, hence the square
is cartesian in all categories.

For the right-hand square, note that we have 
\[
\xymatrix@C=30pt
{
\coprod_{A_1,A_2,s}\scrM^{\mu,\gl}(X,\beta_1,\beta_2,z)
\ar[d]\ar[r] & \coprod_{A_1,A_2,s}\scrM^{\gl}(X,\beta_1,\beta_2,z)
\ar[r]\ar[d]& \scrM(X,\beta,z)\ar[d]\\
\coprod_s\foM^{\mu,\gl,\ev}(\shX,\beta_1,\beta_2,z)\ar[r]&
\coprod_s\foM^{\gl,\ev}(\shX,\beta_1,\beta_2,z)\ar[r]&
\foM^{\ev}(\shX,\beta,z)
}
\]
The first square is cartesian by definition of the spaces in the left-hand
column, while the right-hand square is cartesian as in 
\eqref{eq:Mgl obstruction diagram}. Here, taking the union over all
curve classes is necessary to achieve this, as 
\eqref{eq:Mgl obstruction diagram} was not cartesian. Applying
$\times^{\fs}_{\overline{\scrM}_{0,4}^{\dagger}} B\GG_m^{\ddagger}$
to the outer rectangle then gives the right-hand square of 
\eqref{eq:comparison diagram2} cartesian.

For the properties of $i$, as the morphism $B\GG_m^{\ddagger}\rightarrow
\overline\scrM_{0,4}^{\dagger}$ is a closed embedding of underlying
stacks, finiteness and representability are immediate from Remark
\ref{rem:sat finite representable}.

For the properties of $j$, the morphism 
$\foM^{\gl} (\shX,\beta_1,\beta_2)\rightarrow \foM(\shX,\beta)$ is
representable and finite, see \cite{ACGS18},~Rmk.~3.31.
We may then apply $\times^{\fs}_{\scrP(\shX,s)}
B\GG_m^{\dagger}$ to the morphism to see, using Remark 
\ref{rem:sat finite representable} again, that 
$\foM^{\gl,\ev}(\shX,\beta_1,\beta_2,z)\rightarrow \foM^{\ev}(\shX,\beta,z)$
is finite and representable. 
The morphism $\foM^{\mu,\gl,\ev}(\shX,\beta_1,\beta_2,z)\rightarrow 
\foM^{\gl,\ev}(\shX,\beta_1,\beta_2,z)$ is also finite and
representable (as follows from the construction
of the former stack and Remark~\ref{rem:sat finite representable}), and
thus we may compose the two morphisms
to obtain $\foM^{\mu,\gl,\ev}(\shX,\beta_1,\beta_2,z)
\rightarrow \foM^{\ev}(\shX,\beta,z)$ finite and representable. 
We then apply 
$\times^{\fs}_{\overline{\scrM}_{0,4}^{\dagger}} B\GG_m^{\ddagger}$ to
get $j$, and thus $j$ is similarly finite and representable.
\end{proof}

We will need one local calculation of a fibre product, necessary
for understanding the morphism $i$ of \eqref{eq:comparison diagram2}.
The monoid $Q$ in the lemma below will ultimately be the stalk
of the ghost sheaf of a generic point of $\foM^{\mu,\gl,\ev}(\shX,
\beta_1,\beta_2,z)$ if $\mu$ is chosen sufficiently large. Thus 
the fibre product in the lemma will provide
a model for passing from $\foM^{\mu,\gl,\ev}(\shX,\beta_1,\beta_2,z)$
to $\foM^{\ddagger,\mu,\gl,\ev}(\shX,
\beta_1,\beta_2,z)$. This is where we need to take $\lambda$ sufficiently
large given the choice of $\mu$. Ultimately, one $\lambda$ must be chosen
to work for all possible splittings of $\beta$, hence we do not quantify
precisely how large $\lambda$ needs to be taken: see the last paragraph
of Step 3 of the proof of Theorem~\ref{thm:second comparison}.

\begin{lemma}
\label{lem:local product}
Let $Q$ be a rank two fs monoid with elements $\ell_q,\delta\in Q^{\gp}$
such that $Q$ is rationally generated by $\delta$ and $\ell_q-\mu\delta$.
Let $R,R_{\lambda}$ be as defined in \eqref{eq:R Rlambda def}, and
$J=R\setminus\{0\}$,
$J_{\lambda}=R_{\lambda}\setminus\{0\}$ be the maximal monomial ideals
of $R$ and $R_{\lambda}$ respectively. Let $K\subseteq Q$ be the 
ideal generated by $\delta$ and $\ell_q-\mu\delta$. 
Let $\theta:R\rightarrow Q$ be a homomorphism given by 
\begin{equation}
\label{eq:delta ell}
\delta\mapsto \delta,\quad\quad \ell\mapsto a\ell_q+b\delta,
\end{equation}
for some $a,b\in\QQ$, $a>0$. Then, using the notation of
Definition~\ref{def:APK}, for $\lambda$
sufficiently large, assuming $\mu$ fixed, the projection
\[
(\ul{\shA_{Q,K}\times^{\fs}_{\shA_{R,J}} \shA_{R_{\lambda},J_{\lambda}}})_{\red}
\rightarrow (\ul{\shA_{Q,K}})_{\red}
\]
is an isomorphism of integral stacks. Further, if $a=1$, then
the multiplicities
of $\shA_{Q,K}\times^{\fs}_{\shA_{R,J}} \shA_{R_{\lambda},J_{\lambda}}$
and $\shA_{Q,K}$ coincide.
\end{lemma}

\begin{proof}
We first calculate
$Q_{\lambda}:=Q\oplus_R^{\fs} R_{\lambda}$. 
We write $\theta:R^{\gp}\rightarrow Q^{\gp}$ also for the 
extension of $\theta$ to the groups.
Because $R^{\gp}=R_{\lambda}^{\gp}$,
$Q^{\gp}\oplus_{R^{\gp}}R_{\lambda}^{\gp}\cong Q^{\gp}$,
and thus by Proposition~\ref{prop:fibre product properties}, 
$Q_{\lambda}^{\gp}=Q^{\gp}$ and $Q_{\lambda}$ is the saturated submonoid
of $Q^{\gp}$ rationally generated by $\ell_q-\mu\delta,\delta$ and 
$\theta(\ell-\lambda\delta)
=a\ell_q+(b-\lambda)\delta$. Note that if $\lambda$ is taken sufficiently
large, $\ell_q-\mu\delta$ lies in the saturated submonoid of $Q^{\gp}$
rationally generated by $\delta$ and $\theta(\ell-\lambda\delta)$,
and thus $Q_{\lambda}$ is precisely this latter submonoid.

Next we consider the ideals. Note that $J$ is generated by $\ell$ and
$\delta$ while $J_{\lambda}$ is generated by $\ell-\lambda\delta$ and
$\delta$. Thus,
letting $K'$, $J_{\lambda}'$ be the ideals generated by the images
of $K, J_{\lambda}$ in $Q_{\lambda}$, we see that $K'\cup J_{\lambda}'$
is generated by $\delta, \ell_q-\mu\delta$ and $\theta(\ell-\lambda\delta)$.
In particular, $\Spec \kk[Q_{\lambda}]/(K'\cup J_{\lambda}')$ is finite length.
Of course $\Spec \kk[Q]/K$ is finite length, so the projection
(using Proposition~\ref{prop:fibre product properties}, (4))
\begin{equation}
\label{eq:product description}
\shA_{Q_{\lambda},K'\cup J_{\lambda}'}\cong
\shA_{Q,K}\times^{\fs}_{\shA_{R,J}} \shA_{R_{\lambda},J_{\lambda}}
\rightarrow \shA_{Q,K}
\end{equation}
induces an isomorphism on underlying reduced stacks.

Now assume that $a=1$. By Pick's theorem \cite{Pick}, the length of 
$\kk[Q]/K$ coincides with the area of the parallelogram
with vertices $0, \delta, \ell_q-\mu\delta$ and
$\ell_q-(\mu-1)\delta$. Now with $a=1$, $K'\cup J_{\lambda}'$ is
generated by $\delta$ and $\ell_q+(b-\lambda)\delta$, provided again
that $\lambda$ is sufficiently large. Thus the length of
$\kk[Q_{\lambda}]/(K'\cup J_{\lambda}')$ is the area of the parallelogram
with vertices $0,\delta, \ell_q+(b-\lambda)\delta$, and 
$\ell_q+(b-\lambda+1)\delta$. However, this parallelogram has the
same area as the previous one, and hence the two stacks have the
same multiplicities.\footnote{Note the projection in
\eqref{eq:product description} need not be an isomorphism!}
\end{proof}

\begin{lemma}
\label{lem:first comparison facts}
Suppose $\mu$ is chosen sufficiently large so that $B\GG_m^{2,\mu}
\rightarrow B\GG_m^{\dagger}\times B\GG_m^{\dagger}$ is transverse to
$\Psi$, as in Proposition 
\ref{prop:transverse M2}, (3). Then for $\lambda\gg\mu$,
we have:
\begin{enumerate} 
\item
Both
$\foM^{\mu,\gl,\ev}(\shX,\beta_1,\beta_2,z)$
and
$\foM^{\ddagger,\mu,\gl,\ev}(\shX,\beta_1,\beta_2,z)$ are
pure-dimensional of dimension $0$.
\item The rank of the stalk of the ghost sheaf at each generic point
of 
\[
\hbox{$\foM^{\mu,\gl,\ev}(\shX,\beta_1,\beta_2,z)$
and
$\foM^{\ddagger,\mu,\gl,\ev}(\shX,\beta_1,\beta_2,z)$}
\]
is $2$.
\item The induced morphism on reduced underlying stacks
\[
i_{\red}:\ul{\foM^{\ddagger,\mu,\gl,\ev}(\shX,\beta_1,\beta_2,z)}_{\red}
\rightarrow
\ul{\foM^{\mu,\gl,\ev}(\shX,\beta_1,\beta_2,z)}_{\red}
\]
is an isomorphism on the union of open strata of both spaces. Thus
there is a one-to-one correspondence between irreducible components of
the two spaces.
\end{enumerate}
\end{lemma}

\begin{proof}
{\bf Step 1.}
Statements (1) and (2) for $\foM^{\mu,\gl,\ev}(\shX,\beta_1,\beta_2,z)$
follow immediately from Proposition~\ref{prop:transverse M2}, (1),
bearing in mind that, in that proposition,
$\foM^{\mu}= \foM^{\mu,\gl,\ev}(\shX,\beta_1,\beta_2,z)
\times_{\ul{\shX}}\ul{X}$, hence the shift in dimension.

\medskip

{\bf Step 2. Analysis of $i$ on open strata using Lemma~\ref{lem:local product}.}

\medskip

Let $I=T\setminus \{0\}$, $I_{\mu}=T_{\mu}\setminus\{0\}$ be the
maximal monomial ideals of $T=\NN\delta+\NN\ell_q$, and $T_{\mu}\subseteq
T^{\gp}$ generated by
$\delta$ and $\ell_q-\mu\delta$, see \eqref{eq:T Tmu}.
Note $B\GG_m^{\dagger}\times B\GG_m^{\dagger}=\shA_{T,I}$
and $B\GG_m^{2,\mu}=\shA_{T_{\mu},I_{\mu}}$.

Just as in the proof of Proposition~\ref{prop:transverse M2}, 
the morphism $\Psi:\foM^{\gl,\ev}(\shX,\beta_1,\beta_2,z)\rightarrow
B\GG_m^{\dagger}\times B\GG_m^{\dagger}$ is log smooth, and 
thus so is its base-change 
\begin{equation}
\label{eq:Phimudef}
\Psi_{\mu}:\foM^{\mu,\gl,\ev}(\shX,\beta_1,\beta_2,z)
\rightarrow B\GG_m^{2,\mu}.
\end{equation}
We continue to abuse notation by
also writing $\ell_q, \delta$ for sections of the ghost sheaf
of $\foM^{\mu,\gl,\ev}(\shX,\beta_1,\beta_2,z)$ which are images of
$\ell_q,\delta$ in the ghost sheaf of $B\GG_m^{2,\mu}$ under 
$\overline{\Psi}_{\mu}^{\flat}$.

Now consider an open stratum $S$ of 
$\foM^{\mu,\gl,\ev}(\shX,\beta_1,\beta_2,z)$. Then the ghost sheaf is
locally constant with some stalk $Q$, with $Q$ of rank $2$. In fact,
the ghost sheaf is necessarily constant, as $\ell_q,\delta$ span
$Q^{\gp}\otimes_{\ZZ}\QQ$, and hence there is no possibility of monodromy.
The resulting isomorphism $Q\cong \Gamma(S,\overline{\shM}_S)$
induces a unique strict morphism $S\rightarrow \shA_{Q}$ by 
Proposition~\ref{prop:APrep}. Consider the composition $S\rightarrow 
B\GG_m^{2,\mu}
\hookrightarrow \shA_{T_{\mu}}$. This is determined by the homomorphism
\[
\overline{\Psi}_{\mu}^{\flat}:T_{\mu}\rightarrow Q=\Gamma(S,\overline{\shM}_S),
\]
and hence $S\rightarrow\shA_{T_{\mu}}$ factors as 
$S\rightarrow\shA_Q\rightarrow \shA_{T_{\mu}}$.
Letting $K$ be the ideal generated by the image of $I_{\mu}$
in $Q$, we thus
obtain a factorization of $\Psi_{\mu}|_S$ as
\begin{equation}
\label{eq:Psi mu fact}
S\rightarrow \shA_{Q}\times_{\shA_{T_{\mu}}} \shA_{T_\mu,I_{\mu}}=\shA_{Q,K}
\rightarrow \shA_{T_{\mu},I_{\mu}}=B\GG_m^{2,\mu}.
\end{equation}
Note that $K$ is
the ideal defined in the statement of Lemma~\ref{lem:local product}.
As $\Psi_{\mu}$ is log smooth, so is $\Psi_{\mu}|_S$.
Further, $\shA_{Q,K}\rightarrow\shA_{T_{\mu},I_{\mu}}$
is log \'etale by \cite{ACGS18}, Lem.~B.3. Thus the
composition in \eqref{eq:Psi mu fact} is log smooth and the second morphism
is log \'etale, hence 
$S\rightarrow \shA_{Q,K}$ is log smooth. Since this morphism is
strict, it is also smooth.

Because $S\rightarrow B\GG_m^{2,\mu}$ is integral
by the choice of $\mu$ in Proposition~\ref{prop:transverse M2}, 
the image of the generators $\delta,\ell_q-\mu\delta$ of $I_{\mu}$
rationally generate $Q$, and hence $Q,K$ satisfy the conditions of
Lemma~\ref{lem:local product}. 
Further, the composed morphism
$S\hookrightarrow \foM^{\mu,\gl,\ev}(\shX,\beta_1,\beta_2,z)
\rightarrow\overline{\scrM}_{0,4}^{\dagger}$ induces the stated
morphism $\theta:R\rightarrow Q$ of \eqref{eq:delta ell}
for some choice of $a,b\in\QQ$, $a>0$.

Give the point $y\in\overline{\scrM}_{0,4}$ (determining the morphism
$B\GG_m^{\ddagger}\rightarrow \overline{\scrM}_{0,4}^{\dagger}$)
the induced log structure
from $\overline{\scrM}_{0,4}$, so that $y$ is a standard log point.
Noting that $B\GG_m^{\ddagger}\cong 
(y\times B\GG_m^{\dagger})
\times_{\shA_{R,J}} \shA_{R_{\lambda},J_{\lambda}}$, we have
\begin{align*}
S\times_{\overline{\scrM}_{0,4}^{\dagger}}^{\fs} B\GG_m^{\ddagger}
%\cong {} & S\times^{\fs}_{y\times B\GG_m^{\dagger}}B\GG_m^{\ddagger}\\
\cong {} & S\times^{\fs}_{y\times B\GG_m^{\dagger}} \left(
(y\times B\GG_m^{\dagger})
\times_{\shA_{R,J}} \shA_{R_{\lambda},J_{\lambda}}\right)\\
\cong {} & S\times^{\fs}_{\shA_{R,J}} \shA_{R_{\lambda},J_{\lambda}}.
\end{align*}
Here the first isomorphism holds because $S\rightarrow 
\overline{\scrM}_{0,4}^{\dagger}$ factors through the strict
inclusion $y\times B\GG_m^{\dagger}\hookrightarrow 
\overline{\scrM}_{0,4}^{\dagger}$.
We thus have a diagram cartesian in all categories
\begin{equation}
\label{eq:Sdiagram}
\xymatrix@C=30pt
{
S\times_{\overline{\scrM}_{0,4}^{\dagger}}^{\fs} B\GG_m^{\ddagger}\ar[d]\ar[r]
& S\ar[d]\\
\shA_{Q,K}\times_{\shA_{R,J}}^{\fs} \shA_{R_{\lambda},J_{\lambda}}\ar[r]
&
\shA_{Q,K}
}
\end{equation}
with vertical arrows strict and smooth. Thus the diagram remains
cartesian in the category of underlying stacks after taking the reduction
of all spaces, and thus $(S\times^{\fs}_{\overline{\scrM}_{0,4}^{\dagger}}
B\GG_m^{\ddagger})_{\red}\rightarrow S_{\red}$ is an isomorphism for
$\lambda\gg \mu$ by
Lemma~\ref{lem:local product} as claimed in item (3).

\medskip

{\bf Step 3. $i$ induces a continuous bijection on underlying
topological spaces $|i|:|\foM^{\ddagger,\mu,\gl,\ev}(\shX,\beta_1,
\beta_2,z)|\rightarrow
|\foM^{\mu,\gl,\ev}(\shX,\beta_1,
\beta_2,z)|$.}

\medskip

Recall the morphisms from Remark~\ref{rem:sat finite representable}
\begin{align*}
\foM^{\mu,\gl,\ev}(\shX,\beta_1,\beta_2,z)
\times^{\fs}_{\overline\scrM_{0,4}^{\dagger}}
B\GG_m^{\ddagger}
\mapright{\nu} &
\foM^{\mu,\gl,\ev}(\shX,\beta_1,\beta_2,z)
\times^{\fine}_{\overline\scrM_{0,4}^{\dagger}}B\GG_m^{\ddagger}\\
\mapright{\iota}&
\foM^{\mu,\gl,\ev}(\shX,\beta_1,\beta_2,z)
\times_{\overline\scrM_{0,4}^{\dagger}}
B\GG_m^{\ddagger},
\end{align*}
with $\nu$ finite and surjective and $\iota$
a closed embedding. Since the underlying stack of the last space is
$\ul{\foM^{\mu,\gl,\ev}(\shX,\beta_1,\beta_2,z)}$, the result of Step 2
shows that $\ul{\iota}\circ\ul{\nu}$ is dominant, and
hence $\ul{\iota}$ induces an isomorphism on the underlying
reduced stacks. Thus $\ul{\iota}\circ\ul{\nu}$ is surjective on
geometric points. On the other hand, $\nu$ is the saturation
morphism, and in fact, is also one-to-one on geometric points.
Indeed, if the stalk of the ghost sheaf of a geometric point 
$\bar x \in |\foM^{\mu,\gl,\ev}(\shX,\beta_1,\beta_2,z)|$ is $Q$, then the 
stalk of the ghost sheaf of the corresponding geometric point of
$\foM^{\mu,\gl,\ev}(\shX,\beta_1,\beta_2,z)
\times^{\fine}_{\overline\scrM_{0,4}^{\dagger}}B\GG_m^{\ddagger}$
is $Q_{\lambda}^{\fine}:=Q\oplus^{\fine}_{R} R_{\lambda}$. 
As $R^{\gp}=R_{\lambda}^{\gp}$,
this is a submonoid of $Q^{\gp}$ by Proposition 
\ref{prop:fibre product properties}, (3). Now $Q^{\gp}$ is torsion-free,
and hence so is $Q_{\lambda}$, the saturation of $Q^{\fine}_{\lambda}$.
It follows that $\Spec \kk[Q_{\lambda}]\rightarrow
\Spec \kk[Q_{\lambda}^{\fine}]$ is one-to-one on geometric points, and hence
the fibre of $\ul{\iota}\circ\ul{\nu}$ over $\bar x$ is a single
(possibly non-reduced) point. Thus the projection morphism
$\foM^{\ddagger,\mu,\gl,\ev}(\shX,\beta_1,\beta_2,z)
\rightarrow
\foM^{\mu,\gl,\ev}(\shX,\beta_1,\beta_2,z)$
is injective on geometric points, hence bijective.

\medskip

{\bf Step 4.}
The pure-dimensionality of (1) for $\foM^{\ddagger,\mu,\gl,\ev}(\shX,
\beta_1,\beta_2,z)$ is then immediate from the
pure-dimensionality of $\foM^{\mu,\gl,\ev}(\shX,\beta_1,\beta_2,z)$ and
Step 3, 
and the stalk of the
ghost sheaf at a generic point of this space is of the form 
$Q_{\lambda}$ as in the proof of Lemma~\ref{lem:local product},
which is rank $2$, giving (2).
\end{proof}

\begin{lemma}
Let $Q$ be the stalk at a generic point $\bar\xi$ of either
$\foM^{\mu,\gl,\ev}(\shX,\beta_1,\beta_2,z)$,
$\foM^{\ddagger,\mu,\gl,\ev}(\shX,\beta_1,\beta_2,z)$ or
$\foM_y^{\ddagger,\ev}$, so that $\sigma=\sigma_{\bar\xi}:=\Hom(Q,\RR_{\ge 0})$ 
parameterizes a family
of tropical maps with target $\Sigma(X)$.
Then this family of tropical maps satisfies Assumptions~\ref{ass:tropical}.
\end{lemma} 

\begin{proof}
In all three cases, $Q$ is rank two. Indeed, in the first
two cases, this is Lemma~\ref{lem:first comparison facts}, (2).
For the third case, this follows from log smoothness and integrality
of $\foM^{\ddagger,\ev}_y\rightarrow B\GG_m^{\ddagger}$, along with Proposition 
\ref{prop: realisability}: see the proof of
Proposition~\ref{prop:transverse M2}, (1) for an identical argument.
We now check conditions (1)-(3) of Assumptions~\ref{ass:tropical}. 

(1) is clear in all three cases. The miniversality of (2) is
as follows. The moduli space $\foM^{\gl,\ev}(\shX,\beta_1,\beta_2,z)$
carries the basic log structure,
i.e., the morphism $\foM^{\gl,\ev}(\shX,\beta_1,\beta_2,z)
\rightarrow \foM(\shX,\beta,z)$ is strict.
In particular, the image $\bar\xi'$ of
$\bar\xi$ in this latter moduli space 
satisfies $\sigma_{\bar\xi'}=\sigma_{\bas}$. By Proposition 
\ref{tropicalproduct} and the definition of 
$\foM^{\mu,\gl,\ev}(\shX,\beta_1,\beta_2,z)$, 
$\sigma=\sigma_{\bas}\times_{T^{\vee}_{\RR}}T_{\mu,\RR}^{\vee}$ is then a 
subcone
of $\sigma_{\bas}$ of the same dimension, giving miniversality.
Similar arguments show miniversality for the other two moduli spaces.

Condition (3) holds automatically for $\foM^{\ddagger,\ev}_y$ by construction.
In the case of $\foM^{\mu,\gl,\ev}(\shX,\beta_1,\beta_2,z)$,
$\ell_q$ is a summand of $\ell$ inside $Q$ and $\sigma$ is generated
(over $\RR$) by $\delta^*+\mu\ell_q^*$ and $\ell_q^*$, where
$\delta^*,\ell^*_q$ are the dual basis to $\delta,\ell_q$ of
$Q^{\gp}\otimes_{\ZZ}\RR$. Thus
$\ell_q$ is unbounded for fixed value of $\delta$ on $\sigma$, and
so the same is true for $\ell$. Thus the map 
$(\ell,\delta)$ of Condition (3) has the
claimed properties. Finally, the same is true in the case
of $\foM^{\ddagger,\mu,\gl,\ev}(\shX,\beta_1,\beta_2,z)$ by construction.
\end{proof}

This lemma allows us to make the following definition:

\begin{definition}
We say an irreducible component $\foM$ of 
\[
\hbox{
$\foM^{\mu,\gl,\ev}(\shX,\beta_1,\beta_2,z)$,
$\foM^{\ddagger,\mu,\gl,\ev}(\shX,\beta_1,\beta_2,z)$ or
$\foM^{\ddagger,\ev}_y$}
\]
with generic point $\bar\xi$ is \emph{tail-free} 
(resp.\ \emph{does not have a terminal
tail}, \emph{does not have an internal tail}) if the same is true
of the family of tropical maps parameterized by $\sigma_{\bar\xi}$.
If $\foM$ is not tail-free, we say it has a tail.
\end{definition}

\pagebreak

\subsection{Vanishing contributions from tails}

We will see in the proofs of Theorems~\ref{thm:first comparison} 
and \ref{thm:second comparison}
that it is precisely irreducible components with tails
which may potentially cause trouble for comparing virtual fundamental classes
using both the morphisms $i$ and $j$ in diagram \eqref{eq:comparison diagram2}.
So we now show that irreducible components with tails do not contribute
to the virtual count. We continue to use the notation $i,j,i',j',k_1,k_2,k_3$
for the morphisms of \eqref{eq:comparison diagram2}.

\begin{lemma}
\label{lemma:terminal tail vanishing}
Let $\foM$ be a (reduced) irreducible component of either
\[
\hbox{$\foM^{\mu,\gl,\ev}(\shX,\beta_1,\beta_2,z)$, 
$\foM^{\ddagger,\mu,\gl,\ev}(\shX,\beta_1,\beta_2,z)$, or
$\foM^{\ddagger,\ev}_y$
}
\]
with a terminal tail.  Then $k_1^![\foM]$, $k_2^![\foM]$ or $k_3^![\foM]$ 
respectively vanishes in the rational Chow group.
\end{lemma}

\begin{proof}
The statement follows from an application of
Lemma~\ref{thm:no-tail-lemma}, as in Step 3 of the proof of Lemma
\ref{lemma: invarianceI}, the first step of the independence
of modulus result.
Take $F=\foM$, $\psi:F\rightarrow
\foM^{\ev}(\shX,\beta)$ the tautological morphism, $W=B\GG_m^{\dagger}$
and $g:W\rightarrow\scrP(X,r)$ the standard morphism given by the choice
of $z$ and Proposition~\ref{prop: BGm morphism}. Take $S^{\circ}$ to be the
stratum of $\foM^{\ev}(\shX,\beta)$ containing the image of the
generic point of $F$, so that its closure $S$
contains the image of $\psi$. Thus
condition (1) of Lemma~\ref{thm:no-tail-lemma} is automatic, condition (2)
holds since the image of $\ev_{\shX}:\foM\rightarrow\scrP(X,r)$ is
contained in $z\times B\GG_m \subseteq \ul{Z}^{\circ}_r\times B\GG_m$.
Because $G$ has a terminal tail, (3) holds. 
Condition (4) is clear as
$\foM$ is a component of a moduli space of punctured maps with point constraint
at $z$.  Finally, (5) is immediate:
the transversality statements are trivial because the ghost sheaf on $W$ is
$\NN$. Thus Lemma~\ref{thm:no-tail-lemma} allows us to conclude that
$k_i^![\foM]=0$ in the rational Chow group.
\end{proof}

\begin{lemma}
\label{lemma:internal tail vanishing}
Let $\foM$ be a (reduced) irreducible component of either
\[
\hbox{
$\foM^{\mu,\gl,\ev}(\shX,\beta_1,\beta_2,z)$
or
$\foM^{\ddagger,\mu,\gl,\ev}(\shX,\beta_1,\beta_2,z)$
}
\]
without a terminal tail but with an internal tail. Provided $\mu\gg 0$,
we have $\deg k_1^![\foM]=0$ or $\deg k_2^![\foM]=0$ in the two cases.
\end{lemma}

\begin{proof}
We will reduce to a similar gluing situation as that of 
Theorem~\ref{thm: main gluing theorem}, and in particular use the
gluing setup
of \S\ref{sec:proof of main gluing}. For this purpose, we will make
use of the notation of that section, in particular using the
short-hand $\foM_1,\foM_2,\foM^{\gl},\foM^{\mu,\gl}$ as well
as the maps $p,p',p'',q',q'',\alpha,\alpha'$ and $\alpha''$ of the
diagram \eqref{eq:sub-diagram}.

\medskip

{\bf Step 1. Reduction to the gluing situation of 
\S\ref{sec:proof of main gluing}.}
\medskip

We first observe that if $\foM$ is an irreducible component of
$\foM^{\ddagger,\mu,\gl,\ev}(\shX,\beta_1,\beta_2,z)$ without a
terminal tail but with an interior tail,
then by Lemma~\ref{lem:first comparison facts}, (3), $i_*[\foM]$
is the class of an irreducible component of 
$\foM^{\mu,\gl,\ev}(\shX,\beta_1,\beta_2,z)$. This irreducible component
necessarily also has an internal tail. Note
$i_*' k_2^![\foM]=k_1^! i_* [\foM]$ by \cite{Man},~Thm.~4.1, using
the fact that $i$ is finite and representable by 
Lemma~\ref{lem:finite representable}, hence projective. Thus to show
the desired vanishing, we may assume $\foM$ is an irreducible component
of $\foM^{\mu,\gl,\ev}(\shX,\beta_1,\beta_2,z)$.

Note we have a factorization of $k_1$ given by
\[
\xymatrix@C=30pt
{
\scrM^{\mu,\gl}(X,\beta_1,\beta_2,z)\ar[r]^>>>>{p'\circ\alpha''}&
\foM^{\mu,\gl,\ev(q,x_{\out})}(\shX,\beta_1,\beta_2,z)\ar[r]&
\foM^{\mu,\gl,\ev}(\shX,\beta_1,\beta_2,z),
}
\]
with the second morphism smooth. Hence
we may pull back $\foM$ via the second morphism to obtain
an irreducible component of 
$\foM^{\mu,\gl,\ev(q,x_{\out})}(\shX,\beta_1,\beta_2,z)$.
Using the notation of \eqref{def:Mmu moduli}, we write this
irreducible component still as 
\[
\foM\subseteq \foM^{\mu,\gl,\ev(q,x_{\out})}(\shX,\beta_1,\beta_2,z)
= \foM^{\mu,\gl},
\]
hopefully with no confusion. It is then enough to show
as in \eqref{eq:Muvirt} that $\alpha^!p^![\foM]=0$ in the notation of the 
diagram \eqref{eq:sub-diagram}.
\medskip

{\bf Step 2. Vanishing if $q'|_{\foM}$ does not dominate an
irreducible component of $\foM_2$ or $s\not\in B(\ZZ)$.}
\medskip

Suppose that either (1) $q'(\foM)$ is not dense in an irreducible
component of $\foM_2$ or (2) $s\not\in B(\ZZ)$. We show that 
the desired vanishing holds. Indeed,
\[
\deg \alpha^!p^![\foM] = \deg \alpha^!(q'\circ p')_* (p')^![\foM]
\]
as in \eqref{eq:push-pull}, and
$(q'\circ p')_*(p')^![\foM]$ is supported on the closure of $q'(\foM)$
and hence in case (1) has dimension smaller than $\dim\foM_2$. Now the relative
virtual dimension of $\scrM_2\rightarrow\foM_2$ is 
$\chi((f_2^*\Theta_{X/\kk})(-x_{\out}-x_s'))=A_2\cdot c_1(\Theta_{X/\kk})
-\dim X$. As argued in Step 4 of the
proof of Theorem~\ref{thm: main gluing theorem}, if we are in the 
case that $\pm c_1(\Theta_{X/\kk})$ is nef, then either
one of $\scrM(X,\beta_1)$ or $\scrM(X,\beta_2)$ is empty or
$A_2\cdot c_1(\Theta_{X/\kk})=0$. If we are in the
log Calabi-Yau case, then again either one of
$\scrM(X,\beta_1)$, $\scrM(X,\beta_2)$ is empty or
$A_2\cdot c_1(\Theta_{X/\kk})\le 0$,
with strict inequality if $s\not\in B(\ZZ)$.
Thus in any event, the relative virtual dimension of $\scrM_2\rightarrow
\foM_2$ is at most $-\dim X$, and in case (2) is at most $-\dim X-1$. 
On the other hand, $\foM_2$ is pure-dimensional
of dimension $\dim X$. So we conclude that in any case,
$\alpha^!(q'\circ p')_* (p')^![\foM]$ is of negative dimension
and hence $\deg\alpha^!p^![\foM]=0$.
\medskip

{\bf Step 3. Reduction to the no-tail lemma.}
\medskip

Now assume that $\foM$ dominates the closure of an open stratum
$S_2$ of $\foM_2$ and $s\in B(\ZZ)$ in the log Calabi-Yau case,
so that in any case we may assume that 
$A_i\cdot c_1(\Theta_{X/\kk})=0$. 
We aim to show that $\deg \alpha^!(q'\circ p')_* (p')^![\foM]=0$ by showing
that $(q'\circ p')_* (p')^![\foM]=0$. 
To do so, consider
diagram \eqref{eq:M2 to S2} in our current situation. 
In that diagram
$S^{\mu}$ is an open substack of $\foM^{\mu,\gl}$, and as such determines
an open substack $\foM^{\circ}$ of $\foM$, with $|\foM^{\circ}|=
|\foM|\cap |S^{\mu}|$.
As in \eqref{eq:bunch of equalities}, it is then enough to show that
$(q'\circ p')_* p^![\foM^{\circ}]=0$. We remind the reader,
however, that in \eqref{eq:M2 to S2}, $S_2$ was the union of all open
strata of $\foM_2$, whereas here it is a single open stratum dominated
by $\foM$.

Consider now the diagram
\[
\xymatrix@C=30pt
{
\scrM'_1\times_{\foM^{\mu,\gl}} \foM^{\circ}\ar@{^{(}->}[r] \ar[d]_{p^{\circ}}&
\scrM'_1\times_{\foM^{\mu,\gl}} \foM\ar@{^{(}->}[r]^>>>>>{\iota} \ar[d]_{\bar p}&
\scrM_1'\ar[d]^{p'}\\
\foM^{\circ}\ar[d]_{q^{\circ}} \ar@{^{(}->}[r]&
\foM\ar[d]_{\bar q} \ar@{^{(}->}[r]&
\foM^{\mu,\gl}\ar[d]^{q'}\\
(S_2)_{\red}\ar@{^{(}->}[r]&
(\overline{S}_2)_{\red}\ar@{^{(}->}[r]^{\iota'}&\foM_2
}
\]
Here $(\overline S_2)_{\red}$ denotes the closure of $S_2$ with its
reduced induced stack structure.
All squares but the lower right-hand square are cartesian
in the fs category by definition, while the lower right-hand square
is merely commutative. The left-hand horizontal arrows are open embeddings,
while the right-hand horizontal arrows are closed embeddings. 
As $q'\circ p'$ is proper, so is $q'\circ p'\circ\iota=
\iota'\circ\bar q\circ \bar p$, and hence $\bar q\circ \bar p$ is
proper. By base-change $q^{\circ}\circ p^{\circ}$ is also proper.

Next we claim $q^{\circ}$ is flat. Indeed, it was shown in 
Step 5 of the proof of Theorem~\ref{thm: main gluing theorem} 
(shortly after \eqref{eq:bunch of equalities})
that $q'|_{S^{\mu}}:
S^{\mu}=\foM^{\mu,\gl}\times_{\foM_2} S_2\rightarrow S_2$ is flat.
Further, $\foM^{\mu,\gl}$ is pure-dimensional by Proposition
\ref{prop:transverse M2}, (1). Hence 
the fibres of $q'|_{S^{\mu}}$ are equi-dimensional. Thus the same holds for
$q^{\circ}$ as $\foM^{\circ}$ is an irreducible component
of $\foM^{\mu,\gl}\times_{\foM_2} S_2$.
Further, $(S_2)_{\red}$ is non-singular, 
and $\foM$ is Cohen-Macaulay.
Indeed, since $\Psi_{\mu}:\foM^{\mu,\gl}\rightarrow B\GG_m^{2,\mu}$
is log smooth, see \eqref{eq:Phimudef},
all (reduced) irreducible components
of $\foM^{\mu,\gl}$ are, locally in the smooth topology, isomorphic to
toric varieties and hence Cohen-Macaulay. Thus by
\cite{Mats},~Thm.~23.1, $q^{\circ}$ is flat.

We can now follow the last part of the argument of the proof
of Theorem~\ref{thm: main gluing theorem}. 
Let $s_2\rightarrow S_2$ be a choice of strict morphism from
a standard log point with underlying scheme $\Spec\kk$.
Similarly to \eqref{eq:S2 to p2},
replacing $S_2$ with its reduction
and writing $S_2^{\mu}=S_2\times_{\foM_2}\foM_2^{\mu}$ in place of
$(S_2^{\mu})_{\red} = (S_2)_{\red} \times_{\foM_2} \foM_2^{\mu}$
to simplify notation, we obtain the diagram:
\[
\xymatrix@C=30pt
{
(\scrM'_1\times_{\foM^{\mu,\gl}}\foM^{\circ})\times_{S_2} s_2
\ar[d]_{i_4}\ar[r]^>>>>>{p^{\circ}}&
\foM^{\circ}\times_{S_2} s_2 \ar[r]
\ar@/^2pc/[rr]^{q^{\circ}}
\ar[d]_{i_3}&
S_2^{\mu}\times_{S_2} s_2\ar[r]\ar[d]_{i_2}&s_2\ar[d]_{i_1}\\
\scrM'_1\times_{\foM^{\mu,\gl}}\foM^{\circ}
\ar[r]_{p^{\circ}}& \foM^{\circ}\ar[r]\ar@/_2pc/[rr]_{q^{\circ}}&S_2^{\mu}
\ar[r]&S_2
}
\]
with all squares cartesian in all categories. The same argument as in
\eqref{eq:bernd requested} and
\eqref{eq:a whole bunch of equalities} shows that 
\[
i_1^! (q^{\circ}\circ p^{\circ})_* (p^{\circ})^![\foM^{\circ}]
= (q^{\circ}\circ p^{\circ})_* (p^{\circ})^![\foM^{\circ}\times_{S_2} s_2].
\]
Similarly as in Claim~\ref{claim:claim},(3), proved in the subsequent
Step 5 of the proof of Theorem~\ref{thm: main gluing theorem}, we have
$(q^{\circ}\circ p^{\circ})_* (p^{\circ})^![\foM^{\circ}]$ a multiple of $[S_2]$.
Thus, as $i_1^![S_2]=[s_2]$,
it is sufficient to show, with
\begin{equation}
\label{eq:F}
F=\foM^{\circ}\times_{S_2} s_2,
\end{equation}
that 
\begin{equation}
\label{eq:tailvanishing}
(p^{\circ})^![F]=0.
\end{equation}
As we wish to apply Lemma~\ref{thm:no-tail-lemma}, we require
$\ul{F}$ integral. So to prove \eqref{eq:tailvanishing}, we 
replace $F$ with a reduced irreducible component of $F$.
\medskip

{\bf Step 4. Tropical analysis.}
\medskip

Consider the tropicalization $\Sigma(\pr_2):\Sigma(F)\rightarrow\Sigma(s_2)$
of the projection $\pr_2:F\rightarrow s_2$. As $i_1$ is
strict and $\delta:\foM_2\rightarrow B\GG_m^{\dagger}$ is log smooth, 
$s_2$ is equipped with the standard log point structure,
i.e., has ghost sheaf $\NN$.
From the induced strict morphisms 
\[
F\stackrel{i_3}{\longrightarrow}\foM^{\circ}\hooklongrightarrow\foM
\hooklongrightarrow
\foM^{\mu,\gl}\]
and 
\[
s_2\rightarrow \foM_2,
\]
$\Sigma(F)$ and $\Sigma(s_2)$ both carry families of
tropical maps, the first induced by a punctured map of class
$\beta$ and the second by a punctured map of class $\beta_2$.

Let $\bar\xi$ be the generic point of $F$. Note the stalk of the ghost
sheaf at $\bar\xi$ agrees with the stalk of the ghost sheaf of
$\foM^{\mu,\gl}$ at the generic point of $\foM$. Call this stalk
$Q$, necessarily of rank $2$, and write $\sigma_{\bar\xi}=
\Hom(Q,\RR_{\ge 0})$.
Note $\Sigma(\pr_2): \sigma_{\bar\xi}\rightarrow\sigma_{s_2}=\RR_{\ge 0}$ is
a surjection.

$\Sigma(\pr_2)$ can be described
as follows, bearing in mind that the map $\foM^{\mu,\gl}\rightarrow\foM_2$
is the second component of
the splitting map. A point $m\in \sigma_{\bar\xi}$ determines
a tropical map $h_m:G\rightarrow
\Sigma(X)$, while a point $m'\in \sigma_{s_2}$ determines a tropical
map $h_{m'}:G_2\rightarrow \Sigma(X)$. Here $G_2$
is obtained from $G$ by deleting a point from the interior
of the edge $E_q$, taking the connected
component which contains the edge $E_{\out}$, and replacing the segment
of $E_q$ with an unbounded leg, labelled $E_{x_s'}$. Then
$\Sigma(\pr_2)$ takes
$h_m$ to $h_m|_{G_2}$. This is a slight abuse of notation,
as the map $h_m|_{G_2}$ needs to be extended to the whole of 
$E_{x_s'}$, mapping $E_{x_s'}$ 
to a ray parallel to $\RR_{\ge 0}s$ (keeping in mind
that $s\in \Sigma(X)(\ZZ)$).

Similarly, we have a composition
\[
\psi:F\rightarrow \foM^{\mu,\gl}\rightarrow \foM_1,
\]
where the second morphism is the first component of the splitting
morphism.
Let $\bar\xi'$ be the image of $\bar\xi$
under the morphism $\psi$. Then
$\sigma_{\bar\xi'}$ parameterizes tropical maps of class $\beta_1$,
and $\Sigma(\psi):\sigma_{\bar\xi}\rightarrow \sigma_{\bar\xi'}$
can be described via $h_m\mapsto h_m|_{G_1}$, where
$G_1$ is obtained from $G$ in the same way as $G_2$,
but now taking the connected component containing $E_{x_1}$ and $E_{x_2}$.
This time, the edge $E_{x_s}$ becomes a bounded leg (unless $s=0$)
and is mapped to an edge in $\Sigma(X)$ parallel to $-\RR_{\ge 0}s$.
View $\delta$ as a coordinate on $\sigma_{s_2}$, and as a function
on $\sigma_{\bar\xi}$.
If $r\not=0$, then $h_m|_{G_2}$
varies in a one-dimensional family parameterized by $\delta$,
with $\delta(m)$ determined by \eqref{eq:deltanu}. On the
other hand, if $r=0$, then $h_{m'}:G_2\rightarrow \Sigma(X)$
is in fact independent of $m'$: this is similar to the discussion
of Remark~\ref{rem:delta basic}.
In any event, if we fix $\delta(m)=1$, $h_m|_{G_2}$ is fixed
but $h_m$ itself varies in a one-parameter family. In this family,
$\ell_q$ is unbounded, by construction of $\foM^{\mu,\gl}$. 
In particular, using the notation of 
Definition~\ref{def:common trop}, if $E_q=E_i$ for some $i$, then
$\ell_1,\ldots,\ell_{i-1}$, which yield lengths of edges of
$G_2$, are all proportional to $\delta$, but $\ell_i$
is not proportional to $\delta$. Thus $i$ is the unique index given by
Proposition~\ref{claim:unique}. Since $\foM$ has been assumed to have an
internal tail, it follows from the definition (Definition 
\ref{def:internal tails}) that there are at least two indices $j$
with $\ell_j$ not proportional to $\delta$. Thus there is at least
one index $j>i$ with $\ell_j$ not proportional to $\delta$. In particular,
this implies that $E_i$ is not the last edge in the chain 
$E_1,\ldots,E_{n-1}$, and $\ell_{i+1}+\cdots+\ell_{n-1}$
is also not proportional to $\delta$ (in particular non-zero).
Indeed, since $\delta$ rationally generates a one-dimensional
face of $Q$, a sum of elements of $Q$ is contained in this
face only if all these elements are contained in the face.

We can, however, say more. Let $m\in\sigma_{\bar\xi}$ be such that
$\delta(m)=0$ but $m\not=0$: such can be found using again that
$\rank Q=2$
and $\delta$ rationally generates a one-dimensional face of $Q$. Then
as $\ell_1,\ldots,\ell_{i-1}$ are proportional to $\delta$,
the tropical map $h_m:G\rightarrow \Sigma(X)$
must satisfy $h_m(v_j)=0$ for $j\le i$. In particular,
\begin{equation}
\label{eq:Sigmatgood}
h_m(v_{i+1})\in \RR_{\ge 0}s.
\end{equation}
As $\ell_{i+1}+\cdots+\ell_{n-1}$
is not proportional to $\delta$, necessarily
\begin{equation}
\label{eq:Sigma1tail}
\sum_{j=i+1}^{n-1} \ell_j(m)>0.
\end{equation}

\medskip

{\bf Step 5. Verification of the hypotheses of the no-tail lemma.}
\medskip

Now let $W=S_2^{\mu}\times_{S_2} s_2$. Note that $\ul{W}=B\GG_m$.
Indeed, the underlying stack $\ul{S}_2^{\mu}$ is isomorphic
to $\ul{S}_2\times B\GG_m$ by Lemma~\ref{lem: underlying} and
$s_2\rightarrow S_2$ is strict.
We have a morphism $g:W
\rightarrow \scrP(X,s)$ which is the composition of projection to
$S_2^{\mu}$, the inclusion into $\foM_2^{\mu}$, and the composed
morphism $\foM_2^{\mu}\rightarrow\scrP(X,s)$ in \eqref{eq:big diagram}.
With the first isomorphism by \eqref{eq:bunch of isomorphisms}
and the second isomorphism by associativity of fibre product, we have
\begin{equation}
\label{eq:yet another stupid identity}
W\times^{\fs}_{\scrP(X,s)}\foM_1\cong 
(S_2^{\mu}\times_{S_2} s_2)\times_{\foM^{\mu}_2}\foM^{\mu,\gl}\cong
(S_2^{\mu}
\times_{\foM_2^{\mu}}\foM^{\mu,\gl})
\times_{S_2} s_2.
\end{equation}

We have the diagram
\[
\xymatrix@C=30pt
{
\foM^{\circ} \ar@{^{(}->}[d] \ar[r] & 
S_2^{\mu}\times_{\foM_2^{\mu}}\foM^{\mu,\gl}
\ar@{^{(}->}[d]\ar[r]&S_2^{\mu}\ar[r]\ar@{^{(}->}[d]&S_2\ar@{^{(}->}[d]\\
\foM\ar[r]&\foM^{\mu,\gl}\ar[r]&\foM_2^{\mu}\ar[r]&\foM_2
}
\]
The morphism $\foM^{\circ}\rightarrow 
S_2^{\mu}\times_{\foM_2^{\mu}}\foM^{\mu,\gl}$ along with the
identity map $s_2\rightarrow s_2$ then induces a morphism
$\foM^{\circ}\times_{S_2} s_2\rightarrow 
W\times^{\fs}_{\scrP(X,s)} \foM_1$
via the description \eqref{eq:yet another stupid identity}. From
the description of $F$ as an irreducible component of \eqref{eq:F},
we thus obtain a morphism
\[
F \rightarrow
W\times^{\fs}_{\scrP(X,s)} \foM_1
\]
such that the composition with the projection to $\foM_1$ is given by splitting
the family of curves over $F$ at the node $q$ and taking the connected
component of the domain of class $\beta_1$.

We now need to find a stratum $S^{\circ}\subseteq\foM_1$,
not necessarily open, such that the
data $\psi:F\rightarrow \foM_1$, 
$W$, $g$ and $S^{\circ}$ satisfy the hypotheses of 
Lemma~\ref{thm:no-tail-lemma}, applied with $\foM_1$ playing
the role of $\foM^{\ev}(\shX,
\beta)$ in that lemma.

As in Step 4, with $\bar\xi$ still the generic point of $F$,
let $m\in \sigma_{\bar\xi}$
be a point with $\delta(m)=0$, $m\not=0$, and let $m'=\Sigma(\psi)(m)$,
so that $m'$ corresponds to a tropical map $h_m|_{G_1}$.
Recalling the morphism $\ev_{\shX}:\foM_1\rightarrow\scrP(X,s)$
tropicalizes to evaluation at the vertex of $G_1$ containing the
puncture $x_s$, we have by \eqref{eq:Sigmatgood} that $\Sigma(\ev_{\shX})(m')
\in \RR_{\ge 0}s$. Let $\tau$ be the minimal face of $\sigma_{\bar\xi'}$
containing $m'$, with $\bar\xi'=\psi(\xi)$. Then $\Sigma(\ev_{\shX})(\tau)$ is
contained in the minimal cone of $\Sigma(X)$ containing $s$.
As $\ev_{\shX}$ is log smooth, 
Proposition~\ref{prop: realisability} implies
there exists a stratum $S^{\circ}$ of $\foM_1$ whose closure contains
$\bar\xi'$, hence $\psi(F)$, such that 
the cone of $\Sigma(\foM_1)$ corresponding to $S^{\circ}$ is $\tau$.
Further, we then have $\ev_{\shX}(S^{\circ})\subseteq \ul{Z}_s^{\circ}\times
B\GG_m$. 

We can now show the five conditions necessary to apply
Lemma~\ref{thm:no-tail-lemma}. Conditions (1) and (2) are satisfied
immediately by construction. On the other hand, by
\eqref{eq:Sigma1tail}, $h_{m'}$ does not contract all
edges $E_{i+1},\ldots,E_{n-1}$. Thus if $\phi:\foM_1\rightarrow \Mbf_{0,3}$ is
the forgetful map just remembering the domain, then $\phi(S)
\subseteq D(x_1,x_2\,|\, x_s)$, hence (3). Condition (4) again
follows by construction. 

For condition (5), the first transversality
statement is as in the last paragraph of the proof of Theorem 
\ref{thm: main gluing theorem}. We continue to use the notation 
$\delta', \nu$ of \eqref{eq:def nu}, with $\delta=\nu\delta'$ and
$\delta'$ the generator of the ghost sheaf monoid on $S_2$.

The second transversality comes as follows. By 
Remark~\ref{rem:trop of gluing},
the morphism
$s_2\times B\GG_m^{\dagger}\rightarrow\scrP(X,s)$ tropicalizes
to $\RR_{\ge 0}(\delta')^*\oplus\RR_{\ge 0}\ell_q^*
\rightarrow \sigma\in\Sigma(X)$ given by $a(\delta')^*+b\ell_q^*\mapsto
a \rho+b s$ for some $\rho\in \sigma$. Let $\tau_{\mu}$
be the image of the map $\Sigma(g):\Sigma(W)\rightarrow 
\Sigma(\scrP(X,s))=\Sigma(Z_s)$.
Just as in the last paragraph of Step 5 of the proof of
Theorem~\ref{thm: main gluing theorem},
the tropicalization of $W$ is the
cone generated by $(\delta')^*+\mu\nu\ell_q^*$ and $\ell_q^*$, 
so $\tau_{\mu}$ is the cone in $\sigma$ spanned by
$\rho+\mu\nu s$ and $s$.
Given any other subcone $\sigma'$ of
$\sigma$, we may thus choose $\mu$ sufficiently large so that
$\sigma'\cap \tau_{\mu}$ is a (not necessarily proper) face of
$\tau_{\mu}$.
Thus by taking $\mu$ sufficiently 
large, given any finite collection $\Sigma$ of subcones of 
$\sigma$, we can assume $\tau_{\mu}$ intersects every subcone in this 
collection in a face of $\tau_{\mu}$.\footnote{Note that if $s=0$,
then in fact $\tau_{\mu}$ is a one-dimensional cone, and this
statement then holds for all $\mu$.} 

We may now apply this observation to transversality as follows.
First note that there are only a finite
number of splittings $A_1=A_{11}+A_{12}$
into effective curve classes, and for each choice of $A_{11}$
there are only a finite number of class $\beta_{11}$ of punctured maps
with two punctures $x_{\out}$ and $x$ with $x_{\out}$ of
contact order $-s$. Indeed, the contact order at $x$ is then determined
up to a finite number of choices by Remark~\ref{rem:contact order determined}.
Thus after replacing moduli spaces with finite type approximations,
there are only a finite number of finite type moduli spaces 
$\foM^{\ev}(\shX,\beta_{11})$
for which transversality of $g$ to 
$\ev_{\shX}:\foM^{\ev}(\shX,\beta_{11})\rightarrow\scrP(X,s)$ needs
to be checked. Take $\Sigma$ of the previous paragraph
to be the set of images under $\Sigma(\ev_{\shX})$ of all cones
in $\Sigma(\foM^{\ev}(\shX,\beta_{11}))$ for this finite number
of possible classes $\beta_{11}$. We then choose $\mu$ sufficiently
large as in the previous paragraph.
Then by the transversality condition 
Theorem~\ref{theorem: transversality}, we can guarantee that
$g$ is transverse to all possible $\ev_{\shX}$.

We have now shown that the
hypotheses of Lemma~\ref{thm:no-tail-lemma} hold and we are able
to conclude that $(p^{\circ})^![F]=0$ in the rational Chow group.
As we saw in Step 3, this was sufficient to obtain the desired
$k_1^![\foM]=0$.
\end{proof}

\begin{remark}
\label{rem:wherefore tails}
The reader may wonder why in Lemma~\ref{lemma:terminal tail vanishing}
we showed if $\foM$ was an irreducible component of
of $\foM_y^{\ddagger,\ev}$ with a terminal tail, then $k_3^![\foM]=0$,
but in Lemma~\ref{lemma:internal tail vanishing}, we did not show this
same vanishing if $\foM$ instead had an internal tail.
The reason is that the argument given in the latter lemma relies
subtly on the construction of $\foM^{\mu,\gl,\ev}(\shX,
\beta_1,\beta_2,z)$, and we do not yet know that every irreducible
component of $\foM_y^{\ddagger,\ev}$ with an internal tail is
the image of an irreducible component of
$\foM^{\ddagger,\mu,\gl,\ev}(\shX,\beta_1,\beta_2,z)$. In the case
that the irreducible component of $\foM^{\ddagger,\ev}_y$ does not
have a terminal tail but does have an internal tail, this will be shown
in the course of the proof of Theorem~\ref{thm:second comparison},
see the last paragraph of Step 1 of that proof.
\end{remark}

\subsection{The first comparison}
We now are ready to make the comparisons of virtual degrees
using the morphisms $i$ and $j$ in \eqref{eq:comparison diagram2}.
In this subsection we make the comparison using $i$. This is
the first time internal tails cause difficulties for us:

\begin{lemma}
\label{lem:same multiplicities}
Suppose $\mu$ is chosen sufficiently large so that the
transversality statement of Proposition 
\ref{prop:transverse M2}, (3) holds. For $\lambda\gg\mu$,
the following holds.
Let $\foM'$ be a tail-free irreducible component of 
$\foM^{\ddagger,\mu,\gl,\ev}(\shX,\beta_1,\beta_2,z)$.
Then the multiplicities of $\foM=i(\foM')$ and $\foM'$ in 
$\foM^{\mu,\gl,\ev}(\shX,\beta_1,\beta_2,z)$ and
$\foM^{\ddagger,\mu,\gl,\ev}(\shX,\beta_1,\beta_2,z)$ respectively are the same.
\end{lemma}

\begin{proof}
Following the notation of Definition~\ref{def:common trop},
we have $\ell=\sum_i \ell_i$ and $\ell_q=\ell_i$
for some $i$, and $\ell_q$ and $\delta$ are linearly independent.
It follows from Definition~\ref{def:internal tails} that each
$\ell_j$ for $j\not=i$ is proportional to $\delta$, and hence
$\ell=\ell_q+a\delta$ for some $a\in\QQ$. Note that the
coefficient of $\ell_q$ is $1$, so that the result then follows immediately
from Lemma~\ref{lem:local product} and \eqref{eq:Sdiagram}, with
$S$ in the latter diagram taken to be the open stratum of $\foM$.
\end{proof}

\begin{theorem}
\label{thm:first comparison}
$\deg [\scrM^{\mu,\gl}(X,\beta_1,\beta_2,z)]^{\virt}
=\deg [\scrM^{\ddagger,\mu,\gl}(X,\beta_1,\beta_2,z)]^{\virt}$.
\end{theorem}

\begin{proof}
Lemma~\ref{lem:finite representable} implies $i'$ is finite, and hence
we can apply \cite{Man},~Thm.~4.1 
to see that
\begin{align*}
\deg [\scrM^{\ddagger,\mu,\gl}(X,\beta_1,\beta_2,z)]^{\virt}
= {} &
\deg i'_* [\scrM^{\ddagger,\mu,\gl}(X,\beta_1,\beta_2,z)]^{\virt}
\\
= {} &
\deg k_1^! i_*[\foM^{\ddagger,\mu,\gl,\ev}(\shX,\beta_1,\beta_2,z)].
\end{align*}
Now $i_*[\foM^{\ddagger,\mu,\gl,\ev}(\shX,\beta_1,\beta_2,z)]$
and $[\foM^{\mu,\gl,\ev}(\shX,\beta_1,\beta_2,z)]$
are both linear combinations of irreducible components of
$\foM^{\mu,\gl,\ev}(\shX,\beta_1,\beta_2,z)$. The coefficients
of those irreducible components without tails agree by Lemma 
\ref{lem:same multiplicities}, while those irreducible components
with tails contribute $0$ to the virtual fundamental class, 
by Lemmas~\ref{lemma:terminal tail vanishing} and 
\ref{lemma:internal tail vanishing}, hence the result.
\end{proof}

\subsection{The second comparison}

The following, along with Theorem~\ref{thm:first comparison}, immediately
proves Theorem~\ref{thm:main comparison}. 

\begin{theorem}
\label{thm:second comparison}
\[
\sum_{A_1,A_2,s} \deg[\scrM^{\ddagger,\mu,\gl}(X,\beta_1,\beta_2,z)]^{\virt}
=\deg [\scrM^{\ddagger}_y]^{\virt}
\]
\end{theorem}

\begin{proof}
{\bf Step 1. Sketch.}
\medskip

Fix some choice of $s\in \Sigma(X)(\ZZ)$, leading to classes $\beta_1$,
$\beta_2$ of punctured maps to $\shX$ as usual. Recall that 
types or classes
of punctured maps to $\shX$ do not remember curve classes in 
$X$. By Lemmas~\ref{lem:first comparison facts}, (1)
and \ref{lemma: tau VFC}, (1), provided that one chooses $\lambda$ 
sufficiently large after making a choice of $\mu$ sufficiently large,
$\foM^{\ddagger,\mu,\gl,\ev}(\shX,\beta_1,\beta_2,z)$ 
and $\foM_y^{\ddagger,\ev}$
are both pure dimension zero. As the morphism $j$ is finite
and representable by
Lemma~\ref{lem:finite representable}, if $\foM$ is an irreducible
component of $\foM^{\ddagger,\mu,\gl,\ev}(\shX,\beta_1,\beta_2,z)$,
then $j(\foM)$ is an irreducible component of $\foM^{\ddagger,\ev}_y$.
If $\foM$ has a tail, either terminal or internal, then so does
$j(\foM)$. Thus in this case,
\begin{equation}
\label{eq:boring sequence of equalities}
\deg k_3^!j_*[\foM]=\deg j'_*k_2^![\foM]= \deg k_2^![\foM]=0
\end{equation}
by Lemmas \ref{lemma:terminal tail vanishing}, 
\ref{lemma:internal tail vanishing}, and push-pull for virtual pull-back,
\cite{Man}, Thm.~4.1.

Define $\foU'\subseteq 
\coprod_s \foM^{\ddagger,\mu,\gl,\ev}(\shX,\beta_1,\beta_2,z)$
and $\foU\subseteq\foM^{\ddagger,\ev}_y$ to be the open substacks which are the
unions of open strata which do not have a terminal tail. Also, let
$\foU'_{\tf}\subseteq\foU'$ and $\foU_{\tf}\subseteq\foU$ be the union
of those connected components of $\foU'$, $\foU$ respectively which also
do not have internal tails. Here $\tf$ stands for ``tail-free.'' Because
$j$ is finite and log morphisms take logarithmic strata into 
logarithmic strata, necessarily $j(\foU')\subseteq\foU$ and
$j(\foU'_{\tf})\subseteq \foU_{\tf}$. We write
$\tilde\jmath:\foU'\rightarrow\foU$ for $j|_{\foU'}$. 

We will construct a morphism
$\spl^{\circ}:\foU\rightarrow \foU'$ which is a right inverse to $\tilde\jmath$,
so that $\tilde\jmath$ is surjective. In particular, any irreducible component
$\foM$ of $\foM_y^{\ddagger,\ev}$ with an internal tail but not
a terminal tail is the image of a component
$\foM'$ of $\foM^{\ddagger,\mu,\gl,\ev}(\shX,\beta_1,\beta_2,z)$
with an internal tail. (Note this proves the claim of
Remark~\ref{rem:wherefore tails}). Hence $\deg k_3^![\foM]=0$, by
\eqref{eq:boring sequence of equalities}.
On the other hand, any irreducible component $\foM$ of $\foM_y^{\ddagger,\ev}$
with a terminal tail satisfies $\deg k_3^![\foM]=0$ by Lemma
\ref{lemma:terminal tail vanishing}.

Thus it will be sufficient to show that $\tilde\jmath|_{\foU'_{\tf}}:
\foU'_{\tf}\rightarrow \foU_{\tf}$ is an isomorphism, which we do by
showing that $\spl^{\circ}|_{\foU_{\tf}}$ is its inverse. 
\medskip

{\bf Step 2. The splitting morphism.} 
\medskip

We construct a splitting morphism
$\spl'':\foU\rightarrow \coprod_s \foM^{\gl,\ev}(\shX,\beta_1,\beta_2,z)$.
The key point is showing that the universal family of punctured maps 
$f:\foC/\foU
\rightarrow \shX$ can be split at some node. For this step, replace $\foU$
by a connected component with generic point $\bar\xi$, so that we get a
family of tropical maps $h_m:G\rightarrow\Sigma(X)$ parameterized
by $m\in\sigma_{\bar \xi}$. 
By the definition of $\foU$, this family of maps does not have a terminal tail. 
Let $i$ be the unique index guaranteed
by Proposition~\ref{claim:unique}. Let $Q$ be the stalk of the ghost sheaf
at $\bar\xi$. 
Recall \eqref{eq:R Rlambda def} that the stalk of $B\GG_m^{\ddagger}$ is 
$R_{\lambda}$, generated
by $\ell-\lambda\delta$ and $\delta$.
Then $\delta$ and $\ell_i$ are linearly independent in
$Q^{\gp}$ and $Q$ is rationally generated by $\ell-\lambda\delta, \delta$,
as $\foM^{\ddagger,\ev}_y\rightarrow B\GG_m^{\ddagger}$ is integral by construction.
We can write 
\begin{equation}
\label{eq:another identity we decided to label}
\hbox{$\ell=a\ell_i+b\delta$ for some $a,b\in \QQ$,}
\end{equation}
with $a>0$. Note then that 
\[
\ell_i-\delta={1\over a}(\ell-\lambda\delta)+\left({\lambda-b\over a}-1\right)
\delta.
\]
Thus provided that $\lambda$ is taken sufficiently large so that
$(\lambda-b)/a-1\ge 0$, i.e., $\lambda\ge a+b$, we see that $\ell_i-\delta
\in Q$ and hence $\ell_i$ is in the ideal generated by $\delta$. Now
any section of $\shM_{B\GG_m^{\ddagger}}$ whose image in
$\overline{\shM}_{B\GG_m^{\ddagger}}$ is
$\delta$ has image $0$ under the
structure map $\shM_{B\GG_m^{\ddagger}}\rightarrow \O_{B\GG_m^{\ddagger}}$.
Thus the same is true of any section
of $\shM_{\foU}$ defined in a smooth neighbourhood of
$\bar\xi$ and whose image in $\overline{\shM}_{\foU}$ is $\delta$.
By the observation that $\ell_i$ is in the ideal generated by 
$\delta$, the same vanishing of structure map holds for any
section of $\shM_{\foU}$ lifting $\ell_i$.

Observing that $\ell_i\in Q$ is the smoothing parameter for the
node $q$ corresponding to the edge $E_i$,
we may now apply \cite{ACGS18}, Prop.~5.2, to split the punctured map
$f:\foC/\foU\rightarrow\shX$ at the node $q$,
obtaining two punctured
maps $f_1:\foC_1/\foU\rightarrow\shX$, $f_2:\foC_2/\foU\rightarrow
\shX$ of classes $\beta_1$ and $\beta_2$. Here each curve has three
punctured points, with $\beta_1$ having contact orders $p_1, p_2$ and 
$-s=-u_i$, with $u_i$ the contact order of the edge $E_i$, 
and $\beta_2$ having contact orders $s,p_3$ and $r$.
Note that $s\in \Sigma(X)(\ZZ)$ by Proposition 
\ref{claim:unique}, (3).
Thus we get from the universal property of the glued moduli space
a tautological map $\foU\rightarrow\foM^{\gl,\ev}
(\shX,\beta_1,\beta_2,z)$. This defines $\spl''$.
\medskip

{\bf Step 3. Construction of $\spl:\foU\rightarrow
\coprod_s\foM^{\ddagger,\mu,\gl,\ev}(\shX,\beta_1,\beta_2,z)$.}
\medskip

The diagram of Figure~\ref{Figure71} summarizes the current situation,
with all solid arrows already existing.
\begin{figure}
\[
\xymatrix@C=30pt
{&&&&
\\
&\foU'
\ar@{^{(}->}_{\iota'}[dd] 
\ar@/^1pc/^{\tilde\jmath}[rr]
& & \foU
\ar`u[ul]`[llld]`[ddddl]_{\spl''}[ddddll]
%\ar`r[rd]`[dddddr]^{\eta}[dddddl]
\ar@{-->}@/^8pc/^{\eta}[dddddl]
\ar@{_{(}->}^{\iota}[dd] 
\ar@{-->}@/^1pc/^{\spl^{\circ}}[ll]&
\\
&&&&\\
&\coprod_s\foM^{\ddagger,\mu,\gl,\ev}\ar[r]_{\pi_1}\ar[d]_{\pi_3}
\ar@/^2pc/^{j}[rr]
& B\GG_m^{\ddagger}\ar[d]_{\psi_y} & \foM_y^{\ddagger,\ev}\ar[l]^{\pi_2}\ar[d]^{\pi_4}&\\
&\coprod_s\foM^{\mu,\gl,\ev}\ar[r]^{\Phi'}\ar[ddr]^{\pi_6}\ar[d]_{\pi_5} 
\ar@/_2pc/_{\bar\jmath'}[rr]&
\overline\scrM_{0,4}^{\dagger}& \foM^{\ev}(\shX,\beta,z)\ar[l]_{\Phi}&\\
&\coprod_s\foM^{\gl,\ev}\ar[rdd]_{\Psi}
\ar@/_2.5pc/_{\bar\jmath}[rru]
&&&\\
&&B\GG_m^{2,\mu}\ar[d]^{\epsilon}&&\\
&&B\GG_m^{\dagger}\times B\GG_m^{\dagger}&&\\
&&&&
}
\]
\caption{}
\label{Figure71}
\end{figure}
Here we use short-hand $\foM^{\gl,\ev}$ for $\foM^{\gl,\ev}(\shX,\beta_1,
\beta_2,z)$ etc. There are three fs log cartesian squares, i.e., the
two squares defining $\coprod_s \foM^{\ddagger,\mu,\gl,\ev}$ and
$\foM_y^{\ddagger,\ev}$ and the lower parallelogram defining $\coprod_s 
\foM^{\mu,\gl,\ev}$. The morphism $\tilde\jmath$ was defined
in Step 1 and the morphism $\spl''$ in Step 2. The morphisms
$\bar\jmath$ and $\bar\jmath'$ are the tautological morphisms,
and \emph{except for any commutativity involving $\spl''$},
the diagram is commutative. In particular, 
$\Phi'=\Phi\circ\bar\jmath'$ by construction of $\Phi'$, and $j$ is
obtained from $\bar\jmath'$ by applying the functor
$\times^{\fs}_{\overline{\scrM}_{0,4}^{\dagger}} B\GG_m^{\ddagger}$.
The morphisms $\iota$ and $\iota'$ are the inclusions. At this point,
we do not make any statements about commutativity involving dashed arrows.

In fact, we have a somewhat restricted commutativity involving
$\spl''$, namely,
\begin{equation}
\label{eq:splddcomm}
\spl''\circ \tilde\jmath|_{\foU'_{\tf}} = \pi_5\circ\pi_3\circ
\iota'|_{\foU'_{\tf}}.
\end{equation}
Indeed, let $\foV'\subseteq\foU'_{\tf}$ be a connected component with
universal punctured map $f':\foC'/\foV'\rightarrow\shX$, and with node $q'$
obtained via the gluing. Let $\foV\subseteq\foU_{\tf}$ be the
connected component with $\tilde\jmath(\foV')\subseteq\foV$, with universal
punctured map $f:\foC/\foV\rightarrow\shX$, with node $q$ as chosen in Step 2.
Note that $f'$ is the pull-back of the family of punctured maps $f$ via 
$\tilde\jmath$. Thus if $G$ is the dual graph of $\foC'$, it has
two edges $E_q$, $E_{q'}$ corresponding to these two nodes, and neither
$\ell_q$ nor $\ell_{q'}$ is proportional to $\delta$, the first by
our choice of $q$ using
Proposition~\ref{claim:unique} and the second because $\ell_{q'}$
and $\delta$ rationally generate the group of the stalk of the ghost
sheaf at any point of $\foV'$ by construction.
Thus both edges are of splitting type, and so by the absence of 
internal tails on $\foU_{\tf}$,
$q=q'$. Now the identity $\spl''\circ
\tilde\jmath|_{\foV'} = \pi_5\circ\pi_3\circ\iota'|_{\foV'}$ 
follows from the universal
property of the glued moduli space $\foM^{\gl,\ev}(\shX,\beta_1,\beta_2,z)$.

To define $\eta$, we note that $\Psi\circ\spl''$ factors through
$\epsilon$. Indeed, it is sufficient to check this factorization
at the level of ghost sheaves. The morphism $\Psi\circ\spl''$
is given, at the level
of stalks of ghost sheaves, by $\delta\mapsto \delta$, $\ell_q\mapsto
\ell_i$. Thus $\ell_q-\mu\delta\mapsto \ell_i-\mu\delta\in Q^{\gp}$, 
and we need to ensure that $\ell_i-\mu\delta\in Q$. But recalling
\eqref{eq:another identity we decided to label}, we have
\[
\ell_i-\mu\delta={1\over a}(\ell-\lambda\delta)+\left({\lambda -b\over a}
-\mu\right)\delta,
\]
so provided $\lambda\ge b+\mu a$, we obtain the desired factorization.
Thus we can write $\Psi\circ\spl''=\epsilon\circ\eta$ for a uniquely
determined morphism $\eta$.

The morphisms $\spl''$, $\eta$ give, by the universal property of
fibre product, a morphism 
\[
\spl':\foU\rightarrow\coprod_s \foM^{\mu,\gl,\ev}(\shX,\beta_1,\beta_2,z).
\]
Now  
\[
\Phi'\circ\spl'=\Phi\circ\bar\jmath\circ\pi_5\circ\spl'=\Phi\circ\bar
\jmath \circ\spl''=\Phi\circ\pi_4\circ\iota=\psi_y\circ\pi_2\circ\iota.
\]
Here the third equality holds as 
\begin{equation}
\label{eq:jspl equals pi4iota}
\bar\jmath\circ\spl''=
\pi_4\circ\iota
\end{equation}
are both the tautological morphism
$\foU\rightarrow \foM^{\ev}(\shX,\beta,z)$.
Thus $\spl'$ and $\pi_2\circ\iota$ induce, again by the universal property
of fibre product, a morphism
\[
\spl:\foU\rightarrow \coprod_s \foM^{\ddagger,\mu,\gl,\ev}(\shX,\beta_1,
\beta_2,z).
\]

We pause here to comment on the overall choices of $\mu$ and $\lambda$
for this construction. Having chosen the punctured map class $\beta$,
we make one choice of $\lambda$, as we need $\lambda$ to
define $\foM_y^{\ddagger,\ev}$. This choice has to satisfy a number of 
properties, as given in Theorem~\ref{theorem: independence of modulus}
and Lemma~\ref{lem:local product}, the latter used in the proof of
Lemma~\ref{lem:first comparison facts}, (3). Note the latter choice depends 
on $\mu$. On the other hand, the choice of $\mu$ has to satisfy
a number of properties, given in Proposition~\ref{prop:transverse M2},
(3) and Lemma~\ref{lemma:internal tail vanishing}. However, 
the choice of $\mu$ may vary as $s$ ranges over all points of
$\Sigma(X)(\ZZ)$, which may seem problematic.
But, as
$\scrM^{\ddagger,\mu,\gl}(X,\beta_1,\beta_2,z)$ must be empty
for all but a finite number of $s$ and curve classes $A_1$
and $A_2$ (as $\scrM(X,\beta,z)$ is finite type) we may
in fact shrink all but a finite number of the
$\foM^{\ddagger,\mu,\gl,\ev}(\shX,\beta_1,\beta_2,z)$ to the empty set.
We then find that the parts of the proof imposing constraints on
$\lambda$ and $\mu$ only impose a finite number of constraints.
Thus we may first choose $\mu$ which satisfies all necessary constraints
given the choice of $\beta$, which is fixed, and then similarly choose 
$\lambda\gg\mu$ given the choice of $\mu$.

\medskip

{\bf Step 4. $j\circ\spl=\iota$ and construction of $\spl^{\circ}:
\foU\rightarrow\foU'$.}

\medskip

Indeed, to show $j\circ\spl=\iota$, it is sufficient, by the universal property of
fibre product, to show that
\[
\pi_2\circ j\circ\spl=\pi_2\circ\iota, \quad
\pi_4\circ j \circ \spl = \pi_4 \circ\iota.
\]
For the first equality,
\[
\pi_2\circ j\circ\spl=\pi_1\circ\spl=\pi_2\circ\iota,
\]
(the last equality by the construction of $\spl$)
while for the second,
\[
\pi_4\circ j\circ\spl=\bar\jmath\circ\pi_5\circ\pi_3\circ\spl
=\bar\jmath\circ\pi_5\circ\spl'=\bar\jmath\circ\spl''=\pi_4\circ\iota,
\]
the last equality by \eqref{eq:jspl equals pi4iota}.

We thus see in particular that $\spl$ is injective, and hence the image
of $\spl$ is top-dimensional. Since log morphisms take log strata into
log strata, necessarily $\spl(\foU)\subseteq\foU'$, showing we can
write $\spl=\iota'\circ\spl^{\circ}$ for a morphism $\spl^{\circ}:\foU
\rightarrow\foU'$. By restriction we obtain a morphism
$\spl^{\circ}:\foU_{\tf}\rightarrow \foU'_{\tf}$.
\medskip

{\bf Step 5. $\spl^{\circ}:\foU_{\tf}\rightarrow\foU'_{\tf}$ is inverse to
$\tilde\jmath:\foU'_{\tf}\rightarrow\foU_{\tf}$.}

\medskip

By Step 4, $\iota=j\circ\iota'\circ\spl^{\circ}=\iota\circ\tilde\jmath\circ\spl^{\circ}$,
the latter by commutativity properties of Figure~\ref{Figure71}.
The same remains true after restricting to $\foU_{\tf}$, and thus,
restricting to $\foU_{\tf}$, we have
$\tilde\jmath\circ\spl^{\circ}=\id_{\foU_{\tf}}$.
To show $\spl^{\circ}\circ\tilde\jmath=\id_{\foU'_{\tf}}$, we need to
show
\[
\pi_1\circ\iota'\circ\spl^{\circ}\circ\tilde\jmath=\pi_1\circ\iota',
\quad
\pi_3\circ\iota'\circ\spl^{\circ}\circ\tilde\jmath=\pi_3\circ\iota'.
\]

For the first,
\[
\pi_1\circ\iota'\circ\spl^{\circ}\circ\tilde\jmath
=\pi_1\circ\spl\circ\tilde\jmath=\pi_2\circ\iota\circ\tilde\jmath
=\pi_2\circ j\circ \iota'=\pi_1\circ \iota'.
\]

For the second, we need to check that
\begin{equation}
\label{eq:yabba dabba do}
\pi_5\circ\pi_3\circ\iota'\circ\spl^{\circ}\circ\tilde\jmath=\pi_5\circ\pi_3\circ\iota',
\quad
\pi_6\circ\pi_3\circ\iota'\circ\spl^{\circ}\circ\tilde\jmath=\pi_6\circ\pi_3\circ\iota'.
\end{equation}
For the first equality of \eqref{eq:yabba dabba do}, we have
\[
\pi_5\circ\pi_3\circ\iota'\circ\spl^{\circ}\circ\tilde\jmath=
\pi_5\circ\spl'\circ\tilde\jmath=\spl''\circ\tilde\jmath
=\pi_5\circ\pi_3\circ\iota',
\]
the last equality following from \eqref{eq:splddcomm}.

For the second equality of \eqref{eq:yabba dabba do}, 
we note that $\epsilon:B\GG_m^{2,\mu}\rightarrow 
B\GG_m^{\dagger}\times B\GG_m^{\dagger}$
is a monomorphism in the category of 
fs log stacks, since if $Y$ is any fs log stack, a homomorphism $T_{\mu}
\rightarrow \Gamma(Y,\overline{\shM}_Y)$ is determined by the composition
with $T\rightarrow T_{\mu}$. Thus we may check this equality after
composing on the left by $\epsilon$, and
\[
\epsilon\circ\pi_6\circ\pi_3\circ\iota'\circ\spl^{\circ}\circ\tilde\jmath=
\epsilon\circ\pi_6\circ\spl'\circ\tilde\jmath=\epsilon\circ\eta\circ\tilde
\jmath=\Psi\circ\spl''\circ\tilde\jmath
= \Psi\circ\pi_5\circ\pi_3\circ \iota'=\epsilon\circ\pi_6\circ\pi_3\circ
\iota',
\]
showing the desired identity. Here the
second equality follows from the construction of $\spl'$ and the
fourth equality again follows from \eqref{eq:splddcomm}.

We have thus shown that $\foU_{\tf}$ and $\foU'_{\tf}$
are isomorphic. As explained in Step 1, this complete the proof
of the theorem.
\end{proof}

\section{Remarks on the Frobenius structure conjecture}

In the first arXiv version of \cite{GHKI}, a conjecture was made concerning
the coefficient of $\vartheta_0$ in products of theta functions
$\prod_{i=1}^n \vartheta_{p_i}$. In the language of this paper, we can
make precise the conjecture as follows. We assume we are
in the absolute log Calabi-Yau case.
Fix $p_1,\ldots,p_n\in B(\ZZ)$ and a curve class $A$ in $X$.
Let $\beta$ be the class of punctured map with curve class $A$,
genus $0$, and $n+1$ marked points, $x_1,\ldots,x_n,x_{\out}$, with
contact orders determined by $p_1,\ldots,p_n,0$ respectively. Fix a general
point $z\in X$, to obtain a moduli space $\scrM(X,\beta,z)$,
imposing the constraint $z$ at $x_{\out}$, of virtual
dimension $n-2$, by Proposition~\ref{prop:the invariants}. We then set
\[
N^A_{p_1\ldots p_n0} := \int_{[\scrM(X,\beta,z)]^{\virt}}
\psi_{x_{\out}}^{n-2}.
\]
Here $\psi_{x_{\out}}$ denotes the $\psi$ class associated to the point
$x_{\out}$, i.e., the first Chern class of the cotangent line at $x_{\out}$.
Note this invariant only involves ``classical'' log Gromov-Witten invariants,
i.e., no punctures are involved.

\begin{remark}
It may be more appropriate to rephrase this conjecture in terms 
of ancestor invariants instead of descendent invariants, i.e.,
to take $\psi_{x_{\out}}$ to be the $\psi$ class associated to
$x_{\out}$ pulled back from $\overline\scrM_{0,n+1}$. Travis Mandel
has sketched an argument for us suggesting these in fact agree here.
In any event, the proof of the theorem below in fact shows the equality
of the two invariants for $n=3$.
\end{remark}

The Frobenius structure conjecture can then be partially rephrased 
as 

\begin{conjecture}
\label{conj:frob}
The coefficient of $\vartheta_0$ in the product $\vartheta_{p_1}\cdots
\vartheta_{p_n}$ is
\[
\sum_{A} N^A_{p_1\ldots p_n0} t^A.
\]
\end{conjecture}

\begin{remark}
In fact, the philosophy of the Frobenius structure conjecture is slightly
different. Let us call the above conjecture the \emph{weak Frobenius
structure conjecture}. The \emph{strong Frobenius structure conjecture} then
predicts that there is a unique
choice of structure constants $\alpha_{p_1p_2r}$ determined by the above
values of the constant term for the products $\prod_{i=1}^n\vartheta_{p_i}$.
In particular, this would tell us that the structure constants are determined
by ``classical" log Gromov-Witten invariants. This has the benefit that
frequently such invariants are easier to compute. 

We do not, however,
expect the strong conjecture to hold in complete generality.
Nevertheless, it has been proved in a number of useful situations.
See \cite{M19a},\cite{M19b} for 
some additional discussion on this conjecture. In particular, Mandel
has shown the Frobenius structure conjecture for all cluster varieties,
which includes the case of all log Calabi-Yau surfaces as originally
considered in \cite{GHKI}. More generally, see \cite{KY19} for a proof
of the Frobenius structure conjecture under the hypotheses that
$X\setminus D$ is affine and contains a dense torus.
Crucially, \cite{J22a} uses this fact to make the comparison between
our construction and the construction of Keel and Yu.
\end{remark}

Here we observe our results give

\begin{theorem}
Conjecture~\ref{conj:frob} holds for $n=2$ and $3$.
\end{theorem}

\begin{proof}
The case of $n=2$ is immediate by definition of our structure constants
$\alpha_{p_1p_2r}$. On the other hand, from Theorems 
\ref{theorem: independence of modulus},~\ref{thm: main gluing theorem},
and \ref{thm:main comparison}, it follows that the coefficient of
$t^A\vartheta_0$ in $\vartheta_{p_1}\vartheta_{p_2}\vartheta_{p_3}$
is $\deg[\scrM^{\ddagger}_y]^{\virt}$ for $y\in\overline{\scrM}_{0,4}$ a non-boundary
point. We have a diagram
\[
\xymatrix@C=30pt
{
\scrM^{\ddagger}_y\ar[r]^{i''}\ar[d]_{h'}&\scrM(X,\beta,z)\ar[d]^{h}\\
\foM_y^{\ddagger,\ev}\ar[d]_{\Phi'}\ar[r]_{i'}&\foM^{\ev}(\shX,\beta,z)\ar[d]^{\Phi}\\
y\times B\GG_m^{\dagger}\ar[r]_i&\overline{\scrM}_{0,4}^{\dagger}
}
\]
with all squares cartesian in all categories, as $i$ is strict. Note that
$\Phi$ is flat over a neighbourhood of the image of $i$. Indeed,
$\Phi$ is an integral morphism as the ghost sheaf over
$\overline{\scrM}_{0,4}^{\dagger}$ is $\NN$ in a neighbourhood
of the image of $i$. Further, $\Phi$ is log smooth, and log smooth
and integral morphisms are flat.
As $i$ is also a closed embedding, it is finite and we have 
\[
\deg [\scrM^{\ddagger}_y]^{\virt} = \deg i_*'' [\scrM^{\ddagger}_y]^{\virt}
= \deg h^!i'_*[\foM_y^{\ddagger,\ev}]
=\deg h^!\big(c_1(\Phi^*\shL)\cap [\foM^{\ev}(\shX,\beta,z)]\big),
\]
where $\shL$ is the pull-back of $\O_{\PP^1}(1)$ on $\overline{\scrM}_{0,4}
\cong\PP^1$ to $\overline\scrM_{0,4}^{\dagger}$.
That 
\[
\deg h^!(c_1(\Phi^*\shL)\cap [\foM^{\ev}(\shX,\beta,z]))
=
\int_{[\scrM(X,\beta,z)]^{\virt}} \psi_{x_{\out}}
\]
now follows exactly as in
Step 3 of the proof of Lemma~\ref{lemma: invarianceI}.
\end{proof}

\begin{remark}
After the initial submission of this paper, Conjecture \ref{conj:frob}
was proved by Johnston, see \cite{J22a},~Thm.~1.4.
\end{remark}

\appendix
\section{Tools of the trade}

We collect in this appendix a number of technical aspects of log geometry
which are needed in the main text. It may be skipped
by the casual reader, and referred to only when necessary.

\subsection{Charts for fibre products of fs log stacks}

Here we give several results useful for understanding local models for
fs fibre products of log stacks, and recall some notation.

We fix a ground field $\kk$ of characteristic $0$, as usual, 
and all schemes and stacks are defined over $\Spec\kk$.

We have the following standard definition:

\begin{definition}
\label{def:APK}
Given $P$ a fine monoid and $K\subseteq P$ a monoid ideal, we define
\begin{align*}
A_P:= \Spec\kk[P],\quad\quad\quad&A_{P,K}:=\Spec(\kk[P]/K)\\
\shA_P:=[A_P/\Spec\kk[P^{\gp}]],\quad\quad\quad&
\shA_{P,K}=[A_{P,K}/\Spec\kk[P^{\gp}]].
\end{align*}
Here both stacks carry a canonical log structure coming from
$P$, and the second stack carries a canonical idealized log structure
induced by the monoid ideal $K$, see \cite{Ogus}, III,\S1.3.
\end{definition}

\begin{proposition}
\label{prop:APrep}
Given $P$ a fine (fs) monoid, $\shA_P$ represents the functor
$T\mapsto \Hom(P,\Gamma(T,\overline{\shM}_T))$ in the category of
fine (fs) logarithmic stacks, and $A_P$ similarly represents the functor
$T\mapsto \Hom(P,\Gamma(T,\shM_T))$.
In particular, the set of two-isomorphism classes of log morphisms
$T\rightarrow\shA_P$ is in bijection with $\Hom(P,
\Gamma(T,\overline{\shM}_T))$.
\end{proposition}

\begin{proof}
The statement for $\shA_P$ is \cite{OlssonENS},~Prop.~5.17, while
the statement for $A_P$ is obvious.
\end{proof}

\begin{proposition}
\label{prop:APIrep}
Given $P$ an fs monoid and $K\subseteq P$ a monoid ideal, then $\shA_{P,K}$
represents the functor which associates to an fs log scheme $T$ the
set
\[
\left\{\varphi\in\Hom(P,\Gamma(T,\overline\shM_T))\,\Big |\,
\substack{\hbox{for any $p\in K$, 
any local lift $s\in\shM_T$ }\\
\hbox{of $\varphi(p)\in\overline\shM_T$ satisfies $\alpha_T(s)=0$}}\right\}.
\] 
Similarly, $A_{P,K}$ represents the functor which associates to an fs
log scheme $T$ the set
\[
\{\varphi\in\Hom(P,\Gamma(T,\shM_T))\, |\,
\hbox{for any $p\in K$,
$\varphi(p)\in\shM_T$
satisfies $\alpha_T(\varphi(p))=0$}\}.
\] 
\end{proposition}

\begin{proof}
By Proposition \ref{prop:APrep}, giving a $\varphi$ in the first set above
induces a morphism $T\rightarrow\shA_{P}$. The additional
condition then guarantees that this morphism factors through the closed
substack defined by the ideal $K$, i.e., $\shA_{P,K}$. Conversely, given
a morphism $T\rightarrow\shA_{P,K}$, we first obtain a composition
$T\rightarrow\shA_{P,K}\rightarrow\shA_P$, which is induced by some
$\varphi\in \Hom(P,\Gamma(T,\overline\shM_T))$ by Proposition~\ref{prop:APrep}.
Further, the condition on $\varphi$ is implied by the factorization of this
morphism through $\shA_{P,K}$. The second statement is similar.
\end{proof}

\begin{remark}
\label{rem:APghost}
If $P$ is a fine monoid, define $\overline{P} =P/P^{\times}$. Then 
$\shA_{P}\cong \shA_{\overline P}$.
\end{remark}

We recall that we write $\times^{\fine}$ and $\times^{\fs}$ for
the fibre products in the categories of fine log schemes and fs log
schemes respectively. We also write $\oplus^{\fine}$ and $\oplus^{\fs}$ for
the pushouts in the category of fine and fs monoids respectively.

\begin{proposition}
\label{prop:fibre product properties}
Given fine monoids $Q,P_1$ and $P_2$ with morphisms
$\theta_i:Q\rightarrow P_i$ inducing morphisms $\shA_{P_i}\rightarrow
\shA_Q$, $A_{P_i}\rightarrow A_Q$, we have
\begin{enumerate}
\item 
$\shA_{P_1}\times^{\fine}_{\shA_Q}\shA_{P_2}\cong 
\shA_{P_1\oplus^{\fine}_Q P_2}$ and
$A_{P_1}\times^{\fine}_{A_Q}A_{P_2}\cong 
A_{P_1\oplus^{\fine}_Q P_2}$.
\item If $P_1,P_2$ and $Q$ are also saturated, then
$\shA_{P_1}\times^{\fs}_{\shA_Q}\shA_{P_2}\cong \shA_{P_1\oplus^{\fs}_Q P_2}$
and
$A_{P_1}\times^{\fs}_{A_Q}A_{P_2}\cong A_{P_1\oplus^{\fs}_Q P_2}$.
\item $P_1\oplus^{\fine}_QP_2$ is the image of the canonical
map $P_1\oplus P_2\rightarrow P_1^{\gp}\oplus_{Q^{\gp}} P_2^{\gp}$,
while if $P_1,P_2,Q$ are saturated, then 
$P_1\oplus^{\fs}_QP_2$ is the saturation of $P_1\oplus^{\fine}_Q P_2$.
\item If in addition, $I\subseteq Q$, $J_i\subseteq P_i$ are monoid
ideals with $\theta_i(I)\subseteq J_i$, then
\[
\shA_{P_1,J_1}\times^{\fs}_{\shA_{Q,I}} \shA_{P_2,J_2}
\cong \shA_{P_1\oplus^{\fs}_Q P_2, J_1'\cup J_2'}
\]
and
\[
A_{P_1,J_1}\times^{\fs}_{A_{Q,I}} A_{P_2,J_2}
\cong A_{P_1\oplus^{\fs}_Q P_2, J_1'\cup J_2'},
\]
where $J_i'\subseteq P_1\oplus^{\fs}_Q P_2$ is the ideal generated
by the image of $J_i$. The same isomorphisms hold if instead
the fibre product is taken in the category of idealized fs log stacks.
\end{enumerate}
\end{proposition}

\begin{proof}
Statements (1) and (2) follow immediately 
from Proposition~\ref{prop:APrep},
while (4) follows in the same way from Proposition~\ref{prop:APIrep}.
The third statement is 
\cite{Ogus},~I,Prop.~1.3.4 in the fine case, and is immediate from this in 
the fs case. 
\end{proof}

\subsection{Log fibre dimension}
\label{sec:log fibre dim}
Again we fix a ground field $\kk$ of characteristic $0$, 
and all schemes and stacks are defined over $\Spec\kk$.

We first recall the notion of ordinary fibre dimension, using
the notation of \cite{EGAIV}, \S13. If $\ul{f}:\ul{X}\rightarrow \ul{Y}$
is a morphism of ordinary schemes, and $x\in X$ is a scheme-theoretic
point, $\bar x\rightarrow x$ a choice of geometric point, we write
\begin{equation}
\label{eq:ordinary fibre dim}
\dim_x \ul{f}^{-1}(\ul{f}(x))=\dim \O_{X_y,\bar x}+\mathrm{trdeg}_{\kappa(y)}
\kappa(x)
\end{equation}
for the fibre dimension of $\ul{f}$ at the point $x$.

The following definition, from \cite{AbbesSaito}, pg.~41, generalizes the
above notion to the logarithmic case and allows one 
to state a logarithmic analogue of invariance of fibre dimension of flat 
morphisms. 

\begin{definition}
Let $f:X\rightarrow Y$ be a morphism of log schemes with $\ul{f}$
locally of finite presentation and $x\in X$. Let $y=f(x)$, and let
$X_y$ be the scheme-theoretic fibre with the induced log structure,
and $\kappa(x)$, $\kappa(y)$ the residue fields at $x$ and $y$ respectively.
Choose geometric points $\bar x\rightarrow x$, $\bar y\rightarrow y$.
Then we define
\begin{align}
\label{eq:log fibre dim def}
\begin{split}
&\dim_x^{\log}f^{-1}(f(x)) := \\
&\dim \O_{X_y,\bar x}/\langle \alpha_X(\shM_{X,\bar x}
\setminus \O_{X,\bar x}^{\times})\rangle + \mathrm{trdeg}_{\kappa(y)}
\kappa(x)
+\rank \overline{\shM}^{\gp}_{X,\bar x} -
\rank \overline{\shM}^{\gp}_{Y,\bar y}.
\end{split}
\end{align}
Here we remark that for a fine log scheme $X$, $\O_{X,\bar x}/
\langle \alpha_X(\shM_{X,\bar x}\setminus \O_{X,\bar x}^{\times})\rangle$
is the (Henselian) local ring at $\bar x$ of the 
log stratum of $X$ containing $\bar x$.\footnote{We note that
this formula is not precisely the one given in \cite{AbbesSaito}, where
the last two terms are replaced with $\rank \overline{\shM}^{\gp}_{X,\bar x}/
\overline{\shM}^{\gp}_{Y,\bar y}$, which is not the same if
$\bar f^{\flat}:
\overline{\shM}^{\gp}_{Y,\bar y}\rightarrow \overline{\shM}^{\gp}_{X,\bar x}$
is not injective. However, Lemma 3.10, 2 of \cite{AbbesSaito} fails
with their definition, and the formula we give here makes that statement
true.}

\end{definition}

The following, which states that log flat morphisms are 
flat strata-wise,
should be standard, but we have been unable to find a reference:

\begin{lemma}
\label{lemma: flat strata}
Let $f:X\rightarrow Y$ be a log flat morphism between fs log schemes
with $\ul{f}$ locally of finite presentation, and $x\in X$, $y=f(x)$.
Then
\[
A_{\bar x}:=\O_{X,\bar x}/\langle \alpha_X(\shM_{X,\bar x}\setminus\O_{X,\bar x}^{\times})
\rangle
\]
is flat over 
\[
B_{\bar y}:=\O_{Y,\bar y}/\langle \alpha_Y(\shM_{Y,\bar y}\setminus\O_{Y,\bar y}^{\times})
\rangle.
\] 
\end{lemma}

\begin{proof}
First, suppose $\theta:Q\rightarrow P$ is an injective homomorphism
of fs monoids, and take $X=A_P$, $Y=A_Q$, and $f:X\rightarrow Y$
induced by $\theta$. Take $x\in A_P$ such that $\overline{\shM}_{X,\bar x}
=\overline P:=P/P^{\times}$. By localizing $Q$ along a face if necessary,
we may assume that $f(x)=y\in Y$ with 
$\overline{\shM}_{Y,\bar y}=\overline Q=Q/Q^{\times}$. 
If $X^+:=\Spec\kk[P]/\langle
P\setminus P^{\times}\rangle$, 
$Y^+:=\Spec\kk[Q]/\langle Q\setminus Q^{\times}\rangle$
viewed as closed subschemes of $X, Y$ respectively, then
$x\in X^+$, $y\in Y^+$, and $A_{\bar x}=\O_{X^+,\bar x}$,
$B_{\bar x}=\O_{Y^+,\bar y}$. On the other hand,
$X^+\cong \Spec \kk[P^{\times}]$, $Y^+\cong \Spec\kk[Q^{\times}]$,
and $\theta:Q^{\times}\rightarrow P^{\times}$ is an injective map of
abelian groups, so $X^+$ is flat over $Y^+$. This shows the claim in
this case.

In general, note that the claim is local for $X$ and $Y$ in the
(strict) fppf topology. We may then apply \cite{Ogus},~IV,Lem.~4.1.3
and assume that there is a $\theta:Q\rightarrow P$ as above and a
commutative diagram
\[
\xymatrix@C=15pt
{
X\ar[r]\ar[d]&A_P\ar[d]\\
Y\ar[r]&A_Q
}
\]
with horizontal arrows strict,
$A_P\rightarrow A_Q$ induced by $\theta$, and the induced
morphism $X\rightarrow Y\times_{A_Q} A_P$ flat. Then the desired
flatness follows from this diagram and the case already considered
in the first paragraph.
\end{proof}

\begin{proposition}
\label{prop: flat fibre dim}
Let $f:X\rightarrow Y$ be a log flat 
morphism between fs log schemes, with $\ul{f}$
locally of finite presentation. Then:
\begin{enumerate}
\item $\dim_x^{\log} f^{-1}(f(x))$ is a locally constant function of $x$.
\item If $f$ is  $\QQ$-integral (see \cite{Ogus}, I,Def.~4.7.4)  
in a neighbourhood of $x$ 
or if $f$ is strict, then 
\[
\dim_x^{\log} f^{-1}(f(x)) = \dim_x \ul{f}^{-1}(\ul{f}(x)).
\]
\item If $g:Y\rightarrow Z$ is also log flat, and $x\in X$,
$y=f(x)$, $h=g\circ f$, then
\[
\dim_x^{\log} h^{-1}(h(x))=\dim_x^{\log} f^{-1}(f(x))+
\dim_y^{\log} g^{-1}(g(y)).
\]
\item Given an injective homomorphism of fs monoids $\theta:Q\rightarrow P$,
the induced morphism $f:A_P\rightarrow A_Q$ has log fibre dimension
$\rank P^{\gp}-\rank Q^{\gp}$.
\item If $Y'$ is an fs log scheme with a morphism $Y'\rightarrow Y$,
$f':X'=X\times^{\fs}_Y Y'\rightarrow Y'$ the base-change of $f$, 
$\pi:X'\rightarrow X$ the projection, then
\[
\dim^{\log}_{x'} (f')^{-1}(f'(x'))=\dim_{\pi(x')}^{\log} f^{-1}(f(\pi(x')))
\]
for $x'\in X'$.
\end{enumerate}
\end{proposition}

\begin{proof}
(1) is \cite{AbbesSaito}, Lem.~3.10,~2.

The first case of (2) is the implication $(1)\Rightarrow (7)$ of
\cite{Ogus},~III,Thm.~2.2.7, comparing
the formula given there with \eqref{eq:log fibre dim def} and
\eqref{eq:ordinary fibre dim}. We also need to apply 
\cite{Ogus},~IV,Prop.~4.1.9 to
verify the needed hypothesis of solidity of the idealized fibre.
If instead $f$ is strict, the claim follows immediately from the
definition, bearing in mind that 
$\alpha_Y(s)|_{\bar y}=0$ for any 
$s\in \shM_{Y,\bar y}\setminus\O_{Y,\bar y}^{\times}$, and
the ideal $\langle \alpha_X(\shM_{X,\bar x}\setminus \O_{X,\bar x}^{\times})
\rangle$ in $\O_{X_y,\bar x}$ is generated by pull-backs of such
$\alpha_Y(s)$ in the strict case.

(3) From the definitions, we need to show that
\begin{align*}
\dim \O_{X_z,\bar x}/\langle \alpha_X(\shM_{X,\bar x}\setminus\O_{X,\bar x}^{\times})\rangle
= {} &\\
\dim \O_{X_y,\bar x}/\langle\alpha_X(\shM_{X,\bar x}&\setminus\O_{X,\bar x}^{\times})\rangle
+
\dim \O_{Y_z,\bar y}/\langle\alpha_Y(\shM_{Y,\bar y}\setminus\O_{Y,\bar y}^{\times})\rangle.
\end{align*}
This follows from the flatness statement of Lemma~\ref{lemma: flat strata}
and additivity of fibre dimension for ordinary flat morphisms.

(4) First note that $f$ is log flat, so by (1) we may compute the log fibre
dimension at any point $x\in \Spec \kk[P^{\gp}]\subseteq A_P$. Since
$f$ maps $\Spec\kk[P^{\gp}]$ into $\Spec\kk[Q^{\gp}]$, and the log structures
are trivial on these open sets, the log fibre dimension at $x$ coincides
with the ordinary fibre dimension by (2), which is clearly 
$\rank P^{\gp}-\rank Q^{\gp}$.

(5) Let $x=\pi(x')$, $y=f(x)$, $y'=f'(x')$. 

We first consider the
case that $Y'\rightarrow Y$ is strict. As $X'_{y'}\rightarrow X_y$ is
a strict base-change via $y'\rightarrow y$, strata of $X'_{y'}$ are,
locally, base-changes of strata of $X_y$ by $y'\rightarrow y$. Thus by
\cite{stacks}, Tag 02FY and the definition of log fibre dimension, the
result is immediate.

Now consider $Y'\rightarrow Y$ arbitrary.
To prove the desired result, by (1)
we can assume that $x$ and $y$ are closed points. By the additivity statement
(3) for ordinary fibre dimension for compositions of (strict) flat morphisms,
we may replace $X$ and $Y$ with fppf neighbourhoods of $x$ and $y$. 
In particular, by \cite{Ogus}, IV,Lem.~4.1.3,
we may assume the existence of a chart
$Y\rightarrow A_Q$ inducing an isomorphism
$Q\cong \overline{\shM}_{Y,\bar y}$ and a flat chart for $f$,
i.e., an injective map $\theta:Q\rightarrow P$
and a chart $X\rightarrow A_P$ such that $g:X\rightarrow Y\times_{A_Q}A_P$
is flat. We thus obtain a diagram
\[
\xymatrix@C=30pt
{
X'\ar[d]_{g'}\ar[r]&X\ar[d]^g&\\
Y'\times_{A_Q} A_P\ar[d]_{h'}\ar[r]&Y\times_{A_Q} A_P\ar[r]\ar[d]^h&A_P\ar[d]\\
Y'\ar[r]&Y\ar[r]&A_Q
}
\]
Then $g$, $g'$ are strict flat, and hence 
the ordinary and the log fibre dimensions of $g$ (resp.\ $g'$) agree
at $x$ (resp. $x'$).
Further, by \cite{stacks}, Tag 02FY again, $g$ and $g'$ have the
same ordinary fibre dimensions.
Thus by (3) of the proposition, it is enough to show that 
$h'$ has the same log fibre dimension as $h$.

First, as 
$Y\rightarrow A_Q$ is strict, by the case of the result already
shown, the log fibre dimension of $h$ is the same as that of 
$A_P\rightarrow A_Q$, which by (4) is $\rank P^{\gp}-\rank Q^{\gp}$.

Second, again by the strict case, we may replace $Y'$ with an fppf
neighbourhood of $y'$ and assume given a chart $Y'\rightarrow A_{Q'}$.
Thus we have a diagram
\[
\xymatrix@C=30pt
{
Y'\times_{A_Q} A_P\ar[d]_{h'}\ar[r]&A_{P'}\ar[r]\ar[d]&A_P\ar[d]\\
Y'\ar[r] & A_{Q'}\ar[r]& A_Q
}
\]
with the left-hand square cartesian in all categories and the right-hand
square cartesian in the fs log category. In particular,
$P'=P\oplus^{\fs}_Q Q'$ by Proposition~\ref{prop:fibre product properties},
(2). Thus, again by the strict case already proved, 
the log fibre dimension of $h'$ agrees with that of
$A_{P'}\rightarrow A_{Q'}$. Now since $\theta:Q\rightarrow P$
is injective, so is the induced homomorphism $Q'\rightarrow P'$,
as the same is true at the group level. Thus the log fibre dimension
of $A_{P'}\rightarrow A_{Q'}$ is
\begin{align*}
\rank (P')^{\gp}-\rank (Q')^{\gp} = {} & \big(\rank P^{\gp}+\rank (Q')^{\gp}
-\rank Q^{\gp}\big) -\rank (Q')^{\gp}\\
= {} & \rank P^{\gp} - \rank Q^{\gp}.
\end{align*}
This shows invariance of log fibre dimension under fs base-change.
\end{proof}

While the above statements are all given for morphisms of schemes,
they will often be applied to morphisms of stacks. It is easy to
adapt the above statements using suitable smooth charts, and this will
in general be done without comment.

Here is an immediate application of the above results.

\begin{proposition}
\label{prop:log etale dimensions}
Let $f:X\rightarrow Y$ be a morphism of fs log stacks.
\begin{enumerate}
\item If $f$ is log \'etale, then $f$ has log fibre dimension $0$.
\item If $f$ is integral, log \'etale and $Y$ pure-dimensional,
then $\dim X=\dim Y$.
\item If $X$ and $Y$ carry coherent idealized structures and
$f$ is integral and idealized log \'etale (but not necessarily log
\'etale), then $\dim X \le \dim Y$.
\end{enumerate}
\end{proposition}

\begin{proof}
For (1),
choosing $\bar x\in |X|$, $\bar y= f(\bar x)\in |Y|$, we have strict \'etale
neighbourhoods $X'$, $Y'$ of $\bar x$, $\bar y$ respectively
with $X'\rightarrow \shA_{P}\times_{\shA_Q} Y'$ strict \'etale 
by \cite{ACGS18}, Prop.~B.4.
An ordinary \'etale morphism of stacks always has fibre dimension
$0$. Indeed, this follows immediately by passing to a schematic
\'etale chart for the morphism using the definition of an
\'etale morphism of stacks given in \cite{stacks}, Tag 0CIL (see
Remark~\ref{rem:etale stack}).
As $X'\rightarrow \shA_{P}\times_{\shA_Q} Y'$ is strict \'etale, it
is thus \'etale, and hence
has log fibre dimension
$0$ by Proposition~\ref{prop: flat fibre dim}, (2). 
On the other hand, the composition $\shA_P\times_{\shA_Q} Y'\rightarrow
Y'$ has the same log fibre dimension as $\shA_P\rightarrow\shA_Q$,
which has log fibre dimension $0$. Indeed, we pass to the smooth
chart $A_P\rightarrow A_Q$ for the morphism $\shA_P\rightarrow\shA_Q$,
and conclude by Proposition~\ref{prop: flat fibre dim}, (4). Thus
$X'\rightarrow Y'$ has log fibre dimension $0$ by Proposition
\ref{prop: flat fibre dim}, (3) applied to the composition
$X'\rightarrow \shA_P\times_{\shA_Q} X'\rightarrow Y'$. 
We conclude that $f:X\rightarrow Y$ also has log fibre dimension $0$.

For (2) as integral log \'etale morphisms are log flat 
and flat, 
by Proposition~\ref{prop: flat fibre dim}, (2) the result follows 
immediately from (1) and the dimension formula for flat morphisms.

For (3), again by \cite{ACGS18}, Prop.~B.4, we can 
assume the existence of \'etale neighbourhoods $X'$ and $Y'$ of $X$ and $Y$
with a morphism $X'\rightarrow Y'$ factoring as
$X'\rightarrow Y'\times_{\shA_{Q,J}} \shA_{P,K}\rightarrow Y'$ and the
first morphism \'etale.
Now there is a commutative diagram with the left-hand square cartesian:
\[
\xymatrix@C=30pt
{
Y'\times_{\shA_{Q,J}} \shA_{P,K}\ar[d]\ar[r] &\shA_{P,K}\ar[d]\ar[r]&
\shA_P\ar[d]^{\theta}\\
Y'\ar[r]&\shA_{Q,J}\ar[r]&\shA_Q
}
\]
By integrality, $\theta$ is a flat morphism, necessarily of fibre
dimension $0$, and because the right-hand horizontal arrows are
closed embeddings, $Y'\times_{\shA_{Q,J}} \shA_{P,K}$
is a closed substack of $Y'\times_{\shA_Q}\shA_P$. By flatness of
$\theta$, the latter stack has the same dimension as $Y'$,
and hence $\dim X'\le \dim Y'$, as desired.
\end{proof}

The last observation we make here is often useful for determining
dimensions of various log stacks which appear in this paper.

\begin{proposition}
\label{prop: realisability}
Let $f:X\rightarrow Y$ be a log smooth morphism of fs log stacks,
with $\bar x\in |X|$, $\bar y= f(\bar x)$ geometric points. We then
obtain a morphism of cones $\Sigma(f):\sigma_{\bar x}\rightarrow
\sigma_{\bar y}$. Let $\tau$ be a face of $\sigma_{\bar x}$. Then the following
are equivalent:
\begin{enumerate}
\item
There is a geometric point $\bar x'\in |X|$ which is a generization of $\bar x$
and such that the induced inclusion of faces $\sigma_{\bar x'}\rightarrow
\sigma_{\bar x}$ identifies $\sigma_{\bar x'}$ with $\tau$.
\item
There is a generization $\bar y'\in |Y|$ of $\bar y$ inducing an inclusion
$\sigma_{\bar y'}\subseteq \sigma_{\bar y}$ with $\sigma_{\bar y'}$
the minimal face of $\sigma_{\bar y}$ containing $\Sigma(f)(\tau)$.
\end{enumerate}
\end{proposition}

\begin{proof}
For $(1)\Rightarrow (2)$, note that if $\bar x'$ exists, then
$\bar y'=f(\bar x')$ is a generization of $\bar y$. Thus
$\sigma_{\bar y'}$ is identified with a face of $\sigma_{\bar y}$,
and since $\bar f^{\flat}$ is a local homomorphism of monoids,
$\Sigma(f)(\tau)$ must intersect the interior of $\sigma_{\bar y'}$,
so that $\sigma_{\bar y'}$ is the minimal face of $\sigma_{\bar y}$ containing
$\Sigma(f)(\tau)$.

For the converse, assume (2).
It is sufficient to show (1) holds after replacing $X$ and $Y$ with
schematic smooth neighbourhoods of $\bar x$ and $\bar y$,
so we can assume there is a diagram
\[
\xymatrix@C=30pt
{
X\ar[d]\ar[r] &A_P\ar[d]\\
Y\ar[r] & A_Q
}
\]
with the right-hand vertical arrow induced by an injective monoid homomorphism
$\theta:Q\rightarrow P$, the horizontal morphisms strict,
and $X$ strict smooth over $Y\times_{A_Q} A_P$.
The right-hand vertical arrow tropicalizes to
\[
\theta^t:\sigma_{\bar x}=\Hom(P,\RR_{\ge 0})\rightarrow
\sigma_{\bar y}=\Hom(Q,\RR_{\ge 0}).
\]
It is then standard
toric geometry that the torus orbit of $A_P$ corresponding
to a face $\tau\subseteq\sigma_{\bar x}$ surjects onto the 
torus orbit
of $A_Q$ corresponding to the minimal face $\tau'$ of
$\sigma_{\bar y}$ containing $\theta^t(\tau)$. 

Now let $\bar z\in
Y\times_{A_Q} A_P$ be the image of $\bar x$, necessarily with 
image $\bar y\in Y$. Let $\bar w\in A_P$ be the image of $\bar z$ in
$A_P$ and $\bar s\in A_Q$ be the image of $\bar y$ in $A_Q$.

Suppose that there is a generization
$\bar y'$ of $\bar y$ with $\sigma_{\bar y'}=\tau'$, so that the
image $\bar s'$ of $\bar y'$ in $A_Q$ lies in the stratum corresponding
to $\tau'$. Then by the above mentioned surjectivity of strata,
there is a point $\bar w'\in A_P$
which specializes to $\bar w$,
satisfies $\sigma_{\bar w'}=\tau$,
and the image of $\bar w'$ in $A_Q$ is $\bar s'$. Now
$\bar z\in\cl(\bar y'\times_{\bar s'}\bar w')$, and hence there is a point
$\bar z'$ of $\bar y'\times_{\bar s'}\bar w'$ which specializes
to $\bar z$. Since a smooth morphism is 
locally of finite presentation and flat, hence open, there
is a generization $\bar x'$ of $\bar x$ mapping to $\bar z'$, showing the 
result.
\end{proof}

\pagebreak

\subsection{Log structures on discrete valuation rings}

In this subsection, we let $R$ be a discrete valuation ring with 
field of fractions $K$. Set $T=\Spec R$, with generic point $\xi$ and
closed point $0$. Also let $U=\Spec K$, and let $t\in R$ be a choice
of uniformizing parameter.

\begin{lemma}
\label{lem:log dash}
Let $Q$ be a sharp fs monoid.
The set of isomorphism classes of Zariski fs log structures
$\shM_T$ on $T$ with $\overline{\shM}_{T,0}=Q$ is canonically in bijection 
with the set of pairs $(F,u)$ where $F\subseteq Q$ is a face and $u\in
\Int(F^{\vee})$.

For a log structure $\shM_T$ corresponding to $(F,u)$, we have
$\overline{\shM}_{T,\xi}=Q_K:=(F^{-1}Q)/(F^{-1}Q)^{\times}$.
\end{lemma}

\begin{proof}
Note $u\in \Int(F^{\vee})$ is equivalent to $u:F\rightarrow
\NN$ being a local homomorphism.
Suppose given the data $(F,u)$.
Then the chart
\[
\alpha:Q\rightarrow R, \quad \alpha(q)=
\begin{cases} 
0 & q\in Q\setminus F,\\
t^{u(q)} & q\in F
\end{cases}
\]
defines a log structure $\shM_T$ on $T$. 
Then by
\cite{Ogus}, II,Prop.~2.1.4,
$\overline{\shM}_{T,0}=Q/\alpha^{-1}(R^{\times})=Q$, as 
$u$ local implies that
$u(q)\not=0$ unless $q=0$. 
On the other hand, $\overline\shM_{T,\xi}=Q/\alpha^{-1}(R\setminus 0)=Q/F
= Q_K$.

Conversely, suppose given a Zariski fs log structure $\shM_T$ with
ghost sheaf stalks $Q$, $Q_K$ at $0$ and $\xi$, respectively. Let
$\chi:Q\rightarrow Q_K$ be the generization map, with $\chi^{-1}(0)=F$. 
Now there are (non-canonical) splittings 
\begin{equation}
\label{eq:splittings}
\Gamma(T,\shM_T)=R^{\times} \times Q, \quad
\Gamma(U,\shM_T)=K^{\times} \times Q_K.
\end{equation}
This follows, for example, from the splitting statement in
\cite{Ogus}, III,Thm.~1.2.7, 2, 
The restriction map $\rho:\Gamma(T,\shM_T)\rightarrow\Gamma(U,\shM_T)$
must then take the form
\[
\rho(r,q)= (r \rho'(q), \chi(q))
\]
for some homomorphism $\rho':Q\rightarrow K^{\times}$.
Writing $\nu_R$ for the discrete valuation on $K$, we can write
\[
\rho'(q)=\rho''(q) t^{\nu_R(\rho'(q))}
\]
for some $\rho'':Q\rightarrow R^{\times}$. Set $u=\nu_R\circ\rho':Q\rightarrow
\ZZ$.

Note we can modify the
splitting of $\Gamma(T,\shM_T)$ via an automorphism of the form
$(r,q)\mapsto (r\rho''(q)^{-1},q)$. This has the effect of setting
$\rho''\equiv 1$, so that the restriction map is given by
$\rho(r,q)=(rt^{u(q)},q)$. Furthermore, we can modify the splitting
of $\Gamma(U,\shM_T)$ via an automorphism $(r,q)\mapsto (rt^{\bar u(q)},r)$
for any $\bar u:Q_K\rightarrow\ZZ$. This has the effect of replacing
$u$ with $u+\bar u\circ\chi$, so that the restriction map
is given by $(r,q)\mapsto (rt^{u(q)+\bar u\circ\chi(q)},r)$.
Thus $u$ can be viewed as an element
of $Q^*/Q_K^*\cong F^*$.

Thus the data of $F$ and $u$ determine, up to isomorphism, the
sheaf $\shM_T$. Now let us consider possible structure maps
$\alpha_T:\shM_T\rightarrow\O_T$. Necessarily,
$\alpha_T:\Gamma(U,\shM_T)\rightarrow K$ is given by
\[
(r,q)\mapsto \begin{cases} r & q=0,\\ 0 & q\not=0.\end{cases}
\]
For the structure map $\alpha_T:\Gamma(T,\shM_T)\rightarrow R$
to be compatible with restriction, we thus need for $(r,q)\in R^{\times}\times
Q$ that
\[
\alpha_T(r,q)=\begin{cases} rt^{u(q)} & \chi(q)=0,\\0 & \chi(q)\not=0.
\end{cases}
\]
Note that if $u$ is replaced by $u+\bar u\circ\chi$, the value of
$u(q)$ doesn't change when $\chi(q)=0$, so $\alpha_T$ is well-defined.
In particular, $u$ must be non-negative on $F$. Thus $u\in F^{\vee}$.
Furthermore, the requirement that $\alpha_T^{-1}(R^{\times})\cong
R^{\times}$ then also implies $u(q)>0$ if $q\not=0$ but $\chi(q)=0$.
Thus $u\in \Int(F^\vee)$. This gives the pair $(F,u)$.
\end{proof}

\begin{lemma}
\label{lem:log dash morphism}
Let $\shM_T$, $\shM_T'$ be two log structures on
$\ul{T}$ given by the data $F\subset Q$, $u\in \Int(F^{\vee})$ and
$F'\subset Q'$, $u'\in\Int((F')^{\vee})$ respectively. 

Given local homomorphisms $\varphi,\varphi_K$ 
in a commutative diagram
\[
\xymatrix@C=30pt
{
Q\ar[d]_{\chi}&Q'\ar[l]_{\varphi}\ar[d]^{\chi'}\\
Q_K&Q_K'\ar[l]^{\varphi_K}
}
\]
satisfying $u(\varphi(q))=u'(q)$ for $q\in F'$, there exists
a morphism of log schemes 
$f:(\ul{T},\shM_T)\rightarrow (\ul{T},\shM_T')$
with $\ul{f}$ the identity and $f$
inducing the maps $\varphi$, $\varphi_K$ on stalks of ghost 
sheaves.
\end{lemma}

\begin{proof}
Given the data, choose an extension of $u$ to a homomorphism $u:Q\rightarrow
\NN$. Then $u\circ\varphi$ is an extension of $u'$ to a morphism 
$u':Q'\rightarrow\NN$. We may then choose splittings
for $\shM_T$, $\shM_{T}'$ as in \eqref{eq:splittings} so that the
restriction maps are given by $(r,q)\mapsto (rt^{u(q)},\chi(q))$,
$(r,q)\mapsto  (rt^{u'(q)},\chi'(q))$ for $\shM_T, \shM_T'$. Then define
$f^{\flat}$ on sections by
\begin{align*}
\Gamma(T,\shM_T')\rightarrow \Gamma(T,\shM_T),&
\quad (r,q)\mapsto (r,\varphi(q)),\\
\Gamma(U,\shM_T')\rightarrow \Gamma(U,\shM_T),&
\quad (r,q)\mapsto (r,\varphi_K(q)).
\end{align*}
One then checks compatibility of these maps with restriction
and with $\alpha_T$, $\alpha_{T'}$ to see we obtain a log morphism.
\end{proof}

\subsection{Some results concerning virtual fundamental classes}

We prove several general results about the behaviour of virtual
fundamental classes as operational Chow classes. We feel that there
should be a good reference for these results, but we have been unable to find
one. The proof of the following lemma was provided to us by Andrew
Kresch.

\begin{lemma}
\label{lem:kresch lemma}
Suppose given a cartesian diagram of algebraic stacks over a field of
characteristic zero
\[
\xymatrix@C=30pt
{
M'\ar[r]^{g'}\ar[d]_{f'}& M\ar[d]^f\\
F'\ar[r]_g&F
}
\]
with:
\begin{itemize}
\item
$M$ Deligne-Mumford;
\item $F'$ integral, 
and $F$ and $F'$ stratified by quotient stacks
in the sense of \cite{Kresch}, Def.~3.5.3; 
\item
$f$ is 
equipped with a perfect relative obstruction theory of relative virtual
dimension $p$, and $p+\dim F<0$.
\end{itemize}
Then $[M']^{\virt}=f^!([F'])=0$ in the rational
Chow group of $M'$, where $f^!$ is the virtual pull-back of Manolache
\cite{Man} defined using the given obstruction theory, and the virtual
fundamental class is calculated using the pull-back obstruction 
theory.
\end{lemma}

\begin{proof}
Note that if $g$ were were a lci morphism, then the result
would follow from compatibility of Gysin pull-back and virtual pull-back.
We thus work to reduce to this situation.

If $F''$ denotes the closure of the image of $g$ with its reduced induced
stack structure, the morphism $g$ factors as $F'\rightarrow F''\rightarrow
F$. Pulling back the relative obstruction theory for $f$ to 
$M'':=M\times_F F''\rightarrow F''$, we note that $M''$ is still DM
and $p+\dim F''<0$ as $\dim F''\le \dim F$. Thus the hypotheses
of the lemma still hold if we replace $F$ and $M$ by $F''$ and $M''$. 
So we may make this replacement, and can therefore assume that $F$
is integral and $g$ is dominant.

By \cite{Tem}, Thm~5.1.1, there exists a projective resolution of
singularities $\rho:\widetilde F\rightarrow F$. Let $\widetilde F'$ be
the unique irreducible component of $\widetilde F\times_F F'$ dominating
$F'$ with its reduced induced stack structure. This exists
because $g$ is dominant and $\widetilde F\rightarrow F$ is an isomorphism
on an open set.
In particular, the induced
morphism $\rho':\widetilde F'\rightarrow F'$ is birational and projective.
Let $\tilde g:\widetilde F'\rightarrow \widetilde F$ be the induced
projection. We write $\widetilde M=M\times_F \widetilde F$,
$\widetilde M'=M\times_F \widetilde F'$.

We have a diagram
\[
\xymatrix@C=30pt
{
\widetilde M'\ar[r]\ar[d]
& \widetilde F'\times \widetilde M \ar[r]^{\pr_2}\ar[d] 
& \widetilde M\ar[d]^{\tilde f}\ar[r]&M\ar[d]^f\\
\widetilde F'\ar[r]_{\Gamma_{\tilde g}} &\widetilde F'\times
\widetilde F\ar[r]_{\pr_2}& \widetilde F\ar[r]_{\rho}&F
}
\]
Here $\Gamma_{\tilde g}$ denotes the graph of $\tilde g$.
Note that because $\widetilde F$ is non-singular, $\Gamma_{\tilde g}$
is an lci morphism, as follows from \cite{stacks},
Tag 0CJ2. In this case, \cite{Kresch} defines a Gysin pull-back map
$\Gamma^!_{\tilde g}$.

We now obtain
\[
f^![\widetilde F'] = f^!\Gamma_{\tilde g}^![\widetilde
F'\times\widetilde F]
= \Gamma_{\tilde g}^! f^![\widetilde F'\times\widetilde F]
= \Gamma_{\tilde g}^!f^!\pr_2^*[\widetilde F]
= \Gamma_{\tilde g}^!(\pr_2)^*f^![\widetilde F]
= 0.
\]
Here the first equality follows from the fact that the Gysin pull-back
of the fundamental class is the fundamental class.
The second equality comes from compatibility of virtual
pull-back with Gysin pull-back along regular embeddings, as follows
from \cite{Man}, Thm.~4.3 and the fact that Gysin pull-back along
regular embeddings can be described as a virtual pull-back, see \cite{Man},
Rmk.~3.9.
The third comes from the fact
that $\pr_2$ is flat, and the fourth by compatibility of virtual
pull-back with flat pull-back, \cite{Man} Thm. 4.1, (ii).
The final equality follows from
dimension reasons, as $f^![\widetilde F]$ is a Chow class of
negative degree, and as $\rho$ is projective, hence representable,
$\widetilde M$ is DM and so has no negative degree rational
Chow classes.

We next observe we have a cube with all vertical faces cartesian:
\[
\xymatrix@C=20pt
{
&M'\ar'[d][dd]\ar[rr] & & M
\ar[dd]^>>>>>>>>>>f \\ 
\widetilde M'\ar[ur]^{\rho''}\ar[rr]\ar[dd]
& & \widetilde M\ar[ur]\ar[dd]^>>>>>>>>>>>>>>{\tilde f} & \\
&F'\ar'[r][rr]  & & F
\\
\widetilde F' \ar[rr]\ar[ur]_>>>>{\rho'} & & \widetilde F \ar[ur]_{\rho}
}
\]
By compatibility of virtual pull-back with projective push-forward,
\cite{Man}, Thm.~4.1,~(i),
$0=\rho''_*f^![\widetilde F']=
f^!\rho'_*[\widetilde F']=f^![F']$, as desired. Finally, note that
$f^![F']=(f')^![F']$, where $(f')^!$ is virtual
pull-back defined using the pull-back obstruction theory, by \cite{Man},
Thm.~4.1,~(iii). However by definition $(f')^![F']=
[M']^{\virt}$.
\end{proof}

Lemma~\ref{lem:kresch lemma} is needed for the following theorem, 
which can be interpreted
in a stable map context as saying that if we glue two families of
curves, one of which has negative virtual dimension, then the virtual
fundamental class of the glued family is zero, even if it has positive
virtual dimension.

\begin{theorem}
\label{thm:virtual vanishing}
Suppose given algebraic stacks $M,M_1,M_2$, $F,F_1,F_2$ over a field
of characteristic zero and a cartesian diagram
\[
\xymatrix@C=30pt
{
M\ar[r] \ar[d]_f & M_1\times M_2\ar[d]^{f_1\times f_2}\\
F\ar[r] & F_1\times F_2
}
\]
with:
\begin{itemize}
\item $M,M_1,M_2$ Deligne-Mumford;
\item $F, F_1, F_2$ stratified by quotient stacks;
\item $f_1$, $f_2$ are equipped with perfect relative obstruction
theories, with $f_1$ having relative virtual dimension $p$ and
$p+\dim F_1<0$;
\item $F$ is pure-dimensional.
\end{itemize}
Then $(f_1\times f_2)^![F]=0$ in the rational Chow group of $M$, 
i.e., the virtual fundamental class of $M$ induced by the pull-back relative
obstruction theory vanishes in the rational Chow group of $M$.
\end{theorem}

\begin{proof}
We have the following larger diagram with all squares cartesian:
\[
\xymatrix@C=30pt
{
&&M_2\ar[d]^{f_2}\\
M\ar[r]\ar[d]&M_1\times M_2\ar[ur]\ar[d]&F_2\\
M'\ar[r]\ar[d]&M_1\times F_2\ar[ur]\ar[d]\ar[r]&M_1\ar[d]^{f_1}\\
F\ar[r]&F_1\times F_2\ar[r] &F_1
}
\]
By replacing $F$ by an irreducible component, we can assume $F$ is integral.
By \cite{Man}, Thm.~4.8 and Thm.~4.1,~(iii), $(f_1\times f_2)^!=
f_2^!\circ f_1^!$. However, by Lemma~\ref{lem:kresch lemma},
$f_1^!([F])=0$ in rational Chow, and hence the result.

\end{proof}

\end{document}